\documentclass[12pt,leqno]{amsart}

\usepackage{amsmath,amsthm,amscd,amssymb,eucal,xspace}
\usepackage[PostScript=dvips,balance,midshaft,nohug]{diagrams}
\usepackage{pgf, pgfplots}\usepackage{verbatim}

\usepackage{enumitem,calc,xcolor,tcolorbox,float}
\usepackage[vcentermath]{youngtab}
\usepackage{tikz}\usepackage{tikz-cd}
\usepackage[columns=3, font=footnotesize]{idxlayout}

\makeindex
\usepackage{graphicx}
\swapnumbers\DeclareSymbolFont{bbold}{U}{bbold}{m}{n}
\DeclareSymbolFontAlphabet{\mathbbold}{bbold}
\newarrow{Equals} =====
\newarrow{Implies} ===={=>}
\newarrow{Onto} ----{>>}
\newarrow{Into} C--->
\newarrow{IntoA} {hooka}--->
\newarrow{IntoB} {hookb}--->
\newarrow{Dotsto} ....>
\newarrowhead{sp}{\ }{\ }{\ }{\ }
\newarrowtail{sp}{\ }{\ }{\ }{\ }
\newarrow{Line}----{-}
\newarrow{Eq}{sp}==={sp}
\newtheorem{thm}[subsection]{Theorem}
\newtheorem*{thm*}{Theorem}
\newtheorem{cor}[subsection]{Corollary}
\newtheorem{lem}[subsection]{Lemma}
\newtheorem*{lem*}{Lemma}
\newtheorem{prop}[subsection]{Proposition}
\newtheorem*{prop*}{Proposition}

\theoremstyle{definition}
\newtheorem{defn}[subsection]{Definition}
\numberwithin{equation}{subsection}

\newcommand\T[1]{\rule{0pt}{#1ex}}    
\newcommand\B[1]{\rule[-#1ex]{0pt}{1pt}}
\renewcommand\mod[1]{\ (\mathop{\rm mod}#1)}
\newcommand{\quash}[1]{}

\newcommand{\Sp}{{\rm{Sp}}}
\newcommand{\mb}{\mathbf}
\newcommand{\FF}{\mathbb F}
\newcommand{\DD}{\mathbb D}
\newcommand{\CC}{\mathbb C}\newcommand{\PP}{\mathbb P}
\newcommand{\RR}{\mathbb R}\newcommand{\QQ}{\mathbb Q}

\newcommand{\ZZ}{\mathbb Z}

\newcommand{\Aut}{\rm{Aut}}\newcommand{\cod}{\rm cod}\newcommand{\codim}{\rm cod}

\newcommand{\GL}{\textstyle\mathop{\rm GL}}
\newcommand{\Ind}{{\textstyle\mathop{\rm Ind}}}

\newcommand{\SL}{\textstyle\mathop{\rm SL}}

\newcommand{\Hom}{\mathop{\rm Hom}}\newcommand{\hHom}{{\rm Hom}}
\newcommand{\RHom}{{\rm RHom}}
\newcommand{\End}{\mathop{\rm End}}

\newcommand{\Gal}{\mathop{\rm Gal}}
\newcommand{\Lie}{\mathop{\rm Lie}}

\newcommand{\im}{\rm{Image}}

\renewcommand{\AA}{\mathbb A}
\newcommand{\Ext}{{\rm Ext}}

\newcommand{\coker}{\mathop{\rm coker}}
\newcommand{\Spec}{\mathop{\rm Spec}}

\newcommand{\God}{{\rm God}}
\newcommand{\spt}{\rm{spt}}
\newcommand{\Gr}{\rm{Gr}}
\newcommand{\Tr}{\mathop{\rm Trace}}

\newcommand{\diag}{{\rm diag}}

\renewcommand{\b}{{^{\bullet}}}

\renewcommand{\Im}{\mathop{\sf{Im}}}
\newcommand{\rz}{\textcolor{red}{0}}

\newcommand*\rows{4}
\newcommand*\rowsh{2}
\usetikzlibrary{calc}

\newcommand{\highlight}[2][yellow]{\mathchoice%
  {\colorbox{#1}{$\displaystyle#2$}}%
  {\colorbox{#1}{$\textstyle#2$}}%
  {\colorbox{#1}{$\scriptstyle#2$}}%
  {\colorbox{#1}{$\scriptscriptstyle#2$}}}%

\newcommand{\z}{\text{\sffamily{\bf\textsf{z}}}}

\makeatletter
\newcommand{\uul}[1]{{\mathpalette\uul@{#1}}}
\newcommand{\uul@}[2]{%
  \begingroup
  \sbox\z@{$\m@th#1\underline{#2}$}%
  \dimen@=\dp\z@ \advance\dimen@ -2\uul@dimen{#1}%
  \dp\z@=\dimen@
  \sbox\z@{$\m@th\underline{\box\z@}$}%
  \box\z@
  \endgroup
}
\newcommand\uul@dimen[1]{%
  \fontdimen8
  \ifx#1\displaystyle\textfont\else
  \ifx#1\textstyle\textfont\else
  \ifx#1\scriptstyle\scriptfont\else
  \scriptscriptfont\fi\fi\fi 3
}
\makeatother

\begin{document}
\title{Lecture notes on sheaves and perverse sheaves}
\author{Mark Goresky}

\maketitle
\centerline{version of {\today}}
\setcounter{tocdepth}{1}
{\small{{\tableofcontents}}}

\onecolumn


\subsection*{Acknowledgements}  These notes are an expanded version of a course given by the author at Brown
University, while he was a guest at ICERM in the fall of 2016.  He thanks Brown University and ICERM for their
hospitality.  He thanks George Lusztig for help with some historical comments.
The author is very much indebted to Vidit Nanda for carefully listing many errors
and misprints in an earlier version of this manuscript.
\part{Sheaves}
\begin{section}{Sheaves:  the lightning tour}\end{section}
\subsection{}
Let $R$ be a commutative ring (with $1$). Let $X$ be a topological space. 
\begin{tcolorbox}[colback=yellow!30!white] A presheaf\index{presheaf} of $R$ modules on $X$ is a contravariant functor
\[ S: \text{category of open sets and inclusions} \to \text{category of $R$-modules}.\]\end{tcolorbox}
This is a fancy way to say that $V \subset U$ gives $S(U) \to S(V)$ and $W \subset V \subset U$
gives a commutative diagram
\begin{diagram}[size=2em]
S(U) &          &    \rTo &          & S(V) \\
        & \rdTo &           & \ldTo  &       \\
       &           &  S(W) &           &
\end{diagram}
Also, $U \subset U$ gives $S(U) \to S(U)$ identity and $S(\phi) = 0$.  Elements $s \in S(U)$ are called 
sections\index{section}\index{sheaf!section} 
of $S$ over $U$ (for reasons that will become clear shortly).
If $V \subset U$ and $s \in S(U)$, its image in $S(V)$ is denoted $s|V$ and is called the restriction of the section $s$ to $V$.
A morphism $S \to T$ of presheaves is a natural transformation of functors, that is, a collection of homomorphisms $S(U) \to T(U)$ (for every open set $U \subset X$) which commute with the restriction maps.  This defines a category of presheaves of $R$-modules, and it is an abelian category with kernels and cokernels defined in the obvious manner, for example, the kernel presheaf of a morphism $f:S \to T$ assigns to each open set $U\subset X$ the $R$-module 
$\ker(S(U) \to T(U))$.  

More generally, for any category $\mathcal C$ one may define, in a similar manner, the category of 
presheaves\index{presheaf!category of} on $X$ with coefficients in $\mathcal C$.  If $\mathcal C$ is abelian then so is the category of presheaves with coefficients in $\mathcal C$.  

A presheaf $S$ is a {\em sheaf} if the following {\em sheaf axiom} holds:  Let $\{U_{\alpha}\}_{\alpha \in I}$ be any collection of open subsets of $X$, let $U = \cup_{\alpha \in I}U_{\alpha}$ and let $s_{\alpha} \in S(U_{\alpha})$ be a collection of sections such that 
\[ s_{\alpha}|(U_{\alpha} \cap U_{\beta}) = s_{\beta}|(U_{\alpha}\cap U_{\beta})\]
for all $\alpha, \beta \in I$.  Then there exists a unique section $s \in S(U)$ such that
$s|U_{\alpha} = s_{\alpha}$ for all $\alpha \in I$.  In other words, if you have a bunch of sections over open sets that agree on the intersections of the open sets then they patch together in a unique way to give a section over the union of those open sets. [An essential point in this definition is that the index set $I$ may have infinite cardinality. ]

The {\em category of sheaves} is the {\em full subcategory} of the category of presheaves, whose objects are sheaves.  [This means that $\Hom_{\text{Sh}}(S,T) = \Hom_{\text{preSh}}(S,T)$.]

 \begin{tcolorbox}[colback=yellow!30!white]  The {\em stalk}\index{stalk}\index{sheaf!stalk}
  of a presheaf $S$ at a point $x \in X$ is the $R$-module
\[ S_x = \underset{\overrightarrow{U \ni x}}{\lim}S(U)\] \end{tcolorbox}
This means, in particular, that for any open set $U$ and for any $x\in U$ there is a canonical mapping $S(U) \to S_x$ which we also refer to as ``restriction" and denote by $s \mapsto s|S_x$.

\begin{tcolorbox}[colback=yellow!30!white] The {\em leaf space}\index{leaf space}\index{sheaf!leaf space}
 $LS$ of $S$ is the disjoint union
\[\begin{CD} LS = \coprod_{x \in X} S_x @>{\pi}>> X \end{CD}\]
with a topology that is discrete on each $S_x$ and that makes $\pi$ into a local homeomorphism.
\end{tcolorbox}\noindent
That is, each $U^{\text\small open} \subset X$ and each $s \in S(U)$ defines an open set
\[ U_s = \left\{ (x,t)| x\in U, t\in S_x \text{ and } t = s|S_x\right\}\subset \pi^{-1}(U)\subset LS.\]
Then $\pi:U_s \to U$ is a homeomorphism.

Let $\Gamma(U, LS)$ be the set of continuous
sections of $\pi$ over $U$, that is, the set of continuous mappings $h:U \to LS$ such that $\pi h = \text{identity}$.  The restriction maps of $S$ are compatible, giving $S(U) \to S_x$ for any $U \ni x$ and therefore any $s \in S(U)$ defines a continuous section $h \in \Gamma(U,LS)$ by setting $h(x) = s|S_x$.
\begin{tcolorbox}[colback=cyan!30!white]
\begin{prop}{\rm{(exercise)}}
The presheaf $S$ is a sheaf if and only if the canonical mapping $S(U) \to \Gamma(U,LS)$ is an isomorphism for every open set $U \subset X$.  If $S$ and $T$ are sheaves then there are canonical isomorphisms
\[ \hHom_{\text{Sh}}(S,T) \cong \hHom_X(LS,LT) \cong \hHom_{\text{preSh}}(S,T)\]
where $ \hHom_X(LS,LT)$ denotes the $R$-module of continuous mappings $LS \to LT$ that commute with the projection to $X$ and that consist of $R$-module homomorphisms $S_x \to T_x$ for all $x\in X$.
\end{prop}\end{tcolorbox}

\subsection{Sheafification}\label{subsec-sheafification}\index{sheafification}
An immediate consequence is that if $S$ is a presheaf then we obtain a sheaf $\widehat{S}$ by defining
\[ \widehat{S}(U) = \Gamma(U,LS)\]
to be the $R$-module of continuous sections of the leaf space of $S$ over the open set $U$.  
Then $\widehat{S}$ is called the {\em sheafification}\index{sheafification} of $S$.  The {\em category of sheaves} is the full
subcategory of the category of presheaves whose objects satisfy the above sheaf axiom, in other words,
\[ \hHom_{\text{Sh}}(A,B) = \hHom_{\text{PreSh}}(A,B).\]
To simplify notation, if $S$ is a sheaf we drop the notation $LS$ and we write 
$s \in S(U) = \Gamma(U,S)$ 
\begin{tcolorbox}[colback=cyan!30!white]\begin{prop}{\rm(exercise)}  Sheafification is an exact functor from the category of presheaves to the category of sheaves.  It is left adjoint to the inclusion functor $i:\text{Sheaves} \to \text{Presheaves}$, that is, if $A$ is a presheaf and if $B$ is a sheaf (on $X$) then
\[ \begin{CD}
\hHom_{\text{Sh}}(\widehat{A},B)  @>{\cong}>> \hHom_{\text{preSh}}(A, i(B)). \end{CD}\]
\end{prop}\end{tcolorbox}

Here is an application of this formula.  Following the identity morphisms $B \to B$ through this
series of canonical isomorphisms
\[ \hHom_{\text{Sh}}(B,B) \cong \hHom_{\text{preSh}}(i(B),i(B)) \cong \hHom_{\text{Sh}}
(\widehat{i(B)}, B)\]
gives a canonical isomorphism $\widehat{i(B)} \to B$, that is, if we take a sheaf $B$, look at it as a presheaf, then sheafifity it, the result is canonically isomorphic to the sheaf $B$ that we started with.

\subsection{Caution}  If $A, B$ are sheaves then the set of morphisms $A \to B$ is the same whether we consider $A, B$ to be sheaves or presheaves.  However, care must be taken when considering the kernel, image, or cokernel of such a morphism.  If we consider $f:A \to B$ to be a morphism in the category of presheaves, then $\ker(f)$ is the presheaf which assigns to an open set $U$ the kernel of $f(U):A(U) \to B(U)$, and this turns out to be a sheaf.  But the presheaf 
\[ U \mapsto \text{Image}(A(U) \to B(U))\]
is (usually) not a sheaf, so it is necessary to define the sheaf $\text{Image}(f)$ to be the sheafification of this presheaf.  Similarly for the cokernel.  Consequently, {\em a sheaf mapping $f:A \to B$ is injective (resp. surjective) iff the mapping $f_x:A_x \to B_x$ on stalks is injective (resp. surjective) for all $x\in X$}.   In other words, the abelian category structure on the category of sheaves is most easily understood in terms of the leaf space picture of sheaves.

\subsection{Examples}
\paragraph{(a.)}  Let $j:U \to S^2$ be the inclusion of the open complement $U$ of the north pole $i:N \to S^2$ and
let $j_!(\uul{\ZZ}{}_U)$ be the constant sheaf on $U$ extended by zero over the point $\{N\}$.  Let $i_*(\uul{\ZZ}{}_N)$ be
the skyscraper sheaf on $S^2$ supported at $\{N\}$.  Then $A=j_!(\uul{\ZZ}{}_U) \oplus i_*(\uul{\ZZ}{}_N)$ is a sheaf  on
$S^2$ whose stalk at every point is $\ZZ$ but there are no morphisms between $A$ and the constant sheaf $\uul{\ZZ}$
on $S^2$.\medskip

\paragraph{(b.)}  Let $C^0(0,1)$ be the presheaf that assigns to an open set $U \subset (0,1)$ the vector space of continuous functions $f:U \to \RR$.  Then this is a sheaf because a family of continuous functions defined on open sets that agree on the intersections of those sets clearly patch together to give a continuous function on the union.  Similarly, smooth functions, holomorphic functions, algebraic functions etc. can be naturally interpreted as sheaves.\medskip

\paragraph{(c.)}  Let $S$ be the sheaf on $(0,1)$ whose sections over an open set $U$ are those $C^{\infty}$ functions $f:U \to \RR$ such that $\int_U f(x)^2 dx < \infty$.  This presheaf is not a sheaf because it is possible to patch (infinitely many) $L^2$ functions (defined on smaller and smaller subintervals) together to obtain a function that grows too fast to be square integrable.  In fact, the sheafification of this sheaf is the set of all smooth functions on $(0,1)$.\medskip

\paragraph{(d.)} {\bf Local systems.} \index{local system}
 Let $X$ be a connected topological space with universal cover $\widetilde{X}$.
Let $x_0 \in X$ be a basepoint and let $\pi_1 = \pi_x(X,x_0)$ be the fundamental group of $X$.  This group acts freely on $\widetilde{X}$ (from the right) with quotient $X$.  Let $M$ be an $R$-module and let
$\rho:\pi_1 \to \Aut(M)$ be a homomorphism.  (For example, if $M$ is a vector space over the complex numbers then $\rho:\pi_1\to \GL_n(\CC)$ is a representation of $\pi_1$).  Define
\[\mathcal L = \widetilde{X} \times_{\pi_1} M\]
to be the quotient of $\widetilde{X} \times M$ by the equivalence relation $(yg, m) \sim (y, \rho(g)m)$
for all $y \in \widetilde{X}$, $m \in M$, and $g \in \pi_1$.  The projection $\widetilde X \to X$ passes to a projection $\pi:\mathcal L \to X$ which makes $\mathcal L$ into the leaf space of a sheaf, which is called a  {\em local system}, or {\em bundle of coefficients}.  Its stalk at $x_0$ is canonically isomorphic to $M$ and whose stalk other points $x \in X$ is isomorphic to $M$ but not in a canonical way.  If $U \subset X$ is a contractible subset then there exist trivializations
\[ \pi^{-1}(U) \cong U \times M\]
which identify the leaf space over $U$ with the constant sheaf.  So the sheaf $\mathcal L$ is a {\em locally constant sheaf} and every locally constant sheaf of $R$ modules on a connected space $X$ arises in this way.  If $X$ is a simplicial complex then a simplicial $r$-chain with values in the local system $\mathcal L$ is a finite formal sum $\sum a_i \sigma_i$ where $\sigma_i$ are (oriented) $r$-dimensional simplices and where $a_i \in \pi^{-1}(x_i)$ for some (and hence for any) choice of point $x_i \in \sigma_i$.  So there is a chain group (or module) $C_r(X;\mathcal L)$.  One checks that the boundary map $\partial_r:  C_r(X;\mathcal L) \to C_{r-1}(X; \mathcal L)$ continues to make sense in this setting and so it is possible to define the simplicial homology group $H_r(X; \mathcal L) = \ker(\partial_r)/\im(\partial_{r+1})$.  In other words, locally constant sheaves have homology. \medskip

\paragraph{(e.)} {\bf Flat vector bundles} \index{flat vector bundle}\index{vector bundle, flat}
Let $\pi:\mathcal E \to M$ be a smooth vector bundle (of $\RR$-vector spaces)
 on a smooth manifold.
The sheaf of smooth sections assigns to any open set $U \subset M$ the vector space of smooth
functions $\Gamma(U,\uul{C}^{\infty}(\mathcal E)) = \left\{ s:U \to \mathcal E| \ \pi\circ s(x) = x \text{ for all }\ x \in U \right\}$.  
Its leaf space is infinite dimensional.  Every smooth vector bundle admits a connection.  Suppose
$\mathcal E$ admits a {\em flat} connection $\nabla$.  
A flat section ($\nabla s = 0$) is determined locally by its value at a single point.  So the subsheaf of flat sections
of $\mathcal E$ form a local system which realize $ \mathcal E$ as its leaf space.\medskip

\paragraph{(f.)}  Let $R$ be a commutative ring (with $1$).  There is a topological space, $\Spec(R)$ which consists of all prime ideals in $R$.  The topology on this set was constructed by O. Zariski.  For any subset $E \subset R$ let $V(E) = \{ \mathfrak p| E \subset \mathfrak p\}$ be the set of prime ideals, each of which contains $E$.  These form the closed sets in a basis for a topology, that is, the open sets in this basis are the sets $X - V(E)$.  The topology generated by these open sets is called the Zariski topology.  If $M$ is an $R$-module then it defines a sheaf on this space in the following way. [to be completed]\medskip

\paragraph{(g.)}  Let $K \subset \RR$ be the Cantor set and let $\ZZ_K$ be the constant sheaf (with value equal to the integers, $\ZZ$) on $K$.  Let $j:K \to \RR$ be the inclusion.  Then $j_*(\ZZ_K)$ is a sheaf on $\RR$ (see the definition of $f_*$ below) that is supported on the Cantor set.  So ``bad" sheaves exist on ``good" spaces.\medskip

\paragraph{(h.)}  Fix $r \ge 0$.  For any topological space $Y$ let $C_r(Y;\ZZ)$ be the group of singular $r$-dimensional simplices on $Y$.  (It is the set of finite formal sums of pairs $(\sigma, f)$ where $\sigma$ is an oriented $r$-dimensional simplex and $f:\sigma \to Y$ is a continuous map.)  Now let $X$ be a topological space.  The presheaf of $r$-dimensional singular
\index{cochain, singular}\index{singular cochain} cochains $C^r$ on $X$ assigns to any open set $U \subset X$ the group $C^r(U):=\hHom(C_r(U;\ZZ), \ZZ)$.  If $V \subset U$ then $C_r(V;\ZZ)$ is included in $C_r(U;\ZZ)$ which gives a (surjective) restriction mapping $C^r(U) \to C^r(V)$.  This presheaf is also a sheaf.\medskip

\subsection{Sheaf Hom}\index{Hom!sheaves}\index{sheaf!Hom}
If $A$ is a sheaf on $X$ with leaf space $\pi:LA \to X$ and if $U \subset X$ is an open set let $A|U$ be the restriction of the sheaf $A$ to the subset $U$, in other words, the sheaf on $U$ whose leaf space is $\pi^{-1}(U) \to U$.  In other words, 
If $A, B$ are sheaves of $R$ modules then $\hHom(A,B)$ is again an $R$ module that consists of all sheaf mappings $A \to B$.  However there is an associated presheaf, perhaps we will denote it by ${\mathbf{Hom}}(A,B)$, which assigns to any open set the $R$ module of homomorphism
\[ \hHom_{Sh(U)}(A|U, B|U))\]
of sheaf mappings $A|U \to B|U$.   This presheaf is a sheaf (exercise) for which the group of global sections is
the original module of all sheaf homomorphisms, that is,
\[ \Gamma (X, \mathbf{Hom}(A,B)) = \hHom_{\text{Sh}(X)}(A,B)\]

\subsection{Functoriality}
Let $f:X \to Y$ be a continuous map, let $T$ be a sheaf on $Y$, let $S$ be a sheaf on $X$.  Define $f_*(S)$ to be the presheaf on $Y$ given by
\[ f_*(S)(U) = S(f^{-1}(U)).\]
This presheaf is a sheaf (exercise).  Define $f^*(T)$ to be the sheaf on $X$ whose leaf space is the pull back of the leaf space of $T$, that is,
\[ Lf^*(T) = f^*(LT) = X \times_Y LT = \left\{ (x,\xi) \in X \times LT|\ f(x) = \pi(\xi) \right\}. \]
Then this defines a sheaf, and  the sections of this sheaf are given by
\[ \Gamma(U, f^*T) = \underset{\overrightarrow{V \supset f(U)}}{\lim} \Gamma(V,T).\]
(Although $f(U)$ may fail to be open, we take a limit over open sets containing $f(U)$.)

There is also a pushforward with proper support, $f_!S$ with sections $\Gamma(U, f_!S)$ 
consisting of all sections $s \in \Gamma(f^{-1}(U,S))$ such that the mapping
\[ \text{closure}\{x\in U| s(x) \ne 0\} \to U\]
is proper (that is, the pre-image of every compact set is compact).  It is not so clear what this means, but if $f:X \to Y$ is the inclusion of a subspace then $f_!S(U)$ consists of sections $s \in
\Gamma(U \cap X, S)$ whose support is compact.  This implies, in particular, that $f_!(S)$ vanishes outside $X$, even if $X$ is open.  So, in the case of an inclusion, the sheaf $f_!(S)$ is called the
extension by zero of $S$.
\begin{tcolorbox}[colback=cyan!30!white]\begin{lem} Suppose the space $X$ is locally compact.  If $j:K \subset X$ denotes the inclusion of a closed subset, and if $S$ is a sheaf on $K$
then $j_*(S) \cong j_!(S)$. \end{lem}\end{tcolorbox}
\subsection{Adjunction}\index{adjoint}\index{sheaf!adjunction}
Let $f:X \to Y$ be a continuous mapping, let $A$ be a sheaf on $X$, let $B$ be a sheaf on $Y$.
Then there exist natural sheaf morphisms
\[ f^*f_*(A) \to A \ \text{ and }\ B \to f_*f^*(B).\]
To see this, for the first one, let us consider sections over an open set $U\subset X$.  Then
\[ \Gamma(U, f^*f_*A) = \underset{\overrightarrow{W \supset f(U)}} \lim \gamma(f^{-1}(W),A).\]
If $W \supset f(U)$ then $f^{-1}(W) \supset U$ so we get a mapping from this group to 
$\Gamma(U,A)$.  One verifies that these mappings are compatible when we shrink $U$, and
so this gives a sheaf morphism $f^*f_*(A) \to A$.  For the second morphism, again we look at sections over an open set $V \subset Y$.  If $t$ is a section of $LB$ over $V$ then, pulling it back by $f$
gives a section $f^*(t)$ of the leaf space of $f^*(B)$ over the set $f^{-1}(V)$, in other words, we have defined a map
\[ \Gamma(V,B) \to \Gamma(f^{-1}(V), f^*(B)) = \Gamma(V, f_*f^*(B))\]
which again is compatible with restrictions to smaller open sets.  In other words, this defines a
sheaf morphism $B \to f_*f^*B$.

\begin{tcolorbox}[colback=cyan!30!white]
\begin{prop}\label{prop-hom-adjoint}
The adjunction maps determine a canonical isomorphism
\[ \hHom_{\text{Sh}(X)}(f^*B,A) \cong \hHom_{\text{Sh}(Y)}(B, f_*A)).\]
\end{prop}\end{tcolorbox}
Given $f^*B \to A$, apply $f_*$ and adjunction to obtain $B \to f_*f^*B \to f_*A$.
Given $B \to f_*A$ apply $f^*$ and adjunction to obtain $f^*B \to f^*f_*A \to A$.
This gives maps back and forth between the $\hHom$ groups in the proposition.  We omit the check that they are inverses to each other. 
\bigskip

\begin{section}{Cohomology}\end{section}\index{cohomology}
\subsection{A comment about categories}
In a category\footnote{Saunders MacLane recalled that the first paper \cite{Eilenberg} 
on category theory  was initially rejected by an editor who stated that he had never before seen a paper that was so entirely
devoid of content.} $\mathcal C$ the collection of morphisms $A \to B$ between two objects is assumed to form a set, 
$\hHom_{\mathcal C}(A,B)$ and so we may speak of two morphisms being the same.  However the collection 
of objects do not (in general) form a set and so the statement that ``$B$ is the same object as $A$" does not make sense.  
Rather, ``the morphism $f:A \to B$ is an isomorphism" is the correct way to indicate an identification between two objects.  If 
the set of self-isomorphisms of $A$ is nontrivial, there will be many distinct isomorphisms between $A$ and $B$.

\subsection{Adjoint functors}\index{adjoint!functor}
If $\mathcal C, \mathcal D$ are categories, $f:\mathcal C \to \mathcal D$ and $g:\mathcal D \to \mathcal C$ are functors,
then $f$ and $g$ are said to be adjoint if there are ``natural" isomorphisms 
\[\hHom_{\mathcal C}(A, g(B)) \to \hHom_{\mathcal D}(f(A),B)\] 
for all $A, B$ objects of $\mathcal C, \mathcal D$ respectively (where ``natural" means that these isomorphism are compatible with all morphisms $A_1 \to A_2$ and $B_1 \to B_2$).  We say that $f$ is left adjoint to $g$ and $g$ is right adjoint to $f$.  
\index{left adjoint}\index{right adjoint}

In this case, $f$ essentially determines $g$ and vice versa, that is, given $f$ and two adjoint functors $g_1, g_2$ then there exists a natural equivalence of functors between $g_1, g_2$.  In many cases, the functor $f$ is something simple, and the functor $g$ is something surprising (or vice versa), so it is a fun game to pick your favorite functor and ask whether it has an adjoint.  

 If $H$ is a subgroup of a finite group $G$ then ``restriction" is a more or less obvious functor from the category of representations of $G$ to the category of representations of $H$; its adjoint turns out to be ``induction", a much more subtle construction, and the adjointness property is classically known as the Frobenius reciprocity theorem.    Sheafification (\S \ref{subsec-sheafification})
 is  adjoint to the  inclusion of categories (sheaves $\to$ presheaves).
If $f:X \to Y$ is a continuous mapping then Proposition \ref{prop-hom-adjoint} says that
$f^*$ is left adjoint to $f_*$.  One might ask whether there is an adjoint to the functor $f_!$.  This turns out to be a very subtle question, see \S \ref{def-shriek}.

\subsection{Simplicial sheaves}\label{subsec-simplicialsheaf}\index{simplicial!sheaf}\index{sheaf!simplicial}
This is a ``toy model" of sheaves.  Let $K$ be a (finite, for simplicity) simplicial complex.  Each (closed) simplex $\sigma$ is contained in a naturally defined open set, $St^o(\sigma)$, the open star of $\sigma$.  It has the property that $\sigma < \tau \implies
St^o(\tau) \subset St^o(\sigma)$.  Using these open sets to define a presheaf, and assigning the values of the presheaf to the simplex itself, gives the following definition:

\begin{tcolorbox}[colback=yellow!30!white]\begin{defn}
A simplicial sheaf $S$ (of abelian groups, or $R$-modules, etc.) on $K$ is an assignment of an abelian group $S(\sigma)$ for the interior of each simplex and a {\em restriction} homomorphism $S(\sigma) \to S(\tau)$ 
whenever $\sigma < \tau$, in such a way that whenever
$\sigma < \tau < \omega$ then the resulting triangle of groups and morphisms commutes.\end{defn}\end{tcolorbox}\noindent
 To make some explicit notation, let $j_{\sigma,\tau}:\sigma \to \tau$ denote the inclusion whenever $\sigma < \tau$ and let $S_{\sigma,\tau}:S(\sigma) \to S(\tau)$ be the corresponding homomorphism (so that $S$ becomes a covariant functor from the category of simplices and inclusions to the category of abelian groups).

In this setting there is no distinction between a sheaf and a presheaf.  A simplicial sheaf determines
an actual sheaf whose leaf space\index{leaf space!simplicial sheaf}
 of $S$ is easily constructed as the union
\[ LS = \coprod_{\sigma} \sigma^0 \times S(\sigma) \]
with a topology defined using $S_{\sigma,\tau}$ to glue these pieces together whenever $\sigma < \tau$.
This gives a natural functor from the category of simplicial sheaves to the category of sheaves on $K$ that
are locally constant on the interior of each simplex.

\subsection{Cohomology of simplicial sheaves}  \index{cohomology!simplicial sheaf}
Let $K$ be a finite simplicial complex and let $S$ be a simplicial sheaf
(of abelian groups, or of $R$-modules).  The cohomology $H^*(K,S)$ can be
constructed the same way as ordinary simplicial cohomology $H^*(K)$.  First choose orientations of the simplices.  
(An orientation of a simplex is determined by an ordering of its vertices, two orderings giving the same orientation 
iff they differ by an even permutation.)  The simplest method of orienting all the simplices,  is to choose an 
ordering of the vertices of $K$ and to take the induced ordering on the vertices of each simplex.  Assume this to be done. 

Fix $r \ge 0$.  An $r$-chain with values in the simplicial sheaf $S$ is a function $F$ that assigns, to each (oriented) 
$r$-dimensional simplex $\sigma$ an element of $S(\sigma)$.  The collection of all $r$-chains is denoted $C^r(K;S)$.   
The coboundary $\delta F \in C^{r+1}(K;S)$ is defined as follows.  If $\tau$ is an $r+1$ simplex with 
vertices $v_0,v_1,\cdots, v_{r+1}$ (in ascending order) we write $\tau = \langle v_0,v_1,\cdots, v_{r+1}\rangle$ 
and we denote its $i$-th face by $\partial_i\tau = \langle v_0,\cdots, \widehat{v_i}, \cdots, v_{r+1} \rangle$.  Then
\[ ( \delta F)(\tau) = \sum_{\sigma<\tau} [\tau:\sigma] S_{\sigma,\tau}(F(\sigma))\]
where the sum is over codiimension one faces $\sigma < \tau$.
%
One checks that $\delta \delta F = 0$ (it is the same calculation that is involved in proving that $\partial \partial = 0$
for simplicial homology), so we may define the cohomology, $H^r(K;S) = \ker(\delta)/\Im(\delta)$ to be the cohomology
of the cochain complex
\[ \cdots {\longrightarrow}\ C^{r-1}(K;S) \overset{\delta}{\longrightarrow} C^r(K;S)
\overset{\delta}{\longrightarrow}C^{r+1}(K;S) \ {\longrightarrow} \cdots \]
This combinatorial construction of $H^*(K;S)$  is easily implemented on a computer.  The identification with other
constructions of sheaf cohomology is explained in \S \ref{subsec-simplicialcover} below.

 By reversing the arrows one has the 
analogous notion of a simplicial cosheaf \index{cosheaf} and a similar construction of the {\em homology} of a simplicial cosheaf.

\subsection{Historical interlude}
Many different techniques have been developed for exploring the properties of cohomology of sheaves, the most elegant being the methods associated with the derived category.  Any one of these methods may be used as a ``definition" of the cohomology of a sheaf, and although the historical methods are the most accessible, they are also the most cumbersome.  
We will use the method of injective resolutions.  First we mention a few milestones in the development of
 cohomology and sheaves.

\begin{enumerate}
\item[(1873)] B. Riemann and later, E. Betti, consider the number of ``cuts" of varying dimensions that are needed in order to reduce a space into contractible pieces.
\item[(1892)] H. Poincar\'e, in {\em Analysis Situs} constructs homology of a ``vari\'et\'e" using cycles that consist of the zeroes of smooth functions.
\item[(1898)] P. Heegaard publishes a scathing scriticism of Poincare\'s article for its lack of rigor.
\item[(1900)] H. Poincar\'e publishes {\em Supplement} to his {\em Analysis Situs} in which he essentially describes simplicial homology for a space that has been decomposed into simplices.
\item[(1912)] F. Hausdorff publishes the general definition of a topological space and interprets continuity purely in terms of the open sets.
\item[(1925)] H. Hopf develops the general notion of a chain complex.
\item[(1926)] Alexander, Hopf give precise definition of simplicial complex.
\item[(1928)] H. Hopf, E. Noether describe homology as a group.
\item[(1930)] E. Cartan, G. deRham formalize notion of differential forms, Poincar\'e lemma, de Rham theorem.
\item[(1933)] The drive to develop singular homology theory, with contributions by Dehn, Heegard, Lefschetz, others.
\item[(1934)] E. \v{C}ech develops his approach to cohomology using the open sets in a space; cohomology with coefficients in a ring.
\item[(1935)] H. Whitney develops the abstract theory of differentiable manifolds and their embeddings into Euclidean space.
\item [(1935)] H. Reidemeister develops theory of homology with local coefficients.
\item[(1935-40)] Products in cohomology, modern formulation of Poincar\'e duality, Stiefel Whitney classes, differential forms.  (Until this period, differential forms were ``expressions".)
\item[(1942)] S. Eilenberg and S. MacLane:  Category theory
\item[(1945-46)] J. Leray:  sheaves and their cohomology, spectral sequence of a map
\item[(1946)] S. S. Chern:  Chern classes
\item[(1950)] A. Borel's notes on Leray's theory, vastly increasing its accessibiliy
\item[(1950)] \v{C}ech cohomology of sheaves
\item[(1956)] A. Borel, J. C. Moore, the dual of a complex of sheaves; Borel-Moore homology
\item[(1956)] H. Cartan, S. Eilenberg:  injective resolutions and derived functors
\item[(1957)] A. Grothendieck:  Tohoku paper on homological algebra
\item[(1958)] D. Kan:  notion of adjoint functors
\item[(1961)] J. L. Verdier, derived categories, Verdier duality (published in 1996).

\end{enumerate}

Given the complexity of this history, we will describe several ways to define the cohomology of a sheaf, and leave the proof that they all give the same answer until Proposition \ref{prop-spectralsequence}.  Further details on sheaves, the derived
category and perverse sheaves may be found in any of
the wonderful textbooks on these subjects, including \cite{Dimca, Gelfand, Gelfand2, Iverson, KS, Kiehl, BorelIH, Maxim}.

\subsection {}  
Let $X$ be a topological space and let $S$ be a sheaf (of abelian groups, or of $R$-modules) on $X$.  Let $\sigma \in \Gamma(X,S)$, and consider $\sigma$ to be a section of the leaf space $LS \to X$.  
The {\em support} $spt(\sigma)$ \index{support}
of $\sigma$ is the closure of the set of points $x\in X$ such that $\sigma(x) \ne 0$.  Let $\Gamma_c(X,S)$ denote the group of sections with compact support.  If $K \subset X$ is a closed subset, let $\Gamma_K(S)$ denote the group of sections whose support is contained in $K$.  This may also be identified with the sections $\Gamma(K, j^*(S))$ where $j:K \to X$ denotes the inclusion.

The sheaf $S$ is {\em injective}\index{injective}\index{sheaf!injective}
 if the following holds:  Suppose $f:A\to B$ is an injective morphism of sheaves.  Then any
morphism $h:A \to S$ extends to a morphism $\tilde h:B \to S$.
\begin{diagram}[size=2em]
B &&&\\
\uInto^{f} & \rdDotsto^{\tilde{h}} &\\
A & \rTo_h & S
\end{diagram}
The sheaf $S$ is {\em flabby}
\index{flabby (sheaf)}\index{sheaf!flabby,fine,soft} if $S(U) \to S(V)$ is surjective, for all open subsets $V\subset U$.
The sheaf $S$ is {\em soft (sheaf)} if $\Gamma(X,S) \to \Gamma_K(S)$ is surjective, 
for every closed subset $K \subset X$.
The sheaf $S$ is {\em fine (sheaf)}\index{fine (sheaf)} if, for every 
open cover $X = \cup_{\alpha \in I}U_{\alpha}$ there exists a family of morphisms
$h_i:S \to S$ such that $\spt(h_{\alpha}) \subset U_{\alpha}$ and $\sum_{\alpha}h_{\alpha} = I$.  
(Usually the $h_{\alpha}$ are just a partition of unity with respect to the coefficient ring.)  For completeness we include here the following definition, which  actually requires having previously defined cohomology:
The sheaf $S$ is {\em acyclic} if $H^r(X,S) = 0$ for all $r \ge 1$.

The following fact will not be used 
\begin{tcolorbox}[colback=cyan!30!white]\begin{prop}  For any sheaf $S$ on a locally compact space $X$,
\[ \text{\rm injective }\implies \text{\rm flabby } \implies \text{\rm soft } \implies \text{\rm acyclic }
\text{ and }\ \text{\rm fine }\implies\text{ soft }\implies\text{ acyclic}.\]
\end{prop}\end{tcolorbox}

Of these notions, {\em injective} and {\em acyclic} are categorical, and we will concentrate on them.  However, in order
that a sheaf be injective, it must have certain topological properties and certain algebraic properties.  For example, if the ring $R$ is an integral domain, then the constant sheaf on a single point is injective iff $R$ is a field.  The ring $\ZZ$ is not injective but it has an injective resolution $\ZZ \to \QQ \to \QQ/\ZZ \to 0.$

\begin{tcolorbox}[colback=yellow!30!white]\begin{defn}
An injective resolution\index{resolution!injective}\index{injective!resolution} of a sheaf $S$ on $X$ is an exact sequence
\[ 0 \to S \to I^0 \to I^1 \to \cdots\]
where each $I^r$ is an injective sheaf. \end{defn}\end{tcolorbox}
A abelian category $\mathcal C$ has {\em enough injectives} if every object can be embedded in an injective.  In this case, every object $A$ admits an injective resolution:  just embed $A \to  I^0$ and let $K^0$ be the cokernel of this map.
Then embed $K^0 \to I^1$ and let $K^1$ be the cokenerl of this map.  Then embed $K^1 \to I^2$ and so on.  The resulting sequence $0 \to A \to I^0 \to I^1 \cdots$ is exact.  The category of modules over a commutative ring $R$ has enough injectives, and the category of sheaves of $R$ modules on any topological space $X$ has enough injectives.  
However, there is a canonical and functorial injective resolution of any sheaf, namely
the {\em Godement resolution}.  It will be described later.

{\bf First definition of cohomology}\index{cohomology!of a sheaf}  Let $S$ be a sheaf on a topological space $X$ and let 
$0 \to S \to I^0 \to \cdots$ be the Godement injective resolution of $S$.
Then the cohomology of the complex of global sections
\[ 0 \to \Gamma (X, I^0) \to \Gamma(X, I^1) \to \Gamma(X,I^2) \cdots\]
is called the cohomology of $S$, denoted $H^r(X,S)$.

{\bf Second definition of cohomology}
In fact:  You will get the same answer, up to unique isomorphism, if you use any injective resolution, or any  fine, flabby, soft, or acyclic resolution instead of an injective resolution.  (See \S 4:  injectives and \v{C}ech cohomology)

\bigskip

\section{Complexes of sheaves}\index{complex of sheaves}
\subsection{}  Sheaves tend to occur in complexes.  Familiar examples include
the de Rham complex of differential forms on a smooth manifold $M$, 
\[ \Omega^0(M,\RR) \to \Omega^1(M,\RR) \to \Omega^2(M,\RR) \cdots \]
or the $\bar\partial$ complex on a complex manifold; the sheaf of singular cochains, etc.  Less obvious
examples include the push-forward $f_*(S^{\bullet})$ of a complex by a (continuous, smooth, or algebraic)
mapping $f:X \to Y$.
\subsection{}  A {\em complex}  
 $S^{\bullet}$ (in an abelian category) is a sequence \[ \begin{CD}
\cdots @>>> S^{r-1} @>d^{r-1}>> S^r @>d^r>> S^{r+1} @>d^{r+1}>> \cdots
\end{CD}\]
where $d\circ d = 0$. (We will assume that all of our complexes are bounded from below, that is $S^r = 0$ if
$r$ is sufficiently small, usually if $r<0$).  The {\em cohomology}\index{cohomology!complex of sheaves}
 $H^r(S^{\bullet})$ of the complex is $\ker d^r/ \Im d^{r-1}$.
So the complex $S^{\bullet}$ is {\em exact} (meaning that it forms an exact sequence) iff $H^r(S^{\bullet}) = 0$ for all $r$.
A {\em complex of sheaves} is a complex in which each $S^r$ is a sheaf.  In this case, the {\em cohomology sheaf}
$\mathbf H^r(S^{\bullet})$ is the sheaf  $\ker d/\Im d$. The stalk of the cohomology sheaf coincides with the cohomology
of the stalks (exercise), that is,
\[ \mathbf{H}^r_x(S^{\bullet}) := \mathbf H^r(S^{\bullet})_x \cong H^r(S^{\bullet}_x).\]

A morphism (of complexes in an abelian category)
$S^{\bullet} \to T^{\bullet}$ of complexes is a collection of morphisms that commute with the
differentials.  Two morphisms $f,g:S^{\bullet} \to T^{\bullet}$ are {\em chain homotopic} if there is a
collection of mappings $h:S^{r} \to T^{r-1}$ (for all $r$) so that $d_Th + hd_S = f-g$.
A morphism $\phi:S^{\bullet} \to T^{\bullet}$ is a {\em quasi-isomorphism}\index{quasi-isomorphism} if it induces
isomorphisms on the cohomology objects $H^r(S^{\bullet}) \to H^r(T^{\bullet})$ for all $r$.  
A quasi-isomorphism of complexes of sheaves induces isomorphisms on cohomology
 (see Lemma \ref{prop-spectralsequence}).
 \medskip
\quash{
\paragraph{\bf Example}  For the complex $\mathbf{\Omega}^{\bullet}_M$ of sheaves of differential forms 
on a smooth manifold $M$,  the cohomology sheaves are zero in all degrees except zero,
and  $\mathbf H^0(\mathbf{\Omega}^{\bullet}_M) \cong \uul{\RR}$ is the constant sheaf so
$\uul{\RR} \to \mathbf{\Omega}^{\bullet}_M$ is a quasi-isomorphism.}
\quash{
\subsection{} A {\em double complex} $S^{\bullet \bullet}$ is an array $S^{i,j}$ with horizontal differentials
$d_h:S^{i,j} \to S^{i+1,j}$ and vertical differentials $d_v:S^{i,j} \to S^{i,j+1}$ such that $d_h^2=0$,
$ d_v^2=0$ and $d_vd_h=d_hd_v$.  The {\em associated single complex} is $T\b$ with
$T^m = \oplus_{a+b=m}S^{a,b}$ and $d = \oplus_{a+b=m}d_h + (-1)^bd_v$. 
}
\subsection{Magic Triangles}\index{triangle, distinguished}\index{distinguished triangle}
The {\em mapping cone} $C\b=C(\phi)$ of a morphism $\phi:A\b \to B\b $ is the complex $C^r = A^{r+1} \oplus B^r$
with differential $d_C(a,b) = (d_A(a), (-1)^{\deg(a)}\phi(a) + d_B(b))$.  It is the total complex  of the double complex
\begin{diagram}[size=2em]
{} & & {} \\
\uTo^d && \uTo_d\\
A^2 & \rTo^{\phi} & B^2 \\
\uTo^d &&    \uTo _d\\
A^1 & \rTo^{\phi} & B^1 \\
\uTo^d && \uTo _d\\
A^0 & \rTo^{\phi} & B^0
\end{diagram}
from which we see that there are obvious morphisms $\beta:B^{\bullet} \to C^{\bullet}$ and $\gamma:C^{\bullet} \to
A^{\bullet}[1]$.  It is customary to denote this situation as a triangle of morphisms
\begin{diagram}[size=2em]
  A\b && \rTo^{\phi} && B\b \\
  & \luTo_{[1]} && \ldTo \\
  && C(\phi) &&
  \end{diagram}
 \begin{tcolorbox}[colback=cyan!30!white]

 \begin{lem} \label{lem-triangles}
  If $\phi$ is injective then there is a natural quasi-isomorphism $\coker(\phi) \cong C(\phi)$.
  If $\phi$ is surjective then there is a natural quasi-isomorphism $C(\phi) \cong\ker(\phi)[1]$.  There are
  natural quasi-isomorphisms $A^{\bullet}[1] \cong C(\beta)$ and $B^{\bullet}[1] \cong C(\gamma)$.
  Moreover,
  there is a long exact sequence on cohomology 
  \[
  \cdots\to \mathbf{H}^{r-1}(B\b) \to \mathbf{H}^{r-1}(C\b) \to 
  \mathbf{H}^r(A\b) \to \mathbf{H}^r(B\b) \to \mathbf{H}^r(C\b) \to \cdots.\]
  \end{lem}\end{tcolorbox}
The proof is an exercise, and it works in any abelian category.

\subsection{Double complexes}\label{subsec-doublecomplex}\index{double complex}\index{complex, double}
 A  {\em double complex} is an array $C^{pq}$ with horizontal and vertical
differentials $d_h:C^{pq} \to C^{p+1, q}$ and $d_v:C^{pq} \to C^{p, q+1}$ such that $d_v d_v = 0$, 
$d_h d_h = 0$, and $d_vd_h = d_hd_v$. [Some authors assume instead that $d_vd_h = -d_hd_v$ which eliminates the necessity for a sign in the definition below of the total complex.] For convenience we will always assume that $C^{pq}=0$ unless $p\ge 0$ and $q \ge 0$.
The {\em associated single complex}, or {\em total complex} $T^{\bullet}$ is defined by $T^r = \oplus_{p+q = r} C^{pq}$
with ``total differential" $d_T:T^r \to T^{r+1}$ defined by $d_T(c_{pq}) = (d_h + (-1)^qd_v)c_{pq}$ for $c_{pq} \in A^{pq}$
(that is, change $1/4$ of the signs to obtain that $d_{T}d_T = 0$).

\begin{lem}\label{lem-double-one}
Let $C^{\bullet\bullet}$ be a first quadrant double complex and suppose the rows are exact and that the zeroth
horizontal arrows $d_h^{0q}$ are injections (which is the same as saying that an extra zero may be added
to the left end of each row, without destroying the exactness of the rows).  Then the total complex $T$ is
exact.
\end{lem}  

\begin{proof}
A proper proof involves indices, signs, and an induction that is totally confusing and is best worked out 
in the privacy of your own home.  To see how the argument goes, let us show that the
total complex $T^0 \to T^1 \to T^2 \cdots $ is exact at $T^2$.  Let $x=x_{02} + x_{11} + x_{21} \in T^2$
and suppose that $d_Tx=0$.
\begin{diagram}[size=1.5em]
\cdot \\
\uTo\\
x_{02} & \rTo & \cdot \\
\uTo && \uTo  \\
y_{01} & \rTo& x_{11} & \rTo &\cdot \\
\uTo && \uTo & &\uTo  \\
\cdot & \rTo & y_{10} &\rTo & x_{20} & \rTo & \cdot \\
\end{diagram}
\begin{itemize}
\item Since $d_Tx=0$ we have:
\[ d_v x_{02} = 0;\  d_v x_{11} + d_h x_{02} = 0;\  d_vx_{20} - d_hx_{11} = 0;\  d_h x_{20} = 0.\]
\item Since the bottom row is exact there exists $y_{10}$ so that $d_h y_{10} = x_{20}$.
\item Now consider $x'_{11}=x_{11}-d_vy_{10} $.  Check that $d_h(x'_{11})=0$.
\item Since the first row is exact there exists $y_{01} $so that $d_h y_{01} = x'_{11}$.
\item Now consider $x'_{02} = x_{02} + d_v y_{01}$.  Check that $d_h x'_{02} =0.$
\item But now we are in the left column, so this last operation, $d_h$ was an injective mapping.
This implies that $x'_{02} = 0.$
\item Now just check that we have the right answer, $y=y_{01} + y_{10}\in T^1$.\qedhere
\end{itemize}\end{proof}

\begin{tcolorbox}[colback=cyan!30!white]
\begin{cor}\label{cor-double-one}
Let $\{C^{pq}\}$ be a first quadrant double complex ($p,q \ge 0$) with exact rows.  Let $T^{\bullet}$ denote
the total complex, $T^r = \oplus_{p+q=r}C^{pq}$.
 Let $A^r = \ker(d^{0,r}:C^{0r} \to C^{1r})$ denote the subcomplex of the zeroth column, with its vertical
 differential $d_v$.  Then the
 morphism $A^{\bullet} \to T^{\bullet}$ (given by the inclusion $A^r \to C^{0r} \to T^r$) is a quasi-isomorphism.
 \end{cor}\end{tcolorbox}
 \begin{proof}
 Let us consider the extended double complex obtained by considering the complex $A^{\bullet}$  to be the
 $-1$-st column, that is, $C^{-1,r} = A^r$, with vertical differential $d_v = d_A$ and with horizontal differential
 the inclusion $A^r \to C^{0,r}$.
 \begin{diagram}[size=2em]
 A^2             & \rInto  & C^{02} & \rTo & C^{12} & \rTo & C^{22} & \rTo \cdots \\
 \uTo_{d_A} &         &  \uTo    &        &  \uTo    &        &  \uTo    &             \\
 A^1            & \rInto  & C^{01} & \rTo & C^{11} & \rTo & C^{21} & \rTo \cdots \\
 \uTo_{d_A} &         &  \uTo    &        &  \uTo    &        &  \uTo    &             \\
 A^0             & \rInto  & C^{00} & \rTo & C^{10} & \rTo & C^{20} & \rTo \cdots \\
 \end{diagram}
 Let $S^{\bullet}$ denote the total complex of this extended double complex.  It is precisely the mapping cone
 of the morphism $A^{\bullet} \to T^{\bullet}$ so we have a magic triangle
 \begin{diagram}[size=2em]
 A^{\bullet} & &\rInto & &T^{\bullet}\\
   &\luTo_{[1]}&    &\ldTo&\\
   &&S^{\bullet}&&
 \end{diagram}
 Moreover, the extended double complex has exact rows and the leftmost horizon maps are injective so the previous lemma applies and we conclude that the cohomology of $S^{\bullet}$ vanishes.  By the long exact sequence on cohomology this implies that $H^r(A^{\bullet}) \to H^r(T^{\bullet})$ is an isomorphism for all $r$.  \end{proof}

\subsection{Cohomology of a complex of sheaves}\index{cohomology!complex of sheaves}
We defined the cohomology of a sheaf to be the (global section) cohomology of an injective resolution of the sheaf.
So one might expect an injective resolution of a complex of sheaves $A^{\bullet}$ to be a double complex, and in fact, such a double complex of injective sheves can always be constructed:  If the coefficient ring is
a field, then the simplest way is to use the Godement resolution because it is functorial (see the next section).  For more general coefficient rings, the Godement resolution can be tensored with an injective resolution of the ring.  More generally, a double complex resolution was constructed by Cartan and Eilenberg, and it is known as a Cartan-Eilenberg resolution of the complex.  It is easy to
take an injective resolution of each of the sheaves in the complex, but not so easy to see how to fit them together
so that the ``vertical'' differentials satisfy $d_v^2=0$, but this can be done (see the Stacks Project, or Wikipedia on  Cartan-Eilenberg resolutions).  In any case, let us assume that we have a double complex,
the $r$-th row of which resolves $A^r$. 

Given an injective resolution as a double complex $A^{\bullet} \to I^{\bullet \bullet}$ with horizontal and vertical differentials $d_h, d_v$ respecitvely (and $d_vd_h=d_hd_v$),  form the
associated total complex $T^r = \oplus_{p+q = r}I^{pq}$ and with $d_r = d_h + (-1)^q d_v$.
\begin{lem} The resulting map $A^{\bullet} \to T^{\bullet}$ is a quasi-isomorphism.\end{lem}
This follows immediately from Corollary \ref{cor-double-one}, which says that there is a long exact sequence of
cohomology sheaves,
\[\cdots \to \mathbf H^{r}(A^{\bullet}) \to \mathbf H^r(T^{\bullet}) \to \mathbf H^r(S^{\bullet}) \to \mathbf H^{r+1}(A\b) \to\cdots\]
where the cohomology sheaves of $S^{\bullet}$ vanish.
 
In summary, we have replaced the double complex $I^{pq}$ with a single complex $T^{\bullet}$ so we arrive at the following definition, which works in any abelian category:

\begin{tcolorbox}[colback=yellow!30!white]
\begin{defn} 
\label{def-hypercohomology}\index{hypercohomology}
An injective resolution\index{injective!resolution!of a complex} of a complex  $A^{\bullet}$ is a quasi-isomorphism
$A^{\bullet} \to T^{\bullet}$ where each $T^r$ is an injective object. The cohomology $H^r(X, A^{\bullet})$ 
(also known as the hypercohomology of $A^{\bullet}$)
of a complex of sheaves is defined to be the cohomology of the complex of global sections
\[ \Gamma (X, T^{r-1}) \to \Gamma(X, T^r) \to \Gamma(X, T^{r+1})\]
of any injective resolution $A^{\bullet} \to T^{\bullet}$
\end{defn}\end{tcolorbox}
[A messy technical point:  a complex of injective objects is not necessarily an injective object in
the category of complexes.  So the somewhat ambiguous terminology of ``injective resolution"
could be misleading and some authors refer to these as ``K-injective resolutions". Fortunately
we will not be required to consider injective objects in the category of complexes.]

This fits beautifully with the notion of a resolution of a single sheaf $S$:  it is a quasi-isomorphism,
\[
\begin{CD}
0 @>>> S @>>> 0 @>>>0  @>>>\dots\\
@. @VVV @. @.\\
0 @>>> I^0 @>>> I^1 @>>> I^2 @>>> \cdots
\end{CD}\]
As before, flabby, soft or fine resolutions may be used instead of injective resolutions.
For example, let $M^n$ be a smooth manifold.  Let $x\in M$ and let $U_x$ be a neighborhood of $x$ that is
diffeomorphic to an $n$-dimensional ball.  The Poincar\'e lemma says that if $\xi$ is a closed (i.e. $d\xi = 0$)
differential $r$-form ($r \ge 1$) defined in $U_x$ then there is a differential $r-1$ form $\eta$ so that $d\eta = \xi$.
The sheaf of smooth differential forms is fine, so the Poincar\'e lemma says that this complex of sheaves is a fine
resolution of the constant sheaf,
\[
\begin{CD}
0 @>>> \RR @>>> 0 @.  @.\\
@. @VVV @. @.\\
0 @>>> \Omega_M^0 @>>> \Omega_M^1 @>>> \Omega_M^2 @>>> \cdots
\end{CD}\]
Therefore the cohomology $H^r(M,\RR)$ is canonically isomorphic to the cohomology of the complex of global sections of
$\Omega^*_M$, that is, the de Rham cohomology.

In fact, injective resolutions may be constructed directly without resorting to the Godement
double complex or the Cartan-Eilenberg double complex.

\begin{tcolorbox}[colback=cyan!30!white]
\begin{lem}\label{prop-spectralsequence}
  Let $A^{\bullet} \to B^{\bullet}$ be a quasi-isomorphism of complexes of sheaves on a 
topological space $X$.  Then it
induces an isomorphism on cohomology $H^r(U,A^{\bullet}) \cong H^r(U, B^{\bullet})$ for any open set $U \subseteq X$.\end{lem}\end{tcolorbox}
The morphism $\phi:A\b \to B\b$ induces a mapping on the spectral sequence for cohomology whose $E^2$ page
is
\[ H^p(U;\uul{H}^q(A\b)) \implies H^{p+q}(U;A\b).\]
Since $\phi$ is a quasi-isomorphism it induces an isomorphism on th $E^2$ page and on all other pages, so by the
spectral sequence comparison theorem it induces isomorphism $H^{p+q}(U;A\b) \cong H^{p+q}(U;B\b)$.
\bigskip

\subsection{Examples}
\paragraph{1.}  
Let $X$ be the 2-dimensional simplex.  In the category of simplicial sheaves, suppose that $S$ is a sheaf on $X$ that
assigns the value $\QQ$ to the interior of the $2$-simplex and assigns $0$ to simplices on the boundary.  Find an injective resolution of $S$.  Determine the global sections of each step in the resolution.  Show that the cohomology of $S$ is
$\QQ$ in degree $2$ and is $0$ in all other degrees.
\medskip
\paragraph{2.}  Let $X$ be a triangulation of $S^2$, which may be taken, for example to be the boundary of a
$3$-simplex.  Let $S$ be the constant sheaf on $X$.  Find an injective resolution of $S$ in the category of
simplicial sheaves, and compute the cohomology of its global sections.
\medskip
\paragraph{3.}  Let $X$ be a triangulation of $S^2$ with $10^8$ simplices.  Describe an injective resolution of $X$.
\medskip
\paragraph{4}  Let $X$ be a topological space, let $x_0\in X$ and let $S=S(x_0,\QQ)$ be the presheaf that 
assigns to any open set $U$
\[ S(U) = \begin{cases} \QQ &\text{ if } x_0\in U \\
                                   0 & \text{ else}\end{cases}\]                                
 Show that $S$ is injective and that its leaf space $LS$ is a skyscraper, that is, it consists of a single group $\QQ$
 at the point $x_0$ and zero everywhere else.
\medskip
 \paragraph{5}  In the above example, fix $x_0 \in X$ and let $T^{\bullet}$ be the complex of sheaves $S(x_0,\QQ) \to
 S(x_0, \QQ/\ZZ)$.  Show that this complex is an injective resolution of the skyscraper sheaf that is $\ZZ$ at the point $x_0$.

\begin{section}{Godement and \v{C}ech}\end{section}\index{resolution!Godement}\index{Godement resolution}
These examples show that injective sheaves must be sums of sheaves with tiny support.  This leads one to the following:
\subsection{Godement resolution}\index{resolution!Godement}\index{Godement resolution}
Given a sheaf $A$ on a topological space $X$ it embeds in a flabby sheaf $\God(A)$ with sections
\[ \Gamma(U,\God(A)) = \prod_{x\in U} A_x\]
the product of all the stalks at points in $U$.  It is sometimes called the sheaf of totally discontinuous sections. If we start with the constant sheaf $\uul{\ZZ}$ then a section  $s\in\Gamma(U, \God(\uul{\ZZ}))$ assigns to each point $x\in U$ an integer, without any regard to continuity or compatibility.  It is the sort of sheaf that you definitely do not want to meet in a dark alley.
The Godement
{\em resolution} $\uul{\God}\b(A)$ is obtained by applying this construction to the cokernel of $A \to \God(A)$ and iterating:
 \begin{diagram}[size=2em]
A &\rInto& \God(A) &         & \rTo &        & \God(\coker{}) &        & \rTo    &         & \God(\coker{}) \\
  &       &               &\rdTo&         &\ruInto&                      &\rdTo&           &\ruInto& \\
  &       &               &        &\coker{}&       &                       &        &\coker{}&         & \\
 \end{diagram}
 
 If $A_x$ is
injective for all $x\in X$ (for example, if the coefficient ring $R$ is a field) then this is an injective resolution of $A$.  For many rings there are functorial injective resolutions that can be used, together with the Godement construction to make a double complex, the associated total complex of which is then a canonical injective resolution. 
The Godement resolution is functorial:  a morphism $f:A \to B$ induces a morphism of
complexes $\God(f):\uul{\God}\b(A) \to \uul{\God}\b(B)$ in such a way that  
${\God}(f\circ g) ={\God}(f)\circ\God(g)$.

In summary, injective sheaves are huge, horrible objects and maybe we use them to prove things but never to compute with.
A much more efficient computational tool is the \v{C}ech cohomology.

\subsection{\v{C}ech cohomology of sheaves}\index{Cech cohoology@\v{C}ech cohomology}\index{cohomology!\v{C}ech}
  Let $A$ be a sheaf on $X$.
Let $\mathcal U = \{U_{\alpha}\}_{\alpha \in I}$ be a collection of open sets that
cover $X$. Fix a well ordering on the index set $I$ and for any  finite subset $J \subset I$,
$J = \{ i_0 < i_1 < \cdots < i_q\}$ let
$U_J = \cap_{\alpha \in J}U_{\alpha}$ be the corresponding intersection.  The \v{Cech} cochain complex is
\begin{diagram}
{\check{C}}^r(X,A) &=&  \prod_{|J|=r+1}\Gamma(U_J,A) \\
\dTo_d && \dTo_d \\
{\check{C}}^{r+1}(X,A) &=& \prod_{|K|=r+2}\Gamma(U_K,A)
\end{diagram}
where $d$ is defined as follows:  suppose that $K = \{k_0, k_1, \cdots, k_{r+1}\} \subset I$ and $\sigma \in \check{C}^r(X,A)$.  Then
\[
d\sigma(U_K )= \sum_{i=0}^{r+1} (-1)^i
\sigma(U_{k_0} \cap \cdots \cap \widehat{U}_{k_i} \cap \cdots \cap U_{k_{r+1}})|U_K.\]
For example, if $\sigma \in \check{C}^1$ and if $U = U_0 \cap U_1 \cap U_2$ then
\[ d\sigma(U) = \sigma(U_1\cap U_2)|U - \sigma(U_0 \cap U_1)|U + \sigma(U_0 \cap U_2)|U.\]
  Then one checks that
$d\circ d = 0$ so the cohomology of this resulting complex is defined:
\[ \check{H}^r_{\mathcal U}(X,A) = \ker d/ \Im d.\]
Notice, in particular, that $\check{H}^0_\mathcal U(X,A)$ consists of sections $\sigma_{\alpha}$ over
$U_{\alpha}$ that agree on each intersection $U_{\alpha} \cap U_{\beta}$ so it coincides
with the global sections:  $\check{H}^0_{\mathcal U}(X,A) = \Gamma(X,A)$ for any covering $\mathcal U$.

\begin{tcolorbox}[colback=cyan!30!white]\begin{thm}  \label{thm-opencover}
Suppose $A$ is a sheaf on $X$ and an open cover $\mathcal U$ has the property that  $H^r(U_J,A) = 0$
for every $J \subset I$ and for all $r > 0$.  Then there is a canonical isomorphism $\check{H}^i(X,A)
\cong H^i(X,A)$ for all $i$.  In particular, $H^0(X,A) \cong \Gamma(X, A)$.\end{thm}\end{tcolorbox}
\subsection{Example} \label{subsec-simplicialcover}
 Let $f:|K|\to X$ be a triangulation of a topological space $X$.  This means that  $K$ is a simplicial complex 
with geometric realization $|K|$ and $f$ is a homeomorphism.  Let $A$ be a sheaf on $X$ and suppose that $K$ is
sufficiently fine that $A$ is constant on $f(\sigma)$ for each (closed) simplex $\sigma$ of $K$.
Then $A$ corresponds to a simplicial sheaf $S$ on $K$ in a canonical way.  The open stars of simplices form an open
cover $\mathcal U$ of $X$ that satisfies the condition of Theorem \ref{thm-opencover}.  The \v{C}ech complex
   $\check{C}^{\bullet}_{\mathcal U}(A)$ agrees with the complex $C^{\bullet}(K,S)$ of simplicial cochains (\S
   \ref{subsec-simplicialsheaf}), so they have
   the same cohomology:  $H^*(K,S) = \check{H}^*_{\mathcal U}(X,A).$

\subsection{}
Theorem \ref{thm-opencover} is an incredibly useful result because it says that we can use possibly very few and very large
open sets when calculating sheaf cohomology, and it even tells us how to tailor the open sets
to take advantage of the particular sheaf, whereas the original theorem of \v{C}ech assumed
that all the multi-intersections of the open sets were contractible.   The proof is very simple:

Let $A$ be a sheaf on a topological space and let 
$\mathcal U = \left\{ U_{\alpha} \right\}_{\alpha \in K}$ be an open covering
of $X$.  The \v{C}ech complex  $\check{C}\b(A)= \check{C}_{\mathcal U}\b(A)$ is
\[ \cdots \to  \prod_{|J|=r+1}\Gamma(U_J,A) \to \prod_{|J|=r+2}\Gamma(U_J,A) \to \cdots.\]
Sheafify this construction by defining
\[ \mathbf C^p(A)=\mathbf C^p_{\mathcal U}(A) = \prod_{|J|=p+1}i_*(A|U_J) \to \mathbf C^{p+1}(A) \to\cdots\]
so that $\mathbf C^p(A)(V) = \prod_{|J|=p+1} \Gamma(V\cap U_J, A)$.  This is functorial in $A$.  There is a little
combinatorial argument to show that
\[ 0 \to A \to \mathbf C^0(A) \to \mathbf C^1(A) \to\cdots\]
is exact. (To check exactness of the stalks at a single point $x$, we need to consider the
combinatorics of having $p+1$ open sets whose multi-intersection contains $x$.  Then
exactness comes down to proving that the homology of the $p$-simplex is trivial.)

Now let $A \to I^0 \to I^1 \to I^2 \to \cdots$ be an injective resolution of $A$, and apply the \v{C}ech resolution to each
term in this sequence, which gives a double complex of sheaves:
\begin{diagram}[size=2em]
\mb{C}^0(I^2) & \rTo & \mb{C}^1(I^2) & \rTo & \mb{C}^2(I^2) & \rTo &\\
\uTo       &         &   \uTo     &        &       \uTo  &        &\\
\mb{C}^0(I^1) & \rTo & \mb{C}^1(I^1) & \rTo & \mb{C}^2(I^1) & \rTo & \\
\uTo      &          &    \uTo    &         &   \uTo     &        & \\
\mb{C}^0(I^0)& \rTo  & \mb{C}^1(I^0) & \rTo & \mb{C}^2(I^0)  & \rTo &
\end{diagram}
whose rows are exact.  Let $\mathbf T^{\bullet}$ denote the associated single complex 
(\S \ref{subsec-doublecomplex})  of sheaves.
If we augment the left column with the column $I^0 \to I^1 \to I^2 \to \cdots $ then Corollary 
\ref{cor-double-one} says that the  
resulting map on cohomology sheaves $\mathbf{H}^*(I^{\bullet})\to \mathbf{H}^*(T^{\bullet})$ is an isomorphism,
which is to say that $I^{\bullet} \to T^{\bullet}$ is a quasi-isomorphism.  So it induces an isomorphism on
cohomology, $H^*(I^{\bullet}) = H^*(X,A) \cong H^*(T^{\bullet})$.

On the other hand, let us take global sections to obtain a double complex of groups. The
$r$-th column now reads   (from the bottom up)
\[  \prod_{|J|=r+1}\Gamma(U_J,I^0) \to \prod_{|J|=r+1}\Gamma(U_J, I^1) \to \cdots\]
which is a complex that computes the product of hypercohomology groups 
$\prod_{|J|=r+1}H^*(U_J,A)=0$ by hypothesis.  The kernel of the zeroth vertical map is
exactly the \v{C}ech cochains $\check{C}^r(A)$.  Therefore, if we augment  the bottom row with the complex $\check{C}^0(A) \to \check{C}^1(A) \to \check{C}^2(A) \cdots$ of \v{C}ech cochains then
Corollary \ref{cor-double-one} says that we will obtain a quasi-isomorphism of this complex with the total complex
of this double complex, namely $\Gamma(X, T^{\bullet})$ .  Hence, the cohomology of the
\v{C}ech complex of groups coincides with the cohomology of this total complex, which was
shown above to coincide with the hypercohomology of the sheaf $A$ as computed using
injective resolutions.  \qed

   \quash{
\subsection{} If $A\b$ is a complex of sheaves (bounded below) 
on a topological space $X$ then the \v{C}ech construction applied
to each $A^j$ gives a (first quadrant) double complex $\check{C}\b(A\b)$.  
Let $T\b$ be the associated single complex.  Then there is a canonical isomorphism
\[ H^r(X;A\b) \cong H^r(T\b)\]
and a spectral sequence that converges to this cohomology with
\[ E_2^{pq} \cong H^p(X; \uul{H}^q(A)).\]
}

\section{The sheaf of chains}\label{subsec-chains}\label{subsec-sheafofchains}
\index{sheaf!of chains}\index{chains}

\subsection{The problem}
Most of the complexes of sheaves that were discussed until now have the property that
their cohomology sheaves live only in degree zero.  The sheaf of chains, however, is
a naturally occurring entity with complicated cohomology sheaves.  But
it is not so obvious how to construct a sheaf on $X$ that corresponds to (say) the singular chains.  How do you
define the restriction of a singular simplex that might wrap many times around the whole space?
One might define a presheaf $C_r$ with sections $\Gamma(U,C_r) = C_r(U)$ to be the
group of singular chains on $U$,  with restriction mapping $C_r(U) \to C_r(V)$ that
assigns zero to every singular simplex that is not completely contained in $V$ but this map may fail
preserve the boundary homomorphism (but see \S \ref{subsec-finitetype} below).  

It is necessary to confront the problem of restricting (to $V$) a singular simplex $\sigma \subset U$ that is
not completely contained in $V$.  The solution is to divide $\sigma \cap V$ into (possibly) infinitely
many simplices and to call this the restriction $\sigma|V$.  So the sheaf of chains requires the use of non-compact
chains that may be infinite sums of simplices.  
 If the space $X$ is paracompact then we may restrict to sums that are locally finite.

\subsection{Borel-Moore homology}\index{Borel-Moore homology}\index{Borel-Moore homology}\index{homology, Borel-Moore}
Borel and Moore defined a sheaf $\uul{C}{}^{\bullet}_{BM}$ 
(that is, a complex of sheaves) whose presheaf of sections $\Gamma(U, \uul{C}{}^r_{BM})$
is the ``locally finite $r$ dimensional chains in $U$''.
The cohomology of this sheaf is called the Borel-Moore homology of the space \cite{BorelMoore}.
There are a number of rigorous constructions of this sheaf see \S \ref{subsec-derivedchains} and
\S \ref{subsec-duality}.  For the moment, let us take it on faith
that such a sheaf (or rather, complex of sheaves) exist.  The impatient reader may jump to
\S \ref{subsec-duality} for the general definition, where the sheaf of chains becomes an
essential actor.  The usual (e.g. singular) homology of $X$ is the compact support cohomology:
\[ H^{-i}_c(X;\uul{C}{}^{\bullet}_{BM}) = H_i(X).\]

The stalk cohomology of the sheaf of chains is the local homology:  
$\uul{H}{}^{-r}_x(\uul{C}^{\bullet}_{BM}) \cong
H_r(X, X-x)$.  We use negative degrees in the superscript so that the boundary map on chains raises degree,
as with all complexes of sheaves.

If $X$ is a (topological) manifold of dimension $n$ 
then the top degree local homology group (say, with integral coefficients), $H_n(X,X-x;\ZZ)$, is isomorphic to 
$\ZZ$ but such an isomorphism involves a choice of
local orientation of $X$ near $x$.  This local homology sheaf $\mathcal O= \uul{H}^{-n}({}^{\bar p}{C}{}^{\bullet}_{BM})$ is
the {\em orientation sheaf} of $X$.  \index{sheaf!orientation}\index{orientation sheaf}
An orientation of $X$ (if one exists) is a global section of the orientation sheaf.

\subsection{Finite type}\label{subsec-finitetype}
Let us say that a topological space $X$ has {\em finite type} if it is homeomorphic to
$K - L$ where $K$ is a finite simplicial complex and $L$ is a closed subcomplex.  In this
case $H_r^{BM}(X) \cong H_r(K,L)$ coincides with the relative homology.  Therefore it  can then
be expressed as the homology of a chain complex formed by the simplices that are
contained only in $K$, and by defining the differential so as to ignore all components of the boundary that may lie in $L$.  This gives a simple combinatorial and computable construction of Borel-Moore
homology for spaces of finite type, although it does not describe the Borel-Moore {\em sheaf}.

\section{Homotopy and injectives}
\subsection{Homotopy theory}\index{homotopy}
Two morphisms $f, g:A^{\bullet} \to B^{\bullet}$ of complexes are said to be {\em homotopic} if
there is a collection of mappings $h:A^{r} \to B^{r-1}$ so that $hd_A + d_Bh = f-g$.  This
is an equivalence relation.  Equivalence classes are referred to as homotopy classes of
maps; the set of which is denoted $[A^{\bullet}, B^{\bullet}]$.  Define the complex of abelian
groups (or $R$ modules)
\[ \hHom^n(A^{\bullet}, B^{\bullet}) = \prod_{s}\hHom(A^s, B^{s+n})\]
with differential $df= d_Bf + (-1)^{n+1} fd_A$ where $f:A^s \to B^{s+n}$.  

\begin{tcolorbox}[colback=cyan!30!white]
\begin{lem}  Let $f:A\b \to B\b$ be a morphism of complexes and let $C(f)$ be the cone
of $f$.  For any complex $S\b$ we have a morphism of complexes
$f_*: \hHom\b(S\b, A\b) \to \hHom\b(S\b,B\b)$.  Then there is a canonical isomorphism of complexes
of abelian groups,
\[C(f_*) \cong \hHom\b(S\b, C(f)).\]
For any complex $T\b$ we have a morphism of complexes $f^*:\hHom(B\b, T\b)
\to \hHom(A\b, T\b)$.  Suppose the cohomology of $S\b, T\b$ is bounded, that is,
$H^r(S\b) = H^r(T\b) = 0$ if $|r|$ is sufficiently large.
Then there is a canonical quasi-isomorphism of complexes of abelian groups,
\[ C(f^*)[-1] \cong \hHom\b(C(f), T\b).\]
\end{lem}\end{tcolorbox}
The first statement is obvious because $\Hom(S^s, A^{t+1}\oplus B^t) =
\Hom(S^s, A^{t+1}) \oplus \Hom(S^s, B^t)$.  The second statement is similar.
In particular given $f:A^{\bullet} \to B^{\bullet}$  let $C^{\bullet} = C(f)$ be the cone.  Then 
for any complex $S\b$ there is a long exact sequence
\[ \cdots [S\b, A\b] \to [S\b, B\b] \to [S\b,C\b] \to [S\b, A\b[1]] \to [S\b, B\b[1]] \to\cdots \]
We can do the same with sheaf-Hom.  Recall that $\uul{\hHom}(A,B)(U) =
\hHom_{Sh(U)}(A|U, B|U)$). We obtain a complex of sheaves,
\[\uul{\hHom}^n(A^\b,B^\b) = \prod_s \uul{\hHom}(A^s, B^{s+n})\]
 with the property that 
\[
H^0(X, \uul{\hHom}\b(A^\b,B^\b)) = \Gamma(X, \uul{\hHom}\b(A^\b, B^\b)) =
H^0(\Hom\b(A\b,B\b)) = [A\b,B\b]
\]
 Therefore we have (cf. \S \ref{subsec-Examples7} (4) below): 
 \begin{tcolorbox}[colback=cyan!30!white]
\begin{prop} \label{prop-Homdot}
$ H^n(\hHom\b(A^{\bullet}, B^{\bullet})) = [A^{\bullet}, B^{\bullet}[n]]. $\qed
\end{prop}\end{tcolorbox}

\subsection{}  The bounded {\em homotopy category} $K^b(X)$ of sheaves on $X$ is the category whose objects are
complexes of sheaves whose cohomology sheaves are bounded (meaning that
$H^r(A\b) = 0$ for sufficiently large $r$),
and whose morphisms are homotopy classes of morphisms, that is,
\[ \hHom_{K^b(X)}(A\b, B\b) = [A\b, B\b] = H^0(X, \uul\hHom\b(A\b,B\b)).\]

\subsection{Wonderful properties of injective sheaves}\index{sheaf!injective}
Roughly speaking, when we restrict to injective objects, then quasi-isomorphisms become
homotopy equivalences.  For sheaf theory, this is important because a homotopy of complexes
of sheaves also gives a homotopy on global sections.  In this way, quasi-isomorphisms
of complexes of injective sheaves give isomorphisms on hypercohomology.

\begin{tcolorbox}[colback=cyan!30!white]
\begin{lem}  Let $C^{\bullet}$ be a (bounded below) complex of sheaves and suppose that
the cohomology sheaves $\mathbf H^r(C^{\bullet}) = 0$ for all $r$.  Let $J^{\bullet}$ be a
complex of injective sheaves.  Then any morphism $f:C^{\bullet} \to J^{\bullet}$ is homotopic
to zero, meaning that there exists $h:C^{\bullet} \to J^{\bullet}[-1]$ such that $d_Jh + hd_C = f$.
\end{lem}\end{tcolorbox}

\begin{proof}
 It helps to think about the diagram of complexes:
\begin{diagram}[size=2.5em]
0 & \rTo & C^0 & \rTo^{d^0}  & C^1   &\rTo^{d^1}     & C^2 & \rTo^{d^2}  & C^3    & \rTo \cdots \\
  &         & \dTo^f&\ldTo^{h^1} & \dTo^f&   \ldTo^{h^2} & \dTo^f &  \ldTo^{h^3} & \dTo^f & \\
0 & \rTo & J^0& \rTo_{d^0}  & J^1   &\rTo_{d^1}     & J^2 & \rTo_{d^2}  & J^3    & \rTo \cdots 
\end{diagram}
The first step is easy, since $C^0 \to C^1$ is an injection and since $J^0$ is injective there exists $h^1:C^1 \to J^0$ that
makes the triangle commute, that is, $h^1d^0 = f$.  Now let us define $h^2:C^2 \to J^1$.  Consider the map
$(f-d^0h^1):C^1 \to J^1$.  It vanishes on $\Im(d^0) = \ker(d^1)$ because 
\[(f-d^0h^1)d^0 = fd^0 - d^0h^1d^0 = fd^0-d^0f = 0.\]  Therefore it passes to a vertical mapping  in this diagram:
\begin{diagram}[size=2em]
C^1 & \rTo & C^1/\ker(d^1) & \rInto & C^2 \\
&\rdTo& \dTo_{f-d^0h^1} \\
&&J^1 \end{diagram}
where the second horizontal mapping is an injection.  Since $J^1$ is injective we obtain an extension $h^2: C^2 \to J^1$ such
that $h^2\circ d^1 = f - d^0 h^1$ so that $h^2d^1 + d^0h^1 = f$.  Continuing in this way, the other $h^r$ can be constructed inductively.\end{proof}

Exactly the same argument may be used to prove the following:
\begin{lem}  Let $f:A^{\bullet} \to B^{\bullet}$ be a quasi-isomorphism of (bounded below) complexes.  Then
for any complex $J^{\bullet}$ of injectives, the induced map $[B^{\bullet}, J^{\bullet}] \to
[A^{\bullet}, J^{\bullet}]$ on homotopy classes is an isomorphism. \qed\end{lem}
This result can also be proven by
applying the previous lemma to the cone $C(f)$ and using the long exact sequence on cohomology.

\begin{tcolorbox}[colback=cyan!30!white]
\begin{cor}  The following statements hold.
\begin{enumerate}
\item  Suppose $J\b$ is a complex of injective sheaves and $H^n(J\b)=0$ for all $n$.  Then
$J\b$ is homotopy equivalent to the zero complex.
\item Let $\phi:X^{\bullet} \to Y^{\bullet}$ be a quasi-isomorphism of sheaves of {\em injective} complexes.
Then $\phi$ admits a homotopy inverse $g:Y^{\bullet} \to X^{\bullet}$ (meaning that $g\phi \sim I_X$
and $\phi g\sim I_Y$).
\item Let $A\b \to I\b$ and $B\b \to J\b$ be injective resolutions of complexes $A\b, B\b$.  Then
any morphism $f:A\b \to B\b$ admits a lift $\tilde f: I\b \to J\b$ and any two such lifts are homtopic.
\item Let $f:A^{\bullet} \to B^{\bullet}$ be a quasi-isomorphism of complexes of sheaves.  Then
$f$ induces an isomorphism on hypercohomology $H^r(X,A^{\bullet}) \cong H^r(X,B^{\bullet})$
for all $r$.  \end{enumerate}
\end{cor}\end{tcolorbox}
\begin{proof}  For (1) consider the identity mapping $J\b \to J\b$
For (2) ,  mapping $(X\b \overset{\phi}{\longrightarrow} Y\b)$ to $X\b$ and using the lemma gives an isomorphism
$[Y\b,X\b] \to [X\b,X\b]$, the map given by $f \mapsto f\circ\phi$.  So the identity $X\b \to X\b$
corresponds to some $f$ such that $f \circ \phi \sim Id$, implying that $\phi$ has a left 
homotopy-inverse.
Now consider mapping $ Y\b$ into the triangle $X\b \to Y\b \to C(\phi) \to \cdots$, giving an
exact sequence $\cdots \to [Y\b,X\b] \to [Y\b,Y\b] \to [Y\b, C(\phi)] \to \cdots$.  Since
$C(\phi)$ is injective and its cohomology vanishes, the identity is homotopic to zero, 
hence $[Y\b, C(\phi)] = 0$ so that
$[Y\b,X\b] \cong [Y\b,Y\b]$ with the map given by $g \mapsto \phi \circ g$.  Therefore there
exists $g:Y\b \to X\b$ so that $\phi \circ g \sim Id$ meaning that $\phi$ has a right inverse in
the homotopy category.  If a mapping has both a left inverse and a right inverse then it has
an inverse (in other words, $f$ and $g$ are homotopic, so either of them will behave as a
homotopy inverse to $\phi$).  
For (3), the lemma gives an isomorphism $[A\b, J\b] \to [I\b, J\b]$.
For (4), the hypercohomology is defined in terms of the global sections of an injective
resolution.  So we may assume that $A\b$ and $B\b$ are injective. By the lemma,
the cone $C(f)$ is homotopic to zero.  Let $h$ be such a homotopy.  Now take global sections.
The global sections of
the cone coincides with the cone on the global sections, that is, we have a triangle of groups:
\begin{diagram}[size=2.5em]
\Gamma(X, A\b) && \rTo && \Gamma(X,B\b) \\
 & \luTo_{[1]} && \ldTo & \\
&& \Gamma(X, C(f))
\end{diagram}
 The homotopy $h$ also gives a homotopy on the global sections so that
$H^n(\Gamma(X,C(f))) = 0$.  So the long exact sequence on cohomology implies that
$H^n(X,A\b) \to H^n(X,B\b)$ is an isomorphism.
\end{proof}

\bigskip

\section{The derived category}  \index{derived category}\index{category!derived}
Jean-Louis Verdier \cite{Verdier0, Verdier1}, following a suggestion of Grothendieck, developed the 
notions of the derived category of an abelian category.
There are several different ``models' for the derived category.  The first definition we give is
easy to understand and useful for proofs but the objects themselves are not very natural. 
The second model is less intuitive but the objects occur naturally.

\subsection{The derived category:  first definition}
Let $X$ be a topological space.  The bounded derived category $D^b(X)$ is the
category whose objects are complexes of  injective sheaves whose cohomology sheaves are bounded (meaning that $H^r(A\b) = 0$ for
sufficiently large $r$ and for sufficiently small $r$).  The morphisms are homotopy classes
of morphisms (of complexes of sheaves), so that $D^b(X)$ is the {\em homotopy category}
of (complexes of) injective sheaves.  

If $\mathcal{Sh}(X)$ denotes the category of sheaves on $X$, if $C^b(X)$ denotes the
category of complexes of sheaves with bounded cohomology and if $K^b(X)$ denotes
the homotopy category of (complexes of) sheaves with bounded cohomology
then we have a canonical functors
\[
\begin{CD} \mathcal{Sh}(X) @>>> C^b(X) @>>>K^b(X)\ @>{\text{God}}>{\longleftarrow}> D^b(X)
\end{CD}\] 
that associates to any sheaf $S$ the corresponding
complex concentrated in degree zero, and to any complex $A\b$ its Godement injective resolution.
[This construction makes sense if we replace the category of sheaves with any abelian category $\mathcal C$ provided it has enough injectives.  In this way we define the bounded derived category
$D^b(\mathcal C)$ with functors $\mathcal C \to K^b(\mathcal C) \to D^b(\mathcal C)$. ] 

From the previous section on ``properties of injective sheaves" we therefore conclude:  \begin{itemize}
\item The mapping $K^b(X) \to D^b(X)$ is a functor (that is, a morphism between complexes
determines a morphism in the derived category also).
\item If $A\b \to B\b$ is a quasi-isomorphism of complexes of sheaves then it becomes
an isomorphism in $D^b(X)$.
\item if $A\b$ is a complex of sheaves such that $H^m(A\b)=0$ for all $m$ then $A\b$ is
isomorphic to the zero sheaf.
\end{itemize}

\begin{defn}{}
Let $T:\mathcal{Sh}(X) \to \mathcal{B}$ be a covariant (and additive) functor from the category of sheaves to some other abelian category with enough injectives.   
For the moment, let us also assume that it takes
injectives to injectives.  Define the {\em right derived functor}
$RT: D^b(X) \to D^b(B)$ by $RT(A\b)$ to be the complex $T(I^0) \to T(I^1) \to T(I^2) \to \cdots$  where $A\b \to I\b$ is the canonical (or the chosen)
injective resolution of $A\b$.  Define $R^mT(A)$ (``the $m$-th derived functor", an older terminology)
to be the cohomology object of this complex, $H^m(RT(A\b))$.
\end{defn}
Let $A\b \to I\b$ and $B\b \to I\b$ be the canonical
injective resolutions of $A\b, B\b$.  Then any morphism $f:A\b \to B\b$ has a lift $\tilde{f}:
I\b \to J\b$ that is unique up to homotopy, which is to say that $\tilde{f}$ is a uniquely defined
morphism in the category $D^b(X)$, and we obtain a well defined morphism 
\[RT(f)=T(\tilde f):RT(A\b)=T(I\b) \to T(J\b)= RT(B\b).\]
In other words, the right derived functor of $T$ is obtained by replacing each complex $A\b$ by
its injective resolution $I\b$ and then applying $T$ to that complex.  

\subsection{Examples}\label{subsec-Examples7}\noindent
{1}  
If $f:X \to Y$ is a continuous map between topological spaces, then we show (below)
that $f_*:\mathcal{Sh}(X) \to \mathcal{Sh}(Y)$  takes injectives to injectives.  
Taking $Y = \{pt\}$ we get that the global sections functor $\Gamma$ takes injectives
to injectives. Then, for any complex of sheaves $A\b$,
\[ R^m\Gamma(A\b) = H^m(\Gamma(X, I\b)) \]
where $A\b \to I\b$ is the canonical injective resolution.  
(This is how we defined the hypercohomology of the complex of sheaves $A$ in \S \ref{def-hypercohomology}.) 
\medskip
\paragraph{2}  Let $f:X \to Y$ be a continuous map and let $\uul{\ZZ}$ be the
constant sheaf on $X$.  If $f$ is surjective and its fibers are connected
then  $f_*(\uul{\ZZ})$ is again the constant sheaf, because as a presheaf,
$f_*(\uul{\ZZ}(U) = \uul{\ZZ}(f^{-1}(U)) = \ZZ$ for any connected open set
$U \subset Y$.   Although we cannot hope to understand $Rf_*(\uul{\ZZ})$ we can
understand its cohomology sheaves:
\[ R^mf_*(\uul{\ZZ})(U) = H^m(\Gamma(f^{-1}(U), I\b))\]
where $I\b$ is an injective resolution (or perhaps the canonical injective resolution) of the
constant sheaf.  But this is exactly the definition of the hypercohomology $H^m(f^{-1}(U),
\uul{\ZZ}) = H^m(f^{-1}(U,\ZZ))$.  If $f$ is proper then the stalk cohomology (of the cohomology
sheaf of $Rf_*(\uul{\ZZ})$) at a point $y \in Y$ is equal to $H^m(f^{-1}(y);\ZZ)$,
the cohomology of the fiber. In other words, the sheaf $\uul{\mathbf H}^m(Rf_*(\uul{\ZZ})) $
is a sheaf on $Y$ which, as you move around in $Y$, displays the cohomology of the fiber.
\medskip
\paragraph{3}  We can also determine the global cohomology of the complex $Rf_*(\uul{\ZZ})$,
for it is the cohomology of the global sections $\Gamma(X, I\b)$, that is, the cohomology of $X$.
More generally, the same argument shows that:  for any complex of sheaves $A\b$ on $X$,
the complex $Rf_*(A\b)$ is a sheaf on $Y$ whose global cohomology is
\[ H^*(Y, Rf_*(A\b)) \cong H^*(X, A\b).\]
This complex of sheaves therefore provides data on $Y$ which allows us to compute the
cohomology of $X$.  It is called the Leray Sheaf (although historically, Leray really considered
only its cohomology sheaves $R^mf_*(A\b)$).  In particular we see that {\em the functor
$Rf_*$ does not change the hypercohomology}.  For $f:X \to \{ \text{pt} \}$, if $S$ is a sheaf
on $X$ then $f_*(S) = \Gamma(X,S)$ is the functor of global sections (or rather, it is a sheaf
on a single point whose value is the global sections), so $R^if_*(A\b) = H^i(X,A\b)$ is the
hypercohomology.
\medskip
\paragraph{4}  The $m$-th derived functor of $\Hom$ is called $\Ext^m$, (cf. Proposition \ref{prop-Homdot})
\index{Ext} i.e., it is the group
\[ \Ext^m(A\b, B\b) = H^m(\RHom (A\b,B\b)) = H^m(\hHom\b(A\b, J\b))= H^0(\hHom\b(A\b,
J\b[m]))\]
where $B\b \to J\b$ is an injective resolution.  (We consider $\Hom(A\b, B\b)$ to be a functor
of the $B\b$ variable and derive it by injectively resolving. 
The same result can be obtained by {\em projectively} resolving $A\b$.)  As before, 
there is a sheaf version of $\Hom$, which also gives a sheaf version of $\Ext$:
\[
\uul{\Ext}^m(A\b, B\b)) = \uul{H}^m(\uul{\RHom}(A\b,B\b))
\]

\medskip
\paragraph{Exercise}  Let $G,H$ be abelian groups,  as complexes in degree zero.
Show that $\Ext^1(G,H)$ coincides with the usual definition of $\Ext_{\ZZ}(G,H)$, 

\subsection{The derived category:  second construction}\index{derived category}
The derived category can be constructed as a sort of quotient category
of the homotopy category $K^b(X)$ of complexes, by inverting quasi-isomorphisms.
Let $E^b(X)$ be the category whose objects are complexes of sheaves on $X$,
and where a morphism $A\b \to B\b$ is an equivalence class of diagrams
\begin{diagram}[size=1.5em]
&& C\b && \\
& \ldTo^{qi} && \rdTo \\
A\b && && B\b
\end{diagram} 
where $C^\b \to A\b$ is a quasi-isomorphism, and where two such morphisms 
$A\b \leftarrow C_1^{\bullet} \to B\b$ and $A\b \leftarrow C_2^{\bullet} \to B\b$ are
considered to be equivalent if there exists a diagram
\begin{diagram}[size=2em]
&&  C_1^{\bullet}  && \\
& \ldTo &\uTo& \rdTo \\
A\b & \lTo^{qi} & C_3^{\bullet} & \rTo & B\b \\
& \luTo & \dTo & \ruTo & \\
&& C_2^{\bullet} &&
\end{diagram}
that is commutative up to homotopy.  (Exercise:  figure out how to compose two morphisms and then
check that the result is well defined with respect to be above equivalence relation.)

\begin{tcolorbox}[colback=cyan!30!white]
\begin{thm}  The natural functor $D^b(X) \to E^b(X)$ is an equivalence of categories.
\end{thm}\end{tcolorbox}

\begin{proof}
To show that a functor $F:\mathcal C \to \mathcal D$
is an equivalence of categories it suffices to show (a) that it is
essentially surjective, meaning that every object in $\mathcal D$ is isomorphic to an object
$F(C)$ for some object $C$ in $\mathcal C$, and (b) that $F$ induces an isomorphism on
$\Hom$ sets.  The first part (a) is clear because we have injective resolutions.  Part (b)
follows immediately from the fact that a quasi-isomorphism of injective complexes is a
homotopy equivalence and has a homotopy inverse.
\end{proof}

This gives a way of referring to elements of the derived category without having to
injectively resolve.  Each complex of sheaves is automatically an object in the
derived category $E^b(X)$.  Here are some applications.

\subsection{$T$-acyclic resolutions}\index{resolution!$T$-acyclic}
A functor $T:\mathcal C \to \mathcal D$ between abelian categories is
{\em exact} if it takes exact sequences to exact sequences.  It is {\em left exact} if 
it preserves kernels, that is, if $f:X \to Y$ and if $Z = \ker(f)$ (so $0 \to Z \to X \to Y$ is
exact) then $T(Z) = \ker(T(f))$ (that is, $0 \to T(Z) \to T(X) \to T(Y)$ is exact).   An object
$X$ is $T$-acyclic if $R^iT(X)=0$ for all $i \ne 0$.  This means:  take an injective
resolution $X \to I\b$, apply $T$, take cohomology, the result should be zero except
possibly in degree zero.

The great advantage of $T$-acyclic objects is that they may be often used in place of injective
objects when computing the derived functors of $T$, that is

\begin{tcolorbox}[colback=cyan!30!white]
\begin{lem}  Let $T$ be a left exact functor from the category of sheaves to some abelian category
with enough injectives.
Let $A\b$ be a complex of sheaves and let $A\b \to X\b$ be a quasi-isomorphism, where
each of the sheaves  $X^r$ is $T$-acyclic.  Then $R^rT(A\b)$ is canonically isomorphic
to the $r$-th cohomology object of the complex
\[ T(X^0) \to T(X^1) \to T(X^2) \to \cdots.\]
If $T$ is exact then there is no need to take a resolution at all:  $RT(A\b)$ is canonically
isomorphic to $T(A\b)$.
\end{lem}\end{tcolorbox}
The proof is the standard double complex argument:  Suppose $T$ is left exact.
Let $I^{\bullet\bullet}$ be a double
complex of injective sheaves, the $r$-th row of which is an injective resolution of $X^r$.
Let $Z\b$ be the total complex of this double complex.  It follows that $A\b \to X\b \to Z\b$
are quasi-isomorphisms and so the complex $Z\b$ is an injective resolution of $A\b$.
Now augment the double complex by attaching $X\b$ to the zeroth column, and
apply the functor $T$ to the augmented complex.  Since each $X^r$ is $T$-acyclic, each
of the rows remains exact except possibly at the zeroth spot.  Since $T$ is left exact,
the rows are also exact at the zeroth spot.  So our lemma says that $T(X\b) \to
T(Z\b)=RT(A\b)$ is also a quasi-isomorphism, which gives an isomorphism between
their cohomology objects.  If the functor $T$ is exact then every object $A^r$ is
$T$-acyclic (exercise), so the original complex $A\b$ may be used as its own
$T$-acyclic resolution.

\subsection{Key exercises}
Show that injective objects are $T$-acyclic for any left exact functor $T$.  
If $f:X \to Y$ is a continuous map,  show that
$f^*:{\mathcal{Sh}}(Y) \to \mathcal {Sh}(X)$ is exact (and so it does not need to be derived).  
Using this and the adjunction formula
$\Hom_{{\mathcal{Sh}}(Y)}(B, f_*(I)) \cong \Hom_{\mathcal{Sh}(X)}(f^*(B),I)$ 
show that $f_*$ is left exact and takes injectives to
injectives. Show that fine, flabby, and soft sheaves
are $\Gamma$-acyclic.  In particular, {\em the cohomology of
a sheaf (or of a complex of sheaves) may be computed with respect to any injective, fine,
flabby, or soft resolution}.

\subsection{More derived functors}\label{subsec-derivedfunctors}\index{derived functor}
This also gives us a way to define ``the'' derived functor $RT:D^b(X) \to D^b(\mathcal C)$
for any left exact functor $T: Sh(X) \to \mathcal{C}$ provided the category $\mathcal C$
has enough injectives, namely, if $A\b$ is a complex of sheaves on $X$, take an
injective resolution $A\b \to I\b$, then apply the functor $T$ to obtain a complex
$T(I\b)$ of objects in the category $\mathcal C$, then injectively resolve this complex
by the usual method of resolving each $T(I^r)$ to obtain a double complex, then
forming the associated total complex.  Let $RT$ denote the resulting complex.  Different
choices of resolutions give isomorphic complexes $RT$.

\subsection{The sheaf of chains in the derived category}\label{subsec-derivedchains}\index{sheaf!of chains}
  Suppose $X$ is a piecewise-linear space, that is, a topological space
together with a family of piecewise-linearly related triangulations by locally finite countable simplicial
complexes. 
 Let $U \subset X$ be an open subset,  let $T$ be a
locally finite triangulation of $U$,  and let $C_r^T(U)$ be the group of $r$-dimensional
simplicial chains with respect to this triangulation.  Then the {\em sheaf of piecewise
linear chains} is the sheaf $\uul{\mathbf C}{}^{\bullet}_{PL}$ with sections
\[ \Gamma (U, \uul{\mathbf C}{}^{-r}_{PL}) = 
\underset{\overrightarrow{T}}{\lim}C_r^T(U)\]
for $r \ge 0$. Such a section is a locally finite but possibly infinite linear combination of
oriented simplices.   (We place the chains in negative degrees so that
the differentials will increase degree; it is a purely formal convention.)   

The sheaf $\uul{\mathbf C}{}^{-r}_{PL}$ is {\em soft}, and the resulting complex
\[ \uul{\mathbf C}{}^0_{PL} \longleftarrow \uul{\mathbf C}{}^{-1}_{PL} \longleftarrow
\uul{\mathbf C}{}^{-2}_{PL} \longleftarrow \cdots\]
is quasi-isomorphic to the sheaf of Borel-Moore chains.  

\quash{
Suppose $X$ is a finite simplicial complex.  If $\sigma$ is a (closed) simplex let 
${\uul{\QQ_{\sigma}}}$
denote the constant sheaf on $\sigma$.  It is injective in the category of simplicial sheaves (with $\QQ$ coefficients)
on $X$.  Every injective simplicial sheaf (of rational vector spaces) is a direct sum of
such {\em elementary injectives}.  The sheaf of chains can be realized as
the (injective) complex 
\[
\bigoplus_{\dim(\sigma)=0}\uul{\QQ}{}_{\sigma} \longleftarrow \bigoplus_{\dim(\sigma)=1}
\uul{\QQ}{}_{\sigma} \longleftarrow \bigoplus_{\dim(\sigma)=2} \uul{\QQ}{}_{\sigma} 
\longleftarrow \cdots
\]
in degrees $0, -1, -2 ,\cdots$ respectively. The global sections
of this complex equals the usual complex $C_{\bullet}(X)$ of simplicial chains.  If we use
the constant sheaf $\ZZ_{\sigma}$ everywhere, then the resulting sheaves are soft, rather than injective, but they may still be used to compute the homology of $X$.

Now consider the limit over all subdivisions of the simplicial complex $X$.  We define
a ``topological" sheaf on $X$ that is, in some sense, the sheafification of the direct limit
of these sheaves.  To be precise,
}
If the space $X$ has a real analytic (or semi-analytic or subanalytic or $\mathcal O$-minimal)
structure then one similarly has the sheaf of locally finite subanalytic or $\mathcal O$-minimal
chains, which gives another quasi-isomorphic ``incarnation" of the sheaf of chains.
See also \S \ref{subsec-duality}. 

\subsection{The bad news}
The derived category is not an abelian category.  In fact, the homotopy category of complexes
is not an abelian category.  Kernels and cokernels do not make sense in these categories.
The saving grace is that the cone operation still makes sense and in fact, it passes to the
homotopy category.  So we have to replace kernels and cokernels with triangles.

\begin{tcolorbox}[colback=yellow!30!white]
\begin{defn}
A triangle of morphisms
\begin{diagram}[size=2em]
A\b && \rTo && B\b \\
&\luTo_{[1]}&&\ldTo &\\
&& C\b &&
\end{diagram}
in $K^b(X)$ or in $D^b(X)$ is said to be a {\em distinguished triangle}
\index{triangle, distinguished}\index{distinguished triangle}or an {\em exact triangle}
 if it is  homotopy equivalent to a triangle
\begin{diagram}[size=2em]
X\b && \rTo^{\phi} && Y\b  \\
& \luTo_{[1]}  && \ldTo & \\
&& C(\phi) &&
\end{diagram}
where $C(\phi)$ denotes the cone on the morphism $\phi$.
\end{defn}\end{tcolorbox}
This means that there are morphisms between corresponding objects in the triangles such that the
resulting squares commute up to homotopy.
(In the homotopy categories $K(X)$ and $D^b(X)$ homotopy equivalences are
isomorphisms, so people often define a distinguished triangle to be a triangle in
$K(X)$ that is isomorphic to a mapping cone.)

\begin{tcolorbox}[colback=cyan!30!white]
\begin{lem}
The natural functor $K^b(X) \to D^b(X)$ takes distinguished triangles to distinguished triangles.  If
\begin{diagram}[size=2em]
A\b && \rTo && B\b \\
&\luTo_{[1]}&&\ldTo &\\
&& C\b &&\end{diagram}
is a distinguished triangle and if $X\b$ is a complex (bounded from below) then there are distinguished triangles
\renewcommand{\ss}{\scriptstyle}\newcommand{\fs}{\footnotesize}
\begin{diagram}[size=2em]
{\ss{\rm RHom}(X\b, A\b)} && \rTo &&{\ss {\rm RHom}(X\b, B\b)} & & {\ss{\rm RHom}(A\b,X\b)} && \lTo && {\ss\RHom(B\b,X\b)}\\
& \luTo_{[1]} && \ldTo &    &\ \text{ and }\ &  & \rdTo_{[1]} && \ruTo & \\
&& {\ss\RHom(X\b,C\b)} && && && {\ss\RHom(C\b,X\b)} &&
\end{diagram}
\end{lem}\end{tcolorbox}
The proof is the observation (from the last section) that $\Hom$ into a cone is equal to the
cone of the $\Hom$s.  We stress again that the hypercohomology of $\RHom$ is exactly
the group of homomorphisms in the derived category:
$H^0(X,\RHom(A\b,B\b)) = \hHom_{D^b(X)}(A\b, B\b)$.

\subsection{Exact sequence of a pair}\label{subsec-subspaces}\index{exact sequence of a pair}
Let $Z$ be a closed subspace of a topological space $X$, and let $U = X-Z$, 
say $ Z \overset{i}{\longrightarrow} X \overset{j}{\longleftarrow} U.$

If $S$ is a sheaf on $X$ then there is a short exact sequence of sheaves
$0 \to j_!j^*S \to S \to i_*i^*S \to 0$.  The morphisms are obtained by adjunction, and exactness
can be checked stalk by stalk:  If $x \in Z$ then the sequence reads $0 \to 0 \to S_x \to S_x \to 0$.
If $x \in U$ then the sequence reads $0 \to S_x \to S_x \to 0 \to 0$.
Consequently if $A\b$ is a complex of sheaves then there is a distinguished triangle
\begin{tcolorbox}[colback=cyan!30!white]\begin{equation}\label{eqn-pair}
\begin{diagram}[size=2em]
Rj_!j^*(A\b) && \rTo && A\b \\
& \luTo_{[1]} && \ldTo & \\
&& Ri_*i^*(A\b) &&
\end{diagram}\end{equation}\end{tcolorbox}
Recall from Lemma \ref{lem-triangles} (which is an exercise) 
that the cokernel of an injective morphism is quasi-isomorphic to the mapping cone.
If $X$ is compact the triangle gives an exact sequence
 \[\cdots \to H^r(X;A\b) \to H^r(Z;A\b) \to H^{r+1}_c(U;A\b) \to \cdots\]
but observe that the relative cohomology $H^*_c(U;A\b) = H^*(X,Z;A\b)$ is the cohomology of cochains on $X$ that
vanish on $Z$.  If $A\b$ is the sheaf of chains this gives an exact sequence 
$H_r(U) \to H_r(X) \to H_r(X, U) \to \cdots$ because $i^*(\text{chains})$ is the limit
over open sets containing $Z$ of Borel-Moore chains on that open set.

If $S$ is a sheaf on $X$ define $i^!(S)$ to be the restriction to $Z$ of the 
presheaf with sections {\em supported} in $Z$, that is 
\[ i^!(S) = i^*(S^Z) \text{ where } \Gamma(V, S^Z) = \left\{ s \in \Gamma(V,S)|\ \text{spt}(s) \subset Z \right\}.\]
Thus, if $W \subset Z$ is open then
\[ \Gamma(W, i^!(S)) = \underset{\overrightarrow{V \supset W}}{\lim}\Gamma_Z(V,S)\]
(the limit is over open sets $V \subset X$ containing $W$). 
The functor $i^!$ is a right adjoint to the pushforward with compact support $i_!$, that is,
$\Hom_X(i_!A,B) = \Hom_Z(A, i^!B)$.  In fact, $\Hom_Z(A, i^!B)$ consists of mappings of
the leaf space $LA \to LB|Z$ that can be extended by zero to a neighborhood of $Z$ in $X$,
and this is the same as $\Hom_X(i_!A,B)$. We obtain a canonical morphism
$i_!i^!B \to B$.

\begin{prop}\label{prop-pair-sequence}
 Let $A\b$ be a complex of sheaves on $X$.  There is a distinguished triangle, Verdier dual to
 (\ref{eqn-pair}),
 \begin{tcolorbox}[colback=cyan!30!white]
\begin{diagram}[size=2em]
Ri_*i^!(A\b) && \rTo && A\b \\
& \luTo_{[1]} && \ldTo & \\ 
&&Rj_*j^*(A\b)&&
\end{diagram} \end{tcolorbox}
\end{prop}
Consequently the {\em cohomology supported in $Z$} is:
\[ H^j_Z(S\b) := H^j(Z;i^!S\b)= H^j(X;Ri_*i^!S\b) =H^j(X,X-Z;S\b)\]
and the triangle gives the long exact cohomology sequence for the pair $(X,U)$.
For the (Borel-Moore)  sheaf of chains, $i^!C_{-r}$ gives the homology of a regular neighborhood of $Z$ in
$X$ which, for most nice spaces, will be homotopy equivalent to $Z$ itself.  The sheaf
$j^*(C^{-r})$ gives the Borel-Moore homology of $U$, which is the relative homology
$H_r(X,Z)$.  So this triangle gives the long exact sequence for the {\em homology} sequence  of the pair $(X,Z)$.

\subsection*{Caution}  For a complex of sheaves $A\b$ on $X$ and the inclusion $j_x:x \to X$ of a point, the stalk
cohomology of $A\b$ is $H^r_x(A\b) = H^r(j_x^*A\b)$ but the cohomology supported at $x$ is
$H^r_{\{x\}}(A\b) = H^r(j_x^!A\b)$.

\section{Stratifications}\index{stratification}\label{sec-stratifications}


\subsection{}\label{subsec-stratifications}
The plan is to decompose a reasonable space into a locally finite union of smooth manifolds
(called strata) which satisfy the axiom of the frontier:  the closure of each stratum should be
a union of lower dimensional strata.  If $Y \subset \overline{X}$ are strata we write $Y<X$.
But this should be done in a  locally trivial way.  In Whitney's example below,
\begin{figure}[!h]
\includegraphics[width=.3\linewidth]{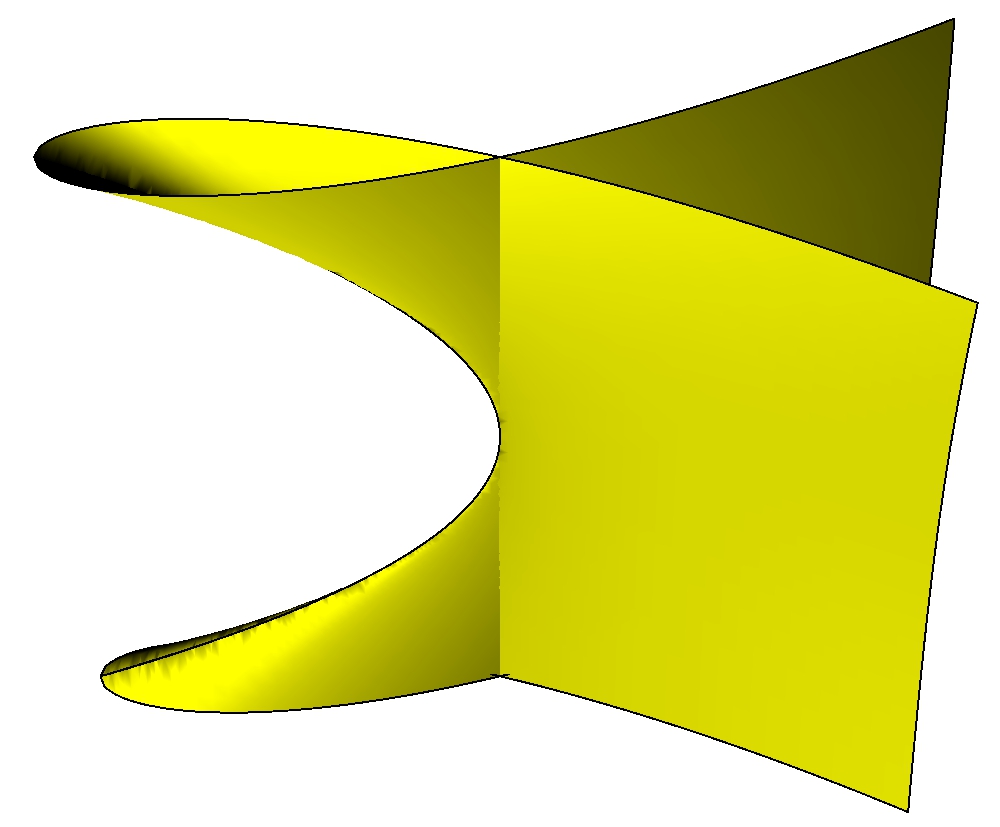}
\caption{Three strata are required.}
\end{figure}
it is not enough to divide this figure into 1- and 2-dimensional strata, even though this gives a decomposition
into smooth manifolds.  If the origin is not treated as another stratum then the stratification
fails to be ``locally trivial".  Whitney proposed a condition that identifies the origin as a
separate stratum in this example.  Let us say that a {stratification} of a closed subset $W$ of some smooth manifold $M$ is a locally finite decomposition $W = \coprod_{\alpha}S_{\alpha}$
into locally closed smooth submanifolds $S_{\alpha}\subset M$
(called strata) so as to satisfy the axiom of the frontier. 

For convenience let us say that strata are connected.
A proper continuous mapping $f:W_1 \to W_2$ between two stratified sets is a {\em stratified mapping}
\index{stratified map} 
if it takes strata to strata submerisvely, that is, if $X \subset W_1$ is a stratum then $f(X)$ is contained in a
stratum $Y$ of $W_2$ and $f:X \to Y$ is a smooth proper submersion.  This implies that $f(X) = Y$ and that
$f^{-1}(Y)$ is a union of strata.


\begin{tcolorbox}[colback=yellow!30!white]
\begin{defn}\index{stratification!Whitney}\index{Whitney condition}
Let $Y \subset \overline{X}$ be strata in a stratification of a closed set $W \subset M$.  
The pair $(X,Y)$ {\em satisfies Whitney's condition B}
at a point $y\in Y$ if the following holds.  Suppose that $x_1, x_2, \cdots \in X$ is a sequence
that converges to $y$, and suppose that $y_1, y_2, \cdots \in Y$ is a sequence that also
converges to $y$.  Suppose that (in some local coordinate system near $y$) the secant lines
$\ell_i = \overline{x_i,y_i}$ converge to some limiting line $\ell$.  Suppose that the tangent
planes $T_{x_i}X$ converge to some limiting plane $\tau$. Then $\ell \subset \tau$.
\end{defn}
\end{tcolorbox}
We say the pair $(X,Y)$ satisfies condition B if it does so at every point $y\in Y$.  
The decomposition into strata is a {\em Whitney stratification} if every pair of strata $Y<X$ 
satisfies condition B at every point in the smaller stratum $Y$.  
If condition B is satisfied, and if the tangent planes $T_{y_i}Y$ also converge to some limiting plane $\eta$ then $\eta \subset \tau$ as well, which Whitney had originally
proposed as an additional condition, which he called Condition A.

It turned out that Whitney's condition B was just the right condition to guarantee that a stratification
is locally trivial, but the verification involved the full development of stratification theory by 
Ren\'e Thom (\cite{Thom, Thom2}) and John Mather (\cite{Mather}).  The problem is that stratifications satisfying Condition B may
still exhibit certain pathologies, such as infinite spirals, so there is a very delicate balance between
proving that local triviality holds while avoiding a host of counterexamples to similar sounding
statements.

\subsection{}\label{subsec-normalslice}
Suppose $W\subset M$ has a stratification that satisfies condition B.  Let $Y$ be a stratum and let $y\in Y$.  Let $N_y \subset M$ be a {\em normal slice},\index{normal slice} that is, a smooth submanifold of dimension
$\dim(N_y) = \dim(M)-\dim(Y)$ that intersects $Y$ transversally in the single point $\{y\}$.  (cf. \S \ref{subsec-transversality})
 It follows from Condition A that
$N$ is also transverse to all strata $Z>Y$ in some neighborhood of $\{y\}$.
Define the {\em link of the stratum $Y$} \index{link}at the point $y \in Y$,
\[  L_Y(y,\epsilon) = (\partial B_{\epsilon}(y))\cap N_y \cap W\]
where $B_{\epsilon}(y)$ is a ball of radius $\epsilon$ (measured in some Riemannian metric
on $M$) centered at the point $y$.

\begin{thm} (R. Thom, J. Mather)\label{thm-ThomMather}
If $\epsilon$ is chosen sufficiently small then 
\begin{enumerate}
\item The closed set $L_Y(y,\epsilon)$ is stratified by its
intersection with the strata $Z$ of $W$ such that $Z>Y$.
\item This stratification satisfies condition B.
\item The stratified homeomorphism type of $L_Y(y,\epsilon)$ is independent of the choice of $N_y, \epsilon$,
and the Riemannian metric so we may denote it by $L_Y(y)$.
\item If the stratum $Y$ is connected then the stratified homeomorphism type of $L_Y(y)$ is
also independent of the point $y$ so we may denote it by $L_Y$.  (So it is notationally convenient to
assume henceforth that all strata are connected.)
\end{enumerate}
Moreover, the point $y$ has a {\em basic neighborhood} $U_y\subset W$ and a stratified
homeomorphism
\[ U_y \cong c^o(L_Y) \times B_1\]
where $c^o(L_Y)$ denotes\footnote{here, we interpret $c^o(L_Y) = \{y\}$ if $L_Y= \phi$ is empty, 
that is, if $y$ is a point in a maximal stratum.} the open cone on $L_Y$ (with its obvious stratification) and where
$B_1$ denotes the open ball of radius $1$ in  $\RR^{\dim(Y)}$.
\end{thm}

This homeomorphism preserves strata in the obvious way:  it takes \begin{enumerate}
\item $\{y\} \times B_1 \to Y \cap U$ with $\{y\} \times \{0\} \to \{y\}$
\item $(L_Y\cap X) \times (0,1) \times B_1 \to X \cap U$ for each stratum $X>Y$
\end{enumerate}
This result says that the set $W$ does not have infinitely many holes or infinitely much
topology as we approach the singular stratum $Y$ and it says that the normal structure
near $Y$ is locally trivial as we move around in $Y$.  In particular, the collection of
links $L_Y(y)$ form the fibers of a stratified fiber bundle over $Y$.
It also implies that (for any $r \ge 0$) the local homology $H_r(W,W-y;\ZZ)$ forms a local 
coefficient system on $Y$ with stalk
\[ H_r(W,W-y;\ZZ) \cong H_{r-\dim(Y)-1}(L_Y;\ZZ).\]

These statements are consequences of the {\em first isotopy lemma} of R. Thom,
\begin{lem}\label{lem-isotopy}\index{isotopy lemma}
Suppose $M,P$ are smooth manifolds and $f:M \to P$ is a smooth proper map.  Let
$W\subset M$ be a  Whitney stratified (closed) subset, and suppose that
 $f|X:X \to P$ is a submersion for each stratum $X$ of $W$.  Then
the mapping $f|W:W \to P$ is locally trivial:  each point $p \in P$ has a neighborhood basis consisting of
open sets $U_p$ for which there exists a stratum preserving homeomorphism
\[ f^{-1}(U_p) \cong U_p \times f^{-1}(p)\]
that is smooth on each stratum.
\end{lem}
\subsection{}
In fact, Thom \cite{Thom} and Mather \cite{Mather} proved that a Whitney stratified $W$ set admits a system of
{\em control data} consisiting of a triple $(T_Z, \pi_Z, \rho_Z)$ for each stratum $Z$,
where $T_Z$ is a neighborhood of $Z$ in $W$, where $\pi_Z:T_Z \to Z$ is a
``tubular projection", $\rho_Z:T_Z \to [0,\epsilon)$ is a ``tubular distance function" so
 that for each stratum $Y>Z$ the following holds:
 \begin{enumerate}
 \item $\pi_Z \pi_Y = \pi_Z$ in $T_Z \cap T_Y$
\item $\rho_Z \pi_Y = \rho_Z$ in $T_Z \cap T_Y$
\item the pair
$(\pi_Z, \rho_Z)|Y \cap TZ : Y \cap T_Z \to Z \times (0, \epsilon)$
 is a smooth submersion.\end{enumerate}
Theorem \ref{thm-ThomMather} then follows by applying the isotopy lemma \ref{lem-isotopy} to the mapping $(\pi_Z,\rho_Z)$.
\subsection{}

\begin{figure}[h!]\begin{center}
\centerline{\begin{picture}(200,80)
\put(40,40){\circle*{6}}
\linethickness{2pt}
\qbezier(40,40)(110,40)(190,40)
\linethickness{.5pt}
\qbezier(40,40)(60,60)(90,60)
\qbezier(90,60)(120,60)(140,60)
\qbezier(140,60)(170,60)(190,40)
\put(190,40){\circle*{6}}
\thicklines
\put(40,40){\circle{40}}
\put(190,40){\circle{40}}
\qbezier(40,40)(60,20)(90,20)
\qbezier(90,20)(120,20)(140,20)
\qbezier(140,20)(170,20)(190,40)
\put(30,45){$Z$}
\put(190,45){$Z'$}
\put(110,45){$Y$}
\put(25,5){$T_Z(\epsilon)$}
\put(100,5){$T_{Y}(\epsilon)$}
\put(180,5){$T_{Z'}(\epsilon)$}
\end{picture}}\caption{Tubular neighborhoods}\end{center}
\end{figure}
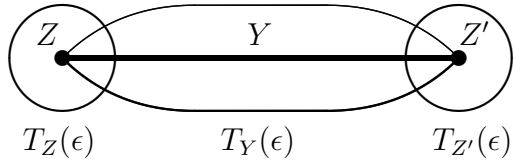
The local triviality statements are proven in \cite{Thom, Mather}, using the control data to
construct {\em controlled vector fields} \index{controlled vector field}\index{vector field, controlled}
whose flow traces out the local triviality of the stratification.
It is a delicate construction:
 a controlled vector field is not necessarily continuous but it has a
continuous flow.  Such a flow may ``spin'' around a small stratum, faster and faster for points
that are close to the small stratum.
\begin{figure}
{\includegraphics[height=2in]{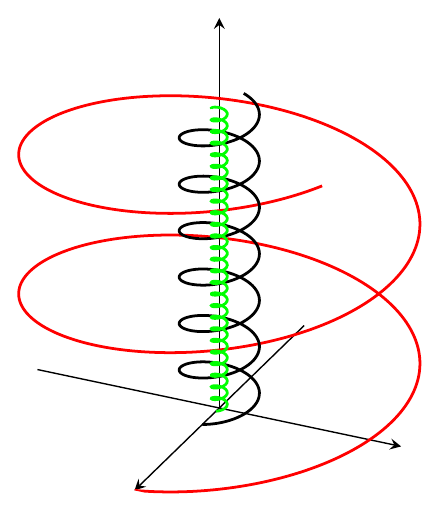}}
\caption{Controlled flow near the $Z$ axis}
\end{figure}
\subsection{}
\begin{figure}[h!]  \centering
\begin{tikzpicture}[scale=.5]
\begin{scope};
\pgfmathsetmacro{\ht}{10};
\draw[thick,red] (0,0) to [out = 0, in = 270] (3,\ht);
\draw[thick,red] (0,0) to [out=0, in = 90] (3,-\ht);
\draw[dotted] (0,0) to [out = 15, in = 270] (.5,\ht);
\draw[dotted] (0,0) to [out =0, in = 270] (5,\ht); 
\draw[dotted](0,0) to [out = -15, in = 90] (.5,-\ht);
\draw[dotted](0,0) to [out = 0, in = 90] (5,-\ht);
\draw [dotted](0,0) circle (4);

\pgfmathsetmacro{\arrowstop}{.9};
\pgfmathsetmacro{\arrowstart}{3}
\draw[dotted] (0,0) circle (8); 
\draw[->](150:\arrowstart) -- (150:\arrowstop);
\draw[->](90:\arrowstart)--(90:\arrowstop);
\draw[->](210:\arrowstart)--(210:\arrowstop);
\draw[->](270:\arrowstart)--(270:\arrowstop);
\end{scope};
\node at (-6,7) (a) {$T_R(2\epsilon)$} {};
\node at (-3,4)(b){$T_R(\epsilon)$};
\node at (6.5,9)(c){$T_S(\epsilon)$};
\node at (6.5,-9)(c){$T_S(\epsilon)$};

\draw [green, fill=green] (0,0) circle (3pt);
\draw[->] (1,10) -- (2,10);\draw[->](5,10)--(4,10);
\draw[->](1,9) -- (2,8.8);\draw[->](5,8.2)  --(4,8.4);

\draw[->](1,-10)--(2,-10);\draw[->](5,-10)--(4,-10);
\draw[->](1,-9) -- (2, -8.8);\draw[->](5, -8.2)  --(4, -8.4);


\draw[->,line width=.5mm,color=blue] (3.15,5.5)--(3.1,4);
\draw[->,line width=.5mm,color=blue] (3.15,-5.5)--(3.1,-4);
\draw[->](1.5,6.5) -- (2.5,4.5);
\draw[->](4.5,6)--(3.7,4.5);
\draw[->](1.5,-6.5) -- (2.5,-4.5);
\draw[->](4.5,-6)--(3.7,-4.5);
\end{tikzpicture}
\caption{{Tubular neighborhoods and retraction for $R<S$}}\label{fig-neighborhoods}
\end{figure}
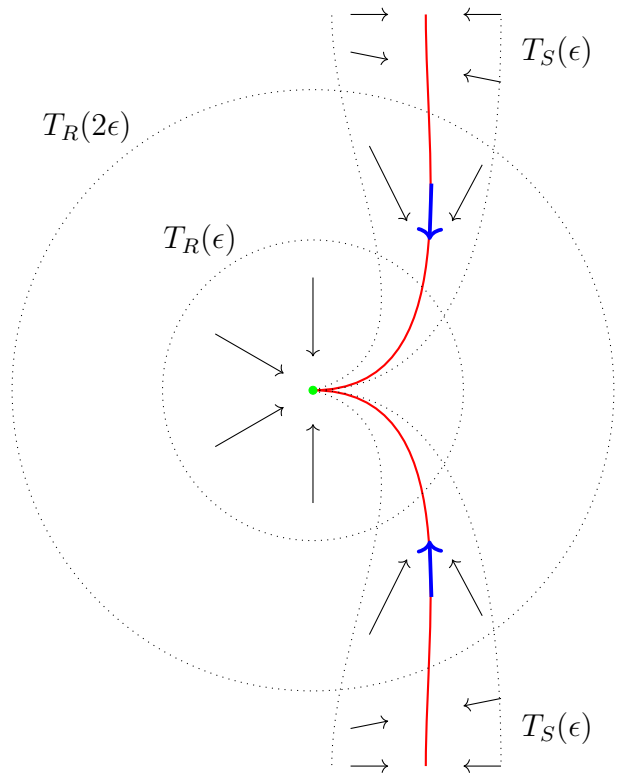

For each stratum $S$ in a Whitney stratified set the union $U_S = \bigcup_{R \le S}T_R$ form a neighborhood
of the closure  $\overline{S}$ and the inclusion $\overline{S} \to U_S$ is a homotopy equivalence.  In fact there exists
homotopy inverse, in fact a deformation retraction, 
$r_1:U_S \to \overline{S}$ that is the time = 1 value of a one parameter family of  maps
$r:U_S\times[0,1] \to U_S$ such that $r_0 =I$ (identity) and each $r_t$ ($t<1$) is stratum preserving and smooth on each stratum,
see Figure \ref{fig-neighborhoods}.
\cite{Hendricks, triangulations, MorseSheaves}
The mapping $r_1$ is a ``weak''  deformation retraction:  its restriction to $\overline{S}$ is the identity outside of a suitable
neighborhood $V_S$ of $\overline{S}-S$ (for example, $V_S = T_R(2\epsilon)$ in Figure \ref{fig-neighborhoods})
 but it is only homotopic to the identity within $V_S \cap \overline{S}$. 

\medskip

\subsection{}
Whitney himself outlined a procedure for proving that any closed subset $W$ of Euclidean 
space defined by analytic equations admits a Whitney stratification.  The idea is to start with the open, nonsingular part $W^0$ of $W$ as the ``top" stratum, and then to look at the set of points
in the singular set $\Sigma=W - W^0$ where condition B {\em fails}.  He proves that this is an analytic subset of codimension two, whose complement in $\Sigma$ is therefore the first singular stratum, $W^1$.  Now, carry both $W^0$ and $W^1$ along, looking at the set of points (in what
remains) where condition B fails, and continue in this way inductively.  Since Whitney's early work, many advances have been made in the subject (see \cite{Wall}).  Theorem \ref{thm-stratifications-exist} below is a partial summary of the work of many people, 
including \cite{Hardt}, \cite{Hironaka} \cite{Hi2}, \cite{Ma2}, \cite{Lo}, \cite{Whitney}, \cite{Thom3}, \cite{Bierstone}, \cite{Hendricks}, \cite{triangulations},
\cite{Johnson}, \cite{Lo2}.

\begin{tcolorbox}[colback=cyan!30!white]
\begin{thm}\label{thm-stratifications-exist}\index{Whitney condition}\index{stratification!Whitney}
The following sets admit Whitney stratifications:  real and complex algebraic varieties,
real and complex analytic varieties, semi-algebraic and semi-analytic varieties, subanaltyic
sets, and sets with o-minimal structure.  Given such an algebraic (resp. analytic etc.) variety $W$ and
a locally finite union $Z$ of algebraic (resp. analytic etc.) subvarieties, the stratification of $W$
can be chosen so that $Z$ is a union of strata.  Given an algebraic (resp. proper analytic etc.) 
mapping $f: W \to W'$ of algebraic (resp. analytic etc.) varieties, it is possible to find algebraic (resp. analytic)
Whitney stratifications of  $W, W'$ so that the mapping $f$ is stratified (see \S \ref{subsec-stratifications}).
\index{stratified map}   Whitney stratified sets can be triangulated by a
triangulation\index{triangulation} that is smooth on each stratum, such that
that the closure of each stratum is a subcomplex of the triangulation.
\end{thm}\end{tcolorbox}
A subanalytic set \index{subanalytic set} is the image under a projection (for example, a linear projection
$\RR^m \to \RR^n)$ of an analytic or a semi-analytic set.  
O-minimal structures\index{O-minimal structure} allow for certain non-analytic functions to be included in the
definition of the set.  Whitney stratifications also make sense for algebraic varieties
defined over fields of finite characteristic.  Given an algebraic mapping $f:W \to W'$
between complex algebraic varieties, it is not generally possible to choose triangulations
of $W, W'$ so that $f$ becomes a simplicial mapping.

\subsection{Pseudomanifolds and Poincar\'e duality}\index{pseudomanifold}
A pseudomanifold of dimension $n$ is a purely $n$ dimensional  (Whitney) stratified space that 
can be triangulated so that every $n-1$ dimensional simplex is a face of exactly two $n$-dimensional simplices.  
This implies that the $n-1$ dimensional simplices can be combined with the $n$-dimensional simplices to 
form an $n$-dimensional manifold, and that the remainder
(hence, the ``singularity set") has dimension $\le n-2$.  If this manifold is orientable then an
orientation\index{orientation}\index{fundamental class}
 defines a fundamental class $[W] \in H_n(W;\ZZ)$.  Cap product with the fundamental
class defines the Poincar\'e duality map $H^r(W;\ZZ) \to H_{n-r}(W;\ZZ)$ which is an isomorphism 
if $W$ is a manifold (or even a homology manifold) but which, in general, is not an isomorphism.

There is a sheaf-theoretic way to say this.  If $W$ is oriented and $n$-dimensional 
then a choice of orientation
determines a sheaf map $s:\uul{\ZZ}_W \to \uul{C}_n^{BM}$ from the constant sheaf to the sheaf of $n$-chains,
cf. \S \ref{subsec-chains}:  Let $U\subset M$ be an open set and choose a triangulation of $U$.  If $U$ is
sufficiently small it is possible to orient all of the $n$ dimensional simplices in $U$ in a compatible way, so that
$ \partial \sum_{\sigma\subset U} \sigma = 0$
where the sum is  over those oriented $n$-dimensional simplices in $U$.  Then set 
\[ s(m) = m \sum_{\sigma \subset U}m\sigma.\]
 (Recall that a PL chains are identified under subdivision.)
For example,  if $W$ is an $n$-dimensional homology manifold, that is, if
 $H_r(W,W-x) = 0$ for all $0 \le r <n$ and $H_n(W, W-x;\ZZ) = \ZZ$  then the composition
\[ \uul{\ZZ}_W[n] \to \uul{C}_n^{BM} \to \uul{C}^{\bullet}_W\]
 is a quasi-isomorphism.  This simple statement is the Poincar\'e duality theorem. \index{Poincar\'e duality}
  For, it says that
this quasi-isomorphism  induces an isomorphism on cohomology, that is,
\[ H^r(W;\ZZ) \cong H_{n-r}^{BM}(W;\ZZ)\]
and an isomorphism on cohomology with compact supports, that is,
\[ H^r_c(W;\ZZ) \cong H_{n-r}(W;\ZZ).\]
[Actually, from this point of view, the deep fact is that $H^i(W;k)$ and $H_i(W;k)$ are dual
over any field $k$, but this is not a fact about manifolds.  Rather, it is a fact about the
sheaf of chains.]

More generally if $W$ is a homology manifold but not necessarily orientable then the {\em orientation sheaf} $\mathcal O_W$ is
the local system whose stalk at $x\in W$ is the top local homology $H_n(W,W-x)$ and the mapping
$O_W \to \uul{C}_n^{BM} \to  \uul{C}^{\bullet}_W$ is a quasi-isomorphism.  So,
for any local coefficient system $\mathcal L$ on $W$ the mapping
$\mathcal L\otimes \mathcal O_W \to \uul{C}^{\bullet}(\mathcal L)$ is a quasi-isomorphism, giving an
isomorphism on cohomology,
\[ H^r(W; \mathcal L\otimes \mathcal O_W) \cong H_{n-r}^{BM}(W;\mathcal L)\]
and on cohomology with compact supports,
\[ H^r_c(W;\mathcal L\otimes \mathcal O_W) \cong H_{n-r}(W;\mathcal L).\]
So this quasi-isomorphism statement includes the  Poincar\'e
duality theorem for orientable and non-orientable manifolds, for non-compact manifolds, and
for manifolds with boundary, and with possibly nontrivial local coefficient systems.

\section{Constructible sheaves and functions}
\subsection{Constructible sheaves}\index{sheaf!constructible}\index{constructible sheaf}
Fix a Whitney stratification of a closed subset $W \subset M$ of some smooth manifold.
A sheaf $S$ (of abelian groups, or of $R$ modules) on $W$ is {\em constructible with 
respect to this stratification} if the restriction
$S|X$ to each stratum $X$ is a locally constant sheaf and the stalks $S_x$ are
finitely generated.  A complex of sheaves $A\b$ on
$W$ is {\em cohomologically constructible} \index{cohomologically constructible} with respect to this stratification if
its cohomology sheaves are bounded (that is, 
$\uul{H}^r(A\b) = 0$ for $|r|$ sufficiently large) and constructible.

If $W$ is an complex algebraic (resp. complex analytic, resp. real algebraic etc.) variety
then a complex of sheaves $A\b$ on $W$ is {\em algebraically construcible} (resp. 
analytically constructible, etc.) if its cohomology sheaves are bounded and constructible with respect to
some algebraic (resp. analytic etc.) Whitney stratification.

In each of these constructibility settings (that is, constructible with respect to a fixed stratification,
or algebraically constructible, etc.) the two constructions of the derived category make sense
(as the homotopy category of injective complexes, or as the category of complexes and equivalence
classes of roofs), which is then referred to as the bounded constructible derived category
and denoted $D^b_c(W)$.

\begin{tcolorbox}[colback=cyan!30!white]
\begin{thm} \label{thm-constructibility}
 Suppose that $W$ is a compact subset of some smooth manifold $M$ and
suppose that $A\b$ is a complex of sheaves that is cohomologically constructible
with respect to some Whitney stratification of $W$. Then the hypercohomology groups
$H^r(W, A\b)$ are finitely generated.  If $U_x$ is a basic neighborhood of $x\in W$
then the stalk cohomology $\uul{H}^r(A\b)_x $ coincides with the cohomology $H^r(U_x , A\b)$
for all $r$ (and so the limit over open sets containing $x$ is essentially constant).  If 
$i:Z \to W$ is a closed union of strata with open complement $j:U \to W$ then $Ri_*i^*(A\b)$
and $Rj_*j^*(A\b)$ are also cohomologically constructible.  If $A\b, B\b$ are cohomologically
constructible then so is $\RHom^{\bullet}(A\b, B\b)$.  If $f:W \to W'$ is a proper stratified mapping
and $A\b$ is CC on $W$ then $Rf_*(A\b)$ is CC on $W'$.
\end{thm}\end{tcolorbox}
\begin{proof}
Let $X$ be the top stratum and let $\Sigma = W-X$ be the singular set.  Let $U$ be the union of the tubular
neighborhoods of the strata in $\Sigma$.  Then $X - (X \cap U)$ is compact, and as $U$ shrinks these form
a sequence of diffeomorphic compact manifolds with corners that exhaust $X$.  If $\mathcal L$ is a local
system on $X$ then (since $W$ is compact)
\[ H^*_c(X, \mathcal L) = H^*(X - (X \cap U), \mathcal L)\]
is finitely generated.  From the spectral sequence for cohomology of a complex, the same holds for
$H^*_c(X, A\b)$ since the cohomology sheaves of $A\b$ are local systems on $X$.

Now consider the exact triangle
\[ Rj_!j^* A\b \to A\b \to Ri_*i^*(A\b) \to \cdots.\]
The cohomology of $Rj_!j^*A\b $ is $H^*_c(X, A\b)$ which is finitely generated as just shown.  The
complex $i^*(A\b)$ is a constructible complex on $\Sigma$, which has smaller dimension, so its
cohomology is finitely generated by induction.  The long exact sequence implies the cohomology
of $A\b$ is finitely generated.  

The stalk cohomology coincides with the cohomology of $U_x$ because the family of these basic
neighborhoods are cofinal in the set of all neighborhoods of $x$ but as they shrink there are stratified
isomorphisms $h:U_x \to U'_x$ with the inverse given by inclusion.  Since the cohomology sheaves of 
$A\b$ are locally constant on each
stratum there is a quasi-isomorphism $h_*(A\b) \to A\b$ which induces isomorphisms on cohomology.
In other words, $H^*(U_x, A\b)$ is independent of the choices, so the limit stabilizes.  
Constructibility of $Ri_*i^*A\b$ is obvious but constructibility of $Rj_*j^*A\b$ takes some work.
Here is the key point:

\begin{tcolorbox}[colback=cyan!30!white]
\begin{lem}
Let $A\b$ be a cohomologically constructible complex of sheaves on $W$.  
Let $Z \subset W$ be a closed subset with
complement $j:V \to W$.  Let $X$ be the largest stratum of $Z$.  Then the stalk cohomology
at $x \in X$ of $Rj_*j^*(A\b)$ is
\begin{equation}\label{eqn-jstar} \uul{H}{}^i(Rj_*j^*A\b)_x \cong  H^{i}(L_X(x),  A\b).\end{equation}
\end{lem}  \end{tcolorbox}
\noindent
Here, $L_X(x)$ denotes the link of the stratum $X$ in $W$ at the point $x$, cf \S \ref{subsec-normalslice}. 
The Lemma follows from the fact that the stalk cohomology is 
$\uul{H}^i(Rj_*j^*A\b)_x = H^i(j^{-1}(U_x \cap V); A\b) $ 
where $U_x \cong c^o(L_X(x)) \times B^{\dim(X)}$ is a basic open neighborhood.
 Since $X$ is the largest stratum of $Z$, we have:
\[U_x \cap V \cong L_X(x) \times (0,1) \times B^{\dim(X)},\]
and the cohomology sheaves of $A\b$ are constant in the Euclidean directions of this product.

Since $L_X(x)$ is compact this cohomology is finitely generated.  It is locally constant as $x\in X$
varies because the same is true of $L_X(x)$ and of the cohomology sheaves of $A\b|L_X(x)$.

The constructibility of $\RHom\b(A\b,B\b)$ follows from the spectral sequence for its stalk  cohomology in
terms of the stalk cohomology of $A\b$ and $B\b$.  
The last  statement in Theorem \ref{thm-constructibility} follows from the isotopy lemma, Lemma
\ref{lem-isotopy}.
\end{proof}

In summary, the derived category $D^b_c(W)$ of complexes whose cohomology sheaves are
bounded and constructibe forms a ``paradise", in the words of Verdier, who had assured us
(when we were writing IH II) that such a category, in which all these operations made sense,
and was closed under pullback, proper push forward, $\Hom$ and Verdier duality, did not exist.

\subsection{Attaching sheaves}
Let us examine the triangle for $Rj_*j^*A\b$ for $i \le \cod(X)-1$ and its stalk cohomology: 
\begin{equation}\label{eqn-beta}\begin{diagram}[size=2em]
Ri_*i^!(A\b)    && \rTo && A\b & & H^{i-\dim(X)}_c(U_x;A\b) && \rTo && \uul{H}^i(A\b)_x \\
& \luTo_{[1]}^{\beta} && \ldTo_{\alpha} &    &\ \text{ and }\ &  & \luTo_{[1]}^{\beta} && \ldTo_{\alpha} & \\
&& Rj_*j^*(A\b) && && && H^i(L_x, A\b) && 
\end{diagram}\end{equation}
The {\em attaching map} $\alpha$ goes from information ($\uul{H}^i(A\b)_x)$ living on the small stratum 
to information ($H^i(L_x,A\b)$) living completely in the larger strata and so it represents the degree to which
the sheaf $A\b$ is ``glued" across the strata.

\paragraph{\bf Exercise}  Suppose $W = X<U$ consists of two strata.  Let $B\b, C\b$ be sheaves on $X$ and on $U$
respectively and let $A\b = Ri_*(B\b) \oplus Rj_!(C\b)$ so that $A\b$ consists of just these two sheaves with no
relation between them.  Show that the attaching homomorphism $\alpha$ is zero.  For example, if $B\b = \ZZ_X$
and if $C\b = \ZZ_U$ then $A\b$ is a sheaf whose stalk at each point is $\ZZ$ however it is not the constant sheaf.
Show that if $A\b = \ZZ_W$ is the constant sheaf then the attaching homomorphism $\alpha$ is injective.

\subsection{Euler characteristic and constructible functions}\label{subsec-constructible}
\index{Euler characteristic}\index{constructible function}\index{function, constructible}
Let $k$ be a field.  In this section all sheaves are sheaves of $k$-vector spaces and all cohomology is
taken with coefficients in $k$.

Suppose $Y$ is a (closed) subcomplex of a finite simplicial complex $X$ and $A\b$ is a complex of sheaves of
$k$-vector spaces on $X$,
constructible with respect to the stratification by interiors of simplices.  For such a sheaf $A\b$ its stalk
Euler characteristic at a point $x \in X$ is
\[ \chi_x(A\b) = \sum_{i\ge 0}(-1)^i \dim H^i_x(A\b) = \sum_{i \ge 0}(-1)^i \dim (\uul{H}^i(A\b)_x)\]
and its global Euler characteristic is
\[\chi(X;A\b) = \sum_{i \ge 0} \dim H^i(X;A\b) = \sum_{\sigma}(-1)^{\dim{\sigma}}\chi_{\hat{\sigma}}(A\b)\]
(sum over simplices $\sigma$ in $X$, where $\hat\sigma$ denotes the barycenter) from which  the long exact sequence
\S \ref{subsec-subspaces} for the pair $(X,Y)$ implies that 
\begin{equation}\label{eqn-additivity}
\chi_c(U;A\b)=\chi(X,Y;A\b) = \sum_{\sigma^o \subset U}
(-1)^{\dim(\sigma)}\chi(\uul{H}{}^*(A\b)_{\hat{\sigma}}) = \sum_{\sigma^o \subset U} \chi_c(A\b|\sigma^o).\end{equation}
The sum is over simplices $\sigma$ whose interior $\sigma^o$ is contained in $U$.   Here, 
$\chi_c$ denotes the Euler characteristic with compact supports,
$\chi_c(S) = \sum_i \dim H^i_c(S;k)$.  It coincides with the Borel-Moore homology Euler characteristic. 
Thus equation (\ref{eqn-additivity}) says that {\em $\chi_c$
is  additive}:  it is the sum over the interiors $\sigma^o$ of simplices of the Euler characteristic with compact
supports of each $\sigma^o$.

As in the previous chapter, fix a Whitney stratification of a closed subset $W \subset M$.  (For convenience, as mentioned in
\S \ref{thm-ThomMather}, we assume the strata are connected.)   A function $f:W \to k$ is
{\em constructible} if it is constant on each stratum.  The Euler characteristic of a constructible function $f$
is defined to be $\chi(W;f) = \sum_S\chi_c(S)f(x_S)$ where the sum is over strata $S$ of $W$, and $x_S\in S$.

If $A\b \in D^b_c(W)$ then its stalk Euler characteristic is a constructible function $x\mapsto \chi_x(A\b)$
which we denote by $\chi_{\square}(A\b)$.
If the total cohomology $H^*(W;A\b)$ is finite dimensional then (using either a spectral
sequence argument or the above commments about additivity applied to a triangulation of $W$), its 
 Euler characteristic $\chi(W;A\b)$ is given by
\[ \chi(W;A\b) = \sum_i(-1)^i \dim H^i(W;A\b)= \chi(W;\chi_{\square}(A\b)) = \sum_S \chi_c(S) \chi_{\square}(A\b).\]

\subsection{Calculus of constructible functions}  \index{constructible function}
If $\pi:W \to Y$ is a proper, weakly stratified mapping\index{weakly stratified map}  (meaning that it takes each
stratum of $W$ submersively to a stratum of $Y$) and if $f:W \to k$ is a constructible function then its
{\em pushforward} $\pi_*(f)$ is the constructible function $\pi_*(f)(y) = \chi_c(\pi^{-1}(y);f)$.  Consequently, if $A\b \in D^b_c(W)$
then $\chi_{\square}(R\pi_*A\b) = \pi_*\chi_{\square}(A\b)$.  If $\pi:X \to \{pt\}$ it is common to write (\cite{Viro, Ghrist})
\[ \pi_*(f) = \int_Xf(x) d\chi(x).\]
If $W$ is compact and $h:W \to \RR$ is a Morse function (cf. \S \ref{subsec-Morsefunctions} below) with
critical points $p_i$ ($1 \le i \le r$) and critical values  $v_1<\cdots< v_r$ then 
$\int_Xf(x)d\chi(x)$ may be computed using Morse theory: The constructible function $h_*(f)$ vanishes at 
$t=-\infty$ and it  only changes at the critical values of $h$.  By Theorem \ref{thm-SMT}
below, the contribution at each
critical value $v_i$ is localized near the critical point $p_i \in W$ and is given by the Euler characteristic
$\chi_i=\chi(V_+,V_-)$ of the local Morse data (\ref{eqn-localMorse}) so that for any $t \in (v_{s}, v_{s+1})$,  
\[ h(t) = \sum_{i=1}^s \chi_i \text{ with } h(\infty) = \int_Xf(x)d\chi(x).\]

A less obvious fact (\cite{KS, Schapira3}) is that Verdier duality \index{Verdier duality}\index{duality!Verdier}
(see \S \ref{subsec-duality}) passes to
constructible functions:  If $f$ is a constructible function on $W$ then there is a dual function $D(f)$ with the
property that for any sheaf $A\b \in D^b_c(W)$ with Verdier dual ${\bf D}(A\b)$,
\begin{equation}\label{eqn-dualfunction}
D(\chi_{\square}(A\b)) = \chi_{\square}({\mathbf D}(A\b)).\end{equation}  
The dual of a general constructible function $f$ is given by
\[ D(f)(x) = \lim_{\epsilon\to 0} \chi_c(B_{\epsilon}^o(x) \cap W; A\b) = \chi(j_x^!;A\b)\]
where $j_x:x \to W$ is the inclusion, and $B_{\epsilon}^o(x)$ denotes an open ball of radius $\epsilon$ centered at $x$.   
As in Theorem
\ref{thm-ThomMather} this limit stabilizes for $\epsilon$ sufficiently small.  Equation (\ref{eqn-dualfunction})
then follows from the fact that $j_x^*({\bf D}(A\b)) = {\bf D}(j_x^!A\b)$.  (Consequently $DD(f) = f$ for any
constructible function $f$, and $D\pi_*f = \pi_*D(f)$ for any $\pi:W \to Y$ as above.) 

\subsection{Applications of constructible functions}
The $k$-th {\em Hadwiger invariant} or
{\em intrinsic volume} \index{invariant, Hadwiger}\index{Hadwiger invariant}\index{intrinsic volume}
(see, for example, \cite{Rota}) of a compact subanalytic
subset $W \subset \RR^n$ is given (\cite{Ghrist2,Viro}) by
\[ \mu_k(W) = c(n,k)\int_{\Gr_{n-k}(\RR^n)}\int_{\RR^n/P} (\pi_P)_*({\bf 1}_W) dx_PdP\]
where $c(n,k)\in\RR$ are specific constants (see \cite{Ghrist2}), ${\bf 1}_W$ is the characteristic function of $W$,
$\pi_P:\RR^n \to \RR^n/P$ is the projection, the Grassmannian
$\Gr_{n-k}(\RR^n)$ is provided with the volume one invariant measure $dP$ and $dx_P$ denotes Lebesgue
measure on $\RR^n/P$.  The intrinsic volume $\mu_k(f)$ of
a constructible function $f$ is similarly defined. (See \cite{Ghrist} and references therein).
    
In \cite{ChernClasses}\index{Chern class} the Chern class $c_*(f) \in H_*(X)$ of a constructible function $f$ 
on a complex algebraic
variety $X$ is constructed so as to have the property that $c_*(\pi_*f) = \pi_*(c_*(f))$ for any proper
algebraic map $\pi:X \to Y$.  The class $c_*({\bf 1}_X)$ was later (\cite{Brasselet}) identified with the Schwartz Chern class
(\cite{Schwartz1, Schwartz2}).

\section{Sheaves and Morse theory} 

\subsection{Conormal vectors}\index{conormal vector}
Throughout this section $W$ denotes a Whitney stratified closed subset of a smooth manifold $M$.
 Let $X$ be a
stratum of $W$ and let $p \in X$.  A cotangent vector $\xi \in T_pM$ is
said to be {\em conormal} to $X$ if its restriction vanishes:  $\xi|T_pX = 0$. 
The collection of all conormal vectors to $X$ in $M$ is denoted $T^*_XM$.  It is a
smooth conical Lagrangian locally closed submanifold of $T^*M$.  An orientation of $M$ induces \cite{VilonenSchmidt}
an orientation of every $T^*_XM$ (whether or not $X$ is orientable).

A subspace $\tau \subset T_pM$ will be said to be a {\em limit
of tanget spaces from $W$} if there is a stratum $Y > X$ ($Y \ne X$) 
and a sequence of points $y_i \in Y$, $y_i \to p$ such that
the tangent spaces $T_{y_i}Y$ converge to $\tau$.   A conormal vector $\xi \in T^*_XM$ at $p$
is {\em nondegenerate} if $\xi(\tau) \ne 0$ for every limit $\tau \subset T_pM$
of tangent spaces from larger strata $Y>X$.  The set of nondegenerate \index{conormal vector}
conormal vectors is denoted $\Lambda_X$.  Evidently,
\begin{equation}\label{eqn-nondegenconormal}
\Lambda_X = T^*_XM - \bigcup_{Y>X} \overline{T^*_YM}\end{equation}
where the union is over all strata $Y>X$ (including the case $Y = M-W$ because
$T^*_MM$ is the zero section, and elements of $\Lambda_X$ are necessarily nonzero).  
If $M$ is an analytic manifold and $W$ is a subanalytically stratified subanalytic subset then each
$\Lambda_X$ has finitely many connected components.

\subsection{Morse functions}\label{subsec-Morsefunctions}\index{Morse!function}
 Let $f:M \to \RR$ be a smooth function and  let $W_{\le a} = W \cap f^{-1}((-\infty,a])$ and similarly for $W_{[a,b]}$, etc.
The main idea in Morse theory is to build the cohomology $H^*(W)$ from the long exact cohomology sequences of
the pairs $(W_{\le v+\epsilon},W_{\le v-\epsilon})$ associated to ``critical values'' $v$ of $f|W$.  The function $f$ is
{\em perfect} if the connecting homomorphisms in these long exact sequences vanish.
Since $W$ is a singular space these notions need to be made precise.

 A {\em critical point of $f|W$} is a point $p\in X$ in some
stratum $X$ such that $df(p)|T_pX = 0$, that is, $\lambda=df(p) \in T^*_XM$.  
(In particular, every zero dimensional stratum
is a critical point.)  The number $v = f(p)$ is said to be a critical value.  It is an {\em isolated critical value} of $f|W$ if
no other critical point $q\in W$ of $f|W$ has $v = f(q)$.

We say that $f$ is a Morse function
for $W$ (cf.~\cite{Lazzeri}) if \begin{itemize}
\item its restriction to $W$ is proper 
\item $f|X$ has isolated nondegenerate critical points for each stratum $X$,
\item at each critical point $p\in X$ the covector $\lambda = df(p) \in \Lambda_X$ is nondegenerate,
that is, $df(p)(\tau) \ne 0$ for every limit of tangent spaces $\tau\subset T_pM$
from larger strata $Y>X$.
\end{itemize}
The set of Morse functions is open and dense in the space of proper smooth mappings with the
Whitney $C^{\infty}$ topology (\cite{Pignoni, Orro}).

\subsection{}  Suppose $p\in X \subset W$ is an isolated nondegenerate critical point of $f:M \to \RR$ as above.
Let $\widetilde{N} \subset M$ be a smooth submanifold intersecting $X$ transversally at the single
point $\{p\}$.  Then $N =\widetilde{N}\cap W$ is a {\em normal slice} to the stratum $X$ at the point $p$ \index{normal slice}
(cf. \S \ref{subsec-normalslice}).
It inherits a stratification from that of $M$ and, at least in some neighborhood of $\{p\}$.

Suppose also that $A\b$ is a complex of sheaves that is 
 (cohomologically) constructible with respect
to the stratification of $W$.  (Or equivalently, $A\b$ is a complex of sheaves on $M$, constructible with respect
to the stratification defined by $W$ and supported on $W$.)   Stratified Morse theory ``reduces'' the problem of computing
the change in cohomology $H^*(W_{\le t};A\b)$ as we pass a critical value $t=v$
to the case of a (germ of a) stratified space $H^*(N_{\le t};A\b)$ with a zero dimensional stratum.  Theorem
\ref{theorem-MorseSheaf} is a consequence of Theorem \ref{thm-SMT} below which is proven
in \cite{SMT0,SMT} using stratum preserving deformations arising from various applications of
Thom's first isotopy lemma \ref{lem-isotopy}.

\begin{tcolorbox}[colback=cyan!30!white]
\begin{thm}{}\cite{SMT} \label{theorem-MorseSheaf}
 Let $f:M \to \RR$ be a smooth function and
suppose $f|W$ is proper.   Suppose $X\subset W$ is a stratum and that $p \in X$ is an isolated nondegenerate
critical point of $f$  (in the above, stratified sense) with isolated critical value $v = f(p)\in (a,b)$.  Let
$\lambda$ denote the  Morse index of $f|X$ at $p$.
 Suppose there are
no additional critical values of $f|W$ in the interval $[a,b]$.  If $0<\delta<<\epsilon$
are chosen sufficiently small (cf. Lemma \ref{lem-epsilon}) 
  then there is a natural isomorphism of {\em Morse groups}\index{Morse!group}
\begin{equation}\label{eqn-Morsegroup} H^r(W_{\le b},W_{\le a}; A\b) \cong H^{r-\lambda}
\left(N_{[v-\delta,v+\delta]}, N_{v-\delta};A\b|N
\right).\end{equation}
For fixed complex $A\b$ the Morse groups
 \[ M^t(\xi) = H^t\left(N_{[v-\delta,v+\delta]}, N_{v-\delta};A\b|N\right)\] 
 depend only on the nondegenerate covector $\xi = df(p)\in T^*_{X,p}M$.
\end{thm}\end{tcolorbox}
For any stratum $X$ and any given nondegenerate covector $\xi \in  T^*_{X,p}M$ it is possible to find a 
smooth function $f:M \to \RR$
with a nondegenerate critical point at $p$ such that $f(p)=0$, $df(p)=\xi$, with no other critical points of critical
value $0$ and so that the Hessian of $f|X$ at $p$ is
positive definite.  Then $\lambda = 0$ so for such a function, and $\epsilon>0$ sufficiently small,
\[ M^t(\xi) = H^t(W_{\le \epsilon},W_{\le -\epsilon}).\]
In the classical smooth Morse theory the Morse group is nonzero only in the single degree $r = \lambda$ where
it is the integers.  In the singular case the Morse group (\ref{eqn-Morsegroup}) may be nonzero in many degrees.
However if $W$ is a complex analytic variety and $A\b$ is a perverse sheaf then (cf. Theorem \ref{thm-onedegree}) 
the Morse group again lives in a single degree $r= c+\lambda$ where $c =\cod_{\CC}(X)$ is the complex
codimension of the stratum containing the critical point $p$.

There are several situations in which the hypothesis ``isolated nondegenerate critical point'' may be relaxed.  See
Theorems \ref{thm-Lefschetz}, \ref{thm-Kirwan}.  The meaning of $0<\delta << \epsilon$ is explained in the following:

\begin{lem}\label{lem-epsilon}
  Fix a Riemannian metric on $M$ and let $B_{\epsilon}(p)$ denote the closed ball of radius $\epsilon$ centered at $p$.
Given $p \in X \subset W \subset M \overset{f}{\longrightarrow}\RR$ and $ N$ as above, there exists $\epsilon_0$ so that for all $\epsilon \le \epsilon_0$ the following holds:
\begin{enumerate}
\item[{\rm({\color{red}*})}] the boundary sphere $\partial B_{\epsilon}(p)$ is transverse to every stratum
of $W$ and of $T \cap W$.
\end{enumerate}
Given such a choice $\epsilon$ there exists $\delta_0=\delta_0(\epsilon) >0$ such that for all 
$0<\delta \le \delta_0$ the following holds:
\begin{enumerate}
\item[{\rm({\color{red}**})}]  $f|N$  has no critical points on
any stratum of 
$N\cap f^{-1}[v-\delta,v+\delta]$
other than $\{p\}$.  
\end{enumerate}

\end{lem}
 For such a choice of $\epsilon,\delta$ we write  $ 0 < \delta << \epsilon$.  The set of possible choices for
$\epsilon,\delta$ will be an open region in the $(\epsilon,\delta)$ plane like the following
shaded area:

\begin{figure}[H]\centering
\label{eqn-graph}
\begin{tikzpicture}[scale = 0.6]
   \begin{axis}[ xmin=0, xmax=1, ymin=0, ymax=.8,
    axis lines = left, label style={font=\Large},
     ylabel={$\delta$}, 
every axis x label/.style={
    at={(ticklabel* cs:1.05)},
    anchor=west,},
ticks=none, ylabel style = {rotate = -90},
     scale=1, restrict y to domain=-1:1]
     \addplot[domain = 0:0.9, fill=gray] {x^3}\closedcycle;
   \end{axis}\node at (7.5,0) {$\epsilon$};
\end{tikzpicture}\caption{{$\delta << \epsilon$ region}}\label{fig-graph} \end{figure}
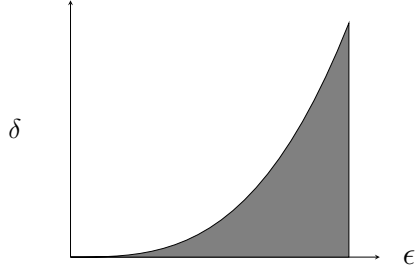
In this situation the main theorem of stratified Morse theory (\cite{SMT}) says:
\begin{tcolorbox}[colback=cyan!30!white]
\begin{thm}\label{thm-SMT} With $0< \delta << \epsilon$ chosen as above, the stratified homeomorphism
type of the pair  
$\left(N_{[v-\delta,v+\delta]}, N_{v-\delta}\right)$ is independent of $\epsilon, \delta$ and
there is a stratum preserving
homeomorphism, smooth on each stratum, between $B_{\epsilon}(p)\cap W_{\le v+\delta}$ and the adjunction space
\begin{equation}\label{eqn-Morse-NxT}
 \left[B_{\epsilon}(p)\cap W_{\le v-\delta}\right] 
\cup (D^{\lambda}, \partial D^{\lambda})
\times D^{s-\lambda} \times 
\left(N_{[v-\delta,v+\delta]}, N_{v-\delta}\right) \end{equation}
\end{thm}\end{tcolorbox}
\noindent
in which case we say that the {\em local Morse data} 
\begin{equation}\label{eqn-localMorse}
(V_+,V_-)=B_{\epsilon}(p)\cap(W_{ [v-\delta,v+\delta]}, W_{v-\delta})\end{equation}
 is homeomorphic to the
product of the {\em tangential Morse data} $(D^{\lambda}, \partial D^{\lambda}) \times D^{s-\lambda}$ and the
{\em normal Morse data} $(N_{[v-\delta,v+\delta]}, N_{v-\delta})$.  (Here, $s = \dim(X)$.)  Moreover:
\begin{tcolorbox}[colback=cyan!30!white]
If $f^{-1}([-\delta,\delta])$ contains no critical points (on any stratum)
other than $p$ then the intersection with $B_{\epsilon}(p)$ may be
removed in the above statement.
\end{tcolorbox}
\subsection{}\label{subsec-independent}
 The stratified homeomorphism type of each of these stratified spaces is indepenent of the choice of
Riemannian metric, $\epsilon$ and $\delta$ (provided $0<\delta << \epsilon$).  To some degree it is
even independent of the Morse function $f$.  If $g:M \to \RR$ is a another smooth
function with a  nondegenerate critical point at $p$ with the same Morse index $\lambda$, and if the covector $dg(p)$ lies in the
same connected component of the set of nondegenerate covectors at $p$ as the covector $df(p)$ then
the Morse data for $g$ at $p$ is stratified-homeomorphic to the Morse data for $f$ at $p$ and in particular
(under the same hypotheses as in Theorem \ref{theorem-MorseSheaf}) the cohomology $H^r(W_{\le b}, W_{\le a};A\b)$ has
the same description when $W_{\le t}$ is defined using either the function $f$ or the function $g$.  If $df(p)$ is a
degenerate covector then Morse data for $g$ may differ from Morse data for $f$ even if $df(p)=dg(p)$.

\subsection{Sheaf theoretic expression}\label{subsec-Rgamma}
Kashiwara and Schapira  \cite{KS} \S 5.1, \S 5.4 and Sch\'{u}ermann \cite{Schuermann} prefer a
 sheaf-theoretic expression for the Morse group.\index{Morse!group}  Let 
\[ x_0 \in X \subset W \subset M \overset{f}{\longrightarrow} \RR\]
as in Theorem \ref{theorem-MorseSheaf} above
 and suppose $f(x_0) = 0$ is an isolated critical value of $f|W$.  Let
$A\b\in D^b_c(W)$ be a constructible complex of sheaves.  Set
\[ Z = \left\{ x \in W| f(x)\ge 0\right\}\] with inclusion $i:Z \to W$ and let
$S_Z^{\bullet}=i_Z^!A\b=R\Gamma_ZA\b$ denote the sheaf obtained from $A\b$ 
with sections supported in $Z$, cf. \S \ref{subsec-subspaces}. 
Let $U = B_{\epsilon}(x_0)\cap W$ be a basic neighborhood of the critical point $x_0$.
If $a<0<b$ and $[a,b]$ contains no critical values other than $0$ then for $0<\delta << \epsilon$
Thom's first isotopy lemma (Lemma \ref{lem-isotopy}) gives isomorphisms of the
Morse groups:
\begin{align*}
 H^r(W_{\le b}, W_{\le a}; A\b) &\cong H^r(U_{\le \delta}, U_{\le -\delta}; A\b)\\
&\cong H^r(U_{\le \delta}, U_{<0}; A\b)\\
&\cong H^r(U_{\le \delta}; i_Z^!A\b)\\
&\cong \uul{H}{}^{\mathbf r}_{x_0}(i_Z^!A\b) = \uul{\mathbf H}{}^{r}_{x_0}(R\Gamma_ZA\b) \end{align*}
since the stalk cohomology is the limit as $\epsilon,\delta \to 0$ provided they
satisfy ({\color{red}*})  and ({\color{red}**}) of Lemma \ref{lem-epsilon}.
If we apply the main theorem in 
stratified Morse theory, this Morse group is identified with
\[M^{r-\lambda}(\xi;A\b) = \uul{\mathbf{ H}}{}^{{r-\lambda}}_{x_0}(i_{Z \cap N}^!(A\b|N))\]
where $\xi = df(x_0)\in \Lambda_X\subset T^*_XM$ and where $N=T\cap W\cap B_{\epsilon}(p)$
 denotes the normal slice to the stratum $X$.
Except for the shift $\lambda$ (which comes from the tangential Morse data), this expression depends only on the nondegenerate covector $\xi=df(x_0)$, cf. \S \ref{subsec-independent}.  (If
 a covector $\xi\in T^*_XM$ is degenerate then the Morse group $M^t(\xi;A\b)$ is not necessarily well defined.) 
If we choose the function $f$ so that the Hessian
of $f|X$ is positive definite at $x_0$ then $\lambda = 0$ and so the Morse group is canonically identified with the derived
functor of sections with support:
 \begin{tcolorbox}[colback=cyan!30!white]
\[ M^t(\xi;A\b) = \uul{\mathbf {H}}{}^t_{x_0}(i_Z^!A\b) = H^t(R\Gamma_ZA^{\bullet}_{x_0}).\]
\end{tcolorbox}

\subsection{Singular support} \index{singular support} \index{support, singular}
 Suppose $W \subset M$ as above, and
$A\b \in D^b_c(W)$ is constructible with respect to a given stratification of $W$.
From the preceding paragraph we see that Morse groups $M^*(\xi;A\b)$ for  $A\b$ are assigned to every
nondegenerate covector $\xi \in \bigcup_X \Lambda_X$ (union over strata $X$ of $W$). 
 Such a covector is not in the singular support
$SS(A\b)$ if all the Morse groups vanish:  $M^t(\xi;A\b) = 0, \forall t$. {\em  Thus, $SS(A\b)$ is the closure of
the set of nondegenerate covectors $\xi \in \bigcup_XT_XM$ for which some Morse group is nonzero. }

If $\xi\in T_p^*M$ is a nondegenerate covector for $W$, the Morse groups $M^*(\xi,A\b)$  form the 
stalk cohomology of a complex of sheaves $\mu_M(A\b)$, the {\em microlocalization}\index{micrlocalization}
 of $A\b$, on $T^*M$, which
is supported on $SS(A\b)$ \cite{KS}. 
 
If a complex of sheaves $A\b$ is not necessarily constructible then the singular support may be defined
in the same manner (cf. \cite{KS}) but it does not necessarily have the simple interpretation described here.

\subsection{Characteristic cycle}\index{characteristic cycle}\index{cycle, characteristic}
  Suppose that $M^n$ is an oriented analytic manifold and $W\subset M$ is a subanalytically
stratified subanalytic subset.  Let $A\b \in D^b_c(W)$.  In \cite{IndexTheorem} M. Kashiwara associates a
Borel-Moore chain $CC(A\b)$ in $C_n^{BM}T^*M$ such that $\partial CC(A\b) = 0$.  It is called the characteristic
cycle of $A\b$.  It can be described Morse-theoretically (see \cite{VilonenSchmidt}) as follows.  
Let $X$ be a stratum of $W$ and let $\Lambda_X \subset T^*_XM$ be the collection of nondegenerate
covectors for $X$.  It is a finite union $\Lambda_X = \cup_{\alpha}\Lambda_{X,\alpha}$ of connected
components.  Each component has dimension $n$ and a canonical orientation, hence a fundamental
class in Borel Moore homology:
\[ \left[ \Lambda_{X,\alpha} \right]   \in C_n^{BM}(\Lambda_{X,\alpha}) \to C_n^{BM}T^*M.\]
For  each $\alpha$  the Morse groups $M^t(\xi;A\b)$ are independent of the choice
of $\xi \in \Lambda_{X,\alpha}$ so the Euler characteristic $\chi_{X,\alpha} = \sum(-1)^t \dim M^t(\xi;A\b)$ is
well defined.  Then the characteristic cycle is the sum over strata $X$ with multiplicity $\chi_{X,\alpha}$ of
these fundamental classes:
\[ CC(A\b) = \sum_X \sum_{\alpha} \chi_{X,\alpha}\left[\Lambda_{X,\alpha}\right] .\]
Kashiwara proves \cite{IndexTheorem} that $\partial CC(A\b) = 0$ and the {\em Kashiwara index theorem}
\index{Kashiwara index theorem}\index{index theorem} \cite{IndexTheorem, KS, VilonenSchmidt}
says if $W$ is compact then the Euler characteristic
of $A\b$ is the intersection number
\[ \chi(W;A\b)= T^*_MM \cap CC(A\b)\]
of the characteristic cycle with the zero section.

\subsection{Hyperbolic Lefschetz numbers} \label{subsec-hyperbolic}
\index{hyperbolic!self map}\index{Lefschetz!number}

Let $f:W \to W$ be a subanalytic self map
defined on a subanalytically stratified subanalytic set $W$. A fixed point $x$ of $f$ is {\em contracting} if
all nearby points are mapped closer to $x$, and it is expanding if nearby points are mapped farther away
from $x$.  In some cases\footnote{Hecke correspondences are hyperbolic \cite{TTF}.
The self map $\CC \to \CC$  given by $z \mapsto z/(1+z)$ $(z \ne -1)$ has
a non-hyperbolic isolated fixed point at $z=0$. } it is possible to identify contracting and expanding directions:
\begin{tcolorbox}[colback=yellow!30!white]
  A connected component $V$ of the fixed
point set of $f$ is said to be {\em hyperbolic}
if there is a neighborhood $U\subset W$ of $V$ and
a subanalytic mapping $r=(r_1,r_2):U \to \RR_{\ge 0} \times \RR_{\ge 0}$ 
so that $r^{-1}(0) = V$ and so that
$r_1(f(x)) \ge r_1(x)$ and $r_2(f(x)) \le r_2(x)$ for all $x \in U$.  \end{tcolorbox}
Hyperbolic behavior of $f:W \to W$ is illustrated in the following diagram.
(Flow lines connecting
$r(x)$ and $r(f(x))$ do not exist in general).
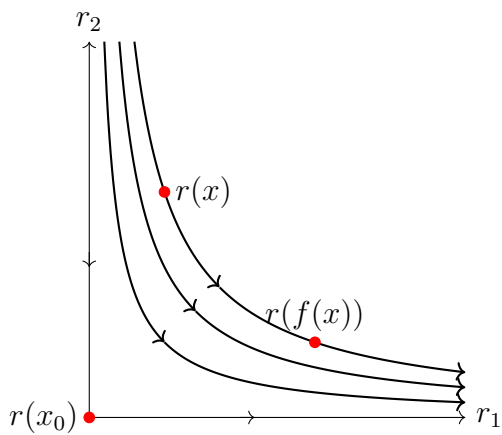
\begin{figure}[H]\centering
\begin{tikzpicture}
      \draw[->] (0,0) -- (5,0) node[right] {$r_1$};
      \draw[->] (0,0) -- (0,5) node[above] {$r_2$};
      \draw[thick,scale=1.0,domain=0.2:1,smooth,variable=\t,black,->] plot({\t},{1/\t});
     \draw[thick,scale=1.0,domain=0.2:1,smooth,variable=\s,black,<-] plot({1/\s},{\s});

      \draw[thick,scale=1.0,domain=0.2:.707,smooth,variable=\t,black,->] plot({2*\t},{1/\t});
     \draw[thick,scale=1.0,domain=0.2:.707,smooth,variable=\s,black,<-] plot({1/\s},{2*\s});
 \draw[thick,scale=1.0,domain=0.2:.577,smooth,variable=\t,black,->] plot({3*\t},{1/\t});
     \draw[thick,scale=1.0,domain=0.2:.577,smooth,variable=\s,black,<-] plot({1/\s},{3*\s});

\draw [->,black] (2,0) -- (2.2,0);
\draw [->,black] (0,2.2) -- (0,2);

\node at  (1,3) (T) {};
\draw [red, fill=red] (T) circle (2pt);
\node[right] at (T) {$r(x)$}; 

\node at (3,1) (S) {};
\draw [red, fill = red] (S) circle (2pt);
\node [above] at (S) {$r(f(x))$};

\node at (0,0) (O) {};
\node [left] at (O) {$r(x_0)$};
\draw[red, fill = red] (O) circle (2pt);
    \end{tikzpicture}
\caption{{Behavior near a hyperbolic fixed point}}\label{fig-graph2} \end{figure}

Suppose $V$ is a hyperbolic fixed point component.  Let $V^+ = r^{-1}(Y)$ and $V^- = r^{-1}(X)$
where $X, Y $ denote the $X$ and $Y$ axes in $\RR_{\ge 0} \times \RR_{\ge 0}$ with inclusions
$j^{\pm}, h^{\pm}$:
\begin{equation}\label{eqn-Vmaps}
\begin{diagram}[size=2em]
V_r & \rTo^{j_r^{\pm}} & V_r^{\pm} & \rTo^{h_r^{\pm}} & W.
\end{diagram}\end{equation}
Let $A\b \in D^b_c(W)$ be a constructible complex of sheaves.
A morphism $\Phi:f^*(A\b) \to A\b$ is called a {\em lift} of $f$ to  $A\b \in D^b_c(W)$.
Such a lift induces a homomorphism $\Phi_*:H^i(W; A\b) \to H^i(W;A\b)$ and defines the Lefschetz number
\[ \text{Lef}(f,A\b) = \sum_{i\in\ZZ} (-1)^i \Tr\left(\Phi_*:H^i \to H^i    \right).\]
 If $V$ is a hyperbolic connected component of the fixed point set define the following {\em hyperbolic
 localizations} or restrictions of $A\b$ to $V$: 
\[A_V^{!*} = (j^+)^!(h^+)^*A\b\ \text{ and }\  A_V^{*!} = (j^-)^*(h^-)^!A\b.\]
Then $H^*(V; A_V^{!*})$ is the cohomology of a neighborhood of $V$ with closed supports in the directions
``flowing'' into $V$ and with compact supports in the directions ``flowing'' away from $V$.

The self map $\Phi$ also induces
self maps $\Phi_V^{!*}$ on $H^i(V; A_V^{!*})$ and $\Phi_V^{*!}$ on $H^i(V; A_V^{*!})$ and
in \cite{Lefschetz} it is proven that
associated {\em local Lefschetz numbers}
$\text{Lef}(\Phi_V^{!*}; A_V^{!*})$ and  $\text{Lef}(\Phi_V^{*!};A_V^{*!})$ are equal\footnote{In many cases the hyperbolic
localizations $A_V^{!*}$ and $A_V^{*!}$ are canonically isomorphic, cf. \cite{Lefschetz, Braden} }. 

 \begin{tcolorbox}[colback=cyan!30!white]
\begin{thm}\cite{Lefschetz}\label{thm-Lefschetz}
Given a self map $f:W \to W$, a complex of sheaves $A\b \in D^b_c(W)$ and a lift $\Phi:f^*(A\b) \to A\b$,  
suppose that $W$ is
compact and that all connected components of the fixed point set are hyperbolic.  Then the global
${\rm Lef}(f, A\b)$ is the sum over connected components of the fixed point set of the local
Lefschetz numbers:
\[ {\rm Lef}(f, A\b) = \sum_V {\rm Lef}(\Phi_V^{!*}; A_V^{!*}) =
\sum_V {\rm Lef}(\Phi_V^{*!}; A_V^{*!}).\]
\end{thm}\end{tcolorbox}
Moreover, each local Lefschetz number ${\rm Lef}(\Phi^{!*}_V)$ 
is the Euler characteristic of a constructible function
$\text{Lef}(\Phi_x, A_V^{!*})$ for $x \in V$ (see \S \ref{subsec-constructible} above). 
Let $V = \coprod V_r$ be a stratification of the fixed point component $V$ so that
the pointwise Lefschetz number $\text{Lef}(\Phi_x, A_V^{!*})$  
is constant on each stratum $V_r$, and call it $L_r(\Phi; A_V^{!*})$. 
If $V$ is compact then (cf. \cite{Lefschetz} \S 11.1),
\[\text{Lef}(\Phi_V^{!*};A_V^{!*}) = \sum_r \chi_c(V_r) L_r(\Phi;A_V^{!*}).\]

\subsection{Torus actions}\label{subsec-revisited}\index{torus action}
The Morse-Bott theory  of critical points for smooth 
manifolds also has various extensions to singular spaces. 
Suppose the torus $T = \CC^*$ acts algebraically on a (possibly singular)
normal projective algebraic variety $W$ and suppose the action extends to some projective
space $M=\PP^m$, containing $W$. The imaginary part of the K\''{a}hler form on $M$ is a nondegenerate symplectic
2-form whose restriction to any complex submanifold of $M$ is also nondegenerate. 
The action of the circle $S^1 = \{e^{i\theta}\} \subset \CC^*$ defines
a Hamiltonian vector field on $M$ with resulting (cf. \S \ref{subsec-momentmap})
moment map $\mu:M \to \RR$ which may be thought of
as a Morse function that increases along the real directions of the $\CC^*$ action.   The critical points
of $\mu$ are exactly the fixed points of the torus action.

Let $W^T = \coprod_r V_r$
denote the fixed point components of the torus action and define
\begin{equation} \label{eqn-Vplus}\begin{aligned}
V_r^- &= \left\{ x \in W|\ \lim_{t\to \infty} t.x \in V_r \right\}\\
    V_r^+ &= \left\{x\in W|\ \lim_{t\to 0} t.x \in V_r \right\}.\end{aligned}\end{equation}
with inclusions as in equation (\ref{eqn-Vmaps}):
\begin{equation}\label{eqn-Vmaps2}
\begin{diagram}[size=2em]
V_r & \rTo^{j_r^{\pm}} & V_r^{\pm} & \rTo^{h_r^{\pm}} & W
\end{diagram}\end{equation}
In \cite{BB} Bialynicki-Birula proves\index{Bialynicki-Birula theorem} that if  $W$ is nonsingular
then each $V_r^{+} \to V$ has the natural structure of an algebraic bundle of affine spaces (of
some dimension $d_r$) and the
moment map $\mu:W \to \RR$ is a $\QQ$-perfect Morse-Bott\footnote{meaning that for each component $V$ of
the fixed point set, the Hessian $\partial^2f/\partial x_i\partial x_j$ is nondegenerate on each normal space 
$T_pM/T_pV$   and the connecting homomorphisms vanish in the long exact sequences of the pairs
$H^i(W_{\le v+\epsilon}, W_{\le v-\epsilon};\QQ)$ where $v = \mu(p)$} function.   
 \begin{tcolorbox}[colback=cyan!30!white]
\begin{thm}  If $W$ is nonsingular then  for each $j\ge 0$ 
there is an isomorphism
\begin{equation}\label{eqn-BB} 
H^j(W;\QQ) \cong \bigoplus_r H^{j-2d_r}(V_r;\QQ)\ \text{ and }\ H_j(W;\QQ) \cong \bigoplus_rH_{j-2d_r}(V_r;\QQ).\end{equation}
\end{thm}\end{tcolorbox}
The original proof involves the Weil conjectures and counting points mod $p$.
In many cases these isomorphisms can be described geometrically.  Suppose $Z \subset V_r$ is an compact $d$-dimensional
algebraic subvariety with fundamental class $[Z] \in H_{2d}(Z)\to H_{2d}(V_r)$.  Let $\pi:V_r^+\to V_r$ be the bundle
projection and let $Z^+ = \pi^{-1}(Z)$.  Then the closure $\overline{Z^+}$ is a $d+d_r$ dimensional variety whose
fundamental class $[\overline{Z^+}] \in H_{2d+2d_r}(W)$ is the image of $[Z]$ under the second isomorphism in
(\ref{eqn-BB}), cf.  \cite{CarrellGoresky}.

Returning to the general case, suppose $W$ is possibly singular.  Each fixed point $p \in W \subset M$ of the torus action is a
critical point for the moment map $\mu:M \to \RR$ but if $p$ is in the singular set of $W$ then the restriction $\mu|W$ 
will fail to be nondegenerate (at $p$)
in the sense of \S \ref{subsec-Morsefunctions} because $d\mu(p)$ kills all limits of tangent vectors.  
So $\mu:W \to \RR$ is not a Morse
function (or a Morse-Bott function) in the sense of \S \ref{subsec-Morsefunctions}.  
However some aspects of Morse theory continue to work in this case.

Each of the sets $V_r^{\pm}$ is algebraic and the projection $V_r^{\pm} \to V_r$ is algebraic. The time $t=1$
map of the torus action is hyperbolic (in the sense of \S \ref{subsec-hyperbolic}) at each fixed point component.
 Suppose $A\b \in D^b_c(W)$ is a complex of sheaves on $W$.  As in \S \ref{subsec-hyperbolic} define
\[A_r^{!*} = (j_r^+)^!(h_r^+)^*A\b\ \text{ and }\  A_r^{*!} = (j_r^-)^*(h_r^-)^!A\b.\]
 If the sheaf $A\b$ is {\em weakly hyperbolic}\footnote{meaning that there exists a local system
 $\mathcal L$ on $T$ and a quasi-isomorphism $\mu^*(A\b) \cong \mathcal L \boxtimes A\b$ where
 $\mu :T \times W \to W$ is the torus action} or if $A\b$ has a lift to the equivariant derived
 category of \cite{BernsteinLuntz} then (\cite{Lefschetz,Braden}) there is a canonical isomorphism
 $A_r^{*!}\cong  A_r^{!*}$ so  (cf.~Theorem \ref{thm-Lefschetz}):
  \begin{tcolorbox}[colback=cyan!30!white]
If $A\b$ is weakly hyperbolic and$W$ is compact then  $\chi(W;A\b) = \sum_r \chi(V_r; A_r^{*!})$.\end{tcolorbox}

The intersection complex $I^{\bf p}C\b$ with arbitrary perversity $\bf p$ on $W$ (see \S ?? and \S ?? below) 
 is weakly hyperbolic in a canonical way.
In \cite{Kirwan} F. Kirwan generalized the Bialynick-Birula theorem (\ref{eqn-BB}) to the case of singular varieties.  She uses
the decomposition theorem (Theorem \ref{sec-decomposition}) to prove the following for the {\em middle} perversity
intersection complex $\uul{IC}\b$:
 \begin{tcolorbox}[colback=cyan!30!white]
\begin{thm}\label{thm-Kirwan} (\cite{Kirwan})  The moment map $\mu$ is a perfect\index{Morse!function!perfect}
 Morse Bott function and it induces a decomposition for all $i$, expressing the intersection
cohomology of $W$ as a sum of locally defined cohomology groups of the fixed point components:
\begin{equation}\label{eqn-IHtorus}
IH^i(W) \cong \bigoplus_r H^i(V_r;\uul{IC}{}_r^{!*}).\end{equation}
\end{thm}\end{tcolorbox}
The result also generalizes to actions of a torus $T = (\CC^*)^m$.  In \cite{Braden} the sheaf
$\uul{IC}{}_r^{!*}$ is shown to be a direct sum of $\uul{IC}^{\bullet}$ sheaves of subvarieties of $V_r$.

\section{Intersection Homology}\index{intersection homology}
\subsection{A motivating example}
Consider $W=\Sigma T^3$, the suspension of the 3-torus with singular points denoted $\{N\}, \{S\}$.  We have natural cycles,
$\text{point}$, $T^1$, $\Sigma T^1$, $T^2$, $\Sigma T^2$, $T^3$, $\Sigma T^3$.  Some
of these hit the singular points, some do not.  The ones that do not are homologous to zero
by a ``homology'' or bounding chain that hits the singular points.  If we restrict cycles and homologies by not
allowing them to hit the singular points, this will change the resulting homology groups.
For $p = 0, 1, 2$ define $C_i^p(W) = \{ \xi \in C_i(W)|\ \xi \cap \{N,S\}  = \phi \text { unless } i \ge  4-p\}$
Here are the resulting homology groups.
\begin{center}
\begin{figure}[!h]
\begin{tabular}{||c|c|c|c||}
\hline
         & $p=0$               & $p=1$             & $p=2$              \\
         \hline 
 $i=4 $& $\Sigma T^{3^{\phantom{M} }}$ & $\Sigma T^3$ &  $\Sigma T^3$  \\
 \hline 
 $i=3$& $0$                  & $\Sigma T^2$ &  $\Sigma T^2$  \\
 \hline
 $i=2$& $T^2 $             &    $0$              &  $\Sigma T^1$   \\
 \hline
 $i=1$ & $T^1$            &    $T^1$           &    $0$               \\
 \hline
 $i=0$ & $\{pt\}$          & $\{pt\}$             &   $\{pt\}$            \\
 \hline
\end{tabular}\qquad
\begin{tabular}{||c|c|c|c||}
\hline
   & $p=0$                & $p=1$              &  $p=2$          \\
   \hline 
$i=4$ &      $T^3$      &           $T^3$    &         $T^3$     \\
\hline 
$i=3$ &          0        &            $T^2$     &           $T^2$    \\
\hline 
$i=2$ &           0       &              0          &          $T^1$      \\
\hline 
$i=1$  &      0          &                0       &             0          \\
\hline 
$i=0$  &         0        &             0         &             0         \\
\hline
\end{tabular}
\caption{Intersection homology and stalk homology of $\Sigma T^3$}
\end{figure}
\end{center}
The larger the number $p$ the more cycles are allowed into the singular points.  If there are more strata we can assign
such numbers to each stratum separately.  Cycles with a large
value of $p$ may get ``locked in'' to the singular set.  Their obstinate and perverse 
refusal to move away from the singular set led to the notion of a ``perversity'' vector.
\bigskip

\subsection{Digression on transversality}\label{subsec-transversality}\index{transversality}
Let $K \subset \RR$ be the Cantor set.  Let $f:\RR \to \RR$ be a smooth ($C^{\infty}$) function that vanishes
precisely on $K$.  Let $A \subset \RR^2$ denote the graph of $f$ and let $B$ denote the $x$-axis.  Then $A, B$ are
smooth submanifolds of $\RR^2$ but their intersection is the Cantor set.  This sort of unruly behavior can be avoided using transversality.

Two submanifolds $A,B \subset M$ of a smooth manifold are said to be transverse at a point of their intersection
$x \in A \cap B$ if $T_xA + T_xB = T_xM$.  If $A$ and $B$ are transverse at every point of their intersection then $A \cap B$ is a smooth submanifold of $M$ of dimension $\dim(A) + \dim(B) -\dim(M)$.  Arbitrary submanifolds $A,B, \subset M$ can
be made to be transverse by moving either one of them, say $A' = \phi_{\epsilon}(A)$ by the flow, for an arbitrarily small time, of a smooth vector field on $M$.  If $V$ is a finite dimensional vector space of vector fields on $M$ which span the tangent space $T_xM$ at every point $x \in M$ then there is an open and dense subset of $V$ consisting of vector fields $v$
such that the time $=1$ flow $\phi_1$ of $v$ takes $A$ to a submanifold $A' = \phi_1(A)$ that is transverse to $B$.

This is a very powerful result.  It says, for example, that two submanifolds of Euclidean space can be made transverse
by an arbitrary small translation.  The proof, due to Marston Morse, is so elegant, that I decided to include it in
 the Appendix \S \ref{appendix-transversality}.

Two Whitney stratified subsets $W_1, W_2 \subset M$ are said to be transverse if each stratum of $W_1$ is transverse to 
each stratum of $W_2$, in which case the intersection $W_1 \cap W_2$ is canonically Whitney stratified.  As in the
preceding paragraph, Whitney stratified sets can be made to be transverse by the application of the flow,
 for an arbitrarily small time, of a similarly chosen smooth vector field on $M$.

\subsection{Intersection homology}
Let $W$ be a compact $n$-dimensional 
Whitney stratified pseudomanifold with strata $S_{\alpha}$
($\alpha$ in some index set $I$, partially ordered by the closure relations between strata with 
$S_0$ being the stratum of dimensioin $n$)  and let $0 \le p_{\alpha}\le \cod(S_{\alpha})-2$ be a collection of  integers
which we refer to as  a perversity $\bar p$.    Define the intersection chains, 

\begin{tcolorbox}[colback=yellow!30!white]
\begin{defn}
\begin{equation}\label{eqn-IC}
IC_i^{\bar p}(W) = \left\{ \xi \in C_i(W)\left|\ \begin{aligned}
&\dim(\xi \cap S_{\alpha}) \le i - \cod(S_{\alpha}) + p_{\alpha}\\
&\dim(\partial \xi \cap S_{\alpha}) \le i - 1 - \cod(S_{\alpha}) + p_{\alpha}
\end{aligned} \text{ for }\alpha>0\right.\right\} \end{equation}
\end{defn}\end{tcolorbox}
Having placed the same restrictions on the chains as on their boundaries, we obtain a chain complex,
in fact a complex of (soft) sheaves $\uul{IC}^{\bar p}$ with resulting cohomology groups $IH_i^{\bar p}(W)$
or $I^{\bar p}H_i(W)$.
(As usual, ``chains" could refer to PL chains, singular chains, subanalytic chains, etc.)  Because $W$ is a
pseudomanifold the singular strata have codimension at least 2.  The condition $p_{\alpha} \le \cod(S_{\alpha})-2$
implies that most of the chain, and most of its boundary are completely contained within the top stratum $S_0$.
So a cycle $(\partial \xi = 0)$ in $IC_i^{\bar p}$ is also a cycle for ordinary homology and similarly a bounding
chain in $IC_i^{\bar p}$ is also a bounding chain for ordinary homology.  So we have a
homomorphism $IH_i^{\bar p}(W) \to H_i(W)$.  Moreover, if $\xi \in IC_i^{\bar p}(W)$ and if $\eta \in
IC_j^{\bar q}(W)$ and if we can arrange that $\xi\cap S_{\alpha}$ and $\eta \cap S_{\alpha}$ are transverse
within each stratum $S_{\alpha}$ then we will have an intersection
\[ \xi \cap \eta \in IC_{i+j-n}^{\bar p + \bar q}(W)\]
which is well defined provided that $p_{\alpha} + q_{\alpha} \le \cod(S_{\alpha}) -2$ for all $\alpha > 0$.

The first problem with this construction is that it is obviously dependent on the stratification.  Moreover, if we
are not careful, large values of $p_{\alpha}$ for small strata $S_{\alpha}<S_{\beta}$ will have the effect of allowing
chains into $S_{\alpha}$ but not into $S_{\beta}$ thereby ``locking" the chain into passing through a small stratum,
resulting in a stratification-dependent homology theory.
This issue can be avoided by requiring that $p_{\alpha}$ depends only on $\cod(S_{\alpha})$ and that
$\beta > \alpha \implies p_{\beta} \ge p_{\alpha}$.  

The second problem involves the effect of refining the stratification.
For a simple case, suppose $W$ consists only of two strata, $S_0$
and $S_c$, the singular stratum having codimension $c \ge 3$, to which we assign a perversity $p_c$.  Now
suppose we refine this stratum by introducing a ``fake" stratum, $S_r$ of codimension $r>c$.  Chains in
$IC^p_i(W)$ may intersect $S_c$ in dimension $\le i -c + p_c$ and for all we know, they may lie completely
in $S_r$, meaning that the chain will have ``perversity" $p_r = p_c + c-r$.  On the other hand if we assume, as
before, that we can arrange for this chain to be transverse to the fake stratum $S_r$ within the stratum $S_c$
then its intersection with $S_r$ will have dimension $\le i-c+p_c -(r-c) = i-r+p_c$ which is to say that it has
``perversity" $p_r = p_c$.  This argument shows (or suggests) that in this case we have natural isomorphisms
between the intersection homology $IH^{p_c}_i(W)$ as computed before the refinement, and the intersection
homology $IH^{p_c,p_r}_i$ after refinement, for any $p_r$ with $p_c \le p_r \le p_c + r-c$, that is,
\[ IH^{p_c, p_c}_i(W) \cong IH^{p_c, p_c+1}_i(W) \cong \cdots \cong IH^{p_c, p_c+r-c}_i(W)\]
In summary, assuming that $p_c \le p_r \le p_c + r-c$ the resulting homology group
$IH^{p_c,p_r}_i$ is unchanged after refinement.  This leads  to the formal definition of intersection homology.

\begin{tcolorbox}[colback=yellow!30!white]
\begin{defn}\index{perversity}
A {\em perversity} is a function $\bar p = (p_2, p_3, \cdots)$ with $p_2 = 0$ and with
$p_c \le p_{c+1} \le p_c +1$.  The complex of sheaves of intersection chains is the complex with sections
\begin{equation}\label{eqn-IC-sheaf1}
\Gamma(U, \uul{IC}{}^{-i}_{\bar p}) =
\left\{ \xi \in C_i(U)\left|\ \begin{aligned}
&\dim(\xi \cap S_{c}) \le i - c + p_{c}\\
&\dim(\partial \xi \cap S_{c}) \le i - 1 - c + p_{c}
\end{aligned} \text{ for } c \ge 2\right.\right\} \end{equation}
where $S_c$ denotes the union of all strata of codimension $c \ge 2$.
\end{defn}\end{tcolorbox}
Intersection homology with coefficients in a local system is defined similarly, however something special
happens in this case.
For any triangulation of a chain  $\xi \in IC^{-i}_{\bar p}$ all of its $i$-dimensional simplices and all of its
$i-1$ dimensional simplices will be completely contained within the top stratum (or ``nonsingular part") of $W$.
So if $\mathcal L$ is a local coefficient system {\em defined only on the top stratum of $W$}, we can still construct the
sheaf of intersection chains $\uul{IC}{}^{\bullet}_{\bar p}(\mathcal L)$ exactly as above.  

Let $\bar 0$ be the perversity $0_c = 0$ and let $\bar t$ be the perversity $t_c = c-2$.
\begin{tcolorbox}[colback=cyan!30!white]
\begin{thm}\label{thm-IH}  Let $W$ be an oriented stratified pseudomanifold of dimension $n$.
For any choice of perversity $\bar p$ equation (\ref{eqn-IC-sheaf1}) defines a complex of soft sheaves 
$\uul{IC}{}^{\bullet}_{\bar p}$ on $W$ and
the following holds.  \begin{enumerate}
\item   The cohomology sheaves $\uul{IH}{}^{-m}_{\bar p}$ and the hypercohomology groups $IH^{\bar p}_i(W)$ are
well defined and are independent of the stratification; 
\item in fact they are topological invariants.    
\item There are canonical maps 
\[H^{n-i}(W) \to IH_i^{\bar p}(W) \to H_i(W)\] 
that factor the Poincar\'e duality
map, 
\item if $\bar p \le \bar q$ then there are also compatible mappings $IH^{\bar p} \to IH^{\bar q}$.  In
sheaf language we have natural maps
\[ \uul{\ZZ}{}_W[n] \to \uul{IC}{}_W^{\bar p} \to \uul{IC}{}_W^{\bar q} \to \uul{C}{}_W^{BM}.\]
\item If the link $L_X$ of each stratum $X$ is connected then for $\bar p = \bar 0$ the first of these maps is a 
quasi-isomorphism, and for $\bar q= \bar t$ the second map is  a quasi-isomorphism.  
\item If $p_c + q_c \le t_c=c-2$ for all $c$ then the 
intersection of {\em transversal} chains determines a pairing 
\[ IH_i^{\bar p}(W) \times IH_j^{\bar q}(W) \to IH_{i+j-n}^{\bar p + \bar q}(W)\]
\item If $\bar p + \bar q = \bar t$ then the resulting pairing
\[ IH_i^{\bar p}(W) \times IH_{n-i}^{\bar q}(W) \to H_0(W) \to \ZZ\]
is nondegenerate over $\QQ$ (or over any field).\end{enumerate}
\end{thm}\end{tcolorbox}

The last statement in Theorem \ref{thm-IH}, Poincar\'e duality, was the big surprise when intersection homology was
discovered for it is a duality statement that applies to singular spaces.  Especially, if
the stratification of $W$ consists only of even codimension strata then there is a ``middle"
choice for $p$, that is, $p_c = (c-2)/2$ for which $IH^{\bar p}(W; k)$ is self-dual for any field $k$.

There is a technical problem with moving chains within a Whitney stratified set $W$, so as to be transverse
within each stratum of $W$.  This can be accomplished with piecewise-linear chains within a piecewise-linear
stratified set $W$, and has recently been accomplished using semi-analytic chains within a semi-analytic
stratified set, but to my knowedge, it has not been accomplished in any other setting.  This is one of
the many problems that is avoided with the use of sheaf theory.  The proof of topological invariance
depends entirely on sheaf theory.  Other results such as the proof of Poincar\'e duality, that can be
established using chain manipulations, are incredibly awkward, requiring a choice of model for
the chains, and delicate manipulations with individual chains.  These constructions are easier, but less
geometric, if they are all made using sheaf theory.  For this purpose we need to identify the quasi-isomorphism
class of the complex of sheaves $\uul{IC}{}^{\bar p}$.

\begin{tcolorbox}[colback=cyan!30!white]
\begin{prop}  Let $W$ be a Whitney stratified pseudomanifold and let $\mathcal E$ be a local
coefficient system defined on the top stratum.
Fix a perversity $\bar p$, and let $x \in X$ be a point in a stratum $X$ of codimension $c\ge 2$.  
Then the stalk
of the intersection homology sheaf at $x$ is
\[\uul{H}^{-i}(\uul{IC}^{\bar p}(\mathcal E))_X = IH_i^{\bar p}(W, W-x;\mathcal E) = \begin{cases}
0 & \text{ if } i < n-p_c\\
IH_{i-n+c-1}(L_X; \mathcal E) &\text{ if } i\ge n-p_c 
\end{cases}
\]
and the stalk cohomology with compact supports is
\[
H^{-i}_c(U_x; \uul{IC}^{\bar p}(\mathcal E)) = IH_i^{\bar p}(U_x) = \begin{cases}
H_i(L_x;\mathcal E) & \text{ if } i \le c-p_c \\
0 & \text{ if } i > c-p_c \end{cases} \]
\end{prop}\end{tcolorbox}
\begin{proof}
Use the local product structure of a neighborhood $U_x \cong c^o(L_x) \times B_{\epsilon}^{n-c}$ and the
K\"unneth formula\footnote{The homology $H_*(B_{\epsilon}, \partial B_{\epsilon})$ is torsion-free so the ${\rm Tor}$
terms vanish.} 
\[ IH^{\bar p}_i(U, \partial U;\mathcal E) \cong IH^{\bar p}_{i-(n-c)}(c(L_X), L_X;\mathcal E). \]
If $\xi \in IC_i^{\bar p}$ and if $(i-n+c) -c + p_c \ge 0$ then  the chain $\xi$ is allowed to hit the
cone point, otherwise it is not.  When it is allowed to hit the cone point, we may assume (using a
homotopy argument) that it locally coincides with the cone over a chain in $L_X$ which satisfies
the same allowability conditions.  Similar remarks apply to $\partial \xi$.  On the other hand, suppose $\xi$
is a compact $i$-dimensional chain in the link $L_X$.  It is the boundary of the cone $c(\xi)$ so the
homology class $[\xi]$ vanishes in the neighborhood $U_x$ provided that cone is allowable, which occurs if 
\[\dim(c(\xi)\cap X) = 0 \le (i+1) -c + p_c\ \text{ that is, if }\ i>c - p_c.\qedhere\]
\end{proof}

Comparing this to the calculation (\ref{eqn-jstar})
of $Rj_*j^*(\uul{IC}^{\bar p})$ where $j:U=W-\overline{S_c} \to W$ is the
inclusion of the open complement of the closure of $S_c$ we see that the intersection homology sheaf on $S_c$ is the {\em truncation}  (\S \ref{subsec-truncation}) of the sheaf $Rj_*(\uul{IC}^{\bar p}|U)$.  
For example, suppose $\dim(W) = 8$ has strata of dimension
$0, 2, 4, 6, 8$ and the perversity is the middle one, $p_c = (c-2)/2$.  Then the stalk cohomology $H^*(\uul{IC}{}^{\bar p}_x)$
of the sheaf  $\uul{IC}^{\bar p}$ is described in Figure \ref{fig-support}, 
where $L^r$ means the $r$-dimensional link of the codimension
$r+1$ stratum and the red zeroes represent homology groups that have been killed by the perversity condition.

\begin{center}
\begin{figure}[!h]\index{support!diagram}
\begin{tabular}{|c||c|c|c|c|c|}
\hline
$i$ & $\cod 0$ & $\cod 2$ & $\cod 4$ & $\cod 6$ & $\cod 8$\\
\hline\hline
0 &                 &                 &                &               &            \rz    \\
-1 &                 &                 &                &               &          \rz   \\
-2 &               &                   &               &                &         \rz   \\
-3 &                &                 &                 &         \rz    &      \rz      \\
-4 &               &                  &                  &       \rz    &      \rz      \\
-5 &               &                  &      \rz        &      \rz      & $IH_4(L^7)$\\
-6 &              &                  &        \rz       & $IH_3(L^5)$&$IH_5(L^7)$\\
-7 &               &        \rz     & $IH_2(L^3)$& $IH_4(L^5)$ & $IH_6(L^7)$\\
-8&   $\ZZ$    &  $ IH_1(L^1)$& $IH_3(L^3)$& $IH_5(L^5)$ & $IH_7(L^7)$\\
\hline 
\end{tabular}\quad \begin{tabular}{||c|c|c|c|c|}
\hline
 $\cod 0$ & $\cod 2$ & $\cod 4$ & $\cod 6$ & $\cod 8$\\
\hline\hline
     $\ZZ$ &$IH_0(L^1)$&$IH_0(L^3)$&$IH_0(L^5)$&$IH_0(L^7)$ \\
                 &   \rz        &$IH_1(L^3)$&$IH_1(L^5)$&$IH_1(L^7)$\\
               &                &    \rz       &$IH_2(L^5)$&$IH_2(L^7)$   \\
                &                 &    \rz      &         \rz    &$IH_3(L^7)$     \\
               &                  &              &       \rz    &      \rz      \\
               &                  &             &      \rz      & \rz         \\
              &                  &              &                &\rz  \\
               &                &               &                 &   \rz \\  
               &                &               &                &  \rz   \\
\hline 

\end{tabular}
\caption{Stalk cohomology and compact support cohomology of $\uul{IC}^{\bullet}$}\label{fig-support}
\end{figure}
\end{center}

This gives an inductive way to construct intersection homology using purely sheaf-theoretic operations, to be
described in the next section.
\bigskip

\section{Truncation}\label{sec-truncation}\index{truncation}

\subsection{Truncation}\label{subsec-truncation}
If $A\b$ is a complex of sheaves define
\[ (\tau_{\le r}A\b)^i = \begin{cases}
0 & \text{ if } i>r\\
\ker(d^r) & \text{ if } i=r\\
A^i & \text{ if } i<r
\end{cases} \]
Then $\uul{H}^i(\tau_{\le r}A\b) = \uul{H}^i(A\b)$ for $i \le r$ and is zero for $i > r$.

During a conversation in October 1976 Pierre Deligne suggested that the following construction might
generate the intersection homology sheaf.
\begin{defn}\index{Deligne's construction}
Let $W$ be a purely $n$ dimensional oriented Whitney stratified pseudomanifold and let
$W_{r}$ denote the union of the strata of dimension $\le r$.  Set $U_k = W - W_{n-k}$
with inclusions 
\[ \begin{CD}
U_2 @>>{j_2}> U_3 @>>{j_3}>  \cdots @>>{j_{n-1}}> U_n @>>{j_n}> U_{n+1}=W \end{CD}\]
 Let $\mathcal E$ be a local coefficient system defined on the top stratum $U=U_2$.  Set 
\[\uul{P}{}^{\bullet} _{\bar p}(\mathcal E)
= \tau_{\le p(n)} Rj_{n*} \cdots \tau_{\le p(3)} Rj_{3*} \tau_{\le p(2)}R{j_2*}\mathcal E.\] 
\end{defn}
The resulting complex of sheaves will have stalk cohomology that is illustrated in Figure 
\ref{fig-stalkcohomology} (in the case of
middle perversity, with $\mathcal E = \ZZ$, for a stratified psuedomanifold $W$ with 
$\dim(W) = 8$ that is stratified with strata of codimension $0$, $2$, $4$, $6$, $8$).  In this
figure we suppress the $\bar p$ on $\uul{P}{}^{\bullet}_{\bar p}$.
\newcommand{\uP}{\uul{P}^{\bullet}}
\begin{center}
\begin{figure}[!h]
\begin{tabular}{|c||c|c|c|c|c|}
\hline
$i$ & $\cod 0$ & $\cod 2$ & $\cod 4$ & $\cod 6$ & $\cod 8$\\
\hline\hline
8 &                 &                 &                &               &                \\
7 &                 &                 &                &               &          \rz   \\
6 &               &                   &               &                &         \rz   \\
5 &                &                 &                 &         \rz    &      \rz      \\
4 &               &                  &                  &       \rz    &      \rz      \\
3 &               &                  &      \rz        &      \rz      & $H^3(L^7,\uP)$\\
2 &              &                  &        \rz       & $H^2(L^5,\uP)$&$H^2(L^7,\uP)$\\
1 &               &        \rz     & $H^1(L^3,\uP)$& $H^1(L^5,\uP)$ & $H^1(L^7,\uP)$\\
0&   $\ZZ$    &  $ H^0(L^1,\uP) $& $H^0(L^3,\uP)$& $H^0(L^5,\uP)$ & $H^0(L^7,\uP)$\\
\hline 

\end{tabular}
\caption{Stalk cohomology of $\uul{P}\b$}\label{fig-stalkcohomology}
\end{figure}\end{center}
We remark, for example, at a point $x \in W$ that lies in a stratum $X$ of codimension $6$, the
stalk cohomology at $x$ equals the cohomology of the link $L^5$ of $X$ with coefficients in the part of the
sheaf $\uul{P}^{\bullet}|U_6$ that has been previously constructed over the strictly {\em larger} strata $Y>X$.

\begin{tcolorbox}[colback=cyan!30!white]
\begin{thm}\label{thm-ICbytruncation}\cite{IH2}
  Let $W$ be an oriented $n$-dimensional stratified pseudomanifold and let $\mathcal E$ be a local coefficient
system on the top stratum.  The orientation map $\mathcal E[n] \to \uul{C}{}^{-n}_{U}(\mathcal E)$ induces an isomorphism
$\uul{P}{}^{\bullet}_{\bar p}(\mathcal E)[n] \cong \uul{IC}{}^{\bullet}_{\bar p}(\mathcal E)$ where $U = U_2$ denotes the largest stratum.  

If the link $L_x$ of every stratum is connected then 
$\uul{P}{}^{\bullet}_{\bar t}(\ZZ)[n] \to \uul{C}^{\bullet}(\uul{\ZZ}{}_U)$ is a
quasi-isomorphism (so that $IH^i_{\bar t}(W) = H_{n-i}(W)$) and $\uul{P}{}^{\bullet}_{\bar 0}(\ZZ) \to \uul{\ZZ}$ is a 
quasi-isomorphism (so that $IH^i_{\bar 0}(W) = H^i(W)$), where $\bar 0_c = 0$ and where $\bar t_c = c-2$ are the
``bottom" and ``top" perversities respectively. 
If $\bar p + \bar q \le \bar t$ (where $t_c = c-2$)
and if $\mathcal E_1 \otimes \mathcal E_2 \to \mathcal E_3$ is a morphism of local systems on $U$ then it  extends canonically to a product
\[
\uul{P}{}^{\bullet}_{\bar p}(\mathcal E_1) \otimes \uul{P}{}^{\bullet}_{\bar q}(\mathcal E_2) 
\to \uul{P}{}^{\bullet}_{\bar p + \bar q}(\mathcal E_3).\]
\end{thm}\end{tcolorbox}

\begin{proof}
For simplicity we discuss the case of constant coefficients.  There are two problems (a) to show that the
orientation map $\uul{\ZZ}{}_U \to \uul{C}{}^{-n}_U$ extends to a (uniquely defined) map in the derived
category $\uul{P}^{\bullet}[n] \to \uul{IC}^{\bullet}$ (for a fixed perversity, which we suppress in the notation)
and (b) to show that this map is a quasi-isomorphism.  These are proven by induction, adding one stratum at a time.
Consider the diagram
\[ \begin{CD}
U_k @>>{j_k}> U_{k+1} @<<{i_k}< X^{n-k} \end{CD}\]
where $X^{n-k}$ is the union of the codimension $k$ strata.  Suppose by induction that we have constructed a
quasi-isomorphism $\uul{P}{}^{\bullet}_k\to \uul{IC}{}^{\bullet}_k$ of sheaves over $U_k$ (where the subscript
$k$ denotes the restriction to $U_k$).  Now compare the two distinguished triangles (writing $i = i_k$ and
$j = j_k$ to simplify notation),
\begin{diagram}[size=2em]
Ri_*i^!(\uul{P}{}^{\bullet}_{k+1}) && \rTo && \uul{P}{}^{\bullet}_{k+1} &&& Ri_*i^!(\uul{IC}{}^{\bullet}_{k+1}) && \rTo && \uul{IC}{}^{\bullet}_{k+1} \\
  & \luTo && \ldTo &&\ \text{ and } && &\luTo && \ldTo &\\
  && Rj_*j^*(\uul{P}{}^{\bullet}_{k+1}) &=&Rj_*(\uul{P}{}^{\bullet}_k )&&& && Rj_*j^*(\uul{IC}{}^{\bullet}_{k+1}) &=&
Rj_*(\uul{IC}{}^{\bullet}_k)
\end{diagram}  
We are actually concerned with the right side of these triangles.
  By induction we have an isomorphism on the bottom row, so we get an isomorphism of the truncations:
\[\uul{P}{}^{\bullet}_{k+1}= \tau_{\le p(k)} Rj_*\uul{P}{}^{\bullet}_k \to \tau_{\le p(k)} Rj_*(\uul{IC}{}^{\bullet}_k)\]
This is the upper right corner of the first triangle and we wish to identify it with the upper right corner of the second triangle.
So it suffices to show that we have an isomorphism (in the derived category),
\[ \uul{IC}{}^{\bullet}_{k+1} \cong \tau_{\le p(k)}Rj_*(\uul{IC}{}^{\bullet}_k).\]
But this is exactly what the local calculation says:  the stalk of the intersection cohomology is the truncation of
the intersection cohomology of the link.

In fact, the formula $\uul{P}{}^{\bullet}_{k+1} = \tau_{\le p(k)}Rj_*\uul{P}{}^{\bullet}_k$ implies that the
attaching morphism $\uul{P}{}^{\bullet}_{k+1} \to Rj_*j^*\uul{P}{}^{\bullet}_{k+1}$ is an isomorphism in 
degrees $r \le p(k)$, or equivalently, that $H^r(i^!\uul{P}{}^{\bullet}) = 0$ for $r \le p(k)+1$.  This is the
same as saying that for any $x \in X^{n-k}$,
\[ H^r_c(U_x; \uul{P}{}^{\bullet}) = H^r(i_x^!\uul{P}{}^{\bullet})= 0 \text{ for } r <p(k)+2+(n-k)=n-q(k)\]
where $q(k)=k-2-p(k)$ is the complementary perversity,
$i_x: \{x\} \to W$ is the inclusion and $U_x$ is a basic open neighborhood of $x$ in $W$.  (See also
Proposition \ref{prop-attachingequivalent}.)

The construction of the pairing is similar.  Start with the multiplication
\[ \uul{\ZZ}{}_{U_2} \otimes \uul{\ZZ}{}_{U_2} \to \uul{\ZZ}{}_{U_2}.\]
Now apply $\tau_{\le p(2)}Rj_*$.  The truncation of a tensor product is not simply the tensor product
of the truncations, there are a lot of cross terms.  By examining the effect on the stalk cohomology one
eventually arrives at a pairing
\[ \uul{P}{}^{\bullet}_{\bar p} \otimes \uul{P}{}^{\bullet}_{\bar q} \to \uul{P}{}^{\bullet}_{\bar p + \bar q}.\qedhere\]
\end{proof}

The Poincar\'e duality theorem for intersection cohomology, translated
 into cohomology indexing says the following:  if $\bar p + \bar q = \bar t$ then the resulting pairing 
\[ IH^i_{\bar p}(W) \times IH^{n-i}_{\bar q}(W) \to IH^{\bar t}_n(W) \to \ZZ\]
is nondegenerate when tensored with any field.  In the next section this will be expressed in a sheaf theoretic way.

\subsection{The dual of a complex of sheaves}\label{subsec-duality}\index{dual!sheaf}\index{sheaf!dual}
Suppose $X$ is a paracompact Hausdorff space and $R$ is a commutative ring with finite cohomological
dimension. Let $A\b$ be a complex of sheaves of $R$-modules on $X$.  A. Borel and J.~C.~Moore 
defined (\cite{BorelMoore}) its dual
$\mathbf{D}(A)^{\bullet}$ as follows.  First, choose a soft or flabby resolution $A\b \cong S\b$.  Then
$\mathbf{D}(A)^{\bullet}$ is the sheafification of the complex of presheaves
\[\mathbf{D}(A)^{-j}( U) = {\rm RHom}_{R-mod}(\Gamma_c(U, S^j), R) = \hHom_{R-mod}(\Gamma_c(U, S^j), I\b)\]
where $I\b$ is an injective resolution\footnote{For $R = \ZZ$ take $I\b$ to be $\QQ \to \QQ/\ZZ$ and
more generally if $R$ is an integral domain take $F \to F/R$ where $F$ is the fraction field of $R$.} of 
(of finite length) of the ring $R$. 
This means that $\mathbf{D}^{\bullet}(U)$ is the single complex obtained from the double complex
$\Hom_{R-mod}^{\bullet}(\Gamma_c(U, S^{\bullet}), I^{\bullet})$ by adding along the diagonals in the usual way. 
If $R$ is a field this is just 
\[\hHom_{R-mod}(\Gamma_c(U, S^j), R).\]
 For example when $R = \RR$, the sheaf of {\em currents} on a smooth manifold is the dual of the sheaf of 
 smooth differential forms.
For  $R = \ZZ$, Borel and Moore proved that there are exact cohomology sequences
\[
0 \to \Ext(H^{n+1}_c(X, S\b), \ZZ) \to H^{-n}(X, \mathbf D(S\b)) \to \Hom(H^n_c(X, S\b), \ZZ) \to 0\]
in analogy with the universal coefficient theorem for cohomology:
\[0 \to \Ext(H_{n-1}(X,\ZZ),\ZZ) \to H^n(X,\ZZ) \to \Hom(H_n(X,\ZZ), \ZZ) \to 0.\]

The Borel-Moore sheaf of chains $\uul{C}{}^{\bullet}_{BM}$ 
is the dual of the constant sheaf but first we must replace the constant sheaf
by the (quasi-isomorphic) flabby sheaf  of singular cochains.  So the sheaf of chains is the sheafification of
the complex of presheaves
\[ C^{-j}_{BM}(U) = {\rm RHom}(C^j_c(U;R), R)\]
where $C^j_c(U;R)$ is the subgroup of the group of singular cochains $\Hom_R(C^{sing}_j(U;R), R)$ with compact support. 
 In short, any reasonable
injective model of the complex $\uul{C}{}^{\bullet}_{BM}$ is a mess since we have ``dualized'' the (cosheaf of)
singular chains twice.  For most applications we do not require a precise
injective model for this sheaf:  it is enough to know that one exists.

The Borel-Moore homology (in degree $j$) is the 
cohomology (in degree $-j$) of the complex $\uul{C}{}^{\bullet}_{BM}$.  
We have previously seen that this is the homology theory of
locally finite chains.  The problem that Borel and Moore failed to resolve is that the double dual of 
$S\b$ does not equal $S\b$; rather, it appears to be something far more complicated.

\begin{figure}[!h]
\includegraphics[width=.6\linewidth]{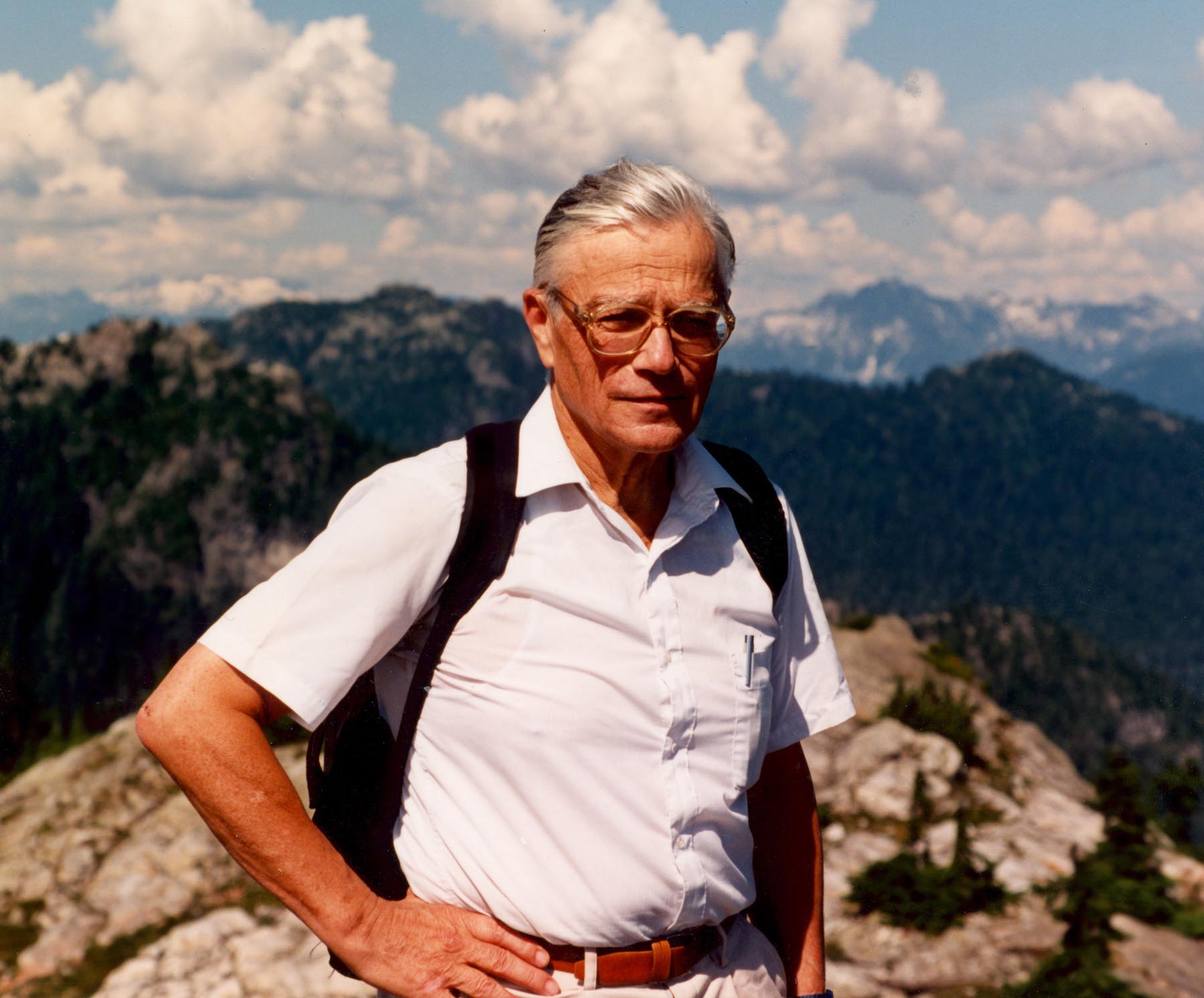}
\caption{Armand Borel (photo by Bill Casselman)}
\end{figure}

\subsection{Dualizing sheaf}\index{dualizing sheaf}\index{sheaf!dualizing}
Later, Verdier\cite{Verdierduality} discovered that the Borel-Moore dual  could be interpreted sheaf theoretically:
\[ \mathbf D(S\b) = \uul{\RHom}^{\bullet}(S\b, \mathbf D^{\bullet})\]
where $\mathbf D^{\bullet}$ is a particular, universal sheaf called the {\em dualizing complex}.  He then showed
that there is a canonical quasi-isomorphism in the derived category,
\[ \mathbf D \mathbf D(S\b) \cong S\b.\]
So double duality is restored.\footnote{a point which A. Borel was particularly sensitive to and which he 
carefully explained in his lectures \cite{BorelIH}}  And what is this magic dualizing complex? 
There is no choice:   for sheaves of abelian groups, 
for example, we have:
\[ \mathbf D^{\bullet} = \uul{\RHom}^{\bullet}(\uul{\ZZ}, \mathbf D^{\bullet}) = \mathbf D(\uul{\ZZ}) \]
is the Borel-Moore dual of the constant sheaf, so it is the sheaf of chains!  
More precisely, it is the quasi-isomorphism class of the complex of
sheaves $\uul C{}^{\bullet}_{BM}$ of Borel-Moore chains.

\subsection{Definition}\index{Verdier duality}\index{duality, Verdier}
  A pairing $S\b \otimes T\b \to \mathbf{D}\b$ of sheaves on a Whitney stratified space
$X$  is said to be a {\em Verdier dual pairing} if the resulting
morphism $S\b \to \uul{\RHom}^{\bullet}(T\b, \mathbf{D}\b)= \mathbf D(T\b)$ is an isomorphism in $D^b_c(X)$.  
In particular, this means that if $R$ is an integral domain then for any open set $U \subset X$ there
is a short exact sequence,
\[\begin{CD}
0 @>>> {\rm Ext}^1_R(H^{-i+1}_c(U, S\b), R) @>>> H^i(U, \mathbf D(S\b)) @>>> 
\hHom_R(H^{-i}_c(U, S\b), R) @>>> 0.
\end{CD}
\]
\begin{tcolorbox}[colback=yellow!30!white]
\begin{defn}\label{def-shriek}
If $f:X \to Y$ is a continuous map and $S\b$ is a complex of sheaves on $Y$ 
define $f^!(S\b) = \mathbf D_X f^* \mathbf D_Y(S\b)$.
\end{defn}\end{tcolorbox}
If $f:X \to \{pt\}$ is the map to a point then the dualizing sheaf is $\mathbf D_X^{\bullet} = f^!(\ZZ)$.

\begin{thm} [Verdier duality]  Let $f:X\to Y$ be a stratified mapping between Whitney stratified spaces.  
Let $A\b$,  $B\b$ and $C\b$ be constructible sheaves of abelian groups on $X$, $Y$ and $Y$ respectively.  
Then $f^*$, $f^!$, $Rf_*$ and $Rf_!$ take distinguished triangles to distinguished triangles.
There are canonical isomorphisms in $D^b_c(X)$ as follows:
\begin{enumerate}
\item $\mathbf D \mathbf D(A\b) \cong A\b$
\item $\mathbf{D}^{\bullet}_X \cong f^! \mathbf{D}^{\bullet}_Y$
\item $f^!(B\b) = \mathbf{D}_Xf^*\mathbf{D}_Y(B\b)$
\item $Rf_!(A\b) = \mathbf{D}_Y Rf_* \mathbf{D}_X(A\b)$
\item[] So $f^!$ is the dual of $f^*$ and $Rf_!$ is the dual of $Rf_*$.
\item $f^!\uul{\RHom}^{\bullet}(B\b,C\b) \cong \uul{\RHom}^{\bullet}(f^*(B\b), f^!(C\b))$
\item $Rf_*(\uul{\RHom}^{\bullet}(A\b, f^!B\b)) \cong \RHom^{\bullet}(Rf_!A\b, B\b)$  [Verdier duality theorem]
\item[] This says that $Rf_!$ and $f^!$ are adjoint, just as $Rf_*$ and $f^*$ are.
\item $Rf_*\uul{\RHom}^{\bullet}(f^*B\b, A\b) \cong \uul{\RHom}^{\bullet}(B\b, Rf_*A\b)$
\item $Rf_!\uul{\RHom}^{\bullet}(A\b, f^!B\b) \cong \uul{\RHom}^{\bullet}(Rf_!A\b, B\b)$
\item If $f:X \to Y$ is the inclusion of an open subset then $f^!(B\b) \cong f^*(B\b)$.
\item If $f:X \to Y$ is the inclusion of a closed subset then $Rf_!(A\b) \cong Rf_*(A\b)$.
\item If $f:X \to Y$ is the inclusion of an oriented submanifold in another, and if $B\b$ is
cohomologically locally constant on $Y$ then $ f^!(B\b) \cong f^*(B\b)[\dim(Y)-\dim(X)]$. 
\end{enumerate}\end{thm}

\paragraph{\bf Exercise} Verify (1) for the category of simplicial sheaves using the canonical model for the (simplicial)
sheaf of chains.  It comes down to a statement about
the second barycentric subdivision.
\paragraph{\bf Remark}  Dualizing complexes exist in many other categories and Verdier duality is now recognied
as the natural categorical statement of duality.

Finally we can state the sheaf theoretic statement of duality for intersection cohomology.  Let
$K$ be a field.
 \begin{tcolorbox}[colback=cyan!30!white]
\begin{thm}\cite{IH1, IH2}\label{thm-IHperfect}
Let $X$ be a Whitney stratified space and let $\bar p + \bar q = \bar t$ be complementary perversities.  
 Let
$\mathcal E_1 \otimes \mathcal E_2 \to \uul{K}$ be a dual pairing of local systems (of $K$-vector spaces)
defined over the top stratum of $X$.  Then the resulting pairing 
\[ \uul{IC}{}^{\bullet}_{\bar p}(\mathcal E_1)[n] \otimes \uul{IC}{}^{\bullet}_{\bar q}(\mathcal E_2) 
\to \uul{IC}{}^{\bullet}_{\bar t}(K) \to \mathbf{D}^{\bullet}_X\]
is a Verdier dual pairing.  
\end{thm}\end{tcolorbox}
The proof of Theorem \ref{thm-IHperfect} is by induction on the strata, as before, adding one stratum at a time, 
using the long exact sequences, duality, and the above formal properties.

\subsection{The middle perversity}  \index{perversity!middle}
Suppose the codimension of every stratum $X$ of the stratified space $W$ is even.
Then there is a middle perversity $\bar{m}(c) = (c-2)/2$ and $\uul{IC}{}^{\bullet}_{\bar{m}}(K)$ is self dual for any field $K$.  
If odd codimension
strata are present there is an upper middle $\bar{m}^+(c) = \lceil{\frac{c-2}{2}} \rceil$ and lower middle
$\bar{m}^-(c) = \lfloor \frac{c-2}{2}\rfloor$ perversity with a canonical morphism
\[ \Phi: \uul{IC}{}^{\bullet}_{\bar{m}^-}(K)  \to
\uul{IC}{}^{\bullet}_{\bar{m}^+}(K)\]
between the dual sheaves $\uul{IC}{}^{\bullet}_{\bar{m}^+}(K)$ and
$\uul{IC}{}^{\bullet}_{\bar{m}^-}(K)$.  The morphism $\Phi$ is a quasi-isomorphism if and only if $W$ is a
{\em $K$-Witt space} (\cite{Siegel, GoreskySiegel}), meaning that for each stratum
$X$ of odd codimension $c$ the stalk cohomology vanishes in the middle dimension:
\[ IH_{\bar{m}^-}^{(c-1)/2}(L_X;K) = 0\]
where $L_X$ denotes the link of the stratum $X$.  For such a space 
\[\uul{IC}{}^{\bullet}_{\bar{m}^-}(K)\cong \uul{IC}{}^{\bullet}_{\bar{m}^+}(K)\] is self dual.

\subsection{Orientation sheaf}  (See also \S \ref{subsec-chains}.)\index{sheaf!orientation}\index{orientation sheaf}
Suppose $X$ is Whitney stratified and $n$-dimensional, with largest stratum $X^0$.
The orientation sheaf $\mathcal O$ on $X^0$ is the local system whose stalk at each point $x \in X^0$ is 
the local homology $H_n(X,X-x;\ZZ) = H_n(X^0, X^0-x;\ZZ)$.   If it is placed in degree $-n$ then it becomes (quasi) isomorphic to the 
dualizing sheaf $\omega=\DD(X^0)$,
that is $\omega = \mathcal O[n]$.  Then $\PP(\omega)$ is the intersection {\em homology} sheaf; it's cohomology is:
\[ H^{-r}(\PP(\omega)) = IH_r(X;\ZZ)\]
for any $r \ge 0$. 
Let $\underline{\ZZ}$ be the constant sheaf on $X^0$, placed in degree $0$.  Then
\[ H^s(\PP(\uul{\ZZ})) = IH^s(X;\ZZ)\]
is the intersection {\em cohomology}.  (Truncations also need to be shifted by $n$ which
accounts for the difference between [16] and [90].)
The canonical pairings of sheaves on $X^0$,
\[ \uul{\ZZ}\otimes \uul{\ZZ} \to \uul{\ZZ}\ \text{ and }\
\uul{\ZZ} \otimes \omega \to \omega \quad\]
induce (with appropriately chosen perversities) cup products 
$IH^*\otimes IH^* \to IH^*$ in intersection cohomology and cap products 
$IH^* \otimes H_* \to H_*$.  These constructions work over the integers or more generally over any commutative ring of finite cohomological dimension.

There is always an isomorphism 
$\uul{\ZZ}[n]\otimes \ZZ/(2)   \to \omega \otimes \ZZ/(2)  $ of sheaves on $X^0$. 
An {\em orientation} of $X^0$, if one exists, is an isomorphism 
$\uul{\ZZ}[n] \to \omega$.  So an orientation induces an isomorphism $\PP(\uul{\ZZ})[n] \to \PP(\omega)$ hence a Poincar\'e duality isomorphism  $IH^*(X) \to IH_*(X)$.
In this case the product $\uul{\ZZ} \otimes \uul{\ZZ} \to \uul{\ZZ}$ on $X^0$ becomes $\omega\otimes\omega \to \omega[n]$
which induces (for appropriate choice of perversity) the (geometrically defined) original intersection  product
in intersection {\em homology}.
The proof of the following is in \S \ref{ProofOfEquivalence}.

\begin{tcolorbox}[colback=cyan!30!white]
\begin{thm}\label{thm-equivalence}\cite{IH2}\index{Deligne's construction}
Let $W$ be an $n$-dimensional Whitney stratified set with biggest stratum $U=U_2$ and let $\bar p$ be a pervesity.  
Then Deligne's construction 
\[ \mathcal E \mapsto \uul{P}^{\bullet}(\mathcal E) =\tau_{\le p(n)} Rj_{n*} \cdots \tau_{\le p(3)} 
Rj_{3*} \tau_{\le p(2)}R{j_2*}\mathcal E\] 
defines an equivalence of categories between the category of local systems of $K$-vector spaces 
($K$ a field) on the nonsingular part $U=U_2$, and the full subcategory of $D^b_c(W)$ consisting
of ``IC sheaves", that is, complexes of sheaves $A\b$, constructible with respect to the given
stratification, such that the following conditions hold 
\begin{enumerate}
\item $A\b | U_2 \cong \mathcal E$ is isomorphic to a local coefficient system
\item $\uul{H}^r(A\b) = 0$ for $r< 0$
\item $H^r(i_x^*A\b) = 0$ for $r > p(c)$\  (``support condition'')\index{support!condition}
\item $H^r(i_x^!A\b) = 0$ for $r < n-q(c)$\  $\phantom{{^{\strut}}}$ (``cosupport condition'')\index{cosupport condition}
\end{enumerate}
for all points $x \in W-U_2$, where $i_x: \{x\} \to W$ is the inclusion of the point and $c$ denotes the (real)
codimension of the stratum containing $x$ and where $q(c) = c-2-p(c)$ is the complementary perversity.
\end{thm}\end{tcolorbox}
If $\mathcal E$ is a local system on $U_2$ and if $A\b$ is a constructible complex of
sheaves that satisfies the  conditions (1)-(4), Theorem \ref{thm-equivalence} says that there is a canonical isomorphism $A\b \cong 
\uul{IC}{}^{\bullet}_{\bar p}(\mathcal E)$.  Since the category of $IC$ sheaves is a full
subcategory of the derived category, the theorem also says that
\[ \RHom(\uul{IC}{}^{\bullet}(\mathcal E_1), \uul{IC}{}^{\bullet}(\mathcal E_2)) \cong \Hom(\mathcal E_1, \mathcal E_2).\]
whenever $\mathcal E_1, \mathcal E_2$ are local systems on $U_2$.
If $\mathcal E$ is an indecomposable  local system (which is not isomorphic to a direct sum of two nontrivial local systems)
then $\uul{IC}{}^{\bullet}_{\bar p}(\mathcal E)$ is an indecomposable complex of sheaves (and is not isomorphic to a
direct sum of two nontrivial complexes of sheaves).

For a perversity $\bar p$ let $p^{-1}(t) = \min \{c| \ p(c) \ge t\}$ and $p^{-1}(t) = \infty$ if $t>p(n)$.
We can reformulate these conditions (2,3,4) in a way that does not refer to a particular stratification as follows: 
\begin{enumerate}
\item[(S1)] $\dim\{x \in W| \ H^r(i_x^*A\b) \ne 0 \} \le n - p^{-1}(t)$ for all $t>0$
\item[(S2)] $\dim \{ x \in W| \ H^r(i_x^!A\b) \ne 0 \} \le n-q^{-1}(n-t)$  for all $t<n$. $\phantom{G^{\strut}}$
\end{enumerate} 
As above, the condition (S2) is the Verdier dual of condition (S1) and may be expressed as
\begin{enumerate}
\item[(S2')] $\dim \{ x \in W| \ H^r(i_x^*\mathbf D(A\b)) \ne 0 \} \le n-q^{-1}(t)$ for all $t>0$.
\end{enumerate}

\part{Perverse sheaves}

\section{Perverse sheaves}\index{sheaf!perverse}\index{perverse!sheaf}
\subsection{}
Let $W$ be a Whitney stratified space with a given stratification.  
We have two notions of the constructible derived category: (1)
as complexes of sheaves that are cohomologically constructible with respect to the given stratification
(2) as complexes of sheaves that are cohomologically constructible with respect to some stratification.
In order to reduce the total number of words in these notes, we shall simply refer to ``the constructible
derived category", meaning either one of these two possibilities.


The category of perverse sheaves is defined by relaxing the support and cosupport 
conditions (\S \ref{thm-equivalence})  for the IC sheaf by one.  Here is the precise
definition in the case of middle perversity:

\subsection{Definition} \label{subsec-defperverse}
 Let $W$ be a $n$-dimensional Whitney stratified (or stratifiable) space that can be stratified with strata of even
dimension.  Let $K$ be a field.
A {\em middle perversity} perverse sheaf on  $W$ is a complex of sheaves $A\b$ in the bounded constructible derived category
$D^b_c(W)$ of $K$-vector spaces, with the following property.
If $S$ is a stratum of dimension $d$  let $j_S:S\to W$ and $j_x:\{x\} \to W$.
Then
\begin{align*}
H^i(j_x^*A\b) = 0 &\text{ for all}\ i>(n-d)/2\\
H^i(j_x^!A\b)=0 &\text{ for all}\ i<(n+d)/2\\
\intertext{Since $H^a(j_x^!A\b) \cong H^{a-d}(j_S^!A\b)_x$ the second condition
may be expressed as:}
H^i(j_S^!A\b) = 0 &\text{ for all}\ i<(n-d)/2
\end{align*}
Equivalently, $A\b$ is middle perverse if and only if for all $r \ge 0$,
\begin{enumerate}
\item[(P1)] $\dim\{x \in W| \ H^r(i_x^*A\b) \ne 0 \} \le n-2r$ and
\item[(P2)] $\dim \{ x \in W| \ H^r(i_x^!A\b) \ne 0 \} \le 2r-n$.   $\phantom{G^{\strut}}$
\end{enumerate}

\begin{figure}[!h]
\begin{tabular}{|c||c|c|c|c|c|}
\hline
$i$ & $\cod 0$ & $\cod 2$ & $\cod 4$ & $\cod 6$ & $\cod 8$\\
\hline\hline
8 &        c        &      c          &        c       &     c         &       c        \\
7 &                 &                 &       c        &     c         &         c      \\
6 &                 &                 &                &      c        &        c       \\
5 &                 &                 &                &               &          c     \\
4 &                 &                 &                &               &         \rz       \\
3 &                 &                 &                &               &      x         \\
2 &                 &                 &                &        x      &    x           \\
1 &                 &                 &       x        &       x       &   x            \\
0 &       x        &      x          &     x         &       x      &      x         \\
\hline 
&\multicolumn{5}{|c|}{$\uul{IC}^{\bullet}$ support}\\
\hline
\end{tabular}\qquad \begin{tabular}{|c||c|c|c|c|c|}
\hline
$i$ & $\cod 0$ & $\cod 2$ & $\cod 4$ & $\cod 6$ & $\cod 8$\\
\hline\hline
8 &        c        &       c         &       c        &        c      &       c        \\
7 &                 &         c       &      c         &       c       &       c        \\
6 &                 &                 &         c      &      c        &        c       \\
5 &                 &                 &                &       c       &        c        \\
4 &                 &                 &                &               &        cx        \\
3 &                 &                 &                &         x     &       x        \\
2 &                 &                 &         x      &        x      &         x      \\
1 &                 &        x        &        x       &       x       &      x         \\
0 &       x         &       x         &        x       &     x         &     x           \\
\hline 
&\multicolumn{5}{|c|}{Perverse sheaf support}\\
\hline\end{tabular}
\caption{Stalk and co-stalk cohomology of $\uul{IC}^{\bullet}$ and perverse sheaves}
\label{fig-cospt}
\end{figure}
In these figures, ``x" denotes regions of possibly nontrivial stalk cohomology and ``c" denotes 
regions of possibly nontrivial stalk cohomology with compact support, $H^r(j_x^!A\b)$.
 \begin{tcolorbox}[colback=yellow!30!white]
\begin{defn}
The {\em category of perverse sheaves} (with middle perversity) is the full subcategory of $D^b_c(W)$ 
whose objects are (middle) perverse sheaves.
\end{defn}\end{tcolorbox}
 \begin{tcolorbox}[colback=cyan!30!white]
\begin{thm}\cite{BBDG}\label{abeliancategory}
The category of (middle) perverse sheaves forms an abelian subcategory of the
derived category $D^b_c(W)$ that is preserved  by Verdier duality.
\end{thm}\end{tcolorbox}
Kernels and cokernels in this category can be described, see \S \ref{subsec-perversecohomology}.  If $W$ is an 
algebraic variety, the simple objects are the shifted $IC$ sheaves with irreducible local coefficients of
irreducible subvarieties.

\subsection{Other perversities}  
There is an abelian category of perverse sheaves for any perversity and even for the generalized perversities considered
in \cite{BBDG}, that is, for an arbitrary function $\bar p$ from $\{\text{strata}\} \to \ZZ$.  For the case of
perversity zero and top see \S \ref{subsec-perversityzero}.
\quash{
For any perversity $\bar p$ the category of $\bar p$-perverse sheaves on $W$ is the
full subcategory of the category of constructible sheaves such that for each stratum $j_S:S \to W$,
\begin{align}\begin{split}
H^r(j_S^*(A\b)) = 0 &\text{ for all}\ r \ge ??\\
H^r(j_S^!(A\b)) = 0 &\text{ for all} r \le ??.
\end{split}\end{align}
} 

\subsection{Historical comment} 
 In 1979 David Kazhdan and George Lusztig \cite{Kazhdan}
posed two conjectures involving representations of Hecke algebras, Verma modules,  and a collection of polynomials
which have come to be known as Kazhdan Lusztig polynomials.
They realized that these questions were related to the failure of Poincar\'e duality for Schubert varieties, which
they had analyzed.  Lusztig had heard
MacPherson lecture about intersection homology and he suspected a relation with the KL polynomials. 
Separately, Raoul Bott suggested to Kazhdan that
intersection homology might be relevant and that he and Lusztig should speak with MacPherson,
at the same time telling MacPherson that he should speak with Kazhdan and Lusztig.
But Kazhdan and Lusztig's computations for Schubert varieties
did not agree with an example in a preprint of Goresky-MacPherson.

In the ensuing conversation MacPherson corrected the error in the preprint and
directed them to contact Deligne who had found a new sheaf-theoretic way to define intersection homology, which
worked in the \'etale setting.  Using Deligne's construction, Kazhdan and Lusztig  were then able to 
prove \cite{Kazhdan_duality}  
as they suspected, that the Kazhdan Lusztig polynomials coincide with the intersection cohomology 
local Poincar\'e polynomial of one Schubert variety at a point in another Schubert cell, 
thereby also proving that the coefficients were nonnegative.  This paper may be viewed as
a step in the proof the K.L. conjectures.   
(See also Lusztig's more extensive comments in 
 \cite{Lusztigarxiv1409, Lusztigarxiv1707}.)
 
 The remaining steps in the proof of the Kazhdan Lusztig conjectures were eventually proven 
by A. Beilinson and J. Bernstein \cite{Beilinson}
and independently by J. L. Brylinski and M. Kashiwara \cite{Brylinski}, by making the
connection between intersection homology and the vast existing theory of $\mathcal D$-modules.  
On an algebraic manifold (in this case, the flag manifold) there is a ring $\mathcal D$ (or rather, a sheaf of rings) 
of differential operators (and an extensive literature, see 
\cite{BorelDModules, Maisonobe, Gelfand, Kashiwara, KashiwaraK2, IndexTheorem, KS, KS2, Mebkhout}).  
To each $\mathcal D$-module there corresponds a sheaf of solutions, which is a constructible sheaf.
 Beilinson, Bernstein, Brylinski and Kashiwara  showed that each Verma module can be associated to a 
 certain holonomic $\mathcal D$-module with regular singularities, whose sheaf of solutions turns out to 
be the IC sheaf.  This provided the link between the Kazhdan-Lusztig polynomials and Verma modules.  

However, the category of $\mathcal D$-modules is an abelian category, whereas the (derived) category of constructible sheaves is not abelian, so it was conjectured that there must correspond an abelian subcategory of the derived category that ``receives" the solution sheaves of individual $\mathcal D$-modules. Thus was born the category of perverse sheaves, 
with middle perversity.  On the other hand, intersection homology is a topological invariant, so then the question arose 
as to whether this category of perverse sheaves could be constructed purely topologically, and for other perversities as well.  
The book  \cite{BBDG} completely answers this question, giving a very general setting in which the category of perverse sheaves, an abelian subategory of the derived category, could be constructed.

\section{Examples of perverse sheaves}
\subsection{IC of subvarieties}
\begin{figure}[!h]
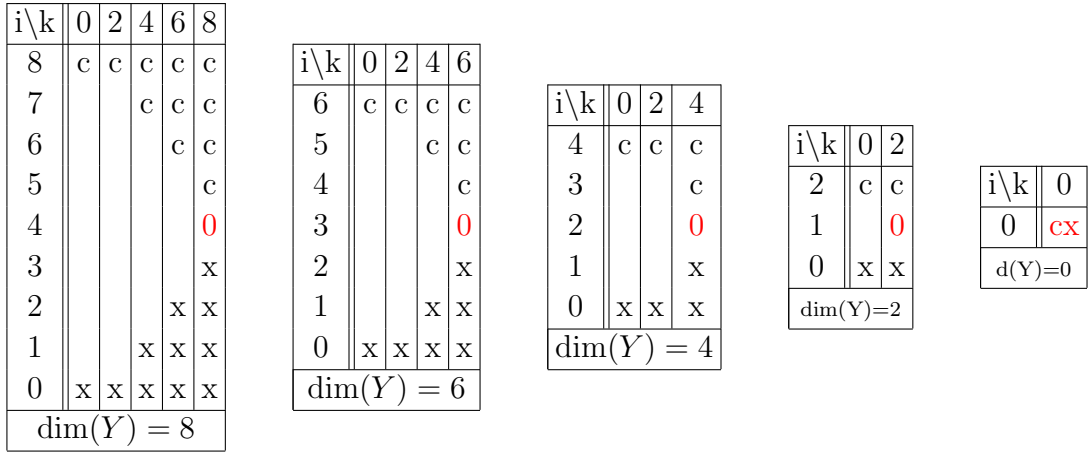

\begin{tabular}{|c||c|c|c|c|c|c|c||c|c|c|c|c|c||c|c|c|c|c||c|c|c|c||c|}
\cline{1-6}
i$\backslash$k&0& 2& 4& 6&8 & \multicolumn{17}{}{} \\
\cline{1-6}\cline{8-12}
8 &c&c&c&c&c&$\phantom{GH}$&i$\backslash$k& 0&2&4&6&\multicolumn{11}{}{} \\
\cline{8-12}\cline{14-17}
7&&&c&c&c&&6&c&c&c&c&$\phantom{GH}$&i$\backslash$k&0&2&4&\multicolumn{6}{}{}\\
\cline{14-17}\cline{19-21}
6&&&&c&c&&5&&&c&c&& 4&c&c&c&$\phantom{GH}$&i$\backslash$k&0&2& \multicolumn{3}{}{}\\
\cline{19-21}\cline{23-24}
5&&&&&c&&4&&&&c&& 3&&&c&&2&c&c&$\phantom{GH}$&i$\backslash$k&0\\
\cline{23-24}
4&&&&&$\rz$&&3&&&&$\rz$&&2&&&$\rz$&&1&&$\rz$&&0&\textcolor{red}{cx}\\
\cline{23-24}
3&&&&&x&&2&&&&x&&1&&&x&&0&x&x&&\multicolumn{2}{|c|}{\tiny{d(Y)=0}}   \\
\cline{19-21}\cline{23-24}
2&&&&x&x&&1&&&x&x&&0&x&x&x&& \multicolumn{3}{|c|}{{\tiny {dim(Y)=2}}}&\multicolumn{2}{}{}  \\
\cline{14-17}\cline{19-21}
1&&&x&x&x&&0&x&x&x&x&& \multicolumn{4}{|c|}{$\dim(Y)=4$}&\multicolumn{6}{}{} \\
\cline{8-12}\cline{14-17}
0&x&x&x&x&x&& \multicolumn{5}{|c|}{$\dim(Y)=6$}&\multicolumn{11}{}{} \\
\cline{1-6}\cline{8-12}
\multicolumn{6}{|c|}{$\dim(Y)=8$} & \multicolumn{17}{}{} \\
\cline{1-6}
 \end{tabular}
\caption{Shifted IC of subvarieties}
\label{fig-perverse-spt}
\end{figure}

As above we consider the middle perversity $\bar m$ and a Whitney stratified space of dimension $n$ with even 
codimension ($= k$ in Figure \ref{fig-perverse-spt}) strata.  Let $Y$ denote the closure of a single stratum, 
$Y^o$.  Let $\mathcal E_Y$ be a local
system on the stratum  $Y^o$.  Then the intersection complex 
$\uul{IC}{}^{\bar m}_Y(\mathcal E_Y)[-\cod(Y)/2]$ is $\bar m$-perverse.  
Here are the support diagrams for an 8 dimensional stratified space with strata of codimension 0, 2, 4, 6, 8 where,
as above, ``x" denotes possibly nonzero stalk cohomology and ``c" denotes possibly nonzero stalk cohomology with compact support.  Adding these up gives the support diagram (Figure \ref{fig-perverse-spt}) for a perverse sheaf. 

\subsection{}  Let $Y^o$ be a stratum of $W$ (which is stratified by even dimensional
strata).  Let $Y$ be its closure with inclusion $j_Y:Y \to W$.  It is stratified by even dimensional strata.
Let $A\b$ be a perverse sheaf on $Y$.  Then $Rj_*(A\b)[-\cod(Y)/2]$ is a perverse sheaf on $W$.

\subsection{Logarithmic perversity}\label{subsec-logarithmic}\index{perversity!logarithmic}
Because the support conditions for (middle) perverse sheaves are relaxed slightly from 
those for $\uul{IC}^{\bullet}$, there are several other perversities for which intersection cohomology forms
a (middle) perverse sheaf.  These include the {\em logarithmic} perversity $\bar{\ell}$, given by
$\bar\ell(k) = k/2=\bar m(k)+1$ and its Verdier dual, the {\em sublogarithmic} perversity, $\bar s$ given by $\bar s(k) = \bar m(k) -1$.
See Figure \ref{fig-logarithmic}.

\begin{figure}[!h]
\begin{tabular}{|c||c|c|c|c|c|}
\hline
$i$ & $\cod 0$ & $\cod 2$ & $\cod 4$ & $\cod 6$ & $\cod 8$\\
\hline\hline
8 &        c        &      c          &        c       &     c         &       c        \\
7 &                 &                 &                &     c         &         c      \\
6 &                 &                 &                &              &        c       \\
5 &                 &                 &                &               &         \rz    \\
4 &                 &                 &                &               &         x       \\
3 &                 &                 &                &         x      &      x         \\
2 &                 &                 &       x       &        x      &    x           \\
1 &                 &                 &       x        &       x       &   x            \\
0 &       x        &      x          &     x         &       x      &      x         \\
\hline 
&\multicolumn{5}{|c|}{$\uul{IC}{}^{\bullet}_{\bar{\ell}}$ support}\\
\hline
\end{tabular}\qquad \begin{tabular}{|c||c|c|c|c|c|}
\hline
$i$ & $\cod 0$ & $\cod 2$ & $\cod 4$ & $\cod 6$ & $\cod 8$\\
\hline\hline
8 &        c        &       c         &       c        &        c      &       c        \\
7 &                 &         c       &      c         &       c       &       c        \\
6 &                 &                  &         c      &       c       &        c       \\
5 &                 &                  &                 &       c       &        c        \\
4 &                 &                  &                 &               &        c       \\
3 &                 &                  &                 &              &       \rz        \\
2 &                 &                  &                 &              &         x      \\
1 &                 &                  &                 &       x       &      x         \\
0 &       x         &       x         &        x       &     x         &     x           \\
\hline 
&\multicolumn{5}{|c|}{$\uul{IC}{}^{\bullet}_{\bar{s}}$ support}\\
\hline\end{tabular}
\caption{Support/cosupport of logarithmic and sublogarithmic IC sheaf}
\label{fig-logarithmic}
\end{figure}

\subsection{Hyperplane complements}\index{hyperplane complement} (see also \cite{Kapranov})
Let $\{H_1,H_2,\cdots, H_r\}$ be a collection of complex affine hyperplanes in $W = \CC^n$.  Stratify
 $W$ according to the multi-intersections of the hyperplanes.  The largest stratum is
 \[ W^o = W - \bigcup_{j=1}^rH_j\]
 and it may have a highly nontrivial fundamental group.  let $\mathcal E$ be a local coefficient system on this
 hyperplane complement.  Then $\uul{IC}{}^{\bullet}_{\bar s}(\mathcal E)$, $\uul{IC}{}^{\bullet}_{\bar m}
 (\mathcal E)$ and $\uul{IC}{}^{\bullet}_{\bar{\ell}}(\mathcal E)$ are perverse sheaves on $W = \CC^n$.
 These are surprisingly complicated objects, and even the case of middle perversity, when the
 hyperplanes are the coordinate hyperplanes, has been extensively studied.  Notice, in this case, that
 the space $W = \CC^n$ is nonsingular, the hyperplane complement $W^o$ is nonsingular, and the 
 sheaf $\uul{IC}(\mathcal E)$ is constructible (with respect to this chosen stratification) but to analyze
 this sheaf we are forced to consider the singularities of the multi-intersections of the hyperplanes.
 
 In the simplest case, $(\CC, \{0\})$ the category of perverse sheaves is equivalent to the category of
 representations of the following quiver
  \begin{diagram}
 {\bullet} & \pile{\rTo^{\alpha} \\ \lTo_{\beta} } & \bullet
 \end{diagram}
 where $I-\alpha\beta$ and $I-\beta\alpha$ are invertible.  
 
 For $\CC^2,\ xy=0$ (the coordinate axes) the perverse category is equivalent to the category of representations of
 the quiver   
 \begin{diagram}
 {\bullet} & \pile{\rTo \\ \lTo} & {\bullet}\\
 \uTo\dTo && \uTo\dTo \\
 \bullet & \pile{\rTo \\ \lTo} & \bullet
 \end{diagram}
 with similar conditions ($I-\alpha\beta$ and $I-\beta\alpha$ invertible)
 on each of the horizontal and vertical pairs, such that all possible ways around the
 outside of the square commute.
 
\subsection{Small and semismall maps}\index{small map}\index{semismall map}\label{subsec-smallmaps}
Let $M$ be a compact complex algebraic manifold and let $\pi:M \to W$ be an algebraic mapping.
Then $\pi$ is said to be {\em semismall} \cite{BorhoMacPherson} if
\[ \cod_W (\{ x \in W|\ \dim \pi^{-1}(x) \ge k\}) \ge 2k.\]
  In other words, if the map has been stratified then for each stratum $S \subset W$
the dimension of the fiber over $S$ is  $\le \frac{1}{2}$ the codimension of $S$. The map is
{\em small} \cite{IH2} if, for each singular stratum $S$, $\dim \pi^{-1}(x) < \frac{1}{2}\cod(S)$ (for all $x \in S$).

If $\pi$ is small then $R\pi_*(\uul{\QQ})$ is a self dual sheaf on $W$ whose support satisfies the
support conditions of (middle) intersection cohomology.  It follows from the axiomatic characterization
that there is a canonical isomorphism (in $D^b_c(W)$), $R\pi_*(\uul{\QQ}) \cong \uul{IC}{}_{\bar m}(W)$.
In other words, the intersection cohomology of $W$ is canonically isomorphic to the ordinary cohomology
of $M$.

If $\pi$ is semi-small then $R\pi_*(\uul{\QQ})$ is (middle) perverse.

Let $W = \left\{ P \subset \CC^4 |\ \dim(P)=2, \dim(P\cap \CC^2) \ge 1 \right\}$ be the singular Schubert variety in the
Grassmannian of 2-planes in 4-space.  It has a singularity when $P = \CC^2$. A resolution of singularities is
$\widetilde{W}= \ \left\{ (P,L)|\ P \in W,\ \text{ and } L \subset P \cap \CC^2 \subset \CC^4 \right\}.$  Then 
$\pi:\widetilde{W}\to W$ is a small map so $R\pi_*(\QQ) \cong \uul{IC}{}^{\bullet}_W$ hence 
$IH^*(W) \cong H^*(\widetilde{W})$.

\quash{
\makeatletter
\renewcommand*\env@matrix[1][*\c@MaxMatrixCols c]{%
  \hskip -\arraycolsep
  \let\@ifnextchar\new@ifnextchar
  \array{#1}}
\makeatother
}

\subsection{Sheaves on $\mathbb{P}^1$} (See \cite{MacPhersonVilonen}.)
Let us stratify $\mathbb{P}^1$ with a single zero dimensional stratum, $N$ (the north pole, say).  The support
diagram for middle perversity sheaves is the following:
\begin{center}
\begin{tabular}{|c||c|c|}
\hline
$i\backslash\cod$ & 0 &  2 \\
\hline
2 & c & c \\
1 &  & cx \\
0 & x & x\\
\hline
\end{tabular}\end{center}
The skyscraper sheaf supported at the point, $\uul{\QQ}{}_N[-1]$ is perverse.  We also have the following:
\[
\uul{\QQ}{}_N[-1]:  \text{  \begin{tabular}{|c|||c|c|}\hline 2&\phantom{c}& \\1& &cx  \\ 0& &  \\ \hline \end{tabular}} \qquad
\uul{\QQ}{}_{\mathbb P^1}: \text{\begin{tabular}{|c|||c|c|}\hline 2& c&c \\1& &  \\0& x & x \\ \hline \end{tabular}} \qquad
j_!(\uul{\QQ}{}_U):  \text{\begin{tabular}{|c|||c|c|}\hline 2& c & c\\ 1&& c \\0& x & \\ \hline\end{tabular}}\qquad j_*(\uul{\QQ}{}_U):
\text{\begin{tabular}{|c|||c|c|}\hline 2&c& \\1& & x \\0& x & x \\ \hline \end{tabular}  }
\]
The first sheaf is self dual.  The second sheaf is self dual.  The third and fourth sheaves are dual to each other.
It turns out that there is one more indecomposable perverse sheaf on this space, which is not an IC sheaf, and its
support diagram is the full diagram.  It is self dual.  Here is how to construct it.  Take a closed disk and put the
constant sheaf on the interior, zero on the boundary, except for one point (or even one segment).  Then map
this disk to the 2-sphere, collapsing the boundary to the N pole, and push this sheaf forward.  

If we started with zero on the boundary and pushed forward we would get the sheaf $Rj_!(\QQ_U)$.  If we started with
the full constant sheaf on the disk and pushed forward we would get the sheaf $Rj_*(\QQ)$.  This new sheaf has both
stalk cohomology and compact support stalk cohomology in degree 1, at the singular point.  Verdier duality switches
these two types of boundary conditions, so when we have a mixed boundary condition as in this case, we obtain a self dual sheaf.

In this case the category of perverse sheaves is equivalent to the category of representations of the quiver
\begin{diagram}
\bullet& \pile{\rTo^{\alpha}\\ \lTo_{\beta}} & \bullet \end{diagram}
where $\alpha \beta = \beta \alpha = I$.  There are five indecomposable objects, one of which is has $\QQ \oplus \QQ$
on one of the vertices of the graph.

\subsection{Homological stratifications}  \index{stratification!homological}
Although we have assumed the space $W$ is Whitney stratified, all the
preceding arguments remain valid for spaces with a homological stratification.  
Homological stratifications were used to prove (\cite{IH2}) the topological invariance of
intersection homology and are also required in the \'etale setting. 

Suppose $W$ is a locally compact
Hausdorff space with a locally finite decomposition into locally closed subsets called ``strata'' which satisfy
the condition of the frontier:  the closure of each stratum is a union of strata.  If $T \subset \overline{S}$ write
$T<S$.  
\begin{tcolorbox}[colback=yellow!30!white]
Such a decomposition is a homological stratification of $W$ provided
\begin{itemize}
\item the strata are topological manifolds
\item There is an open dense stratum $U \subset W$ so that $W-U$ has codimension $\ge 2$
\item For each stratum $j:S \to W$ and for every (finite dimensional) local coefficient system $\mathcal E$ 
on $S$ and for every $i$, the
sheaf $R^ij_*(\mathcal E) = \uul{H}^i(Rj_*(\mathcal E))$ is constructible, that is, locally constant 
and finite dimensional on every stratum $T<S$.
\end{itemize}
\end{tcolorbox}  
If we restrict the coefficients to lie in a field $k$ (so that the local coefficient systems are systems of
$k$-vector spaces) then condition (1) may be replaced by the weaker assumption that each stratum $S$ is a
$k$-homology manifold\footnote{meaning that $H_r(S, S-x; k) = 0$ for $r \ne \dim(S)$ and $H_{\dim(S)}(S, S-x;k) \cong k$
for all $x \in S$.}, as required for arguments involving Poincar\'e-Verdier duality.

\subsection{Perversity zero}\label{subsec-perversityzero}\index{perversity!zero} (See also \S \ref{subsec-t=0}.)
  Let $W$ be a stratified pseudomanifold of dimension $n$ (with a fixed stratification).
The category of Perverse sheaves on $W$ with perversity zero, constructible with respect to this stratification, is
equivalent to the category of sheaves on $W$  that are constructible with respect to this stratification, that is,
 sheaves whose restriction
to each stratum is locally trivial.  In this case, the ``abelian subcategory"
defined by the perversity condition simply coincides with the abelian category structure of the category of (constructible)
sheaves.  Here is the support diagram for $IC^{\bar 0}$ of a six dimensional pseudomanifold, where ``x'' denotes
possibly nonzero cohomology and ``c'' denotes possibly nonzero compact support cohomology:

\begin{figure}[!h]
\begin{tabular}{|c||c|c|c|c|c|c|c|}
\hline
$i$ & $\cod 0$ & $\cod 1$ & $\cod 2$ & $\cod 3$ & $\cod 4$ &$\cod 5$ & $\cod 6$\\
\hline\hline
6 &        c        &      c          &        c       &     c         &       c   &c     &c     \\
5 &                 &                 &               &     c         &         c    &c      &c     \\
4 &                 &                 &                &               &         c    &c     &c   \\
3 &                 &                 &                 &               &            &c    &c   \\
2 &                 &                 &                 &              &            &         &c  \\
1 &                 &                 &                 &              &            &          &\rz          \\
0 &       x        &      x          &     x         &       x      &      x        &x    &x    \\
\hline 
&\multicolumn{7}{|c|}{$\uul{IC}{}^{\bullet}_{\bar{0}}$ support}\\
\hline
\end{tabular}
\end{figure}
The support diagram (Figure \ref{fig-0-spt})  for the constant sheaf is the same if the link of every stratum
is connected.  In general the ``c'' in the bottom of each column in the following 
support diagram is given by the reduced cohomology $h^0(L)$ of the link of
the stratum.

\begin{figure}[!h]
\begin{tabular}{|c||c|c|c|c|c|c|c|}
\hline
$i$ & $\cod 0$ & $\cod 1$ & $\cod 2$ & $\cod 3$ & $\cod 4$ &$\cod 5$ & $\cod 6$\\
\hline\hline
6 &        c        &      c          &        c       &     c         &       c   &c     &c     \\
5 &                 &                 &        c       &     c         &         c    &c      &c     \\
4 &                 &                 &                &       c        &         c    &c     &c   \\
3 &                 &                 &                 &               &       c     &c    &c   \\
2 &                 &                 &                 &              &            &   c      &c  \\
1 &                 &                 &                 &              &            &          &c         \\
0 &       x        &      x          &     x         &       x      &      x        &x    &x    \\
\hline 
&\multicolumn{7}{|c|}{constant sheaf support}\\
\hline
\end{tabular}\qquad \begin{tabular}{|c||c|c|c|c|c|c|c|}
\hline
$i$ & $\cod 0$ & $\cod 1$ & $\cod 2$ & $\cod 3$ & $\cod 4$ & $\cod 5$   & $\cod 6$  \\
\hline\hline
6 &        c        &       c         &       c        &        c      &       c     &c     &c   \\
5 &                 &        c        &      c         &       c       &       c     &c     &c   \\
4 &                 &                  &       c         &       c       &        c    &c     &c   \\
3 &                 &                  &                 &       c       &     c         &c     &c       \\
2 &                 &                  &                 &              &       c       &   c    &c      \\
1 &                 &                  &                 &              &              &     c   &   c       \\
0 &       x         &       x         &        x       &     x         &     x     & x      &  xc      \\
\hline 
&\multicolumn{7}{|c|}{perverse sheaf support}\\
\hline\end{tabular}
\caption{Support/cosupport for constant sheaf and perversity $\bar0$ sheaf}
\label{fig-0-spt}
\end{figure}

If $j:X \to W$ is the inclusion of a stratum and $E$ is a local system on $X$ then $j_!E$ is a perversity zero sheaf on $W$.
 The second table in Figure \ref{fig-0-spt} may be obtained by ``summing'' support diagrams for sheaves of the form
$j_!(E)$ of strata of dimension 0, 1, 2, 3, 4, 5 and 6 (with no shift).  A ``pseudomanifold'' has no strata of codiimension one but
in these diagrams we have included the possibility of codimension one strata to allow for pseudomanifolds with
boundary, and to allow for stratification by the simplices of a triangulation.

\subsection{Top perversity}  \index{perversity!top}
For $\bar p(S) =\bar t(S)= \codim(S) -2$ the support and cosupport diagram for a
6 dimensional stratified space are shown
in Figure \ref{fig-top-IC} and \ref{fig-top-spt}.  We have shifted the degrees by 6 so that the dualizing
sheaf (the sheaf of Borel Moore chains) is $\bar t$ perverse.

\begin{figure}[!h]
\begin{tabular}{|c||c|c|c|c|c|c|c|}
\hline
$i$ & $\cod 0$ & $\cod 1$ & $\cod 2$ & $\cod 3$ & $\cod 4$ &$\cod 5$ & $\cod 6$\\
\hline\hline
0 &        c        &      c          &        c       &     c         &       c   &c      &c     \\
-1 &                 &                 &               &                &               &         &  0     \\
-2 &                 &                 &                &               &               &         &  x    \\
-3 &                 &                 &                 &               &              &    x    &  x   \\
-4 &                 &                 &                 &              &        x     &    x    &  x  \\
-5 &                 &                 &                 &       x      &        x     &   x     &  x   \\
-6 &       x        &      x          &     x         &       x      &      x        &   x    &  x    \\
\hline 
&\multicolumn{7}{|c|}{$\uul{IC}{}^{\bullet}_{\bar{t}}$ support}\\
\hline
\end{tabular}\qquad
\caption{Support/cosupport for $\bar p = \bar t$ sheaf}\label{fig-top-IC}
\end{figure}

 For the sheaf of chains,
 the uppermost ``x'' in each column Figure \ref{fig-top-spt}  represents the reduced homology $h_0(L)$ of the link
 of the stratum.  If $j:X \to W$ is a stratum and $E$ is a local system on $X$ then $Rj_*(E)[\dim(X)]$ is a sheaf
 on $W$ with top perversity.

\begin{figure}[!h]
\begin{tabular}{|c||c|c|c|c|c|c|c|}
\hline
$i$ & $\cod 0$ & $\cod 1$ & $\cod 2$ & $\cod 3$ & $\cod 4$ &$\cod 5$ & $\cod 6$\\
\hline\hline
0 &        c        &      c          &        c       &     c         &       c   &c      &c     \\
-1 &                 &                 &               &                &               &         &  x     \\
-2 &                 &                 &                &                &               &     x    &  x    \\
-3 &                 &                 &                 &               &         x     &    x    &  x   \\
-4 &                 &                 &                 &        x      &        x     &    x    &  x  \\
-5 &                 &                 &      x          &       x      &        x     &   x     &  x   \\
-6 &       x        &      x          &     x         &       x      &      x        &   x    &  x    \\
\hline 
&\multicolumn{7}{|c|}{$\uul{\bf \DD}^{\bullet}$ support}\\
\hline
\end{tabular}\qquad
\begin{tabular}{|c||c|c|c|c|c|c|c|}
\hline
$i$ & $\cod 0$ & $\cod 1$ & $\cod 2$ & $\cod 3$ & $\cod 4$ &$\cod 5$ & $\cod 6$\\
\hline\hline
0 &        c        &      c          &        c       &     c         &       c   &c      &cx     \\
-1 &                 &                 &               &                &               &      x   &   x    \\
-2 &                 &                 &                &               &         x     &    x     &  x    \\
-3 &                 &                 &                 &        x      &       x       &    x    &  x   \\
-4 &                 &                 &       x         &       x       &        x     &    x    &  x  \\
-5 &                 &      x          &      x          &       x      &        x     &   x     &  x   \\
-6 &       x        &      x          &     x         &       x      &      x        &   x    &  x    \\
\hline 
&\multicolumn{7}{|c|}{perverse sheaf support}\\
\hline
\end{tabular}\qquad
\caption{Support/cosupport for the Borel-Moore sheaf and perversity $\bar t$ sheaf}
\label{fig-top-spt}
\end{figure}
Suppose $A\b$ is a top-perverse sheaf with on an $n$-dimensional pseudomanifold $W$.  Then
$A\b$ may be interpreted as a {\em cosheaf}, \index{cosheaf} namely, the covariant functor on the
category of open sets,
\begin{equation}\label{eqn-cosheaf}
U \mapsto H^0_c(U;A\b)\end{equation}
which assigns  the group homomorphism
$H^0_c(U;A\b) \to H^0_c(V;A\b)$ to each inclusion $U \subset V$ in a functorial way.  (Recall that
we have introduced a shift by $n$, the dimension of the space.  With no shift, this cosheaf would be 
$H^n_c(U;A\b)$.)  The ``co-stalk'' at  a point $j_x:\{x\}\in W$ is the limit over neighborhoods $U_x$
of $x$:
\[H^0(j_x^!A\b) = \underset{\longleftarrow}{\lim}H^0_c(U_x;A\b).\]  
The support conditions (see Figure \ref{fig-top-spt}) imply that the co-stalk cohomology vanishes except
in this single degree and the cosheaf (\ref{eqn-cosheaf}) determines $A\b$ up to
quasi-isomorphism.

\subsection{BBDG numbering system}
In their book \cite{BBDG} the authors modified the indexing system for cohomology in a way that vastly reduces the amount of notation and arithmetic involving indices.  Although the new system is extremely simple, it is deceptively so, because it takes us one step further away from any intuition concerning perverse sheaves.  The new system
works best in the case of a complex algebraic (or analytic) variety $W$, stratified with complex algebraic (or analytic)
strata, and counted according to their complex dimensions.  The idea is to shift all degrees by 
$\dim_{\CC}(W) = \dim(W)/2$ because cohomology is symmetric about this point.  So the support conditions look like this:

\begin{center}\begin{figure}[H]
\begin{tabular}{|c||c||c|c|c|c|c|}
\hline
new & old & $\cod 0$ & $\cod 2$ & $\cod 4$ & $\cod 6$ & $\cod 8$\\
\hline\hline
4& 8 &        c        &       c         &       c        &        c      &       c        \\
3& 7 &                 &         c       &      c         &       c       &       c        \\
2& 6 &                 &                  &         c      &       c       &        c       \\
1& 5 &                 &                  &                 &       c       &        c        \\
0& 4 &                 &                  &                 &               &        cx      \\
-1& 3 &                 &                  &                &     x         &      x       \\
-2& 2 &                 &                  &      x        &      x        &      x      \\
-3&  1 &                 &        x       &       x       &       x       &      x         \\
-4&  0 &       x         &       x      &        x       &     x         &     x           \\
\hline 
&&\multicolumn{5}{|c|}{Perverse sheaf support}\\
\hline\end{tabular}\caption{Deligne's numbering for middle perversity}\label{fig-newnumbers}\end{figure}\end{center}
If $S$ is a stratum of dimension $d$ and codimension $c$ ($d+c = n = \dim(W)$) let $j_S:S\to W$ and $j_x:\{x\} \to W$.
Then subtracing $n/2$ from the indices in \S \ref{subsec-defperverse} gives support conditions for a middle
perverse sheaf $A\b$ in the new indexing scheme:
\begin{align}\begin{split}\label{eqn-newspt}
H^i_{new}(j_S^*A\b) = 0 &\text{ for all}\ i > -d/2\\
H^i_{new}(j_S^!A\b)=0 &\text{ for all}\ i< -d/2
\end{split}\end{align}

Changing the degree of cohomology implies a corresponding change in the perversity function.
The authors of \cite{BBDG} chose to further simplify the numerology by removing the ``$-2$'' that
occurs in IC calculations.  They define the middle perversity to take values $p(S) = - \dim(S)/2$
for every even dimensional stratum $S$ so that equations (\ref{eqn-newspt}) become
\begin{align*}\begin{split}
H^i_{new}(j_S^*A\b) = 0 &\text{ for all}\ i > p(S)\\
H^i_{new}(j_S^!A\b) = 0 &\text{ for all}\ i < p(S).
\end{split}\end{align*} 
The $IC$ sheaf is defined by replacing these strict inequalities by $\le$ and $\ge$.

\subsection{}\label{subsec-numbering}
More generally the authors of \cite{BBDG} consider a perversity to be a function of {\em dimension}:
\begin{tcolorbox}[colback=yellow!30!white]
 A perversity is an integer valued function $p:\ZZ_{\ge 0} \to \ZZ_{\le 0}$ with $p(0) = 0$ and
\begin{equation}\label{eqn-p-inequalities}
p(d) \ge p(d+1) \ge p(d)-1.\end{equation}\end{tcolorbox}\noindent
(See diagram in \S \ref{subsec-perverse-dimension}.)
The category of perverse sheaves is defined to be
(abelian) full subcategory of $D^b_c(W)$ consisting of complexes $A\b$ so that for any stratum $S$ (except
the largest stratum),
\begin{itemize}
\item $H^n(i_x^*(A\b)) = 0$ for all $n > p(S)$
\item $H^n(i_x^!(A\b)) = 0$ for all $n < p(S)+\dim(S)$
\end{itemize}
for some (and hence any) $x \in S$ while the IC sheaf satisfies 
\begin{itemize}
\item $H^n(i_x^*(\uul{IC}\b)) = 0$ for all $n \ge p(S)$
\item $H^n(i_x^!(\uul{IC}\b)) = 0$ for all $n \le p(S) + \dim(S)$.
\end{itemize}
Each perversity involves its own ``shift'':  for a space $W$ of dimension $n$ the stalk cohomology of the $IC$
sheaf at points in the top stratum is nonzero in degree $p(n)$.  This is consistent with Deligne's numbering
system in which the middle perversity is $p(S) = -\dim(S)/2$.

\section{ $t$-structures and Perverse cohomology}\index{t structure@$t$-structure}
\noindent
\begin{tcolorbox}[colback=yellow!30!white]{\begin{defn}
An {\em indecomposable} object $A$ in an abelian category is one that cannot be expressed nontrivially as
a direct sum $A = B \oplus C$.  A {\em simple} object $A$ is one that has no nontrivial subobjects $B \to A$ (where
the morphism is a monomorphism).  An object is {\em semisimple} if it is a direct sum of simple objects.  An object
is {\em Artinian} if descending chains stabilize and is {\em Noetherian} if ascending chains stabilize.
A category is Artinian (resp. Noetherian) if every object is.  Each object in an Artinian Noetherian category 
can be  expressed as a finite iterated sequence of extensions of simple objects. \end{defn} }\end{tcolorbox}

\subsection{The perversity zero $t$ structure: ordinary sheaves}\label{subsec-t=0}
\index{perversity!zero}\index{t structure@$t$-structure!perversity zero}
Let $W$ be a stratified space.  The category $Sh_c(W)$ of (ordinary) sheaves on $W$ that are constructible with respect
to this stratification is Artinian and Noetherian:   If $S$ is a constructible sheaf on 
$W$ there is a largest stratum
$j:X \to W$ so that the local system $j^*(S)$ is nonzero.  So there is an exact sequence
\[ 0 \to j_!j^*S \to S \to \uul{\coker} \to 0\]
and the cokernel is supported on smaller strata.  Continuing by induction we conclude that
$S$ is an iterated extension of sheaves of the type $j_!(A)$ where $A$ is a local system on a single
stratum, which is itself an iterated extension of simple local systems.

Now let $A^{\bullet}$ be a complex of sheaves.  We have truncation functors
\begin{equation}\label{eqn-tau}
 \begin{diagram}[size=2em]
A\b &=                &(\cdots &\rTo^{d^{r-2}}&A^{r-1}&\rTo^{d^{r-1}}&A^r&\rTo^{d^r}&A^{r+1}&\rTo^{d^{r+1}}\cdots)\\
\tau_{\le r}A\b & = &(\cdots &\rTo & A^{r-1} & \rTo &\uul{\ker}(d^r) & \rTo &0 & \rTo\cdots)\\
\tau^{\ge r}A\b& = &(\cdots & \rTo& 0  & \rTo &\uul{\coker}(d^{r-1}) & \rTo &A^{r+1} & \rTo\cdots)
\end{diagram}\end{equation}
Then there is a short exact sequence $0 \to \tau_{\le 0} A\b \to A\b \to \tau^{\ge 1} A\b \to 0$, and
the cohomology sheaf of $A\b$ is given by\index{cohomology!sheaf}\index{sheaf!cohomology!and $t$ structure}
\begin{tcolorbox}[colback=cyan!30!white]\[\uul{H}^r(A\b) = \tau_{\le r}\left(\tau^{\ge r}A\b\right) = 
\tau^{\ge r}\left(\tau_{\le r} A\b\right)\]\end{tcolorbox}

To summarize, let $Sh_c(W)$ be the category of (ordinary) sheaves on $W$ that are constructible with respect to this
stratification.  Then the the following holds:

\begin{tcolorbox}[colback=cyan!30!white]
\begin{thm}  The cohomology functor $\uul{H}^r:D^b_c(W) \to Sh_c(W)$ is given by
$\tau_{\le r}\circ\tau^{\ge r}$ and also by $\tau^{\ge r}\circ\tau_{\le r}$.  The functor $\uul{H}^0$ restricts to
an equivalence of categories
between $Sh_c(W)$ and the full subcategory of $D^b_c(W)$ whose objects are complexs $A\b$ such that
$\uul{H}^r(A\b) = 0$ for $r \ne 0$.  This category is Artinian and Noetherian
and its simple objects are the sheaves $j_{!}(\mathcal E)$
where $\mathcal E$ is a simple local system on a single stratum $j:X\to W$. \end{thm}\end{tcolorbox}


\subsection{}  The simple observations in the preceding paragraph reflect a general principal.
Fix a perversity $\bar p$.  
Let $\mathcal P(W)$ denote the category of $\bar p$-perverse sheaves on $W$
that are constructible with respect to a given stratification.  There are truncation functors, cf.~\S \ref{subsec-construction},
\[ {}^{\bar p}\tau_{\le r} \ \text{ and }\ {}^{\bar p}\tau^{\ge r}: D^b_c(W) \to D^b_c(W)\]
which are cohomological, that is, they take distinguished triangles to exact sequences, and satisfy
\begin{enumerate}
\item[\textcolor{magenta}{(T1)}]
${}^{\bar p}\tau_{\le r}(A\b) = ({}^{\bar p}\tau_{\le 0}(A\b[r]))[-r]$.
\end{enumerate}
 From this, define the {\em perverse cohomology}\index{cohomology!perverse}\index{perverse!cohomology}
\[ {}^{\bar p}H^r(A\b) = {}^{\bar p}\tau_{\le r} \left( {}^{\bar p}\tau^{\ge r} A\b\right).\]
Then ${}^{\bar p}H^r: D^b_c(W) \to \mathcal P(W)$ and $A\b \in \mathcal P(W)$ if and only 
if ${}^{\bar p}H^r(A\b) = 0$ for all $r \ne 0$.  In this case ${}^{\bar p}H^0(A\b) = A\b$. The
same method of proof (cf. \cite{BBDG} 4.3.1) gives:

\begin{tcolorbox}[colback=cyan!30!white]
\begin{thm}  The cohomology functor ${}^{\bar p}{H}^r:D^b_c(W) \to \mathcal P(W)$ is given by
${}^{\bar p}\tau_{\le r}\circ {}^{\bar p}\tau^{\ge r}$ and also by ${}^{\bar p}\tau^{\ge r}\circ {}^{\bar p}\tau_{\le r}$.  
The functor ${}^{\bar p}{H}^0$ restricts to an equivalence of categories
between $\mathcal P(W)$ and the subcategory of $D^b_c(W)$ whose objects are complex $A\b$ such that
${}^{\bar p}{H}^r(A\b) = 0$ for $r \ne 0$.  This category is Artinian and Noetherian and its simple objects are the sheaves 
$Rj_*(\uul{IC}{}^{\bullet}_{\bar p}(\mathcal E))$ where $\mathcal E$ is a simple local system on a single
stratum $X$ and where $j:\bar X \to W$ is the inclusion.
 \end{thm}\end{tcolorbox}
\noindent
In particular, a {\em semisimple} perverse sheaf is one which is a direct sum of (appropriately shifted)
intersection cohomology sheaves of (closures of) strata.


\subsection{$t$ structures}\index{truncation!and $t$ structure}
The truncation functors are determined by the support and cosupport conditions.
At this point it becames  essential to shift to Deligne's numbering scheme (see \S \ref{subsec-numbering}),
but for simplicity we consider only the case of middle perversity $p(S) = -\dim(S)/2$.

Define the full subcategory $D^b_c(W)_{\le 0}\subset D^b_c(W)$ to 
consist of constructible complexes $A\b$  that satisfy the support condition, that is,
\[ \dim_{\CC} \{ x \in W|\ H^i(j_x^*A\b) \ne 0 \} \le -i\]
and define $D^b_c(W)^{\ge 0}$ to consist of complexes $A\b$ that satisfy the cosupport condition,
\[ \dim_{\CC} \{ x \in W| \ H^i(j_x^!A\b) \ne 0\} \le i.\]  
Then $\mathcal P(W) = D^b_c(W)_{\le 0} \cap D^b_c(W)^{\ge 0}$.  Define $D^b_c(W)_{\le m}=D^b_c(W)_{\le 0}[-m]$ 
and $D^b_c(W)^{\ge m}= D^b_c(W)^{\ge 0}[-m]$.   In \S \ref{subsec-construction} we show:
\begin{enumerate}
\item[\textcolor{magenta}{(T2)}] $D^b_c(W)_{\le 0} \subset D^b_c(W)_{\le 1}$ and $D^b_c(W)^{\ge 0} \supset
D^b_c(W)^{\ge 1}$.
\item[\textcolor{magenta}{(T3)}] $\Hom_{D^b_c(W)}(A\b, B\b) = 0$ for all $A\b \in D^b_c(W)_{\le 0}$ and
$B \in D^b_c(W)^{\ge 1}$\end{enumerate}
\begin{enumerate}
\item[\textcolor{magenta}{(T4)}]For any complex $X\b$ there is a distinguished triangle
\[\begin{diagram}[size=2em]
A\b && \rTo && X\b \\
& \luTo_{[1]} && \ldTo &\\
&&B\b&& \end{diagram}\]
where $A\b \in D^b_c(W)^{\le 0}$ and where $B\b \in D^b_c(W)^{\ge 1}$.
\end{enumerate} 
This determines $A\b$ and  $B\b$ up to unique isomorphism in $D^b_c(W)$ and it also determines
the {\em $\bar p$-perverse truncation functors}\index{t structure@$t$-structure}\index{truncation}
\begin{tcolorbox}[colback=yellow!30!white]
\[{}^{\bar p}\tau_{\le 0}X\b = A\b\qquad \text{ and }\qquad  {}^{\bar p}\tau^{\ge 1}X\b = B\b.\]
\end{tcolorbox}

\subsection{Proof of \textcolor{magenta}{(T2)}, \textcolor{magenta}{(T3)}, \textcolor{magenta}{(T4)}}\label{subsec-construction} 
Statement \textcolor{magenta}{(T2)} is clear and \textcolor{magenta}{(T3)} comes from Lemma \ref{lem-lifting}.  
Statement \textcolor{magenta}{(T4)} is the essential technical step in the theory: for any $X\b$ in $D^b_c(W)$
we must construct its ``truncations'' $A\b$ in $D^b_c(W)^{\le 0}$ and $B\b$ in $D^b_c(W)^{\ge 1}$. 
If $W$ is a manifold and if $X\b$ is cohomologically locally constant on $W$ then take $A\b = \tau_{\le 0}X\b$
and $B\b = \tau^{\ge 1}X\b$ as in (\ref{eqn-tau}), that is, ${}^{\bar p}\tau = \tau$ is the usual truncation.  
Otherwise, assume by induction that we have defined  these truncation functors on an open union $U$ of
strata, and consider the addition of a single stratum $S$.  We may assume that $W = U \cup S$.
\[ \begin{CD}
U @>>{j}> W @<<{i}< S
\end{CD}\]
Given a complex $X\b$ on $W$ let $Y\b[-1]$ be the mapping cone of $X\b \to Rj_*{}^{\bar p}\tau^{\ge 1}j^*X\b$ and let
$A\b[-1]$ be the mapping cone of $Y\b \to Ri_* {}^{\bar p}\tau^{\ge 1}i^*Y\b$. 
\[ \begin{diagram}[size=2em]
&&&&  X\b & & \rTo && Rj_*{}^{\bar p}\tau^{\ge 1}j^*X\b \\
&&&\ruDotsto&  & \luTo&& \ldTo_{[1]} &\\
&&A\b&&\rTo  && Y\b && \\
&&&\luTo_{[1]}&&\ldTo&&&\\
&&&&Ri_*{}^{\bar p}\tau^{\ge 1}i^*Y\b&&&&
\end{diagram} \]
This gives the desired morphism $A\b \to X\b$ as the composition $A\b \to Y\b \to X\b$.  
Then $B\b$ is the third term in the triangle defined by $A\b \to X\b$.  
There is a lot to check that $A\b, B\b$ have the desired properties.
 This argument is in \cite{BBDG} p. 48.

For example, suppose we are in the situation of middle perversity when there are two strata, the smaller
stratum of even codimension. Let $X\b = Rj_*(\QQ)$.  Then ${}^{\bar p}\tau^{\ge 1}j^*X\b = 0$ so
$Y\b = X\b$.  Now the stalk cohomology of $i^*(Y\b)$ equals the cohomology of the link (with no shift).  
So ${}^{p}\tau^{\ge 1}i^*(Y\b)$ is the
cohomology of the link in degrees (strictly) above the middle.  Therefore the mapping cone, which is $A\b$, keeps the cohomology
of the link in degrees $\le$ the middle:  it is the $\uul{IC}{}_{\bar \ell}^{\bullet}$ (logarithmic perversity, in the
``classical'' numbering scheme) sheaf, which is perverse (see \S \ref{subsec-logarithmic} ).  
We end up with $0 \to A\b \to X\b \to B\b \to 0$ where $A\b$ is perverse 
and where $B\b$ is in $D^b_c(W)^{\ge 1}$.  

If instead, we start with $X\b = Rj_!(\QQ)$ (which is already in
$D^b_c(W)^{\le 0}$) then $X\b = Y\b$ as before but $i^*Y\b = 0$ so $A\b = Y\b  = X\b$ .  

\subsection{Perverse cohomology}\label{subsec-perversecohomology}\index{perverse!cohomology}\index{cohomology!perverse}
 If $\phi:A\b \to B\b$ is a morphism (in $D^b_c(W)$) between
two perverse sheaves, we may consider it to  be a morphism of complexes $A\b \to B\b$.  As such, it has
a kernel and a cokernel.  These are unlikely to be perverse, and moreover, they may change if we choose
different (but quasi-isomorphic) representative complexes for $A\b, B\b$:  as a morphism in the derived
category, $\phi$ does not have a kernel or cokernel.  It only has a mapping cone $M\b$.

However, the kernel and cokernel of $\phi$
{\em in the category of perverse sheaves} are again perverse sheaves, and are well defined as elements of the
derived category.  Moreover, various constructions from the theory of abelian cartegories can be implemented.
For example, suppose $A^{\bullet}_0 \overset{d}{\to} A^{\bullet}_1 \overset{d}{\to} A^{\bullet}_2 \overset{d}{\to} \cdots$ 
is a complex of perverse sheaves, that is, a complex such that $d\circ d = 0$ in the derived category.  
Then $\ker(d)/\Im(d)$ makes sense
as a perverse sheaf, so we obtain the {\em perverse cohomology} ${}^{\bar p}H^r(A^{\bullet}_{\bullet})$ of such a complex.

There is a beautiful way to see how the mapping cone $M\b$ of a morphism $\phi:A\b \to B\b$ of perverse sheaves
is broken into a (perverse) kernel and cokernel.  The long exact sequence on perverse cohomology for the 
triangle $A\b \to B\b \to M\b$ reads as follows:
\[
\begin{diagram}[size=2em]
\cdots & \rTo &{}^{\bar p}H^{-1}(B\b) &\rTo& {}^{\bar p}H^{-1}(M\b) &\rTo & {}^{\bar p}H^0(A\b)
& \rTo & {}^{\bar p}H^0(B\b) & \rTo & {}^{\bar p}H^0(M\b) & \rTo & {}^{\bar p}H^1(A\b) & \rTo & \cdots\\
&  &||&&||&&||&&||&&||&&||&\\
&&0&\rTo&{}^p\ker(\Phi)&\rTo&A\b &\rTo^{\Phi}& B\b &\rTo& {}^p\coker(\Phi) &\rTo& 0
\end{diagram}\]
Thus, the perverse kernel and cokernel are precisely the perverse cohomology sheaves of $M\b$.  (The same
statement holds in the category of sheaves:  the kernel and cokernel of a sheaf morphism $A \to B$ are
isomorphic to the cohomology sheaves of the mapping cone.)

Perverse versions of other functors can be defined by using ${}^{\bar p}\uul{H}^0$ to ``project" 
the result into the category $\mathcal P(W)$.  For example, if $j:U \to X$ and if $A\b$ is perverse on $U$
 then ${}^{\bar p}j_*A\b = {}^{\bar p}H^0(Rj_*A\b)$
and ${}^{\bar p}j_!A\b = {}^{\bar p}H^0(Rj_!A\b)$.  If we start with a local system $\mathcal E$ on the nonsingular
part $U= U_2$ of a stratified space $W$ with even codimension strata then with middle perversity,
\[ {}^{\bar p}j_!(\mathcal E) = \uul{IC}{}^{\bullet}_{\bar s}(\mathcal E) \quad \text{ and }  \quad 
{}^{\bar p}j_*(\mathcal E) = \uul{IC}{}^{\bullet}_{\bar\ell}(\mathcal E)\]
are the logarithmic and sub-logarithmetic intersection complexes.
The (perverse) image of the first in the second is the middle IC sheaf, and it is sometimes referred to as
\[ {}^{\bar p}j_{!*}(\mathcal E) = \uul{IC}{}^{\bullet}_{\bar m}(\mathcal E).\]

\subsection{}  More generally, a $t$-structure on a triangulated category $\mathcal D$ is defined to be a pair of strictly
full subcategories
$\mathcal D_{\le 0}$ and $\mathcal D^{\ge 0}$ satisfying \textcolor{magenta}{ (T1), (T2), (T3), (T4)} above, where
$\mathcal D^{\ge m} = \mathcal D^{\ge 0}[-m]$ and $\mathcal D_{\le m} = \mathcal D_{\le 0}[-m]$.  It is proven in [BBD] that 
under these hypotheses the {\em heart}\index{t structure@$t$-structure!heart}
 $P = \mathcal D_{\le 0} \cap \mathcal{D}^{\ge 0}$ is an abelian full subcategory.

\section{Algebraic varieties and the decomposition theorem}
\subsection{Lefschetz theorems}\index{Lefschetz!theorem}
Suppose $W \subset \CC\PP^N$ is a complex projective algebraic variety of complex dimension $n$. 
 Let $L^j \subset \CC\PP^n$ be a codimension $j$ linear subspace.
Let $Y^j = L^j \cap W$.  If $L^j$ is transverse to each stratum of a Whitney stratification of $W$ then there are natural
morphisms $\alpha:IH^r(W;\QQ) \to IH^r(Y^j;\QQ)$ and $\beta:IH^s(Y^j;\QQ) \to IH^{s+2j}(W;\QQ)$.
\begin{tcolorbox}[colback=cyan!30!white]
\begin{thm}
If $L^1$ is transverse to $W$ then the restriction mapping $IH^r(W;\ZZ) \to IH^r(Y^1;\ZZ)$ is an isomorphism for 
$r \le n-2$ and is an injection for $r = n-1$.  If $j \ge 1$ and $L^j$ is transverse to $W$ then the composition 
$L^j:IH^{n-j}(W;\QQ) \to IH^{n-j}(Y^j;\QQ) \to IH^{n+j}(W;\QQ)$ is an isomorphism.
\end{thm}\end{tcolorbox}
These maps are illustrated in the following diagram.

\[
\begin{diagram}[size=2em]
IH^0\ &\ IH^1\ &\ IH^2 \ &\ IH^3\  &\  IH^4\  &\ IH^5\  &\  IH^6\  &\ IH^7\ &\ IH^8\ &\ IH^9\ &\ IH^{10}  \\
IP^0 &\rTo& L(IP^0) &\rTo& L^2(IP^0) &\rTo & L^3(IP^0) &\rTo& L^4(IP^0) & \rTo &L^5(IP^0)\\
&\phantom{L^2\cdot}IP^1 &\rTo & L(IP^1) & \rTo & L^2(IP^1) & \rTo & 
L^3(IP^1) & \rTo & L^4(IP^1) & \\
&& IP^2 & \rTo& L(IP^2) & \rTo & L^2(IP^2) & \rTo & L^3(IP^2)  &&\\
&&& IP^3 &\rTo &L(IP^3) &\rTo& L^2( IP^3) &&&\\
&&&& IP^4 &\rTo& L(IP^4) &&&&\\
&&&&&IP^5&&&&& 
\end{diagram}\]
where (for $j \le n$) the {\em primitive intersection cohomology} $IP^j \subset IH^j$ is the kernel of 
$\cdot L^{n-j+1}$.  (It may also be identified with the cokernel of
$\cdot L:IH^{j-2} \to IH^j$.)  The resulting decomposition is called the
{\em Lefschetz decomposition} for $r \le n$, $IH^r \cong \oplus_{j=0}^{\lfloor r/2\rfloor}L^j\cdot IP^{r-2j}$.
The combination of Poincar\'e duality 
\[IH^{n+r}(W;\QQ) \cong \Hom(IH^{n-r}(W;\QQ),\QQ)\] and the Lefschetz
isomorphism $L^r:IH^{n-r}(W;\QQ) \cong IH^{n+r}(W;\QQ)$ induces a nondegenerate bilinear pairing on $IH^{n-r}(W;\QQ)$.
With respect to this pairing the Lefschetz decomposition is orthogonal and its signature is described by the Hodge Riemann
bilinear relations.

\subsection{Hodge theory and purity} \index{Hodge decomposition}
 Let $W$ be a complex projective algebraic variety.  Then there is a natural
decomposition $IH^r(W;\CC) \cong \oplus_{a+b=r}IH^{a,b}(W)$ with complex conjugate
$\overline{IH^{a,b}} \cong H^{b,a}$.  The
Lefschetz operator $\cap H^1$ induces $IH^{a,b} \to IH^{a+1,b+1}$.

Let $X$ be a projective algebraic variety defined over a finite field $k$ with $q = p^t$ elements.    Then the \'etale
intersection cohomology $IH^s_{ \acute{e}t}(X(\bar k); \QQ_{\ell})$ carries an action of $\Gal(\bar k/k)$ which is
topologically generated by the Frobenius $Fr$.  The eigenvalues of $Fr$ on $IH^s_{\acute{e}t}(X)$ have absolute
value equal to $\sqrt{q}^s$.

\subsection{Modular varieties} \index{modular variety}
 Let $G$ be a semisimple algebraic group defined over $\QQ$ of Hermitian type,
with associated symmetric space $D = G(\RR)/K$ where $K\subset G(\RR)$ is a maximal compact subgroup.
Let $\Gamma \subset G(\QQ)$ be a torsion free arithmetic group.  The space $X = \Gamma \backslash D$
is a Hermitian locally symmetric space and it admits a natural ``invariant" Riemannian metric.  Let $\mathcal E$
be a finite dimensional metrized local system on $X$.  The $L^2$ cohomology $H^r_{(2)}(X;\mathcal E)$ is defined
to be the cohomology of the complex of differential forms
\[\Omega_{(2)}^j(X) = \left\{ \omega \in \Omega^j(X;\CC)|\ \int_X \omega \wedge * \omega < \infty,\ 
\int_X (d\omega) \wedge *(d\omega) < \infty \right\}.\]
W. Baily and A. Borel \cite{BailyBorel} showed how to compactify the space $X$ so as to obtain a complex
algebraic variety $\overline{X}$.  \index{Baily-Borel compactification}
Their construction was modified by G. Shimura \cite{Shimura} and
P. Deligne \cite{travauxdeShimura} to obtain a variety that is defined over the rational numbers or over
a number field.  The complex $\Omega_{(2)}^{\bullet}(X)$ 
 may be interpreted as the global sections of a complex of sheaves $\uul{\Omega}^r$
on the compactification $\overline{X}$, whose sections over an open set 
$U \subset \bar X$ consist of differential forms $\omega \in \Omega^r(U \cap X)$ such that $\omega, d\omega$
are square integrable. It is true but not obviouis that this is a fine sheaf.  (A cutoff function may
introduce large derivatives, see \cite{Zucker}.)  Since this sheaf is fine its cohomology equals the cohomology
of its global sections, which are the $L^2$ differential forms on $X$.
 The following conjecture of S. Zucker (\cite{Zucker}) \index{Zucker conjecture}\index{conjecture!Zucker}
 was proven by E. Looijenga \cite{Looijenga}
 and independently by L. Saper and M. Stern \cite{SaperStern}.

\begin{tcolorbox}[colback=cyan!30!white]
\begin{thm}There is a natural quasi-isomorphism of sheaves $\uul{\Omega}{}^{\bullet}_{(2)}(X) \cong
\uul{IC}{}^{\bullet}_{\overline{X}}$ which induces an
isomorphism $H^r_{(2)}(X) \cong IH^r(\overline{X};\CC)$.
\end{thm}\end{tcolorbox}
Since the complex of sheaves $\uul{\Omega}{}^{\bullet}_{(2)}$ is self dual, the proof consists of showing that it satisfies the
support conditions of intersection cohomology.  In \cite{Looijenga} this is accomplished by applying the
decomposition theorem (below) to a resolution of singularities.

Very roughly speaking, this result provides a link beween the (infinite dimensional) representation theory of $G(\RR)$ 
(and modular forms) and the \'etale intersection cohomology of $\overline{X}$, which admits an action of a
certain Galois group.

\subsection{Morse theory again}
Suppose $W\subset M$ is a Whitney stratified complex algebraic or complex analytic subvariety of a complex manifold
and that the strata closures are also complex analytic.  The proof of the following theorem is in Appendix
\S \ref{sec-Morseperverse}.\index{Morse!theory!and perverse sheaves}
\begin{tcolorbox}[colback=cyan!30!white]
  \begin{thm}\label{thm-onedegree}\index{perverse!sheaf!and Morse theory}\index{Morse!theory!and perverse sheaves}
  Let $A\b$ be a constructible sheaf on  $W\subset M$.  Then $A\b$ is (middle) perverse
  if and only if for each stratum $X$, each point $x\in X$ and each nondegenerate covector 
$\xi \in T^*_{X,x}M$ the Morse group
  $M^t(W,A\b,\xi)$ vanishes except possibly in the single degree $t = \cod_W(X)$. 
  \end{thm}\end{tcolorbox}
Consequently if $f:M \to \RR$ is a  a Morse function in the  sense
of \S \ref{subsec-Morsefunctions} and  $x\in X\subset W$ is a nondegenerate critical point
of $f$ then the Morse group
\[  H^i(W_{\le v+\epsilon}, W_{\le v-\epsilon};A\b)\]
is nonzero in at most a single degree, $i = \cod_{\CC}(X) + \lambda$ where $\lambda$ is the Morse index of
$f|X$.
The proof of this and other Morse-theoretic properties of perverse sheaves is discussed  in Appendix
\ref{sec-Morseperverse} following \cite{SMT}.

\subsection{Decomposition theorem}  \index{decomposition theorem}\label{sec-decomposition}
This incredibly useful theorem was first formulated and proven in \cite{BBDG}.  An accessible proof,
valid also in the complex analytic context may be found in \cite{deCataldo}.
Let $f:X \to Y$ be a proper complex projective algebraic map with $X$ nonsingular.
The decomposition theorem says that $Rf_*(\uul{IC}{}^{\bullet}_X)$ breaks into a direct sum of intersection
complexes of subvarieties of $Y$, with coefficients in various local systems, and with shifts.  In many cases
this statement is already enough to determine the constituent sum.   More precisely,
\begin{tcolorbox}[colback=cyan!30!white] \begin{enumerate}
\item $Rf_*(\uul{IC}{}^{\bullet}_X) \cong \bigoplus_i {}^{\bar p}\uul{H}^i(Rf_*(\uul{IC}{}^{\bullet}_X))[-i]$ 
(This says that the push forward sheaf is a direct sum of perverse sheaves; these are its own perverse
cohomology sheaves, shifted.)
\item Each ${}^{\bar p}\!\uul{H}^i(Rf_*(\uul{IC}{}^{\bullet}_X))$ is a semisimple perverse sheaf.  (This
says that it is a direct sum of $\uul{IC}^{\bullet}$ sheaves of strata.  In particular, each summand enjoys
all the remarkable properties of intersection cohomology that were described in the previous section.)
\item The big summands come in pairs, ${}^{\bar p}\!\uul{H}^i(Rf_*(\uul{IC}{}^{\bullet}_X)) \cong
{}^{\bar p}\!\uul{H}^{-i}(\mathbf{D}(Rf_*(\uul{IC}{}^{\bullet}_X)))$ (This is because $\uul{IC}{}^{\bullet}_X$ is
self-dual, hence so is its pushforward.)
\item If $\eta$ is the class of a hyperplane on $X$ then for all $r$, 
\[ \cdot \eta^r: {}^{\bar p}\!\uul{H}^{-r}( Rf_*(\uul{IC}{}^{\bullet}_X)) \to 
{}^{\bar p}\!\uul{H}^r(Rf_*(\uul{IC}{}^{\bullet}_X))\]
is an isomorphism. (This is the relative hard Lefschetz theorem.)\index{hard Lefschetz theorem}\index{Lefschetz!theorem!hard}
\item If $L$ is the class of a hyperplane on $Y$ then for all $s$ and all $r$,
\[ \cdot L^s: H^{-s}(Y, {}^{\bar p}\!\uul{H}^r(Rf_*(\uul{IC}{}^{\bullet}_X)))
\to H^s(Y,{}^{\bar p}\!\uul{H}^r(Rf_*(\uul{IC}{}^{\bullet}_X)))\]
is an isomorphism.  (This is just the statement that each summand satisfies hard Lefschetz.)
\end{enumerate} \end{tcolorbox}
Moreover, if the mapping $f$ can be stratified then the resulting perverse sheaves are constructible with
respect to the resulting stratification of $Y$.  An example is given in \S \ref{subsec-semismall}.
\subsection{}
The hard Lefschetz isomorphisms give rise to the Lefschetz decomposition into primitive pieces as above.
The combination of the Lefschetz isomorphism and the Poincar\'e duality isomorphism gives a  nondegenerate
bilinear form on each $H^r(Rf_*(\uul{IC}^{\bullet})$.  The Lefschetz decomposition is orthogonal with respect
to this pairing, and the signature of the components is given by the Hodge Rieman bilinear relations.  
\subsection{} If $f:X \to Y$ is a proper projective algebraic map, recall that the $i$-th cohomology
sheaf of $Rf_*(\uul{IC}{}^{\bullet}_X)$ is the constructible sheaf
\[  R^i f_*(\uul{IC}{}^{\bullet}_X) = \uul{H}^i(Rf_*(\uul{IC}{}^{\bullet}_X)\]
whose stalk at a point $y \in Y$ is the cohomology $H^i(f^{-1}(y); \uul{IC}{}^{\bullet}_X | f^{-1}(y))$.  Let
$U \subset Y$ be the nonsingular part.  Then the {\em invariant cycle theorem}
\index{invariant cycle theorem} says that
\begin{tcolorbox}[colback=cyan!30!white]
\begin{enumerate}
\item[(6)]The restriction map
\[ IH^i(X) \to H^0(U,R^if_*(\uul{IC}{}^{\bullet}_X)) = \Gamma(U, R^if_*(\uul{IC}{}^{\bullet}_X))\]
is surjective. \end{enumerate}\end{tcolorbox}

 Let $y \in Y$ and let $B_y$ be the intersection of $Y$ with a small ball around $y$.  Let $y_0 \in
B_y \cap U$.  For simplicity assume $X$ is nonsingular so that $\uul{IC}{}^{\bullet}_X = \QQ_X[n]$.  Then the
{\em local invariant cycle theorem} says\index{invariant cycle theorem}
\begin{tcolorbox}[colback=cyan!30!white]
\begin{enumerate}
\item[(7)]the restriction mapping to the monodromy invariants
\[ H^i(f^{-1}(B_y, \QQ) = H^i(f^{-1}(y), \QQ) \to H^0(B_y, R^if_*(\QQ_X[n]))
\cong H^i(f^{-1}(y_0))^{\pi_1(U \cap B_y)})\]
\noindent is surjective. \end{enumerate}\end{tcolorbox} 

\subsection{}
Let us suppose that $X$ is nonsingular with complex dimension $n$,
so that (in Deligne's numbering system) $\uul{IC}{}^{\bullet}_X \cong
\QQ_X[n]$.  The decomposition theorem contains two hard Lefschetz theorems and they work 
against each other to limit the types of terms that can appear in this decomposition.  Let $\eta \in
H^2(X)$ denote a hyperplane class and let $L \in H^2(Y)$ denote a hyperplane class.
Statement (4) says that, for each $j \ge 0$  the cup product with $\eta^j$ induces an isomorphism
\[ H^r(Y; {}^{\bar p}H^{-j}(Rf_*(\QQ[n]))) \cong H^{r+2j}(Y; {}^{\bar p}H^j(Rf_*(\QQ[n]))) \text{ for all } r.\]
Statement (2) says that ${}^{\bar p}H^j(Rf_*(\QQ[n]))$ is a direct sum of intersection cohomology sheaves, each of
which satisfies hard Lefschetz (with respect to $L$) so that, for any $t \ge 0$ and for all $j$, the
 cup product with $L^t$ induces an isomorphism
\[
H^{r-t}(Y; {}^{\bar p}H^{j}(Rf_*(\QQ[n]))) \cong H^{r+t}(Y; {}^{\bar p}H^j(Rf_*(\QQ[n]))).\]

\subsection{}
Suppose $\pi:X \to Y$ is a resolution of singularities.  The decomposition theorem says that
$R\pi_*(\QQ[n])$ is a direct sum of intersection cohomology sheaves of subvarieties.  The stalk
cohomology of this sheaf, at any nonsingular point $y\in Y$ is $H^*(\pi^{-1}(y);\QQ[n])$ which is $\QQ$ in
degree $-n$.  So the sheaf $\uul{IC}{}^{\bullet}_Y$ is one of the summands, that is:  the intersection
cohomology of $Y$ appears as a summand in the cohomology of any resolution.

\subsection{}
Suppose $X, Y$ are nonsingular and $f:X \to Y$ is an algebraic fiber bundle.  Then
$Rf_*(\QQ_X[n])$ decomposes into a direct sum of perverse sheaves on $Y$,  each of
which is therefore a local system on $Y$, that is, 
\[ H^r(X;\QQ) \cong \oplus_{i+j=r} H^i(Y;H^j(F))\]
where $H^j(F)$ denotes the cohomology of the fiber, thought of as a local system on $Y$.
In other words, the Leray spectral sequence for this map degenerates (an old theorem of Deligne)
and hard Lefschetz applies both to $H^j(F)$ and to $H^i(Y)$.

\subsection{Three proofs}
The first and original proof is in [BBD] and uses reduction to varieties in characteristic $p >0$, purity of
Frobenius, and Deligne's proof of the Weil conjectures.  The second proof is due to Morihiko Saito, who
developed a theory of mixed Hodge modules in order to extend the proof to certain analytic settings.
The third proof is due to deCataldo and Migliorini, who used classical Hodge theory.  Their proof works
in the complex analytic setting and some people feel it is the most accessible of the three.
 
\bigskip

\section{Cohomology of toric varieties}\index{cohomology!toric varieties}\index{toric variety}
\subsection{}
In 1915 Emmy Noether proved that if a Hamiltonian system is preserved by an 1-parameter infinitesimal symmetry
(that is to say, by the action of a Lie group) then a certain corresponding ``conjugate" function, or ``first integral"
is preserved under the time evolution of the system.  Time invariance gives rise to conservation of energy.
Translation invariance gives rise to conservation of momentum.  Rotation invariance gives rise to conservation
of angular momentum.  
\begin{figure}[!h]
\includegraphics[width=.4\linewidth]{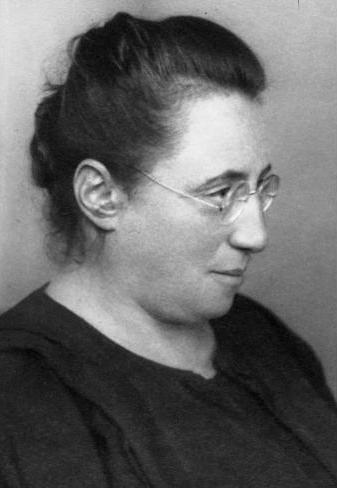}
\caption{Emmy Noether}
\end{figure}

Today, this is known as the {\em moment map}: \index{moment map} Suppose $(M,\omega)$ is a symplectic
manifold, and suppose a compact lie group $G$ acts on $M$ and preserves the symplectic form.  
 The infinitesimal action of $G$ in the direction of $V$ 
defines a vector field $X$ on $M$.  Contract this with the symplectic form to obtain a 1-form
$\theta = \iota_X(\omega)$.  It follows that $d\theta=0$.  If the action of $G$ is {\em Hamiltonian} then in
fact, $\theta = df$ for some smooth function $f:M \to \RR$ (defined up to a constant).  This is the conserved
quantity.

\subsection{Moment map}\label{subsec-momentmap} 
 In summary, if the action of a compact Lie group $G$ on a smooth manifold $M$ is Hamiltonian then there 
 exists a {\em moment map}, 
that is, a smooth mapping $\mu:M \to \mathfrak g^*$
so that for each $X \in \mathfrak g$ the differential of the function $p \mapsto\langle\mu(p), X\rangle $ 
equals $\iota_X(\omega)$.

For example consider the Fubini Study metric $h(z,w)=\sum{dz_i\wedge d\bar z_i}$ on projective space.  The real and
imaginary part, $h = R + i\omega$ are respectively, positive definite and sympletic.  
Fix $a_0, a_1, \cdots, a_n \in \ZZ$.  If $\lambda \in \CC^{\times}$
acts on $\CC^{n+1}$ by 
\[ \lambda \cdot (z_0, z_1, \cdots, z_n) = (\lambda^{a_0}z_0, \lambda^{a_1}z_1,\cdots, \lambda^{a_n}z_n)\]
then, restricting the action to $(S^1)$, the resulting moment map $\mu:\CC\PP^n \to \RR$ is
\[ \mu([z_0:z_1:\cdots:z_n]) = \frac{a_0|z_0|^2 + a_1|z_1|^2 + \cdots + a_n |z_n|^2}{(|z_0|^2 + \cdots + |z_n|^2)}.\]
If $(\lambda_0, \lambda_1,\cdots,\lambda_n)\in (\CC^{\times})^n$ acts on $\CC^n$ by
\[ (\lambda_0,\cdots,\lambda_n)\cdot (z_0,\cdots, z_n) = (\lambda_0z_0, \cdots,\lambda_nz_n)\]
then, restricting the action to $(S^1)^n$ the restulting moment map $\mu:\CC\PP^n \to \RR^n$
\[ \mu([z_0:\cdots:z_n]) = \frac{(|z_0|^2, |z_1|^2, \cdots, |z_n|^2)}{(|z_0|^2+ \cdots + |z_n|^2)}\]
and it is the standard simplex contained in the hyperplane $x_0 + \cdots + x_n = 1$.  
These actions are Hamiltonian and the moment map collapses orbits of $(S^1)^n$.

\subsection{}
Now let $X \subset\CC\PP^N$ be an $n$-dimensional subvariety on which a torus $T=(\CC^{\times})^n$ acts with finitely
many orbits.  In this case the action extends to a linear action on projective space of the sort described above and
the moment map image (for the action of $(S^1)^n$), $\mu(X) \subset \mu(\CC\PP^n)$ is convex.  
In fact, the {\em convexity theorem}\index{conexity theorem} of Atiyah, Kostant, Guillemin, 
Sternberg (\cite{Atiyah, Guillemin})states that
\begin{tcolorbox}[colback=cyan!30!white]
\begin{thm}The moment map image $\mu(X)$ is the convex hull of the images $\mu(x_i)$ of the $T$-fixed points in $X$.
The image of each $k$-dimensional $T$-orbit is a single $k$-dimensional face of this polyhedron.
\end{thm}\end{tcolorbox}
It turns out, moreover, that the toric variety is a rational homology manifold if and only $\mu(X)$ is a simple polytope,
meaning that each vertex is adjacent to exactly $n$ edges.

Algebraic geometers prefer a presentation of a toric variety from a {\em fan}, a collection of homogeneous cones in
Euclidean space.  From a fan one constructs a convex polynedron by intersecting the fan with a ball centered at the origin,
and then flattening the faces.  The resulting convex polyhedron is the dual of the moment map polyhedron.  If the
moment map polyhedron is {\em simple} then the fan-polyhedron is {\em simplicial}, meaning that the faces are
simplices.

\begin{tcolorbox}[colback=yellow!30!white]
\begin{defn}  If $Y$ is a complex algebraic variety define the intersection cohomology Poincar\'e polynomial
\[ h(Y,t) = h_0 + h_1t + h_2t^2 + \cdots + h_nt^n\]
where $h_r = {\rm rank}\ IH^r(Y;\QQ)$.
If $y \in Y$ define the local Poincar\'e polynomial 
$h_y(Y,t) = \sum_{r \ge 0} {\rm rank}\left(\uul{H}^r(\uul{IC}{}^{\bullet})_y\right)t^r.$
\end{defn}\end{tcolorbox}
If $Y$ is defined over $\FF_q$ we use the same notation for the Poincar\'e polynomial
of the \'etale intersection cohomology.  

\subsection{Counting points}
There is a very general approach to understanding the cohomology and intersection cohomology 
of an $n$-dimensional algebraic
variety defined over a finite field $\FF_q$, provided its odd degree cohomology groups vanish.  The variety $Y/ \FF_q$
is said to be {\em pure}\index{cohomology!pure}\index{pure cohomology}\index{Frobenius eigenvalues}
 if the eigenvalues of Frobenius on $H^r(Y; \QQ_{\ell})$ have absolute value $\sqrt{q}^r$ with
respect to any embedding into the complex numbers.  The Weil conjectures (proven by Grothendieck \cite{}
and Deligne \cite{})\index{conjecture!Weil}\index{Weil conjectures} say that
\[ \sum_{r=0}^{2n} (-1)^r{\rm Tr}(Fr_q: H^r(Y) \to H^r(Y)) = |Y(\FF_q)|\]
the right hand side being the number of points that are fixed by the Frobenius morphism.  The intersection
cohomology of any projective algebraic variety is pure.  If the variety $Y$ is also nonsingular (so that
$IH^*(Y) = Y^*(Y)$ and {\em Tate} (which means that the eigenvalues on $H^r$ are in fact equal to
$(\sqrt{q})^r$) then this gives
\[ h(Y, \sqrt{q}) =  \sum_{r=0}^n {\rm rank}H^{2r}(Y)q^r = |Y(\FF_q)|.\]
For example, if such a variety $Y$ is defined over the integers, is nonsingular and $Y(\CC)$ has an algebraic cell decomposition with $m_r$ cells of (complex) dimension $r$ then $h_{2r+1}=0$ and $h_{2r}=m_r$ accounts for
$m_r q^r$ points over $\FF_q$.
In the case of a nonsingular toric variety whose moment map image is a convex polyhedron with $f_r$ faces of
dimension $r$ this gives
\[   h(Y, \sqrt{q})=
\sum_{s=0}^n {\rm rank}H^{2s}(Y;\CC) q^s = \sum_{r=0}^nf_r(q-1)^r\]
since each $r$-dimensional orbit is itself (isomorphic to) a torus of dimension $r$.  The hard Lefschetz
theorem says $h_{2s-2} \le h_{2s}$ for $2s \le n$ which in turn gives inequalities between the
face numbers, as observed by Stanley in 1980.

\subsection{}
If we wish to use intersection cohomology rather than ordinary cohomology in the Weil conjectures then the
formula becomes
\[ \sum_{r=0}^{2n} (-1)^r{\rm Tr}(Fr_q: IH^r(Y) \to IH^r(Y)) = |Y(\FF_q)|_{\text{mult}}\]
where each point $y \in Y(\FF_q)$ is counted with a multiplicity
equal to the (alternating sum of) trace of Frobenius acting on the stalk of the intersection
cohomology at $y\in Y(\FF_q)$.  If this is pure and if the stalk cohomology vanishes in odd degrees, then this 
multiplicity equals the Poincar\'e polynomial $h_y(Y,\sqrt{q})$ of the stalk of the intersection cohomology.
In conclusion, if the intersection cohomology of $Y$ is Tate and vanishes in odd degrees then
\begin{tcolorbox}[colback=cyan!30!white]
\begin{equation}\label{eqn-global}
\sum_{s=0}^n {\rm Tr}(Fr_q|(IH^{2s}(Y))) = \sum_{s=0}^n \text{rank}\ IH^{2s}(Y)q^s
= h(Y, \sqrt{q}) = \sum_{y \in Y(\FF_q)} h_y(Y,\sqrt{q}).\end{equation}\end{tcolorbox}
Let us now try to determine these multiplicities $h_y$.  If $f(x) = a_0 + a_1x + \cdots + a_nx^n$ 
define the truncation $\tau_{\le r}f$ to be
the polynomial $a_0 + \cdots + a_rx^r$ consisting of those terms of degree $\le r$.

\begin{tcolorbox}[colback=cyan!30!white]
\begin{lem}
Let $Z \subset \CC P^{N-1}$ be a projective algebraic variety of dimension $d$,
with intersection cohomology Poincar\'e polynomial
\[ g(t)=g_0 + g_1 t+ \cdots + g_{2d} t^{2d} = \sum_{r=0}^{2d} \dim(IH^r(Z))t^r.\]
Then the stalk of the intersection cohomology of the complex cone $Y={\rm cone}_{\CC}(Z)\subset \CC P^N$
at the cone point $y\in Y$ has Poincar\'e polynomial 
\begin{equation}\label{eqn-local}
h_y(Y,t)=\tau_{\le d} \left(g(t)(1-t^2)\right)\end{equation}
\end{lem}\end{tcolorbox}
\begin{proof}
The complex projective space $\CC \PP^N$ is the complex cone over $\CC \PP^{N-1}$.  In fact, if we remove the
cone point then what remains is a line bundle $E \to \CC \PP^{N-1}$ whose first Chern class 
$c^1(E) \in H^2(\CC \PP^{N-1})$ is the class of a hyperplane section.  This is to say that there exists a section of this
bundle that vanishes precisely on a hyperplane; it may be taken to be
\[ s([z_0:...z_{N-1}]) = [z_0:...:z_{N-1}, \Sigma_j a_j z_j] \in \CC \PP^N\]
for any fixed choice (not all zero) of $a_0, a_1, \cdots, a_{N-1} \in \CC$.  The vanishing of the last coordinate
is a hyperplane in $\CC \PP^{N-1}$.  So this class may be used as a hard Lefschetz class.

If $Z \subset \CC \PP^{N-1}$ is a projective algebraic variety
 then ${\rm cone}_{\CC}(Z) \subset \CC \PP^{N}$ is a singular
variety and the link $L$ of the cone point can be identified with the sphere bundle of this line bundle 
$E_Z\to Z$.
The Gysin sequence becomes
\[ \begin{CD}
IH^{i-2}(Z)@>{\cup c^1}>>IH^i(Z) @>>> IH^i(L) @>>> IH^{i-1}(Z) @>{\cup c^1}>> IH^{i+1}(Z) @>>> 
\end{CD}\] 
It follows from the hard Lefschetz theorem for $Z$ that $IH^i(L) $ is the primitive part of the
intersection cohomology of $Z$ for $i \le d$, that is, 
\[ IH^i(L) \cong IP^i(Z) = \coker(\cdot c^1(E): IH^{i-2}(Z) \to IH^i(Z))\]
for $i \le \dim(Z)$, and hence its Poincar\'e polynomial is given by 
\[ g_0 + g_1t + (g_2-g_0)t^2 + (g_3-g_1)t^3 + (g_4-g_2)t^4 + \cdots + (g_d-g_{d-2})t^d =
\tau_{\le d}g(t)(1-t^2).\qedhere\]  
\end{proof}

\subsection{Some geometry}  Let $\mu:Y \to P\subset \RR^m$ be the moment map corresponding to the action 
of a torus $T\cong (\CC^{\times})^m$ on a toric variety $Y$.
If $F$ is a face of $P$, the link of $F$ can be realized as another
conex polyhedron, $L_F = P \cap V$ where $V\subset \RR^N$ is a linear subspace such that
$\dim(V) + \dim(F) = N-1$, which passes near $F$ and through $P$.  (For example, $V$ may be taken to lie
completely in the plane $F^{\perp}$.)  In fact, $L_F$ is the moment map image of a sub-toric variety
$Y_F$ on which a sub-torus $T_F$ acts.

\subsection{}  In the case of a toric variety $Y$, \index{toric variety}
a given face $F$ corresponds to a stratum $S_F$ of the toric variety.
The link of this stratum is therefore isomorphic to a circle bundle over a toric variety whose moment map
image is the link $L_F$ of the face $F$.  Let $h(Y_F,t)$ be the intersection cohomology Poincar\'e polynomial
of this ``link" toric variety.  Then equations (\ref{eqn-local}) and (\ref{eqn-global}) give:
\begin{tcolorbox}[colback=cyan!30!white]
\begin{thm} The $IH$ Poincar\'e polynomial of $Y$ is
\[h(Y,t)= \sum_{F}(t^2-1)^{\dim(F)}. \tau_{\le n-\dim(F)}\left( (1-t^2)h(Y_F,t) \right)\]
\end{thm}\end{tcolorbox}

\subsection{}
In particular, the intersection cohomology only depends on the combinatorics of the moment map image $P=\mu(Y)$,
and moreoer, the functions $h(Y_F,t)$ may be determined (inductively) from the moment map
images $ L_F=\mu_F(Y_F)$.  The hard Lefschetz theorem (which says that $h_{2r} \ge h_{2r-2}$ for all $r \le \dim(Y)$)
then implies a collection of inequalities among the numbers of chains of faces.

\subsection{Remarks.}
This formula simplifies if $P$ is a simple polyhedron, to:  
\[h(Y,t) = \sum_F(t^2 -1)^{\dim (F)} = f(t^2-1)\] 
where $f(s) = f_0 + f_1s + \cdots + f_d s^d$ and $f_j$ is the number of faces of dimension $j$.  The polytopes considered
here are always rational, meaning that the vertices are rational points in $\RR^d$, or equivalently, the faces are
the kernels of linear maps $\RR^d \to \RR$ with rational coefficients.  Any simple (or simplicial) polytope can be
perturbed by moving the faces (resp. the vertices) so as to make them rational.  Therefore the inequalities
arising from hard Lefschetz apply to all simple polytopes.  However a general polytope cannot necessarily
be perturbed into a rational polytope with the same face combinatorics.  The Egyptian pyramid, for example, has a
square face.  Lifting one of the vertices on this face, an arbitrarily small amount, will force the face to ``break".  In
order to prove that the inequalities arising from hard Lefschetz for intersection cohomology can be applied to
any polytope it was necessary to construct something like intersection cohomology in the non-rational case.
This was constructed in \cite{Barthel} and the Hard Lefschetz property was provven in \cite{Karu}.

\subsection{}
There is another way to prove this result using the decomposition theorem (which does not involve passing to
varieties over a finite field).  The singularities of the toric variety $Y$ can be resolved by a sequence of steps,
each of which is toric with moment maps that correspond to `cutting off the faces" that are singular.  For example,
the Egyptian pyramid has a single singular point.  The singularity is resolved by a mapping $\pi:\tilde Y \to Y$
as illustrated in this diagram:

\begin{figure}[!h]
\begin{tikzpicture}[scale=.3]
\coordinate (0) at (0,0); \coordinate (T) at (5,20);
\coordinate (A) at (10,0); \coordinate (B) at (15,3);
\coordinate (C) at (5,3);
\coordinate (I0) at (3,12);\coordinate (IA) at (7,12);
\coordinate (IB) at (9, 13.2);\coordinate (IC) at (5, 13.2);

\draw (0) -- (A) -- (B);
\draw [red, dotted] (0) -- (C) -- (B);
\draw[red, dotted] (C) -- (IC);

\draw (A) -- (IA);  \draw (B) -- (IB);

\draw (0) -- (I0) -- (IA) -- (IB) -- (IC) -- (I0);

\begin{scope} [shift = {(25,0 )}, rotate = 0]
\draw (0,0) -- (10,0) -- (15,3) -- (5,20) -- (0,0);
\draw (5,20) --(10,0); \draw[red, dotted] (5,20) -- (5,3);
\draw[red, dotted] (15,3) -- (5,3) -- (0,0);

\end{scope}

\draw [-latex] (13,14) -- (25,19); 
\draw [-latex] (16.5,1.5) -- (21.5,1.5);
\end{tikzpicture}
\caption{Moment map of a resolution}\label{fig-resolution}
\end{figure}
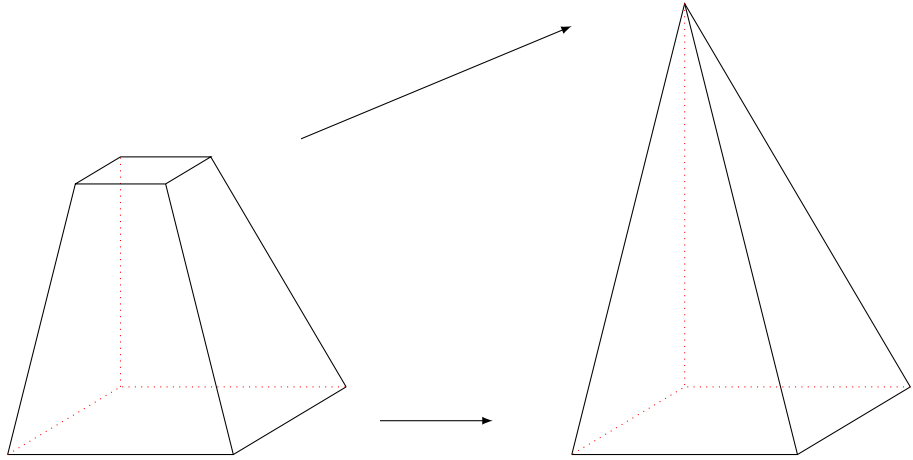

\medskip\noindent

Let us examine the decomposition theorem for this mapping.  The mapping is an isomorphism everywhere except
over the singular point $y \in Y$ and $\pi^{-1}(y) \cong \PP^1 \times \PP^1$.  The stalk cohomology of the
pushforward $R\pi_*(\uul{\QQ}{}_{\tilde Y})$ is $(\QQ, 0, \QQ \oplus \QQ, 0, \QQ)$.  Put this into the support diagram for
a 3 dimensional variety, see Figure \ref{fig-Rpistar}.  (In this figure, the degree $i$ is the ``usual" cohomology 
degree notation and the degree $j$ is the ``perverse degree" notation.)

\begin{figure}[!h]
\begin{tabular}{|c||c||c|c|c|c|c|}
\hline
$i$ & $j$& $\cod 0$ & $\cod 2$ & $\cod 4$ & $\cod 6$ & $H^*(\pi^{-1}(y))$ \\
\hline\hline
6 & 3&        c      &        c        &      c        &      c        &               \\
5 &  2&               &                 &       c        &      c         &               \\
4 &  1&               &                 &                &        c       &        $\QQ$        \\
3 &  0&               &                 &                &      \rz      &         0      \\
2 & -1&                &                 &                &        x      &     $ \QQ\oplus\QQ$         \\
1 &  -2&               &                 &       x        &       x       &       0        \\
0 & -3&      x        &      x          &     x         &       x      &       $ \QQ $      \\
\hline 
&&\multicolumn{4}{|c|}{$\uul{IC}^{\bullet}$ support}&\\
\hline
\end{tabular}
\caption{Support of $R\pi_*(\uul{\QQ})$}\label{fig-Rpistar}
\end{figure}

From the decomposition theorem we know that one term will be $\uul{IC}{}^{\bullet}_Y$ and that there are 
additional terms supported at the singular point $y$.  From the support condition it is clear that
$\QQ[3]$  (on the bottom row) is part of the $IC$ sheaf.  It is not so clear how much of the $(\QQ \oplus \QQ)[1]$ 
belongs to $\uul{IC}{}^{\bullet}_Y$ and how much belongs to the other terms.  However the $\QQ[-1]$ (at the top
of the column; in degree $j=1$) is definitely not part of $IC$.  By Poincar\'e and especially by Hard Lefschetz, 
it must be paired with one copy of $\QQ$ in degree $-1$.  So this leaves $\QQ[3] \oplus \QQ[1]$ (in degrees
$-3$ and $1$ respectively) for the IC sheaf.  A closer inspection of this argument shows that these two
terms constitute the primitive cohomology of the fiber, as we saw earlier.

Thus, the decomposition theorem singles out the primitive cohomology of the fiber as belonging to the IC sheaf.
Now, assuming by induction that the formula holds for $IH^*(\tilde Y)= H^*(\tilde Y)$ (which is less singular that $Y$)
and knowing how these terms decompose, it is easy to conclude that the formula must also hold
for $IH^*(Y)$.

\bigskip

\section{Springer Representations}\index{Springer!representation}
\subsection{The flag manifold}
Let $G = SL_n(\CC)$.  It acts transitively on the set $\mathcal F$
of complete flags $0 \subset F^1 \subset \cdots \subset F^{n-1}
\subset \CC^n$ and the stabilizer of the standard flag is the ``standared" Borel subgroup $B$ 
 of (determinant $=1$) upper triangular matrices, giving an isomorphism $\mathcal F \cong G/B$.  The Lie algebras are
$\mathfrak g$ (matrices with trace $=0$) and $\mathfrak b = $ upper triangular matrices with trace $=0$.  If $x\in G$ and
$x Bx^{-1} = B$ then $x \in B$.  So we may identify $\mathcal F$ with the set $\mathcal B$ of all subgroups of $G$ that are
conjugate to $B$ or equivalently to the set  of all subalgebras of $\mathfrak g$ that are conjugate to
$\mathfrak b$, that is, the {\em variety of Borel subalgebras} of $\mathfrak g$.
\begin{tcolorbox}[colback=yellow!30!white]
\begin{defn}\label{defn-Springer}
Let $\mathcal N$ be the set of all nilpotent elements in $\mathfrak g$.  Define 
\[
\begin{CD}\widetilde{\mathcal N} &=& \left\{ (x\in \mathcal N, A \in \mathcal B)|\ x \in \Lie(A) \right\}@>{\phi}>>\mathcal B\\
 &&@V{\pi}VV&& \\
 && \mathcal N &&\end{CD}\]
\end{defn}\end{tcolorbox}
The mapping $\pi$ is proper and its fibers $\mathcal B_x = \pi^{-1}(x)$ are called {\em Springer fibers}.  
\index{Springer!fiber}
In a remarkable series of papers [Springer 1976, 1978], T. A. Springer constructed an action of the 
symmetric group $W$ on the cohomology of each Springer fiber $\mathcal B_x$, even though $W$ does not actually act
on $\mathcal B_x$.

\begin{figure}[!h]
\includegraphics[width=.9\linewidth]{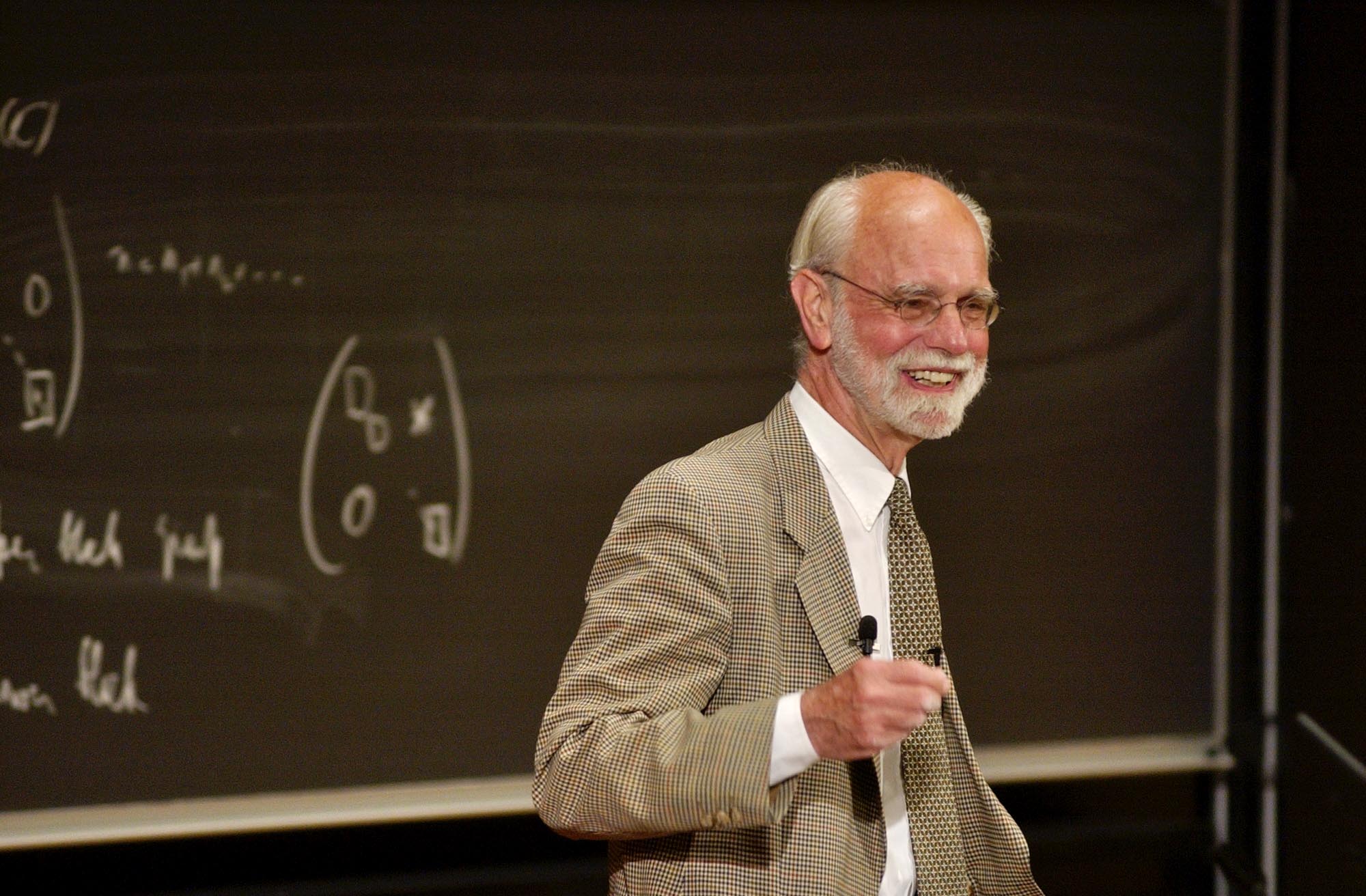}
\caption{Tonny Springer}
\end{figure}

Let $A$ be the subgroup that preserves a flag $F_A=(0=A^0 \subset A^1 \subset \cdots \subset A^n = \CC^n)$ then
the following are equivalent:
\begin{enumerate}
\item $(x,A) \in \widetilde{\mathcal N}$ 
\item $x \in \Lie(A)$
\item  $\exp(x)$ preserves the flag $F_A$
\item the vectorfield $V_x$ (defined by $x$) on the flag manifold $\mathcal F$ vanishes on $F_A$ 
\item  $xA^j \subset A^{j-1}$ for $1\le j \le n$.
\end{enumerate}
So the Springer fiber $\pi^{-1}(x)$ is the zero set of the vectorfield $V_x$; it 
is the set of all flags that are preserved by $x$ and is often referred to as the {\em variety of fixed flags}.
For the subregular nilpotent $x\in \mathfrak g$ the Springer fiber turns out to be a string of $n-1$ copies of
$\PP^1$, each joined to the next at a single point.  [picture]

\begin{lem}  The mapping $\phi:\widetilde {\mathcal N} \to \mathcal B$ 
identifies $\widetilde{\mathcal N}$ with the cotangent bundle to the flag manifold. \end{lem}
\begin{proof}
The tangent space at the identity to $\mathcal F$ is $T_I(G/B) = \mathfrak g/\mathfrak b$.  So its dual space is
\[ T^*_I(G/B) = \left\{ \phi:\mathfrak g \to \CC |\ \phi(\mathfrak b) = 0\right\}.\]
The canonical inner product $\langle, \rangle: \mathfrak g \times \mathfrak g \to \CC$ given by
$\langle x, y \rangle = \Tr(xy)$ is symmetric and nondegenerate.  Using this to identify $\mathfrak g^*$ with $\mathfrak g$
gives
\[ T^*_I(G/B) \cong \left\{ x \in \mathfrak g|\ \langle x, \mathfrak b \rangle = 0 \right\}=\mathfrak n\]
is the algebra of strictly upper triangular matrices, that is, the nilradical of $\mathfrak b$.  So for each Borel subgroup
$A \subset G$, the cotangent space $T^*_A(G/B) \cong \mathfrak n(A)$ is naturally isomorphic to the nilradical of
$\Lie(A)$.  But this is exactly the fiber, $\phi^{-1}(A)$.
\end{proof}
\subsection{}
The group $G$ acts on everything in the diagram (\ref{defn-Springer}).  It acts transitively on $\mathcal B$ and it
acts with finitely many orbits on $\mathcal N$, each of which is a {\em nilpotent conjugacy class}.
\index{nilpotent conjugacy class}\index{conjugacy class!nilpotent}  These
form a Whitney stratification of $\mathcal N$ by complex algebraic strata.  It follows from
Jordan normal form that each nilpotent conjugacy class corresponds to a partition
$\lambda_1 \ge \lambda_2 \ge \cdots \ge \lambda_r$ with $\sum \lambda_i = n$ or equivalently to a Young frame
\[
\yng(5,2,1)\]
Stratum closure relations correspond to refinement of partitions with the largest stratum corresponding to the
case of a single Jordan block ($\lambda_1 = n$) and the smallest stratum corresponding to $0 \in \mathfrak N$,
which is the partition $1+1+1\cdots +1 = n$.

\subsection{}  The {\em Grothendieck simultaneous resolution} is the pair
\[
\begin{CD}
 \widetilde{\mathfrak g} = \left\{ (x,A) \in \mathfrak g \times \mathcal B|\ x \in A\right\}
 @>{\widehat\pi}>> \mathfrak g \end{CD}\]
  Given $(x,A) \in \widetilde{\mathfrak g}$ choose $h\in G$ which conjugates $A$ into the standard Borel subgroup $B$.
  It conjugates $x$ into an element $x' \in \Lie(B)$ and the diagonal entries $\alpha= \alpha(x) \in \mathfrak t$
  are well defined where $\mathfrak t$ is the set of diagonal matrices with trace $=0$.  On the other hand, the
  characteristic polynomial $ch(x)$ of $x$ is determined by the diagonal matrix $\alpha$ but is independent of the order of
  the entries.  The set of possible characteristic polynomials forms a vector space, with coordinates given by
  the coefficients of the characteristic polynomial, which provides an example of a remarkable theorem of Chevalley 
  that the quotient $\mathfrak t/W$ is
  again an affine space.  The nilpotent elements $\mathcal N$ in $\mathfrak g$ map to zero in $\mathfrak t/W$.
   In summary we have a diagram
  \[\begin{diagram}[size=2em]
  \widetilde{\mathcal N} & \rInto & \widetilde{\mathfrak g} & \rTo^{\alpha} & \mathfrak t \\
    \dTo^{\pi} && \dTo_{\widehat\pi} && \dTo_p \\
    \mathcal N & \rInto & \mathfrak g & \rTo_{ch} & \mathfrak t/W
  \end{diagram}
  \]
  \begin{tcolorbox}[colback=cyan!30!white]
  \begin{thm}[Grothendieck, Lusztig, Slowdowy]
  The map $\widehat{\pi}:\widetilde{\mathfrak g} \to \mathfrak g$ is small.\index{small map}  
  The map  $\pi:\widetilde{\mathcal N} \to \mathcal N$ is semi-small.  For each $A\in \mathfrak t$
  the map $\widehat{\pi}:\alpha^{-1}(A) \to ch^{-1}(p(A))$ is a resolution of singularities.
  \end{thm}\end{tcolorbox}

\subsection{Adjoint quotient}\index{adjoint!quotient}
There is another way to view the map $ch$.  Each $x \in \mathfrak g$ has a unique Jordan decomposition
$x = x_s + x_n$ into commuting semisimple and nilpotent elements.  Then $x_s$ is conjugate to an element
of $\mathfrak t$, and the resulting element is well defined up to the action of $W$.  The quotient
$\mathfrak t/W$ turns out to be isomorphic to the geometric invariant theory quotient $\mathfrak g//G$ and
the map $ch:\mathfrak g \to \mathfrak t/W$ is called the adjoint quotient map.\index{braid group}

The vector space $\mathfrak t$ consists of diagonal matrices $A = \diag(a_1,\cdots, a_n)$
with trace zero.  The {\em reflecting hyperplanes} are the subspaces 
$H_{ij }= \left\{A| a_i = a_j\right\}$ where two entries coincide and
they are permuted by the action of $W$.  Their image in $\mathfrak t/W$ is the discriminant variety ${Disc}$
consisting of all (characteristic) polynomials with multiple roots.  The complement of the set $\cup_{i \ne j}H_{ij}$
is sometimes called the {\em configuration space of $n$ ordered points in $\CC$}; 
its fundamental group is the {\em colored braid group}.  The
complement of the discriminant variety in $\mathfrak t/W$ is the configuration space of $n$ unordered points, and
its fundamental group is the braid group.

Suppose $x \in \mathfrak g^{rs}\subset \mathfrak g$ is regular and semisimple, 
meaning that its eigenspaces $E_1, E_2, \cdots, E_n$
are distinct and form a basis of $\CC^n$.  Then the flag $E_1 \subset E_1\oplus E_2 \subset E_1\oplus E_2 \oplus E_3
\cdots$ is fixed by $x$ and every fixed flag has this form, for some ordering of the eigenspaces.  Therefore there
are $n!$ fixed flags and the symmetric group permutes them according to the regular representation.

\subsection{Springer's representation}\index{Springer!representation}
Since $\widehat{\pi}$ is a small map, we have a canonical isomorphism 
\[R\widehat{\pi}_*(\uul{\QQ}{}_{\widetilde{\mathfrak g}}[n])
\cong \uul{IC}^{\bullet}(\mathfrak g;\mathcal E),\] 
the intersection cohomology sheaf on $\mathfrak g$, constructible with respect
to a stratification of $\widehat{\pi}$, with coefficients in the local system $\mathcal E$ over the regular semisimple
elements whose fiber at $x \in \mathfrak g^{rs}$ is the direct sum $\oplus_F \QQ_F$ of copies of $\QQ$, one for each
fixed flag.  The symmetric group $W$ acts on $\mathcal E$ which induces an action on $\uul{IC}^{\bullet}(\mathcal E)$
and therefore also on the stalk cohomology at each point $y \in \mathfrak g$, that is, on
\[ \uul{H}^r(R\widehat{\pi}_*(\uul{\QQ}{}_{\widetilde{\mathfrak g}})_y = H^r(\widehat{\pi}^{-1}(y)) \cong 
\uul{IH}^r(\mathfrak g;\mathcal E)_y.\]
For $y \in \mathcal N$ this action of $W$ on $H^*(\pi^{-1}(y))$ turns out to coincide with Springer's representation.
The decomposition theorem for the semismall map $\pi$ 
provides an enormous amount of information about these representations.

\subsection{Decomposition theorem for semismall maps}
\label{subsec-semismall}\index{decomposition theorem!semismall}
Recall that a proper algebraic morphism $f:X \to Y$ is semismall if it can be stratified so that
\[ 2\dim_{\CC}f^{-1}(y) \le \codim(S)\]
for each stratum $S \subset Y$, where $y \in S$.  This implies that $d=\dim(X) = \dim(Y)$ and, if $X$ is
nonsingular, that $Rf_*(\QQ_X)[d]$ is perverse on $Y$.  A stratum $S \subset Y$ is said to be {\em relevant}
if $2d = \codim(S)$ where $d = \dim f^{-1}(y)$.  In this case, the top degree cohomology $H^{2d}(f^{-1}(y)$ forms a local
system $L_S$ on the stratum $S$.  The following result is due to W. Borho and R. MacPherson\cite{BorhoMacPherson}.

\begin{tcolorbox}[colback=cyan!30!white]
\begin{prop}
Suppose $f:X \to Y$ is semismall and $X$ is nonsingular of complex dimension $d$.  
Then the decomposition theorem has the following special form:
\[ Rf_*(\QQ_x)[d] \cong \bigoplus_{S} \uul{IC}{}^{\bullet}_{\bar S}(L_S)\]
where the sum is over those strata $S$ that are relevant (with no shifts, if we use Deligne's numbering).  
In particular, the endomorphism algebra of this
sheaf $\End(Rf_*(\QQ_X)[d]) \cong \bigoplus_S \End_S(L_S)$ is isomorphic to the direct sum of the endomorphism
algebras of the individual local systems $L_S$.
\end{prop}\end{tcolorbox}

\begin{proof}
The top stratum, $Y^o$ is always relevant.  If no other strata is relevant then the map is {\em small}. 
By \S \ref{subsec-smallmaps} this implies
$Rf_*(\QQ_X[d]) \cong \uul{IC}{}^{\bullet}_{Y}(L_{Y^o})$.   So there is only one term in the decompsition theorem.  
Now suppose the next relevant stratum
has codimension $c$ so that the fiber over points in this stratum has (complex) dimension $d=c/2$ (and 
in particular, the complex codimension $c$ is even).
One term in the decomposition theorem is $\uul{IC}{}^{\bullet}_Y(L_{Y^o})$.  Consider the support diagram (e.g. if $c=4$):

\begin{figure}[!h]
\begin{tabular}{|c||c||c|c|c|c|c|c|}
\hline
$i$ & $j$& $\cod_{\CC} 0$ & $\cod_{\CC} 1$ & $\cod_{\CC} 2$ & $\cod_{\CC} 3$ & $\cod_{\CC}4$ & $H^*(f^{-1}(y))$ \\
\hline\hline
8 & 4 &      c        &     c          &      c        &       c        &  c   &    \\
7 & 3 &               &                 &     c          &     c         &  c   &    \\
6 & 2&               &                 &                &      c         &   c        &    \\
5 & 1&               &                 &                &                 &    c      &     \\
4 &  0&               &                 &                &                &  \rz    &       $H^4(f^{-1}(y))$   \\
3 &  -1&               &                 &                &               &    x    &    $H^3(f^{-1}(y))$     \\
2 & -2&                &                 &                &        x      &   x   &      $H^2(f^{-1}(y)) $         \\
1 &  -3&               &                 &       x        &       x       &  x     &    $ H^1(f^{-1}(y))$       \\
0 & -4&      x        &      x          &     x         &       x      &  x    &     $ H^0(f^{-1}(y))$      \\
\hline 
&&\multicolumn{5}{|c|}{$\uul{IC}^{\bullet}$ support}&\\
\hline
\end{tabular}\end{figure}
From this diagram we can see that a new summand must be added to the decomposition, and it is the
local system $H^4(f^{-1}(y)) = H^{2d}(f^{-1}(y) = L_S$ arising from the top cohomology of the fiber, that is,
from the irreducible components of the fiber.  This local system on $S$ gives rise to the summand
$\uul{IC}{}^{\bullet}_{\bar S}(L_S)$ over the whole of the closure $\bar S\subset Y$. 
Let $U_S = \cup_{T \ge S}$ be the open set consisting of $S$ together with the strata larger than $S$.  We
have constructed an isomorphism of $Rf_*(\QQ_X)$ with 
\[ Rf_*(\QQ_X[d]) \cong \uul{IC}{}^{\bullet}_Y(L_{Y^o}) \oplus \uul{IC}{}^{\bullet}_{\bar S}(L_S)\]
over the open set $U_S$.  The exact same method as in Proof No.1 shows that this isomorphism
extends uniquely to an isomorphism  over the larger open set that contains additional (smaller) strata 
until we come to the next relevant stratum.
Continuing in this way by induction gives the desired decomposition.

Finally, if $L_R, L_S$ are local systems on distinct strata $R, S$ of $Y$ then
\[ \hHom_{D^b_c(Y)}(\uul{IC}{}^{\bullet}_{\bar R}(L_R), \uul{IC}^{\bullet}{}_{\bar S}(L_S)) =0\]
which implies that the endomorphism algebra of this direct sum decomposes into a direct sum of
endomorphism algebras.
\end{proof}

\subsection{Some conclusions}  Let $d = \dim(G/B) = n(n-1)/2$.
For the $SL_n$ adjoint quotient $R\pi_*(\QQ_{\widetilde{N}}[d]) \cong \oplus_S \uul{IC}^{\bullet}(\overline{S}; L_S)$
\begin{enumerate}
\item Every stratum $S$ is relevant:  for $x \in S$, $2\dim(\mathcal B_x) = \cod(S)$. 

\item The odd cohomology of each Springer fiber $H^{2r+1}(\mathcal B_x)=0$ vanishes.
\item Let $S$ be a stratum corresponding to some partition (or Young frame) $\lambda=\lambda(S)$.  Then the Springer
action on the top cohomology $H^{\codim(S)}(\mathcal B_x)$ is the irreducible representation 
$\rho_{\lambda}$ corresponding
to $\lambda$ (modulo possible notational normalization involving transpose of the partition and tensoring with the
sign representation).  [Note:  this representation is not (cannot) be realized via permutations of the components
of $\mathcal B_x$, but the components give, nevertheless, a basis for the representation.]
\item Every irreducible representation of $W$ occurs in this decomposition, and it occurs with multiplicity one.
\item The local systems occurring in the decomposition theorem ($\SL_n$ case only!) are all trivial. 
\item Putting these facts together, let $S_{\lambda}$ denote the stratum in $\mathcal N$ corresponding to the
partition $\lambda$.  Then, using ``classical" degree indices, the decomposition theorem in this case becomes:
\[R\pi_*(\QQ_{\widetilde{N}}) \cong
\oplus_{\lambda} \uul{IC}^{\bullet}(\overline{S}_{\lambda})[-2d_{\lambda}]\otimes V_{\lambda}\] 
where $V_{\lambda}$ is the (space of the) irreducible representation $\rho_{\lambda}$ of $W$ and
$d_{\lambda}$ is the complex dimension of the stratum $S_{\lambda}$. 
\item  Applying stalk cohomology at a point $x \in S_{\mu}$ to this formula gives:
\[ H^i(\mathcal B_x) \cong \oplus_{\lambda \ge \mu} \uul{IH}{}^{i-2d_{\lambda}}_y
(\overline{S}_{\lambda}) \otimes V_{\lambda}.\]
Consequently, if an irreducible representation $\rho_{\lambda}$ of $W$ occurs in $H^*(\mathcal B_x)$
then the stratum $S_{\mu}$ containing $x$ is in the closure of the stratum $S_{\lambda}$ corresponding to $\rho$.
\item More generally, suppose  $R<S$ are strata corresponding to partitions $\mu, \lambda$ respectively
and let $x\in R$.   Then $\dim \uul{IH}{}^{i-2d_{\lambda}}(\overline{S})_x$ is the
multiplicity of the representation $\rho_{\lambda}$ in the cohomology $H^{i}(\mathcal B_x)$ 
of the Springer fiber $\mathcal B_x$. 
In fact the Poincar\'e polynomial of these multiplicities 
\[P_{\lambda,\mu}(t) = \sum_i \text{mult}(\rho_{\lambda}, H^{2i}(\mathcal B_{x}))t^i = \sum_i \text{rank}
(\uul{IH}^{2i}(\overline{S})_x) t^i\]
turns out to be the Kostka-Foulkes polynomial.
\item For $x =0 \in \mathcal N$ the Springer fiber $\mathcal B_x = G/B$ is the full flag variety and the representation
of $W$ is the regular representation.  Moreover, the full endomorphism algebra 
\[ \End(R\pi_*(\QQ_{\widetilde{\mathcal N}}) \cong \CC[W]\]
is isomorphic to the full group-algebra of the Weyl group, with its regular representation.
\end{enumerate}

Even for the $W$ action on the full flag manifold $\mathcal B \cong G/B$ these results are startling.  In this case it
had been shown by Borel and Leray that the action of $W$ was the regular representation, but which irreducible
factors appeared in which degrees of cohomology had appeared to be a total mystery.  Many computations were
done by hand and the result appeared to be random.  The above conclusions explain that the multiplicity of
each representation $\rho_{\lambda}$ in $H^i(G/B)$ is given by the rank of the local intersection cohomology, 
in degree $i$, at the origin $o \in \mathcal N$ of the stratum $S$ that corresponds to the partition $\lambda$.

\begin{center}
\begin{figure}[H]
\begin{tabular}{|c|c|c|}
\hline
degee & rank & Young \\
\hline
$H^{12}(G/B)$ & 1 & \T{2.5}\tiny{\yng(4)} \\
\hline
$H^{10}(G/B)$  & 3 &\T{4.5} \tiny{\yng(3,1)}\\
\hline
$H^8(G/B)$ & 5 & \T{4}\tiny{\yng(2,2)} $\oplus$ \tiny{\yng(3,1)}\\
\hline
$H^6(G/B)$ & 6 & \T{6}\tiny{\yng(2,1,1)} $\oplus$ \tiny{\yng(3,1)}\\
\hline
$H^4(G/B)$ & 5 & \T{6} \tiny{\yng(1,1,1)} $\oplus$ \tiny{\yng(2,2)}\\
\hline
$H^2(G/B)$ & 3 & \T{6} \tiny{\yng(1,1,1)}\\
\hline
$H^0(G/B)$ & 1 & \T{7.5} \tiny{\yng(1,1,1,1)}\\
\hline
\end{tabular}
\caption{Springer representations for $SL_4$.
Each rep occurs as often as its dimension} \end{figure}\end{center}



\bigskip

\section{Iwahori Hecke Algebra}\index{Hecke algebra}\index{Iwahori Hecke algebra}
\subsection{}  The (Iwahori) Hecke algebra of the Weyl group of an algebraic group is
usually defined using generators and relations (see Proposition \ref{prop-Hecke-definition})
without motivation or explanation for the mysterious formula $(T_s-q)(T_s+1) = 0$.  In this
section we explain the geometric nature of this equation.
\subsection{}
Let $G$ be a finite group.  The convolution\index{convolution}
 of two functions $f, f':G \to \CC$ is the function
\[ (f*f')(x) = \frac{1}{|G|}\sum_{h\in G}f(xh^{-1})f'(h) = \frac{1}{|G|}\sum_{a \in G}f(a)f'(a^{-1}x).\]
this product is associative.  If $H\subset G$ is a subgroup the Hecke algebra is
\[ \mathcal H(G,H) = \left\{ \phi:G \to \CC|\ \phi(kgk') = \phi(g) \text{ for all } k, k' \in H\right\}\]
with algebra structure given by convolution.  It is the convolution algebra of functions on 
the double coset space $H\backslash G/H$.  If $\rho:H \to \GL(V)$ is a representation, the
induced representation is
\[ \Ind_H^G(\rho) = \left\{ \phi:G \to V|\ \phi(hx) = \rho(h)\phi(x)\ \text{ for all } h \in H , x\in G\right\}\]
with action $(g.\phi)(x) = \phi(xg^{-1})$.  Then
\[ \hHom_G(\Ind_H^GV,W) \cong \hHom_H(V, \text{Res}^G_H(W)) \ \text{ and }
\mathcal H(G,H) \cong \hHom_G(\Ind_H^G(\mathbbold 1), \Ind_H^G(\mathbbold 1)).\]
Now let $G = \SL_{n}(\FF_q)$ and $H = B$ the collection of upper triangular matrices (or determinant one).  
Let $W=S_n$ be the symmetric group which may be thought of as acting on the standard basis 
vectors $\{ e_1, \cdots, e_n\}$.   It is generated by the ``simple reflections'' 
$S=\{s_1, \cdots, s_{n-1}\}$ where $s_i$ exchanges  $e_i$ and
$e_{i+1}$.   The {\em length} $\ell(w)$ of an element $w\in W$ is the minimum number of elements 
required to express $w$ as a product of simple reflections, and it is well defined.
The Bruhat decomposition says that $G = \coprod_{w \in W} BwB$.  Each $B$ orbit $BwB/B \subset G/B$
is isomorphic to an affine space of dimension $\ell(w)$. 
\begin{tcolorbox}[colback=yellow!30!white]
\begin{defn}
The Hecke algebra $\mathcal H$ is the algebra of $B$-bi-invariant functions on $G$.
It has a basis consisting of functions
\[\phi_w = {\mathbbold 1}_{BwB}.\]\end{defn}
\end{tcolorbox}
The unit element in $\mathcal H$ is the function 
$\phi_1=\mathbbold 1_B$.  In this algebra we will use the following normalization for convolution of bi-invariant
functions $f, f':G \to \CC$,
\[ (f*f')(x) = \frac{1}{|B|} \sum_{h \in G}f(xh^{-1})f'(h).\]
\begin{lem}  If $s \in S$ is a simple reflection and if $w\in W$ then the following holds:
\begin{tcolorbox}[colback=cyan!30!white]
\[\begin{cases}
\phi_w * \phi_{w'} = \phi_{ww'} &\text{  if } \ell(w) + \ell(w') = \ell(ww')\\
\phi_s * \phi_s = (q-1)\phi_s + q\mathbbold \phi_1 &  \\
\phi_s * \phi_w = (q-1) \phi_{w} + q\phi_{sw} &\text{  if }\ell(sw) < \ell(w)
\end{cases}\]  \end{tcolorbox}\end{lem}
\begin{proof}  (The third equation follows from the second by induction.)
The key nontrivial point (proven in \cite{Iwahori}; see for example, \cite{BumpHecke}) is the following geometric
property of these double cosets:
\begin{align*}
\ell(w) + \ell(w') = \ell(ww') &\implies  (BwB)(Bw'B) = Bww'B\\
\ell(ws) = \ell(w)-1 &\implies (BwB)(BsB) \subset (BwB) \cup (BwsB).\end{align*}
Following \cite{BumpHecke}, for $f \in \mathcal H$ let $\epsilon(f) = \frac{1}{|B|}\sum_{g\in G}f(g)$ so that
$\epsilon(f*f') = \epsilon(f)\epsilon(f')$ and $\epsilon(\phi_w) = q^{\ell(w)}$.  
If $\ell(w)+\ell(w') = \ell(ww')$ then $\phi_w*\phi_{w'}$ is supported on
$Bww'B$ and is $B$ bi-invariant.  Apply epsilon to conclude that $\phi_w*\phi_{w'} = \phi_{ww'}$.  Similarly, 
$\phi_s*\phi_s$ is supported on $(BsB) \cup B$ so it equals $\alpha \phi_s + \beta \phi_1$ for some
$\alpha, \beta \in \CC$.  Apply $\epsilon$ to conclude that $q^2 = \alpha q + \beta$.  Evaluate at $x = I\in G$
to get  $\phi_s*\phi_s(I) = |BsB|/|B| = q = \alpha.0 + \beta.1$ So $\alpha = q-1$. 
\end{proof}
 The same holds for any semisimple algebraic group $G$ defined over $\FF_q$,
where $B$ is a Borel subgroup and $W$ is the Weyl group and $S$ denotes the set of
simple reflections:
\begin{tcolorbox}[colback=cyan!30!white] 
 \begin{prop}\label{prop-Hecke-definition}
The Hecke algebra  $\mathcal H(W,S)$ is the free $\ZZ[q,q^{-1}]$ module
with basis elements $\phi_w$ for $w\in W$ and relations
\begin{itemize}
\item $\phi_s\phi_w = \phi_{sw}$ if $sw>w$
\item $(\phi_s-q)(\phi_s+1)=0$.
\end{itemize}
If $q=1$ this is the group algebra $\ZZ[W]$.\end{prop}\end{tcolorbox}
More generally,  if $(W,S)$ is a Coxeter group with
generators $S$ and resulting order $>$ then
the Hecke algebra of $(W,S)$ is defined to be the $\ZZ[q^{1/2}, q^{-1/2}]$ algebra generated by
symbols $T_w$ and satisfying the relations in the box.  (The reason for introducing $\sqrt{q}$ in
the coefficients will become clear in the next paragraph.)

\subsection{}  Kazhdan and Lusztig \cite{Kazhdan} discovered a mysterious new basis for the Hecke 
algebra that appeared to be closely related to infinite dimensional representations of the Lie algebra $\mathfrak g$ of $\SL_n$.
Each element $\phi_w \in \mathcal H$ is invertible and the algebra $\mathcal H$ admits an involution defined by
\[ \iota(q^{1/2}) = q^{-1/2} \text{ and }\ \iota(\phi_w) = (\phi_{w^{-1}})^{-1}.\]
\begin{tcolorbox}[colback=cyan!30!white]
\begin{thm}[Kazhdan, Lusztig]\label{thm-KL}
  For each $w \in W$ there is a unique element $c_w \in \mathcal H$ and a uniquely determined polynomial
  $P_{yw}$ for $y \le w$ such that: $\iota(c_w) = c_w$;  
  $P_{ww}=1$; $P_{yw}$ (for $y<w$) is a polynomial of degree $\le \frac{1}{2}(\ell(w)-\ell(y)-1)$; and
\[c_w = q^{-\ell(w)/2} \sum_{y \le w} P_{yw}(q) \phi_y\]
\end{thm}\end{tcolorbox}
Existence and uniqueness of $c_w$ and $P_{yw}$  is easily proven by induction.  
Kazhdan and Lusztig conjectured that the coefficients of $P_{yw}$
were nonnegative integers.  They further conjectured that, in the Grothendieck group of Verma modules, 
\[ [L_w] = \sum_{y \le w} (-1)^{\ell(w)-\ell(y)} P_{yw}(1)[M_y]\]
where $M_w$ is the Verma module corresponding to highest weight $-\rho-w(\rho)$ and $L_w$ is its unique
irreducible quotient.  This second conjecture became known as the {\em Kazhdan-Lusztig conjecture} and was
eventually proven by J. L. Brylinski and M. Kashiwara \cite{Brylinski}
and independently by A. Beilinson and D. Bernstein \cite{Beilinson}.  This
circle of ideas became extremely influential in representation theory.  But what exactly is the meaning of $c_w$ 
and $P_{yw}$?  The answer (Theorem \ref{thm-hecke-sheaves} below)  eventually came from intersection cohomology.

\section{Algebra of correspondences}\index{correspondence}
\subsection{}  Let us return to the complex picture with $G = \SL_n(\CC)$, $B$ the Borel subgroup of
upper triangular matrices and $W$ the symmetric group.  let $X = G/B$ be the flag manifold, or equivalently,
the variety of Borel subgroups of $G$.  The group $G$ decomposes as a disjoint union
$G = \coprod_{w]\ in W} BwB$.  It follows that the group $B$ acts on $X$ with finitely many orbits, and these are
indexed by the elements of $W$.  
For each $w \in W$ the {\em Schubert cell} or {\em Bruhat cell} \index{Schubert cell}\index{Bruhat cell}
$X_w = BwB/B \subset G/B$ indexed by $w$ contains
$X_y$ in its closure iff $ y < w $ in the Bruhat order.  Similarly the group $G$ acts on $X \times X$ with finitely
many orbits, each of which contains a unique point $(B,wB)$ (thinking of the standard Borel subgroup $B$ as
being the basepoint in the flag manifold) for some $w \in W$.  It consists of
pairs of flags  $(F_1 \subset F_2 \subset \cdots F_n = \CC^n)$ and $(F'_1 \subset F'_2 \subset \cdots
\subset F'_n = \CC^n)$ that are in relative position $w \in W$, meaning that there exists an ordered basis
$(e_1, e_2, \cdots, e_n)$ of $\CC^n$ so that, for each $1 \le i \le n$ 
\[ \langle e_1, e_2, \cdots, e_i \rangle = F_i\ \text{ and }\ 
\langle e_{w(1)}, e_{w(2)}, \cdots, e_{w(i)} \rangle = F'_i.\]
Equivalently, two flags $F, F'$ are in relative position $w$ if
\[ \dim(F_i \cap F'_j) = | ([1,i] \cap \sigma([1,j]))|.  \] 

If $q_2:X \times X \to X$ denotes projection to the second factor then this orbit, let us denote it by
$\mathcal O_w$ fibers over $X$ with fiber equal to $X_w$; 
in particular it is simply connected.   In other words, each $G$ orbit on $X \times X$ intersects the fiber $X$ 
in a single $B$ orbit.

\subsection{Algebra of correspondences}  The following construction  was discovered independently by
R. MacPherson, G. Lusztig, Brylinski and Kashiwara, Beilinson and Bernsstein and is described
in Springer's article \cite{Springer_Bourbaki}.  It is convenient here to change notation in the Hecke algebra, setting
\[ t = \sqrt{q}.\]
Consider the derived category $D^b_{c,even}(X \times X)$ of sheaves $A^{\bullet}$ on $X\times X$, 
cohomologically constructible with respect to this stratification such that $\uul{H}^i(A^{\bullet})=0$ for $i$ odd.
For for such a sheaf $A^{\bullet}$ let $\uul{H}^{2i}(A\b)_w$ denote the stalk cohomology of $A\b$ at a point in
$\mathcal O_w$ and define
\begin{tcolorbox}[colback=yellow!30!white]
\begin{align*}
&h:D^b_{c,even} \to \mathcal H \\
h(A) &= \sum_{w \in W}\sum_{i \ge 0}\dim(\uul{H}^{i}(A^{\bullet})_w) t^{i}\phi_w\end{align*}\end{tcolorbox}
Note that all elements of $\mathcal H$ obtained in this way have nonnegative coefficients so the image of $h$ ends
up in a sort of ``positive cone" in the Hecke algebra.  Consider the following diagram of correspondences.
\begin{diagram}[size=2em]
&& X \times X \times X &&\\
& \ldTo^{q_{12}}& \dTo_{q_{13}} & \rdTo^{q_{23}} &\\
X \times X && X \times X && X \times X
\end{diagram}

If $A^{\bullet}, B^{\bullet} \in D^b_c(X \times X)$  define their convolution product
\begin{tcolorbox}[colback=yellow!30!white]
\[A\b \circ B\b = Rq_{13!}\left(q_{12}^*(A\b) \otimes q_{23}^*(B\b)\right).\]\end{tcolorbox}

\begin{tcolorbox}[colback=cyan!30!white]
\begin{thm}  \label{thm-hecke-sheaves}
Let $j_w:\mathcal O_w \to X$ and $\bar{j}_w:\overline{\mathcal O_w} \to X$ denote the inclusions.  Then 
\[ h\left(j_{w!}(\uul{\CC}{}_{\mathcal O_w})[ \ell(w)]\right) = \phi_w\ \text{ and }
h(\uul{IC}{}^{\bullet}_w) = c_w\] 
where $\uul{IC}{}^{\bullet}_w = R\bar{j}_{w*}\left(\uul{IC}{}^{\bullet}_{\overline{\mathcal O}_w}\right)$.  
Moreover, $h(\uul{IC}{}^{\bullet}_w \circ \uul{IC}{}^{\bullet}_v) = h(\uul{IC}{}^{\bullet}_w).
h(\uul{IC}{}^{\bullet}_v)$ for all $v,w \in W$. 
 \end{thm}\end{tcolorbox}
In other words, $P_{y,w}$ is the local intersection cohomology Poincar\'e polynomial of $X_w$ at a point in $X_y$
(originally proven in \cite{Kazhdan_duality}).
It vanishes in odd degrees and its coefficients are non-negative.  Using ``classical" indexing for sheaf cohomology
(and $q = t^2$),
\[ P_{yw}(q) =  \sum_{i\ge 0} \dim \uul{IH}{}^{2i}_y(\overline{X}_w)q^{i}.\]
This stalk cohomology vanishes in odd degrees and the highest power of $t$ that can occur here is $t^{\ell(w)-\ell(y)-1}$.
Besides making the essential connection with geometry this result is a ``categorification"
\index{categorification} of the Hecke algebra:  it replaces numbers and coefficients with (cohomology) groups.   
It implies that the coefficients of $P_{yw}$ are non-negative integers.

\subsection{}  \label{subsec-pushforward-calculation}
The mapping $h$ in Theorem \ref{thm-hecke-sheaves} is multiplicative on sums of $\uul{IC}^{\bullet}$ sheaves, and
the resulting elements $c_w$ of the Hecke algebra form a basis of the Hecke algebra.  However, the convolution 
operation on the sheaves $j_{w!}$ (corresponding to elements $\phi_w \in \mathcal H$) does not necessarily agree with
Hecke multiplication, essentially because of the minus signs that arise in the product formula.

Let $s =(1,2) \in S_3$ denote the simple reflection that exchanges $1$ and $2$.  Then $c_s =\phi_1+\phi_s$.  Let us
compute the convolution of this element with itself.  Then $BsB$ consists of pairs of flags
$(E,F)$ in $\CC^3$ such that $F_1 \ne E_1$, $F_1 \subset E_2$, $F_2 = E_2$.  So the union
$(BsB) \cup (B1B)$ consists of pairs of flags $(E,F)$ such that $F_1 \subset E_2$ is arbitrary.  Let us denote
this relationship by $E \overset{1+s}{\longrightarrow} F$.  We can assume that $E$
is the standard flag, so that the flag $F$ is completely determined by $F_1 \subset E_2$.  This
set is isomorphic to $\PP^1$ hence $\uul{IC}{}^{\bullet}_s = \uul{\CC}{}_{\overline{\mathcal O}_s}$.

Now let us compose this element $1+s$ with itself.    The correspondence consists of triples 
\[E \overset{1+s}{\longrightarrow} F \overset{1+s}{\longrightarrow}F'\] 
of flags and it is apparently a $\PP^1$-bundle over $\PP^1$.
Now let us understand the mapping $\pi_{13}:X\times X \times X \to X \times X$, in other words, 
we project the tripe $(E,F,F')$ of the correspondence to the pair $(E,F')$.  
Taking $E$ to be the standard flag, the second flag
$F'$ is determined by $F'_1 \subset E_2$ and apparently all possible subspaces occur. If $F'_1 \subset E_2$ then the
fiber consists of all 1-dimensional subspaces $F_1 \subset E_2$, which is a $\PP^1$, and this is independent of
which point $F'_1$ was chosen.  So the pushforward
(which is also the pushforward with compact supports) of the constant sheaf gives $\uul{\QQ} \oplus 
\uul{\QQ}[2]$ on the closed set $(BsB) \cup (B1B)$.  In other words, this calculation says that
\[ (\phi_1 + \phi_s).(\phi_1 + \phi_s) = (q+1)(\phi_1 + \phi_s)\]
from which it follows that $\phi_s.\phi_s = (q-1)\phi_s + q\phi_1$, which is one of the defining formulas for $\mathcal H$.

 Similarly, the product $(\phi_1+\phi_s).\phi_s$ corresponds to a map of triples $E \overset{1+s}{\longrightarrow}
 F \overset{s}{\longrightarrow}F'$ for which $F_1, F'_1 \subset E_2$ but $F_1 \ne F'_1$.  When we project to the
 pair $E \overset{1+s}{\longrightarrow}F'$ we find that the fiber over each point consists of $\PP^1 - \{pt\}$ 
 so the pushforward with
 compact supports yields $\uul{\QQ}[2]$ which gives the formula $(\phi_1 + \phi_s).\phi_s = q(1+\phi_s)$.
 
 These calculations are made using the ``classical" indexing for intersection cohomology.  If Deligne's indexing is
 used then the Hecke algebra should be redefined (as Lusztig often does) as follows:
 $(T_s-q^{\frac{1}{2}})(T_s + q^{-\frac{1}{2}}) = 0$.

\subsection{}
The proof of Theorem \ref{thm-hecke-sheaves}
is tedious but does not require sophisticated methods; it is completely worked out in
online notes  \cite{Rietsch}
(following an outline of T. A. Springer \cite{Springer_Bourbaki}, 
as communicated to him by R. MacPherson) and \cite{Haines} .   
Here is the outline.  It is obvious from the definitions that $h(Rj_{w!}(\mathcal O_w)) = \phi_w$.  Now consider
$h(\uul{IC}{}^{\bullet}_s)$.  The support condition for intersection cohomology implies that 
$h(\uul{IC}{}_w)$ is a linear combination of
$h(j_{w!}\uul{\CC}{}_w)$ with coefficients that are polynomials which satisfy the degree restriction.

The same argument as in the previous section \S \ref{subsec-pushforward-calculation} shows: 
if $s$ is a simple reflection with corresponding orbit closure
$\overline{\mathcal O}_s \subset X \times X$ then
\[ h(\CC_{\overline{\mathcal O}_s} \circ \uul{IC}{}^{\bullet}_w) = (\phi_s +1). h(\uul{IC}{}^{\bullet}_w) \in \mathcal H\]
although the argument needs to be modified slightly when $sw<w$.  Next, one verifies that $h(\uul{IC}{}^{\bullet}_w)$ is
preserved by $\iota$.
One way to prove this is to consider the Bott-Samelson resolution of the Schubert variety $X_w$.  It is obtained
as a sequence of blowups by simple reflections.  One checks at each stage of the induction that the result is
preserved by $\iota$.  This proves that $h(\uul{IC}{}_w)$ is preserved by $\iota$ and satisfies the combinatorial conditions
that uniquely identify it as  $c_w \in \mathcal H$.

Finally,  again using the Bott Samelsom resolution, the decomposition theorem and induction, one proves that 
$h(\uul{IC}{}_w \circ \uul{IC}{}_{w'}) = c_wc_{w'}$.

\subsection{Digression:  Hecke algebra and modular forms}
Let $G = \SL_n(\RR)$, let $K = O(n)$ and let $\Gamma_0 = \SL_n(\ZZ)$.  
Then $D = G/K$ be the (contractible) symmetric space of positive definite symmetric matrices of determinant one.
Let $X = \Gamma_0 \backslash D$.  This is the moduli space of Riemannian tori.  (For each lattice $L \subset
\RR^n$ of determinant one we get a torus $\RR^n/L$ and an invariant Riemannian metric on it.) 
Each $g \in G_{\QQ} = \SL_n(\QQ)$ gives a correspondence\index{correspondence}
 on this space as follows.  Let $\Gamma' =
\Gamma_0 \cap (g^{-1} \Gamma_0 g)$ and let $X' = \Gamma' \backslash D$.  Then the correspondence
$X' \to X \times X$ is given by $\Gamma' x \mapsto (\Gamma_0x, \Gamma_0 gx)$.  It is well defined and each
of the projections $X' \to X$ is a finite covering.  Moreover, the isomorphism class of this correspondence
depends only on the double coset $\Gamma_0 g \Gamma_0$.  (Replacing $g$ by $\gamma g $ where
$\gamma \in \Gamma_0$ does not change the correspondence.  Replacing $g$ by $g \gamma$ changes the
correspondence but it gives an isomorphic correspondence.)  Therefore points in the double coset space
\[ \Gamma_0 \backslash \SL_n(\QQ) / \Gamma_0\]
may be interpetred as defining correspondences on $X$, which therefore acts on the homology, cohomology,
functions etc. of $X$.  So the same is true of linear combinations of double cosets.  In summary, the Hecke algebra of
compactly supported functions (meaning, functions with finite support) on $\Gamma_0 \backslash \SL_n(\QQ)
/ \Gamma_0$  acts on $H^*(X)$ by correspondences.  The composition of correspondces turns out to
coincide with convolution of compactly supported functions. Such functions are called Hecke operators.

This construction makes more sense in the ad\`elic setting where natural Haar measures can be used 
in order to define the algebra structure and the action without resorting to correspondences.  In this
setting there is an equality
\[ \SL_n(\ZZ) \backslash \SL_n(\RR)/K \cong \SL_n(\QQ) \backslash \SL_m(\AA_{\QQ})/K.\SL_n(\widehat{\ZZ})\]
and the Hecke algebra is the convolution algebra of locally constant functions with compact support
\[ f \in C_c^{\infty}(\SL_n(\QQ) \backslash \SL_n(\AA_f) / \SL_n(\QQ))\]
where $\AA_{\QQ}$ denotes the ad\`eles of $\QQ$ and $\AA_f$ the finite ad\`eles.

\section{The affine theory}
\subsection{Affine Weyl group}\index{affine Weyl group}\index{Weyl group!affine}
The symmetric group $S_n$ has Dynkin diagram : [diagram]
It is generated by simple reflections $s_1, \cdots, s_{n-1}$ with the relations
\begin{enumerate}
\item $s_i^2 = 1$
\item  $s_is_{i+1}s_i = s_{i+1}s_is_{i+1}$ for $(1 \le i \le n-2)$
\item $s_is_j = s_js_i$ if $|i-j|>1$.
\end{enumerate}
  It can be interpreted as acting on $\RR^{n-1}= \{ \sum x_i =0\} \subset \RR^n$ with
$s_i$ acting as reflection across the hyperplane $x_i = x_{i+1}$.  This decomposes $\RR^{n-1}$ into Weyl
chambers, one for each element of $S_n$.

The affine symmetric group $\widetilde{S}_n$ has Dynkin diagram:  
It is generated by simple reflections $s_0, s_1, \cdots, s_n$ with the same relations as $S_n$ and the
additional relations (corresponding to edges $s_0 s_1$ and $s_n s_0$).
It can be interpreted as acting on $\RR^{n-1}$ by adding a reflecting hyperplane to the previous picture.
Then it acts simply transitively on the alcoves.  If we take the fundamental alcove as the basepoint (identity),
then every alcove becomes labeled by a unique element of $\widetilde{S}_n$.

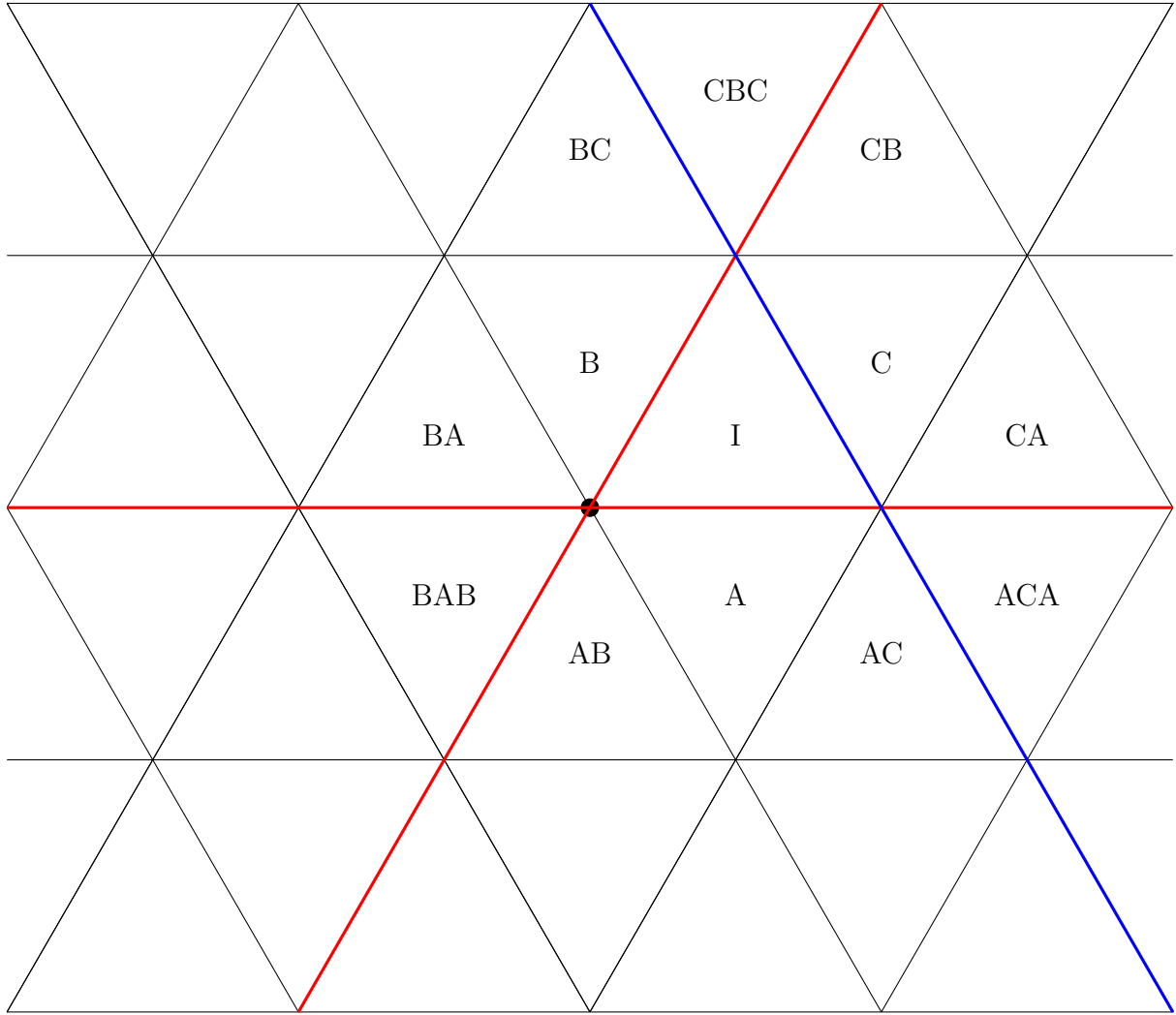
\begin{figure}[!h]
\begin{tikzpicture}[scale=4]
\foreach \r in {0,1,...,{\rowsh}} {
 \draw($(0,{\r*sqrt(3)})$)  -- ($({\r},0)$) ; 
 \draw($(0,{\rowsh*sqrt(3)}) + (\r,0)$) -- ($(\rowsh,0) + (\r,0)$); 
 \draw($(\rowsh,{\rowsh*sqrt(3)})+(\r,0)$) -- ($(\rows,0) + (0, {\r*sqrt(3)})$);
                                          }
 \foreach \r in {0,1,...,\rows} {
  \draw($(0,{\r*sqrt(3)/2})$) -- ($(\rows,{\r*sqrt(3)/2})$) ;  
                                           }
\foreach \r in {0,1,...,\rowsh}{  
 \draw($(\r,0)$) -- ($((\r,0)+({\rows/2}, {\rows*sqrt(3)/2})$);
 \draw($(0,{\r*sqrt(3)})$) -- ($(\rowsh,{\rowsh*sqrt(3)}) - (\r,0)$); 
 \draw($(\rowsh,0)+(\r,0)$) -- ($(\rows,{\rowsh*sqrt(3)}) - (0, {\r*sqrt(3)})$);
};

\draw[fill] (2, {sqrt(3)}) {circle[radius=0.030]};
 
\draw[red, very thick] (0,{sqrt(3)}) -- (4,{sqrt(3)});
\draw[red, very thick] (1,0) -- (3, {2*sqrt(3)});
\draw[blue,  very thick] (2, {2*sqrt(3)}) -- (4,0);

\node at (2.5, {sqrt(3)+.25}) {I};
\node at (2.5, {sqrt(3)-.3}) {A};
\node at (2, {sqrt(3)+.5}) {B};
\node at (2, {sqrt(3)-.5}) {AB};
\node at (1.5, {sqrt(3)-.3}) {BAB};
\node at (1.5, {sqrt(3) + .25}) {BA};

\node at (3, {sqrt(3) +.5}) {C};
\node at (2, {2*sqrt(3) - .5}) {BC};
\node at (2.5, {2*sqrt(3)-.3}) {CBC};
\node at (3, {2*sqrt(3)-.5}) {CB};

\node at (3.5, {sqrt(3)+.25}) {CA};
\node at (3.5, {sqrt(3)-.3}) {ACA};
\node at (3, {sqrt(3)-.5}) {AC};

\end{tikzpicture}

\caption{This diagram ilustrates the affine Weyl group for $\SL_3$.  
The red lines are the reflecting hyperplanes of the finite Weyl group (generated by A and B).  
The blue line is the affine reflection, C.}  
\end{figure}
\subsection{}
The affine symmetric group can also be described as $S_n \ltimes A$ where $A$ is the root lattice of translations,
\[\left\{ (a_1, \cdots, a_n) \in \ZZ^n|\ \sum a_i = 0 \right\}\]
on which $S_n$ acts by permutations.  The group $A$ may also be interpretred as the cocharacter group
of the maximal torus $\mathcal T$ consisting of diagonal matrices of determinant one.
The affine Weyl group can be described as 
\[ \{\omega: \ZZ \to \ZZ| \ \omega(i +n) = \omega(i) +n \}.\]
In this realization each element is determined by its value on $\{1, \cdots, n\}$ and so it may be written as
$[\omega(1),\cdots,\omega(n)]$.  Then elements in the lattice of translations are the elements $[a_1,\cdots,a_n]$
with $\sum a_i = 0$ and they act by addition, that is, $\omega(i) = a_i + i$ for $1 \le i \le n$.  Then every element
$\omega \in \widetilde{S}_n$ can be expressed as a permutation followed by a translation.

\subsection{Affine Hecke algebra}\index{Hecke algebra!affine}\index{affine Hecke algebra}
Let $G_{\QQ_p} = \SL_n(\QQ_p)$, as a locally compact topological group and set $K = G(\ZZ_p) = \SL_m(\ZZ_p)$.
Let $B_{\FF_p}$ be the Borel subgroup of $\SL_n(\ZZ/p\ZZ)$.  The Iwahori subgroup 
\index{Iwahori subgroup} $I_p$ is the preimage
$\pi^{-1}(B_p)$ under the $\mod p$ mapping $\phi:K \to SL_n(\ZZ/p\ZZ)$, that is, it consists of $n \times n$
matrices with entries in $\ZZ_p$, whose diagonal entries are invertible in $\ZZ_p$, and whose lower diagonal
entries are multiples of $p$.  It is compact and open in $G_{\QQ_p}$.

\begin{tcolorbox}[colback=yellow!30!white]
\begin{defn} The Iwahori Hecke algebra is the convolution algebra of locally constant complex valued functions
\[ f\in C_c^{\infty}( I_p \backslash G_{\QQ_p} / I_p) \]
with compact support on $G_{\QQ_p}$ that are bi-invariant under $I_p$. 
\end{defn}\end{tcolorbox}
Haar measure $\mu$ on $G_{\QQ_p}$ is
normalized so that $\mu(I_p)=1$.

\quash{
Then $K$ and $I_p$ are compact subgroups of $G_{\QQ_p}$ so we can define the Hecke
algebras $\mathcal H_K$ (the ``spherical" Hecke algebra) and $\mathcal H_I$ (the ``Iwahori Hecke algebra")
to be the convolution algebra of functions $f:G_{\QQ_p} \to \CC$ that are
locally constant and have compact support and are (respectively) $K$-bi-invariant or $I_p$-bi-invariant. (Typically
we normalize the Haar measure $\mu$ on $G$ so that $\mu(K)=1$ or $\mu(I_p)=1$.)
}

The Bruhat decomposition in this case reads 
$G_{\QQ_p} = \coprod_{w \in W_a} I_pwI_p$ where $W_a$ is the affine Weyl group.
Then $\mathcal H_I$ is generated by characteristic functions $\phi_w$ for $w \in W_a$ and with the same
relations as before:  $\phi_s^2 = (q-1)\phi_s + q\phi_1$.  The Kazhdan Lusztig canonical basis $c_w$ is defined
exactly as before.  The Kazhdan Lusztig theorem works in this context as well and it gives a basis of
$\mathcal H_I$ consisting of elements $c_w$ for $w \in W_a$ the affine Weyl group.

The field $\QQ_p$ is analogous to the field $\FF_q((T))$ of formal Laurent series (meaning formal power series with
finitely many negative powers of $T$ and coefficients in $\FF_q$).  The ring $\ZZ_p$ corresponds to
$\FF_q[[T]]$ (the ring of formal power series).  Reduction modulo $T$ gives a homomorphism
$\phi:\FF_q[[T]] \to \FF_q$ and  the Iwahori subgroup 
\index{Iwahori subgroup} $I_{p((T))} = \phi^{-1}(B_q)$ is defined
similarly.  Then an Iwahori Hecke algebra over $\FF_q$ is defined to be the convolution algebra of
locally constant complex valued functions
\[ f \in C_c^{\infty}(I_{p[[T]]}\backslash G_{\FF_q((T))} / I_{p[[T]]})\]
with compact support that are bi-invariant under the Iwahori subgroup.

All of this has a complex analog following the same procedure as in the finite case.
Instead of the flag manifold over $\CC$ one uses the ``affine flag manifold" $\SL_n(\CC((T)))/I$
where $\CC((T))$ is the field of formal Laurent series (that is, power series with finitely many negative
powers of $T$) and where $I=\phi^{-1}(B)$ is the Iwahori subgroup defined by $\mod T$ reduction,
\[ \phi: \SL_n(\CC[[T]]) \to \SL_n(\CC).\]
 The quotient  $\SL_n{((T))}/I$ is infinite dimensional but it is an increasing limit of finite 
 dimensional complex algebraic
varieties, and each $I$ orbit of an element $w \in W_a$ is a (generalized) ``Schubert cell" or Bruhat cell,
of dimension $\ell(w)$.  The Kazhdan Lusztig
polynomials $P_{yw}$ have non-negative coefficients and they may be interpreted as the 
local intersection cohomology Poincar\'e polynomials of one Schubert cell at a point in another Schubert cell.

However, the sheaf-convolution construction does not work in this setting because the orbits of $I$ on
$X \times X$ (where $X$ denotes the affine flag manifold) have infinite dimension and infinite codimension.
Instead, another approach is needed, which will be described later in the case of the affine Grassmannian.

\subsection{Overview}
There are strong analogies between these constructions over different fields.  The following chart gives some idea
of the parallels between the different cases.

\begin{figure}[!h]
\begin{tabular}{|c|c|c|c|c|c|}
\hline
\tiny{Field}&$\CC$ & $\FF_q$ & $\QQ_p$ & $\CC((T))$ &$\FF_q((T)$\\
\hline
\tiny{Group}&$SL_n(\CC)$ & $\SL_n(\FF_q)$ & $\SL_n(\QQ_p)$ & $\SL_n(\CC((T)))$ & $\SL_n(\FF_q((T)))$ \\
\hline \tiny{symbol} & $G_{\CC}$ & $G_q$ & $G_{\QQ_p}$ & $G_{((T))}$ & $G_{q((T))}$\\
\hline\hline\hline
\tiny{Borel/Iwahori}
 &\T{3.5}\B{2.5}$\left(\begin{smallmatrix} *&*&*&*\\&*&*&*\\&&*&*\\&&&* \end{smallmatrix}\right)$ &  
$\left(\begin{smallmatrix} *&*&*&*\\&*&*&*\\&&*&*\\&&&* \end{smallmatrix}\right)$ 
 &        $\phi^{-1}(B) $ & $\phi^{-1}(B)$& $\phi^{-1}(B_q)$\\
 \hline \tiny{symbol} &$B_{\CC}$ & $B_{\FF_q}$ & $I_p$ & $I_{[[T]]}$ & $I_{q[[T]]}$\\
 \hline \hline
 \tiny{Weyl group} & $S_n$ & $S_n$ & $W_a$ & $W_a$ & $W_a$\\
 \hline
\tiny{ Bruhat decomp }& \T{2.5}\B{1.3}$\coprod_{w\in W} BwB$ & $\coprod_{w\in W} BwB$ & 
 $\coprod_{w\in W_a} IwI$ & $\coprod_{w\in W_a} IwI$ & $\coprod_{w\in W_a} IwI$\\
  \hline \tiny{Flag manifold} & $G_{\CC}/B_{\CC}$ & $G_{\FF_q}/B_{\FF_q}$ & $G_{\QQ_p}/I_p$ 
  & $G_{((T))}/I_{((T))} $ &\B{1.2} $ G_{q((T))}/I_{q[[T]]}   $\\ 
  \hline \hline
  \tiny{Parabolic/Parahoric} & \T{3.5}\B{2.5} $\left(\begin{smallmatrix} *&*&*&*\\ *&*&*&*\\&&*&*\\
 &&*&*\end{smallmatrix}\right)$ & $\left(\begin{smallmatrix} *&*&*&*\\ *&*&*&*\\
 &&*&*\\&&*&*\end{smallmatrix}\right)$ &
 $\SL_n(\ZZ_p)$ & $ \SL_n(\CC[[T]])$ & $ \SL_n(\FF_q[[T]])$    \\
 \hline 
 \tiny{symbol} & $P_{\CC}$ & $P_{\FF_q}$ & $K_p$ & $G_{[[T]]}$ & $\B{1.2} G_{q[[T]]}   $ \\
  \hline\hline
\tiny{ Bruhat decomp }&&&\T{2.5}\B{1.3} $\coprod_{a \in A}K_paK_p$ & $\coprod_{a\in A} G_{[[T]]} a G_{[[T]]} $
& $\coprod_{a\in A} G_{q[[T]]} a G_{q[[T]]}$ \\ 
 \hline
\tiny{ Grassmannian} & $G_{\CC}/P_{\CC}$ & $G_{\FF_q}/P_{\FF_q}$ & $G_{\QQ_p}/K_p$ & $G_{((T))}/G_{[[T]]}$
&\B{1.2} $G_{q((T))}/G_{q[[T]]}$ \\
 \hline
 \end{tabular}\end{figure}

\bigskip

\begin{tcolorbox}[colback=yellow!30!white]
\begin{defn}  The affine Grassmannian is the quotient \index{affine Grassmannian}\index{Grassmannian, affine}
\[X = \SL_n(\CC((T)))/\SL_n(\CC[[T]]).\]\end{defn}   \end{tcolorbox}
If we think of $\CC((T))^n = \cup_{N=0}^{\infty}t^{-N}\CC[[T]]$ then a {\em lattice}
 in $\CC((T))^n$ is a $\CC[[T]]$ submodule $M \subset \CC((T))^n$ (meaning that it is preserved under multiplication
 by $T$) such that
 \[ T^{-N}\CC[[T]]^n \supset M \supset T^N \CC[[T]]^n\]
 for sufficiently large $N$, and
 which satisfies the determinant one conditon, $\wedge^nM = \CC[[T]]$. The affine Grassmannian is the
 set of all such lattices.  In fact the group $\SL_n(\CC((T)))$ acts transitively on the set of such lattices 
 and the stabilizer of the standard lattice
 $\CC[[T]]^n$ is the parahoric subgroup $\SL_n(\CC[[T]])$. 
 
 \[\left(\begin{matrix} T^{-1} &&&\\ 0 & T^{2} && \\ 0&0& T^{-3}& \\ 0 & 0 & 0 & T^2 \end{matrix} \right) K =
 \begin{tabular}{|c|c|c|c|c|}
 \hline
 $T^{-4}$ & & & & \\ 
 \hline
 $T^{-3}$ & && $\bullet$& \\
\hline $T^{-2}$ & && $\bullet$ &\\
\hline$T^{-1}$ & $\bullet$ &&$\bullet$&\\
\hline $T^{0}$ & $\bullet$ && $\bullet$&\\
\hline $T^1$ & $\bullet$ && $\bullet$ &\\
\hline $T^2$ & $\bullet$ & $\bullet$ & $\bullet$&$\bullet$ \\
\hline $T^3$ & $\bullet$ & $\bullet$ & $\bullet$&$\bullet$ \\
\hline  
 \end{tabular}
 \]

\subsection{Its stratification}\label{subsec-xK}
The affine Grassmannian $X$ is an infinite increasing union of projective varieties.  It has two interesting
stratifications.  The first, is by orbits of the Iwahori subgroup $I_{[[T]]}$.  These orbits are indexed by the
group of translations $A$ in affine Weyl group and the orbit $X_a=I_{[[T]]}aK$ is a Schubert cell:  it is an affine space of dimension $\ell(a)$
with $X_y \subset \overline{X}_w$ iff $y < w$ in the Bruhat order on $W_a$.  In other words (setting $I = I_{[[T]]}$ for brevity)
\[ X = \coprod_{a \in A} IaK/K \text{ because } \SL_n(\CC((T))) =  IAK\]
which is the analog of the Iwasawa decomposition of $G$.  The lattice of translations $A$ may be identified with 
the group of all diagonal matrices $\diag(T^{a_1}, T^{a_2}, \cdots T^{a_n})$ such that $\sum_i{a_i} = 0$.
Such an element may be interpreted as a cocharacter of the maximal torus $\mathcal T$ of diagonal matrices, that is,
$A \cong \chi_*(\mathcal T)$.

The second stratification is by orbits of the subgroup $K = \SL_n(\CC[[T]])$.  The (finite dimensional) Bruhat
decomposition $\SL_n(\CC) = \coprod_{w \in W} BwB$ implies that $K = \coprod_{w \in W} IwI$.  So each $K$
orbit is a union of $n!$ Schubert cells.  Let $\mathcal T$ be the torus of diagonal matrices (of determinant
one) in $\SL_n$.  Then $W$ acts on $\mathcal T$ and on its group of cocharacters 
$A = \chi_*(\mathcal T)$, which we have identified
with the lattice of translations in the affine Weyl group $W_a$.  A fundamental domain for this action is the 
{\em positive cone} $A_+ \subset \mathcal T_+$.  So the Bruhat decomposition becomes, in this case:
\[ \SL_n(\CC((T))) = \coprod_{a \in A_+} K a K \ \text{ and }  X = \coprod_{a \in A_+} KaK/K.\]
The strata are no longer cells, but each stratum has the structure of a vector bundle over a nonsingular projective
algebraic variety, so it is simply connected.

For example, in the affine Grassmannian for ${\rm PGL}_3$ consider the $K$ orbit of the lattice represented by 
\[xK=\left( \begin{matrix} T^{-2} &  & \\
0 & 1\ & \\0 & 0& 1 \end{matrix} \right)K =
 \begin{tabular}{|c|c|c|c|}
 \hline
 $T^{-3}$ & &&  \\
\hline $T^{-2}$ & $\bullet$ &  &\\
\hline$T^{-1}$ & $\bullet$ &&\\
\hline $T^{0}$ & $\bullet$ &$\bullet$& $\bullet$\\
\hline $T^1$ & $\bullet$ &$\bullet$& $\bullet$ \\
\hline $T^2$ & $\bullet$ & $\bullet$ & $\bullet$ \\
\hline  
 \end{tabular} \]
 The leading term $T^{-2}$ determines a line in $\CC^3$.  But, acting by elements of $K$ we can also 
 obtain lattices like this: 
 \begin{tabular}{|c|c|c|c|}
 \hline
 $T^{-3}$ & &&  \\
\hline $T^{-2}$ & $\bullet$ &  &\\
\hline$T^{-1}$ & $\bullet$ &$\star$&$\star$\\
\hline $T^{0}$ & $\bullet$ &$\bullet$& $\bullet$\\
\hline $T^1$ & $\bullet$ &$\bullet$& $\bullet$ \\
\hline $T^2$ & $\bullet$ & $\bullet$ & $\bullet$ \\
\hline  
 \end{tabular} which means that the orbit $KxK/K$ has the structure of a two dimensional vector bundle over $\CC P^2$.
  
\subsection{Two more views of the affine Grassmannian}

Let $\CC[T]$ be the ring of polynomials and let $\CC(T)$ be the field of rational functions $p(T)/q(T)$.  There is a
natural map
\[ G(\CC(T))/ G(\CC[T]) \to G(\CC((T)))/ G(\CC[[T]]).\]
It turns out to be an isomorphism.  We can also consider $G(\CC(T))$ to be the {\em loop group} 
\[ LG = \left\{ f:S^1 \to G| f \in \CC(T) \right\}\]
consisting of all mappings which are rational functions.  (Similarly one could consider analytic, smooth, or continuous
functions; the results are homotopy equivalent).  If $LG^+$ denotes mappings that can be extended (holomorphically,
or as a polynomial) over the origin in $\CC$ then the quotient 
\[LG/LG^+ \cong G(\CC(T))/ G(\CC[T])\]
is sometimes referred to (by physicists) as the {\em fundamental homogeneous space}.

 \section{Perverse sheaves on the affine Grassmannian}
\subsection{Spherical Hecke algebra} \index{Hecke algebra!spherical}
As in the previous section we take $G = \SL_n$, although the following is valid for any semisimple algebraic
group defined over the appropriate ring.  The Hecke algebra 
\[\mathcal H(G(\QQ_p)\slash\slash G(\ZZ_p))
\text{ resp. } \mathcal H(G(\FF_q((T)) \slash \slash G(\FF_q[[T]])) \text{ etc.}\]
of locally constant compactly supported bi-invariant functions is called the {\em spherical Hecke algebra},
that is,
\begin{tcolorbox}[colback=yellow!30!white]
\[ \mathcal H(G(\QQ_p)\slash\slash G(\ZZ_p))= C_c^{\infty}(\SL_n(\ZZ_p) \backslash \SL_n(\QQ_p) / \SL_n(\ZZ_p)).\]
\end{tcolorbox}
Recall that the representation ring $R(G)$ of a (complex) reductive group $G$ is isomorphic isomorphic to
the Weyl invariants
\[ R(G) \cong \ZZ[\chi^*(\mathcal T)]^W\]
in the group of characters of a maximal torus $\mathcal T$.  In fact, a fundamental domain for the action of $W$ on
$\chi^*(\mathcal T)$ is given by the positive Weyl chamber, $\chi^*(\mathcal T)_+$.  To such a character
$\lambda \in \chi^*(\mathcal T)_+$ one associates the irreducible representation $V_{\lambda}$
with highest weight $\lambda$.  Its trace is a character of $\mathcal T$.

As a consequence, there are many equivalent ways to view this
Hecke algebra.  \begin{enumerate}
\item By  theorems of Satake and MacDonald, there is a natural isomorphism
\[ \mathcal H(G(\QQ_p)\slash\slash G(\ZZ_p)) \cong \CC[X_*(\mathcal T)]^W\]
of the Hecke algebra with the Weyl invariants in the group algebra of the cocharacter group of the maximal torus.
\item  This in turn may be identified with the Weyl invariants $\CC[X^*(\widehat{\mathcal T})]^W$ in the
characters of the dual torus. 
\item Which, by the adjoint quotient map, is isomorphic to the group of conjugation-invariant polynomial functions
on ${\rm PGL}_n$.  (Recall that we previously identified this as a polynomial algebra, given
by the coefficients of the characteristic polynomial.)
\item This may be identified with $\CC \otimes K({\rm Rep}_{\rm PGL_n})$ (that is, the Grothendieck group of
the category of finite dimensional (rational) representations of ${\rm PGL}_n$) by associating, 
to any representation $\rho$ its character (or trace), which is a Weyl invariant polynomial function. 
\item In fact, these identifications can be made over the integers.   See [Gross] who considers the Hecke algebra
$\mathcal H_G$ of $\ZZ$-valued functions on this double coset and describes isomorphisms
\[\begin{CD} \mathcal H_G @>>> \left(\mathcal H_{\mathcal T}\otimes \ZZ[q^{\frac{1}{2}}, q^{-\frac{1}{2}}]\right)^W
@<<< {\rm Rep}(\widehat{G} )\otimes \ZZ[q^{\frac{1}{2}}, q^{-\frac{1}{2}}].\end{CD}\]
Here $\widehat{G}$ is the {\em Langlands dual group} of the group $G$. (The dual of $\SL_n$ is $PGL_n$.) 
\end{enumerate}

   \subsection{Digression:  Langlands Dual group}
   (from Wikipedia, who took it from \cite{Springer_Bourbaki})\index{Langlands dual group}\index{dual!group}

A {\em root datum}\index{root datum}
 consists of a quadruple $(X^*, \Phi, X_*, \Phi^{\vee})$ where $X^*, X_*$ are free abelian groups
of finite rank together with a perfect pairing $\langle,\rangle:X^* \times X_* \to \ZZ$, where $\Phi \subset X^*$
and $\Phi^{\vee} \subset X_*$ are finite subsets, and where there is a bijection $\Phi \to \Phi^{\vee}$, denoted
$\alpha \mapsto \alpha^{\vee}$, and satisfying the following conditions: \begin{enumerate}
\item $\langle \alpha, \alpha^{\vee} \rangle = 2$ for all $\alpha \in \Phi$
\item The map $x \mapsto x - \langle x, \alpha^{\vee}\rangle \alpha$ takes $\Phi$ to $\Phi$ and 
\item the induced action on $X_*$ takes $\Phi^{\vee}$ to $\Phi^{\vee}$.
\end{enumerate}
If $G$ is a reductive algebraic group over an algebraically closed field then it defines a root datum where
$X^*$ is the lattice of characters of a (split) maximal torus $T$, where $X_*$ is the lattice of cocharacters
of $T$, where $\Phi$ is the set of roots and $\Phi^{\vee}$ is the set of coroots.
\begin{tcolorbox}[colback=cyan!30!white]
A connected reductive algebraic group over an algebraically closed field $K$ is determined up to isomorphism
 (see \cite{Springer_AMS})
by its root datum and every root datum corresponds to such a group. 
\end{tcolorbox}
 Let $G$ be a connected reductive algebraic group over an algegraically closed field, with root datum 
$(X^*, \Phi, X_*, \Phi^{\vee})$.  Then the connected reductive algebraic
group with root datum $(X_*, \Phi^{\vee}, X^*, \Phi)$ is called the Langlands dual group and it is
denoted $G^{\vee}$ or sometimes  ${}^LG$.

Langlands duality switches adjoint groups with simply connected groups.  It takes type $A_n$ to $A_n$ but it
switches types $\Sp(2n)$ with ${\rm SO}(2n+1)$.  It preserves the type ${\rm SO}(2n)$.  A maximal torus
in the dual group $G^{\vee}$ may be identified with the dual of a maximal torus in $G$.

\subsection{Lusztig's character formula}  \index{Lusztig character formula}\index{character formula, Lusztig}
As above let $G = \SL_n(\CC((T)))$, let $K=\SL_n(\CC[[T]])$, let $I\subset K$
be the Iwahori subgroup. 
The affine flag manifold $Y=G/I$ fibers over the affine Grassmannian $X = G/K$ with fiber isomorphic to 
$K/I \cong \SL_n(\CC)/B(\CC)$ the finite dimensional flag
manifold, which is smooth.  Consequently the singularities of $I$-orbit closures in $X$ are the same
as the singularities of $I$ orbit closures in $Y$.
The $K$ orbits on $X$ are indexed by cocharacters in the positive cone.
If $\lambda =\diag(a_1,a_2,\cdots,a_n) \in \ZZ^n$ is in the positive cone (and $\sum a_i=0$) let $x_{\lambda} = 
\diag(T^{a_1}, T^{a_2}, \cdots, T^{a_n}) \in \SL_n(\CC((T)))$.  The $K$ orbit $X_{\lambda}$ corresponding to $a$ is 
$X_{\lambda}=K x_{\lambda}K/K \subset G/K$.  If $\mu \le \lambda$ then the point $x_{\mu}$ lies in the closure
of the stratum $X_{\lambda}$ and the local intersection cohomology Poincar\'e polynomial
\[ \sum_{i \ge 0}\dim(\uul{IH}{}^{2i}_{x_{\mu}}(\overline{X}_{\lambda}) )t^i = P_{\mu,\lambda}(t)\]
is given by the Kazhdan Lusztig polynomial $P_{\mu,\lambda}$ for the affine Weyl group.
Lusztig [Lu] proves
\begin{tcolorbox}[colback=cyan!30!white]
\begin{thm}
Let $\mu \le \lambda \in \chi_*(\mathcal T)_+\cong \chi^*(\mathcal T^*)_+$.  
Let $V_{\lambda}$ be the representation (of ${}^LG(\CC)$) of highest weight  $\lambda$.  It decomposes into
weight spaces $(V_{\lambda})(\mu)$ under the action of the maximal torus.  Then
\[ \dim (V_{\lambda}(\mu)) = P_{\mu\lambda}(1).\] \end{thm}
\end{tcolorbox}
That is, the local intersection cohomology Euler characteristic of the affine Schubert varieties (and of the affine
$K$ orbits) equals the weight multiplicity in the irreducible representation.  (If you wish to add up these polynomials
in order to get the intersection cohomology of the whole orbit closure, then you must do so with a shift
$t^{\ell(\lambda) - \ell(\mu)}$ corresponding to the codimension in $X_{\lambda}$ of the $I$-orbit that contains the
point $x_{\mu}$.)  Consequently Lusztig considers the full
Kazhdan Lusztig polynomial $P_{\mu\lambda}(q)$ be a {\em $q$-analog of the weight multiplicity}.  The individual
coefficients were 
eventually shown (by R. Brylinski, Lusztig, others) to equal the multiplicity of the weight $\mu$ in a certain
layer $V_{\lambda}^r/V_{\lambda}^{r-1}$ of the filtration of $V_{\lambda}$ that is induced by the principal nilpotent
element.

\subsection{Moment map}\index{moment map}
The complex torus $\{(a_1,\cdots,a_n) \in \CC^n|\ \prod_i a_i = 1\} \cong (\CC^{\times})^{n-1}$ acts on $X$ with a
the moment map (for the action of $(S^1)^{n-1}$) $\mu:X \to \mathfrak a^*$.  Each fixed point 
$x_{\lambda}=(T^{\lambda_1}, \cdots, T^{\lambda_n})$ 
corresponds to a cocharacter $\lambda = (\lambda_1,\cdots,\lambda_n)
\in \ZZ^n$ (with $\sum_i \lambda_i = 0$).    The torus action preserves both stratifications and the image of each
stratum closure is a convex polyhedron.  We can put all this information on the same diagram.
Fix $\lambda \in \chi_*(\mathcal T)_+$.  Let $V_{\lambda}$ be the irreducible representation of $\text{\rm PGL}_3$ with
highest weight $\lambda$.  Let $X_{\lambda}$ be the $K$-orbit of the point  $x_{\lambda}$.  Then $\mu(X_{\lambda})$
is the convex polyhedron spanned by the $W$ orbit of the point $\lambda$.  The lattice points $\mu$
inside this polyhedron correspond to the weight spaces $V_{\lambda}(\mu)$.  At each of these points the
Kazhdan Lusztig polynomial $P_{\mu,\lambda}(t)$ gives the dimension of this weight space.

For $\SL_3$ the moment map image gives the triangular lattice which can be interpreted as the weight lattice
for $\text{\rm PGL}_3$.  The moment map image of the first stratum is a hexagon and in fact the first stratum has
the structure of a vector bundle over the flag manifold $F_{1,2}(\CC^3)$.  The closure of this stratum consists of
adding a single point, so it is the Thom space of this bundle.  The local intersection cohomology at the vertex
(which appears at the origin in the diagram) is the primitive cohomology of the flag manifold (that is, $1+t$).

\begin{figure}[!h]

\begin{tikzpicture}[scale=4]
\foreach \r in {0,1,...,{\rowsh}} {
 \draw($(0,{\r*sqrt(3)})$)  -- ($({\r},0)$) ; 
 \draw($(0,{\rowsh*sqrt(3)}) + (\r,0)$) -- ($(\rowsh,0) + (\r,0)$); 
 \draw($(\rowsh,{\rowsh*sqrt(3)})+(\r,0)$) -- ($(\rows,0) + (0, {\r*sqrt(3)})$);
                                          }
 \foreach \r in {0,1,...,\rows} {
  \draw($(0,{\r*sqrt(3)/2})$) -- ($(\rows,{\r*sqrt(3)/2})$) ;  
                                           }
\foreach \r in {0,1,...,\rowsh}{  
 \draw($(\r,0)$) -- ($((\r,0)+({\rows/2}, {\rows*sqrt(3)/2})$);
 \draw($(0,{\r*sqrt(3)})$) -- ($(\rowsh,{\rowsh*sqrt(3)}) - (\r,0)$); 
 \draw($(\rowsh,0)+(\r,0)$) -- ($(\rows,{\rowsh*sqrt(3)}) - (0, {\r*sqrt(3)})$);
                                               }
   \draw[dashed,red] (2,0) -- (2,{2*sqrt(3)});
   \draw[dashed,red] (0, {sqrt(3)/3}) -- (4, {5*sqrt(3)/3});
    \draw[dashed,red]( 4, {sqrt(3)/3}) -- (0, {5*sqrt(3)/3});        

\draw[fill] (3, {sqrt(3)}) circle [radius=.025];   \node[below] at (3,{sqrt(3)}) {$1$};
\draw[fill] (1,{sqrt(3)})  circle [radius=.025];  \node[below] at (1,{sqrt(3)}) {$1$};
\draw[fill] (1.5, {sqrt(3)/2}) circle [radius = .025];  \node[below] at (1.5, {sqrt(3)/2}) {$1$};
 \draw[fill] (2.5, {sqrt(3)/2}) circle [radius = .025]; \node[below] at (2.5, {sqrt(3)/2}) {$1$};
 \draw[fill] (1.5, {3*sqrt(3)/2}) circle [radius = .025]; \node[below] at (1.5, {3*sqrt(3)/2}) {$1$};
 \draw[fill, red] (2.5, {3*sqrt(3)/2}) circle [radius = .025];  \node[below] at (2.5,{3*sqrt(3)/2}) {$1$};
 \draw[fill](2,{sqrt(3)}) circle [radius = .025]; 
 \node[below] at (2,{sqrt(3)}) {$1+t$};

 \draw[blue,very thick] (2.5,{3*sqrt(3)/2}) -- (3,{sqrt(3)});
 \draw[blue,very thick] (3,{sqrt(3)}) -- (2.5, {sqrt(3)/2});
 \draw[blue,very thick] (2.5, {sqrt(3)/2}) --(1.5, {sqrt(3)/2});
 \draw[blue,very thick](1.5,{sqrt(3)/2}) -- (1, {sqrt(3)});
 \draw[blue,very thick](1,{sqrt(3)}) -- (1.5, {3*sqrt(3)/2});
 \draw[blue,very thick](1.5,{3*sqrt(3)/2}) -- (2.5, {3*sqrt(3)/2});

 \node[above right] at (3,{sqrt(3)}) {$(1,0,-1)$};
 \node[above right] at (2.5,{3*sqrt(3)/2}) {$(0,1,-1)$};
 \node[above right] at (1.5,{3*sqrt(3)/2}) {$(0,-1,1)$};
 \node[above right] at (2.5,{sqrt(3)/2}) {$(-1,0,1)$};
 \node[above right] at (1.5,{sqrt(3)/2}) {$(-1,1,0)$};
 \node[above right] at (1, {sqrt(3)}) {$(1,-1,0)$};
                                            
\end{tikzpicture}
\caption{Affine moment map}\label{fig-affinemomentmap}
\end{figure}
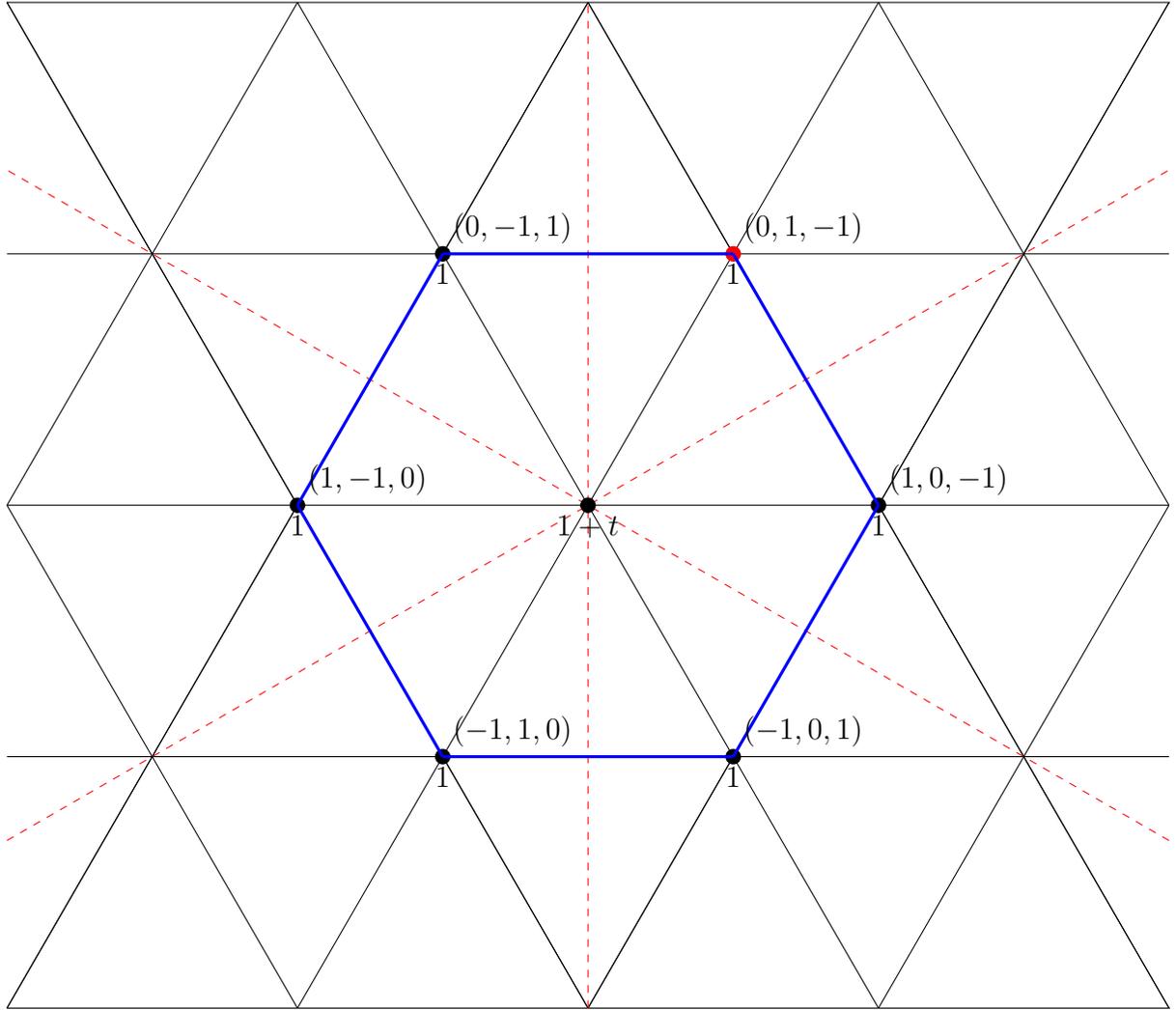

Figure \ref{fig-affinemomentmap} represents the moment map image of the affine Grassmannian for $\SL_3$; it also
equals the weight diagram for $\text{\rm PGL}_3(\CC)$.  The dotted red lines are reflecting hyperplanes for
the Weyl group.  The red dot is a highest weight $\lambda$ for $\text{\rm PGL}_e$; 
the other dots are their Weyl images.  The blue hexagon is the outline of the moment map images; it is also the 
collection of weights in the irreducible representation of highest weight $\lambda$.  The $1$ and $1+t$ beneth
the dots are the Kazhdan Lusztig polynomials.
\medskip

 \subsection{Perverse sheaves on $X$}
 Throughout this section for simplicity let $G = \SL_n(\CC((T)))$, $K = \SL_n(\CC[[T]])$ and $X = G/K$.
 We would like to imitate the construction with the flag manifold, and create a convolution product for
 sheaves on $X \times X$ that are constructible with respect to the stratification by $G$ orbits, that is,
 if $p_{ij}:X\times X \times X \to X \times X$ as before, set
 \[ A\b \circ B\b = Rp_{13*}(p_{12}^*(A\b) \otimes p_{23}^*(B\b)).\]
 Unfortunately the orbits of $G$ on $X \times X$ have infinite dimension and infinite codimension so this
 simply does not make sense.  V. Ginzburg and (later) K. Vilonen and I. Mirkovi\v{c} found a way around this problem.

   Let $\mathcal P(X)$ denote the category of perverse sheaves, constructible with respect to the above orbit
   stratification of $X$.  (Ginzburg shows this is equivalent to the category of $K$-equivariant perverse sheaves
   on $X$.)  Since each stratum is simply connected the local systems associated to these 
   sheaves are trivial.  It turns out that the intersection cohomology sheaves live only in even degrees and this
   implies that every perverse sheaf is isomorphic to a direct sum of $\uul{IC}^{\bullet}$ sheaves of
   stratum closures.
   
   Mirkovi\v{c} and Vilonen define a tensor product structure on $\mathcal P(X)$ as follows.  Consider the diagram
   \[
   \begin{CD}
   X \times X @<{p}<< G \times X @>{q}>> G\times_K X @>{m}>> X
   \end{CD}\]
   Here, $k.(g,x) = (gk,k^{-1}x)$ so that $G \times_K X$ is a bundle over $X$ whose fibers are copies
   of $X$, and $m(g,x) = gx$. If $A\b, B\b \in P(X)$ it turns out that there exists $C\b$ a perverse sheaf
   on $G \times_KX$, constructible with respect to the $G$ orbits on this space, such that 
   \[ q^*(C\b) = p^*(\pi_1^*(A\b) \otimes \pi_2^*(B\b)\]
   where $\pi_1, \pi_2:X \times X \to X$ are the two projections.    Then set $A\b \circ B\b = Rm_*(C\b)$.  
   
\begin{tcolorbox}[colback=cyan!30!white]
   \begin{thm}\label{thm-geomSatake}
If $A\b, B\b \in \mathcal P(X)$ then so is $A\b \circ B\b$.  The functor $h:A\b \mapsto H^*(X;A\b)$
is exact and it induces an equivalence of categories
\[ \mathcal P(X) \sim \text{\rm Rep}({}^LG)\]
which takes $A\b \circ B\b$ to the tensor product $h(A\b)\otimes h(B\b)$
of the associated representations.  If $\lambda \in \chi_*(\mathcal T)_+$ then $h(\uul{IC}^{\bullet}
(\overline{X}_{\lambda})) = V_{\lambda}$ is the irreducible representation of highest weight $\lambda$.   
   \end{thm}\end{tcolorbox}
   
 \subsection{}  Although it sounds intimidating, the convolution product of sheaves is exactly parallel to the previous case
 of the finite (dimensional) flag manifold.  Consider the weight diagram for $\SL_3$ and the moment map image of
 torus fixed points in the affine Grassmannian for ${\rm PGL}_3$.  The coordinate lattices are indicated on 
 Figure \ref{fig-SL3},
 for example, the point $(1,0,0)$ corresponds to the lattice $sK=$\ 
 \begin{tabular}{|c|c|c|c|}
\hline $T^{-2}$ &  &  & \\
\hline$T^{-1}$ & $\bullet$ &&\\
\hline $T^{0}$ & $\bullet$ &$\bullet$& $\bullet$ \\
\hline $T^1$ & $\bullet$ &$\bullet$& $\bullet$ \\
\hline $T^2$ & $\bullet$ & $\bullet$ & $\bullet$ \\
\hline  
 \end{tabular} \
 where $s$ is a certain simple reflection in the affine Weyl group.  So the orbit  of $G(\CC[[T]])$ in $ X \times X$ 
corresponding to this element consists of the set of pairs of lattices
 $L_0\overset{s}{\longrightarrow}L$ that are in relative position $s$, that is,
 \[ \left\{ (L_0, L_1)| \ L_0 \subset L_1 \subset T^{-1}L_0 \text{ and }\ \dim(L_1/L_0)=1 \right\}.\]
If we fix $L_0$ to be the standard lattice, then we see that the orbit $\mathcal O_s \subset X$ is
isomorphic to $\PP^1$ so its $\uul{IC}$ sheaf is the constant sheaf.

\renewcommand*\rows{6}
\renewcommand*\rowsh{3}
\newcommand*\xx{.5}
\newcommand*\yy{sqrt(3)/2}
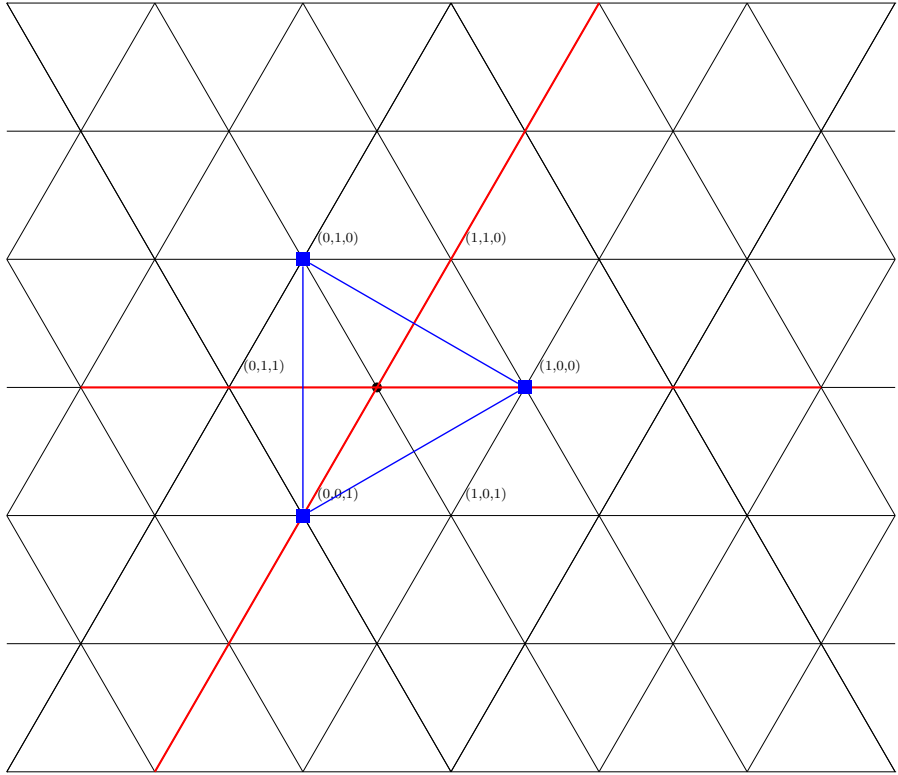
\begin{figure}[!h]
\resizebox{12cm}{!}{
\begin{tikzpicture}[scale=3]
\tikzstyle{dc}   = [circle, minimum width=8pt, draw, inner sep=0pt]

\foreach \r in {0,1,...,{\rowsh}} {
 \draw($(0,{\r*sqrt(3)})$)  -- ($({\r},0)$) ; 
 \draw($(0,{\rowsh*sqrt(3)}) + (\r,{0})$) -- ($(\rowsh,0) + (\r,0)$); 
 \draw($(\rowsh,{\rowsh*sqrt(3)})+(\r,0)$) -- ($(\rows,0) + (0, {\r*sqrt(3)})$);
                                          }
 \foreach \r in {0,1,...,\rows} {
  \draw($(0,{\r*sqrt(3)/2})$) -- ($(\rows,{\r*sqrt(3)/2})$) ;  
                                           }
\foreach \r in {0,1,...,\rowsh}{  
 \draw($(\r,0)$) -- ($((\r,0)+({\rows/2}, {\rows*sqrt(3)/2})$);
 \draw($(0,{\r*sqrt(3)})$) -- ($(\rowsh,{\rowsh*sqrt(3)}) - (\r,0)$); 
 \draw($(\rowsh,0)+(\r,0)$) -- ($(\rows,{\rowsh*sqrt(3)}) - (0, {\r*sqrt(3)})$);
};

\draw[fill] (2+\xx, {sqrt(3)+\yy}) {circle[radius=0.030]};
 
\draw[red, very thick] ({0+\xx},{sqrt(3)+\yy}) -- ({5+\xx},{sqrt(3)+\yy});
\draw[red, very thick] ({.5+\xx},{\yy-sqrt(3)/2}) -- (3.5+\xx, {5*sqrt(3)/2+\yy});

\node[fill, blue, label=above right:{\tiny(1,0,0)}] (n100) at (3+\xx,{sqrt(3)+\yy}) {};
\node[fill, blue, label= above right:{\tiny(0,1,0)}] (n010) at (1.5+\xx, {3*sqrt(3)/2+\yy}) {};
\node[label = above right:{\tiny(1,1,0)}] (n110) at (2.5+\xx, {3*sqrt(3)/2+\yy}) {};
\node[fill, blue, label = above right:{\tiny(0,0,1)}] (n001) at (1.5+\xx, {sqrt(3)/2+\yy}) {};
\node[label = above right:{\tiny(0,1,1)}] (n011) at (1+\xx, {sqrt(3)+\yy}) {};
\node[label = above right:{\tiny(1,0,1)}] (n101) at (2.5+\xx, {sqrt(3)/2+\yy}) {};

\draw[blue, thick] (n100) -- (n010);
\draw[blue,thick] (n010) -- (n001);
\draw[blue,thick] (n001) -- (n100);

\end{tikzpicture}
}
\caption{$SL_3$: Reflecting hyperplanes and
{moment map image of $\mathcal O_{(1,1,0)}$}}\label{fig-SL3}
\end{figure}

Let us consider
 the convolution product of this sheaf with itself.  Thus the total space of the correspondence 
 consists of triples of lattices
 \[C = \left\{ (L_0,L_1,L_2)|\  L_0 \overset{s}{\longrightarrow} L_1 \overset{s}{\longrightarrow}L_2\right\}\]
 in their appropriate relative positions.  Apply $\pi_{13}$ to obtain the correspondence
 \[ \pi_{23}(C)=\left\{ (L_0, L_2)| \ L_0 \subset L_2 \subset T^{-2}L_0 \ \text{ and }\ \dim(L_2/L_0) = 2\right\}.\]
 Again, taking $L_0$ to be the standard lattice, we need to understand the decomposition into orbits
 of $R\pi_{23*}(\uul{\QQ}{}_C)$.  There are two types of such lattices:  the first type  considered in $\S \ref{subsec-xK}$ consists of those lattices in the orbit $xK$ that was, that is, lattices like this:
 \begin{tabular}{|c|c|c|c|}
 \hline
 $T^{-3}$ & &&  \\
\hline $T^{-2}$ & $\bullet$ &  &\\
\hline$T^{-1}$ & $\bullet$ &$\star$&$\star$\\
\hline $T^{0}$ & $\bullet$ &$\bullet$& $\bullet$\\
\hline $T^1$ & $\bullet$ &$\bullet$& $\bullet$ \\
\hline $T^2$ & $\bullet$ & $\bullet$ & $\bullet$ \\
\hline  
\end{tabular}\
which project under the moment map to an image that contains the $W$-translates of $(2,0,0)$.  The second type is
lattices like this:
 \begin{tabular}{|c|c|c|c|}
 \hline
 $T^{-3}$ & &&  \\
\hline $T^{-2}$ &  &  &\\
\hline$T^{-1}$ & $\bullet$ &$\bullet$&\\
\hline $T^{0}$ & $\bullet$ &$\bullet$& $\bullet$\\
\hline $T^1$ & $\bullet$ &$\bullet$& $\bullet$ \\
\hline $T^2$ & $\bullet$ & $\bullet$ & $\bullet$ \\
\hline  
 \end{tabular}\
 which project under the moment map to an image containing $(1,1,0)$.  That is,
 \[ m: C \to \mathcal O_{(2,0,0)} \cup \mathcal O_{(1,1,0)}.\]
 The map $m$ is guaranteed to be semi-small and $Rm_*(\uul{\QQ}{}_C)$ breaks into a direct sum of copies
 of $\uul{IC}$ sheaves of these two strata the multiplicities equal to the number of components of the fiber.
 One checks that the multiplicity equals one.
 
\renewcommand*\rows{6}
\renewcommand*\rowsh{3}
\renewcommand*\xx{.5}
\renewcommand*\yy{sqrt(3)/2}

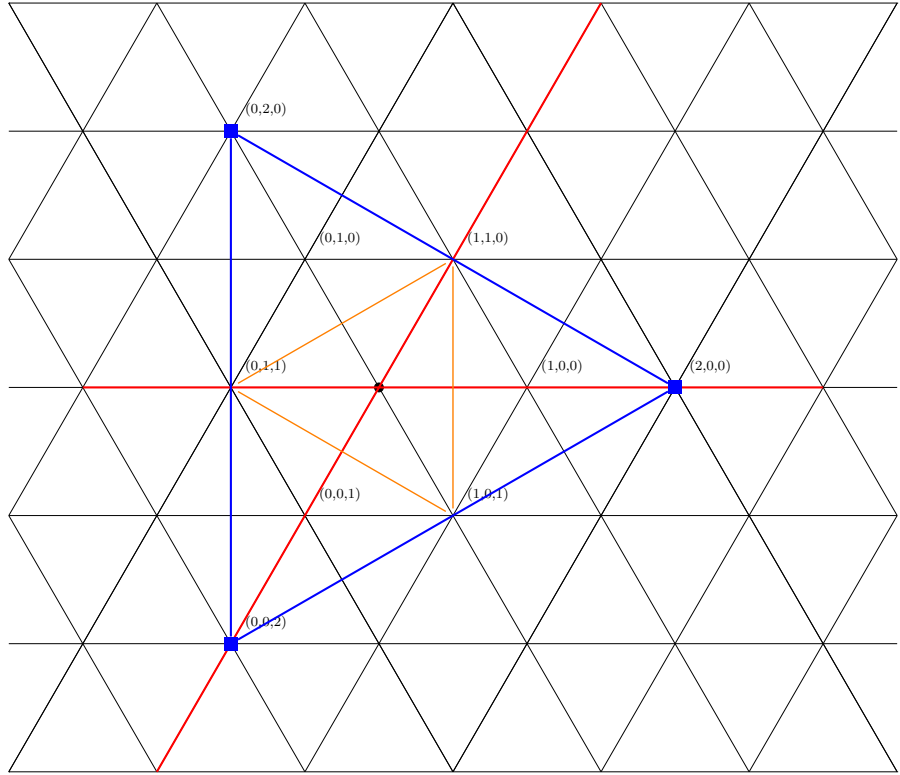
\begin{figure}[!h]
\resizebox{12cm}{!}{
\begin{tikzpicture}[scale=3]
\tikzstyle{dc}   = [circle, minimum width=8pt, draw, inner sep=0pt]

\foreach \r in {0,1,...,{\rowsh}} {
 \draw($(0,{\r*sqrt(3)})$)  -- ($({\r},0)$) ; 
 \draw($(0,{\rowsh*sqrt(3)}) + (\r,{0})$) -- ($(\rowsh,0) + (\r,0)$); 
 \draw($(\rowsh,{\rowsh*sqrt(3)})+(\r,0)$) -- ($(\rows,0) + (0, {\r*sqrt(3)})$);
                                          }
 \foreach \r in {0,1,...,\rows} {
  \draw($(0,{\r*sqrt(3)/2})$) -- ($(\rows,{\r*sqrt(3)/2})$) ;  
                                           }
\foreach \r in {0,1,...,\rowsh}{  
 \draw($(\r,0)$) -- ($((\r,0)+({\rows/2}, {\rows*sqrt(3)/2})$);
 \draw($(0,{\r*sqrt(3)})$) -- ($(\rowsh,{\rowsh*sqrt(3)}) - (\r,0)$); 
 \draw($(\rowsh,0)+(\r,0)$) -- ($(\rows,{\rowsh*sqrt(3)}) - (0, {\r*sqrt(3)})$);
};

\draw[fill] (2+\xx, {sqrt(3)+\yy}) {circle[radius=0.030]};
 
\draw[red, very thick] ({0+\xx},{sqrt(3)+\yy}) -- ({5+\xx},{sqrt(3)+\yy});
\draw[red, very thick] ({.5+\xx},{\yy-sqrt(3)/2}) -- (3.5+\xx, {5*sqrt(3)/2+\yy});

\node[label=above right:{\tiny(1,0,0)}] (n100) at (3+\xx,{sqrt(3)+\yy}) {};
\node[label= above right:{\tiny(0,1,0)}] (n010) at (1.5+\xx, {3*sqrt(3)/2+\yy}) {};
\node[label = above right:{\tiny(1,1,0)}] (n110) at (2.5+\xx, {3*sqrt(3)/2+\yy}) {};
\node[label = above right:{\tiny(0,0,1)}] (n001) at (1.5+\xx, {sqrt(3)/2+\yy}) {};
\node[label = above right:{\tiny(0,1,1)}] (n011) at (1+\xx, {sqrt(3)+\yy}) {};
\node[label = above right:{\tiny(1,0,1)}] (n101) at (2.5+\xx, {sqrt(3)/2+\yy}) {};

\draw[orange, thick] (n110) -- (n011);
\draw[orange,thick] (n011) -- (n101);
\draw[orange,thick] (n101) -- (n110);

\node[fill, blue, label=above right:{\tiny(2,0,0)}] (n200) at (4+\xx, {sqrt(3)+\yy}) {};
\node[fill, blue, label=above right:{\tiny(0,2,0)}] (n020) at (1+\xx, {2*sqrt(3)+\yy}) {};
\node[fill, blue, label=above right:{\tiny(0,0,2)}] (n002) at ({1+\xx}, {0+\yy}) { };

\draw[blue, very thick] (n200) -- (n020);
\draw[blue, very thick] (n020) -- (n002);
\draw[blue, very thick] (n002) -- (n200);

\end{tikzpicture}
} 

\caption{$\mathcal O_{(1,0,0)} \circ \mathcal O_{(1,0,0)} = \mathcal O_{(2,0,0)} + \mathcal O_{(1,1,0)}$}
\end{figure}

 By taking the cohomology of these sheaes we obtain highest weight representations:  $V_{(1,0,0)} = \text{std}$,
 $V_{(1,1,0)} = \text{std}^{\vee}$, and $V_{(2,0,0)} = \wedge^2(\text{std})$ (where std is the standard representation)
 so that
 \[ \text{std} \otimes \text{std} \cong \text{std}^{\vee} \oplus \wedge^2(\text{std}).\]
 
 \subsection{Remarks}
 Theorem \ref{thm-geomSatake} should be regarded as a categorification of Satake's isomorphism.  
\index{Satake isomorphism}   In fact, taking the $K$ group of the
 Grothendieck group on both sides gives
 \[ K(\mathcal P(X)) \cong \mathcal H(G,K) \cong \chi_*(\mathcal T)^W\]
 which is the classical Satake isomorphism.
 
 If we could duplicate the construction in the finite dimensional case we would consider $G$ orbits on
 $X \times X$.  Ginzburg, Mirkovi\v{c} and Vilonen replace $X \times X$ with $G \times_K X$ which is a
 fiber bundle over $X$ with fiber isomorphic to $X$.  They replace the $G$ orbits with the strata $S_{\lambda,\mu}$
 which is a fiber bundle over $X_{\mu}$ with fiber isomorphic to $X_{\lambda}$.
 
 It is totally nonobvious that the convolution of perverse sheaves is perverse.  This depends on the fact
 that the mapping $m:G \times_KX \to X$ is semi-small in a very strong sense.  For each $\lambda, \mu
 \in \chi_*(\mathcal T)_+$ let $S_{\lambda,\mu} = p^{-1}(X_\lambda) \times_{K} X_{\mu}$.   These form a
 stratification of $G \times_KX$.  It turns out
 that the restriction $m:\overline{S}_{\lambda,\mu} \to X$ is semi-small (onto its image, which is a union of
 strata $X_{\tau}$).  This implies that $Rm_*(\mathcal E[{d_{\lambda,\mu}}])$ is perverse, 
 for any locally constant sheaf $\mathcal E$ on $S_{\lambda,\mu}$ (where $d_{\lambda,\mu} =
 \dim(S_{\lambda,\mu})$.

\subsection{Tannakian category}(see \cite{DeligneMilne})\index{Tannakian category}\index{category!Tannakian}
There is a general theorem that a reductive algebraic group (say, over an algebraically closed field) can be
recovered from its cartegory of representations.  More generally if $\mathcal C$ is a tensor category, that is,
an abelian category together
with a ``tensor" structure $(A, B \in \mathcal C \implies A \circ B \in \mathcal C)$ that is commutative and
associative in a functorial way (satisfies the {\em associativity constraint} and the
{\em commutativity constraint}, and if $h:\mathcal C \to \{VS\}$ is a rigid {\em fiber functor} (that is, an 
exact functor to vector spaces such that $h(A \circ B) = h(A)\otimes h(B)$) then the group of automorphisms of $h$ 
(that is, the group of natural transformations $h \to h$) is an algebraic group $G$ whose
category of representations is equivalent to the original category $\mathcal C$.

It turns out that $\mathcal P(X)$ is such a Tannakian category and that $h$ is a rigid fiber functor.  Therefore the
group of automorphisms of $h$ is isomorphic to the Langlands dual group ${}^LG$.


\section{Cellular perverse sheaves}\label{sec-cellular-perverse-sheaves}
\subsection{} The theory of cellular perverse sheaves was developed by R. MacPherson and described
in lectures \cite{MP1,MP2} and
was explained in detail in the thesis \cite{Vybornovthesis2, Vybornovthesis} of M. Vybornov.  
Cellular perverse sheaves interpolate between
cellular sheaves and cellular cosheaves in a remarkable way.   The theory was extended and simplified
in \cite{Polishchuk, Vybornov, Vybornov2}. 

Suppose $K$ is a finite simplicial complex\footnote{Statements in this section hold without modification for any
finite strongly  regular cell complex, meaning that  the closure of each cell is homeomorphic to a
closed ball, and the intersection of two (closed) cells is a cell, or else empty.  Beware that the Schubert cell decomposition
of the flag and Grassmann manifold is not regular.}.  We will be sloppy and identify $K$ with its geometric
realization $|K|$.  Each simplex $\sigma$ is a closed subset of $K$ and we denote its interior by $\sigma^o$ which
we refer to as a {\em cell}.  If $\sigma \subseteq\tau$ we write $\sigma \leq \tau$ and similarly for proper inclusions.
We also write $\sigma \leftrightarrow\tau$ if {\em either} $\sigma \le \tau$ or $\tau \le \sigma$.
The barycenter of a simplex $\sigma$ is denoted $\hat\sigma$. An $r$-simplex $\theta$ in the
barycentric subdivision $K'$ is the span $\theta = \langle \hat\sigma_0, \hat\sigma_1,\cdots, \hat\sigma_r \rangle$
where $\sigma_0 < \sigma_1 < \cdots < \sigma_r$.  If $S$ is a collection of vertices of $K$ then its {\em span} 
\index{span (simplicial)} consists of those simplices $\sigma$ such that every vertex of $\sigma$ is in $S$.

\subsection{Simplicial derived category}  \label{subsec-simplicialderived}\index{derived category!simplicial}
We summarize some results in the thesis of Allen Shepard   
\cite{Shepard} that were developed in a seminar at Brown University in 1977-78.
We consider sheaves of $\QQ$-vector spaces but the coefficients may be taken in any field.

  Let $D^b_c(K)$ denote the bounded derived category of
sheaves on $K$ that are cohomologically constructible with respect to the stratification by cells (i.e. interiors
of simplices).

Recall \S \ref{subsec-simplicialsheaf} that a simplicial sheaf \index{sheaf!simplicial} $A$
assigns to each simplex $\sigma$ a finite dimensional $\QQ$-vector space $A_{\sigma}$ and to each face 
$\sigma < \tau$ a  {\em restriction} homomorphism
$s_{\sigma\tau}:A_{\sigma} \to A_{\tau}$  in a way that is compatible with compositions.
 Morphisms of simplicial sheaves
 are required to be compatible with these restriction mappings.    The category of simplicial sheaves is
 equivalent to the subcategory of $D^b_c(K)$ whose objects have stalk cohomology only in degree $0$ (that is,
 the category  of $\bar 0$-perverse sheaves).
\begin{tcolorbox}[colback=yellow!30!white]\index{Hom!simplicial sheaf}
If $A, B$ are simplicial sheaves, the sheaf $\uul{\Hom}(A,B)$ is the simplicial sheaf
\[ \uul{\Hom}(A,B)(\sigma) = \Hom(A|St(\sigma), B|St(\sigma))\]
where $St(\sigma)$ denotes the open star of $\sigma$.  
 \end{tcolorbox}

So the group of sections $\uul{\Hom}(A,B)(\sigma)$ consists of homomorphisms
$A(\sigma) \to B(\sigma)$ that are compatible with all corestriction maps emanating from $\sigma$.

 Let $C$ be a cell (or stratum) of $K$.  
Let $i:C \to K$ be the inclusion.  The closure of its image is a simplex $\sigma$ and 
$Ri_*(\uul{\QQ}{}_C) = i_*(\uul{\QQ}{}_C) = \uul{\QQ}{}_{\sigma}$ is
the constant sheaf on the (closed) simplex $\sigma$.  To streamline notation we denote this by
$\uul{\QQ}{}_{\sigma*}$ and set $\uul{\QQ}{}_{\sigma!} = i_!(\uul{\QQ}{}_C)$, the extension by zero.
\quash{
  For any simplex $\sigma$ with inclusion $i_{\sigma}:\sigma \to K$ by abuse of notation we denote
  the constant sheaf  $(i_{\sigma})_*\uul{\QQ}{}_{\sigma}$ on $\sigma$  by $\uul{\QQ}{}_{\sigma*} = \uul{\QQ}{}_{\sigma}$.
  If $C = \sigma^o$ is a cell with $j:C \to K$ the inclusion then we denote the constant sheaf on $C$ by
   $\uul{\QQ}{}_{\sigma!} = j_!\uul{\QQ}{}_C$ .}
     A key exercise is to verify the following:
\begin{tcolorbox}[colback=cyan!30!white]
In the category $Sh^{\Delta}_K$ of simplicial sheaves each\index{injective!elementary}\index{elementary injective}
$\uul{\QQ}{}_{\sigma*}$  is injective and every injective sheaf is a direct sum of such {\em elementary} injectives.  \qed
\end{tcolorbox}

 Every simplicial sheaf of  $\QQ$ vector spaces has a canonical injective imbedding,
 \[ A \to I^0(A)= \bigoplus_{\sigma}(A_{\sigma})_{\sigma*} = \bigoplus_{\sigma}A_{\sigma}\otimes \uul{\QQ}{}_{\sigma*}\]
 which, by repeated application to the cokernel gives rise to a canonical injective resolution 
 \[A \to I^0(A) \to I^1(A) \cdots.\]
where $I^k(A)$ is the sum over flags with codimension one steps  (cf.~\cite{Shepard} \S 1.4),
\[  I^k(A) = \bigoplus_{\tau_k<\tau_{k-1}< \cdots< \tau_0} A_{\tau_k}\otimes \uul{\QQ}{}_{\tau_k*}.\]
So its global sections $\Gamma(K, I^k(A))$ may be identified with the simplicial $k$-cochains $C^k(K';A)$
of the first barycentric subdivision $K'$ of $K$.

 Consequently any (bounded below) complex $A\b$ of simplicial sheaves has a canonical injective resolution
 $A\b \overset{\sim}{\longrightarrow} T\b$  where $T\b$ is the single complex associated to the double
 complex $I^p(A^q)$.  This gives canonical models for derived functors such as
 the cohomology (derived functor of global sections),
 \[ H^i(K,A\b) = R^i\Gamma(K,A\b) = H^i(\Gamma(K, T\b))\]
 and
 \[ {\rm R}\uul{\rm Hom}\b(B\b, A\b) = \uul{\Hom}\b(B\b, T\b)\]
 (Recall that $\uul{\Hom}\b(B\b,T\b)$ is defined as the single complex associated to the double complex
 $\uul{\Hom}(B^i, T^j)$.)  The main properties of these derived functors can then be proven directly using
 these canonical models. In \cite{Shepard}, \cite{KS}\S 8.1.11 is proven:

\begin{tcolorbox}[colback=cyan!30!white]
The natural functor (\S \ref{subsec-simplicialsheaf})
$D^b(Sh_K^{\Delta}) \to D^b_c(K)$ from the bounded derived category of simplicial sheaves on $K$ to
the bounded derived category of sheaves that are (cohomologically) constructible with respect to the cell decomposition
of $K$, is an equivalence of categories.   \end{tcolorbox}
 
  The dualizing sheaf $\DD\b$\index{dualizing sheaf!simplicial} has a canonical injective model in the category
  $D^b(Sh_K^{\Delta})$.  \index{dualizing sheaf}  If $j:\sigma<\tau$ is a
 codimension one face, denote the composition $\uul{\QQ}{}_{\tau*} \to j_*j^*\uul{\QQ}{}_{\tau*} =\uul{ \QQ}{}_{\sigma*}$
 by $\partial_{\tau,\sigma}$.
  Choose an ordering of each simplex.   
  Let $[\tau:\sigma] = \pm 1$ denote the incidence number which is $+1$ if the orientation
 of $\sigma$ followed by an inward pointing vector agrees with the orientation of $\tau$. 
 Set $[\tau:\sigma] = 0$ if $\sigma$ is not a codimension one face of $\tau$. Then:
 \[ \DD^{-j} = \bigoplus_{\dim(\sigma) = j} \uul{\QQ}{}_{\sigma*}\]
 with $d:\DD^{-j} \to \DD^{-j+1}$ given by 
 \[\bigoplus_{\dim\tau=j}\bigoplus_{\dim\sigma=j-1}[\tau:\sigma]\partial_{\tau,\sigma}\]
It is independent, up to unique quasi-isomorphism, of the choice of orientation.  (\cite{Shepard})

For any $A\b$ in $D^b(Sh^{\Delta}_K)$ this gives a canonical construction of the dual ${\rm D}A\b = \Hom\b(A\b, \DD\b)$,
and the natural isomorphism ${\rm D D }A\b \cong A\b$ may be proven directly (\cite{Shepard}).  Combinatorially,
${\rm DD} A\b$ corresponds to the sheaf induced by $A\b$ on the second barycentric subdivision of $K$.

 If $f:A\b \to B\b$ is a morphism of injective simplicial sheaves then $\ker(f)$ is not necessarily injective, and so the
truncation functor $\tau_{\le a}$ does not preserve injectives.

\newcommand{\QQQ}{\uul{\QQ}{}}
\newcommand{\pd}{ {}^{\delta}}

\subsection{Local calculations}\label{subsec-singlecell}
Let $x \in \sigma$.  It lies in the
interior $\tau^o$ of some simplex $\tau \le \sigma$.  Let $j_x:\{x\} \to \sigma$.  Then one checks that

\begin{itemize}
\item $H^0(j_x^*(\uul{\QQ}{}_{\sigma*})) = \QQ$ and $H^n(j_x^*(\uul{\QQ}{}_{\sigma*}))=0$ for $n \ne 0$
\item $H^n(j_x^!(\uul{\QQ}{}_{\sigma*})) = 0$ for all $n$
\item $H^n(j_x^*(\uul{\QQ}{}_{\sigma!})) = 0$ for all $n$
\item $H^{\dim(\sigma)}(j_x^!(\uul{\QQ}{}_{\sigma!})) = \QQ$ and $H^n(j_x^!(\uul{\QQ}{}_{\sigma!}))=0$ for  $n \ne \dim(\sigma)$
\end{itemize}


If $A\b, B\b\in D^b(Sh^{\Delta}_K)$ are (constructible) complexes of sheaves, and if $B\b \to I\b$ is an injective
resolution of $B\b$, recall that \index{Ext}
 \[\Ext^0(A\b,B\b)=\hHom_{D^b(Sh_K)}(A\b,B\b) = H^0(K;R\uul{\hHom}\b(A\b,B\b))=H^0(K; \uul{\hHom}\b(A\b, I\b))\] 
 so $\Ext^i(A\b,B\b) = \Ext^0(A\b, B\b[i])$
is the $i$-th cohomology group of the single complex associated to the double complex $\Hom_{D^b(Sh_K)}(A\b, I\b)$.
It is an exercise to calculate these $\Ext$ groups for simplices $\sigma, \tau$, and they turn out to live in
a single degree. 
Set $d(\sigma) = \dim(\sigma)$.  Using the canonical injective resolution
of $\uul{\QQ}{}_{\tau!}$ and noting that $\uul{\QQ}{}_{\tau*}$ is injective, we find:

\begin{align*}
 \Ext^i(\QQQ_{\sigma*}, \QQQ_{\tau*}) &=\begin{cases} \QQ &\text{ if }\ i=0 \ \text{ and }\sigma \ge \tau \\ 
 0 &\text{otherwise}\end{cases}\\ 
 \Ext^i(\QQQ_{\sigma!}, \QQQ_{\tau!}) &= \begin{cases} \QQ &\text{ if }\ i=   d(\tau)-d(\sigma)
 \text{ and }\ \sigma \le \tau \\ 0 &\text{otherwise}\end{cases}\\
\Ext^i(\QQQ_{\sigma!}, \QQQ_{\tau*}) &= \begin{cases} \QQ &\text{ if }\ i=0 \text{ and }\ \sigma = \tau 
\\ 0 &\text{otherwise}\end{cases}\\
 \Ext^i(\QQQ_{\sigma*}, \QQQ_{\tau!})&= \begin{cases} \QQ &\text{ if } i=
 d(\tau) - d(\sigma \cap \tau) \text{ and } \sigma \cap \tau \ne \phi\\
 0 &\text{ otherwise}\end{cases}
\end{align*}

\subsection{Simplicial cosheaves}\index{cosheaf!simplicial}\label{subsec-simplicial-cosheaf}
  There are two ways to view a simplicial cosheaf $R$ on a simplicial complex $K$.  The
first, as a pre-cosheaf that assigns to each simplex $\sigma$ a vector space $R_{\sigma}$ and to each incidence
$\sigma < \tau$ a homomorphism $A_{\tau} \to A_{\sigma}$,in a way that is compatible with compositions.  The second way
involves the {\em dual cells}\index{dual!cell} of $K$.  

Let $K'$ denote the first barycentric subdivision of $K$.  For each simplex $\sigma$ in $K$ the dual cell
$D(\sigma)$ \index{dual cell} \index{cell, dual}
is the union of those simplices $x$ in $K'$ spanned by $\hat\sigma$ and the barycenters $\hat\tau$ of
simplices $\tau > \sigma$.  (Although it is contractible to $\hat\sigma$ the dual cell $D(\sigma)$ is not necessarily
homeomorphic to a disk.)  If $\sigma < \tau$ then $D(\tau)\subset D(\sigma)$ and we write $D(\tau)<D(\sigma)$.

\begin{center}\begin{tikzpicture}[scale = 1.6]
\node at (0,0) (a) {};
\node at (0,2) (b) {};
\node at (2,0) (c) {};
\node at (3,2) (d) {};
\node at (3,-2) (e) {};

\draw [ shade] (0,0) -- (2,0) -- (3,-2) -- (0,0);  

\draw[ fill=gray!80] (0,0) -- (2,0) -- (3,2) -- (0,0);  

\draw[ shade]  (0,0) -- (0,2) -- (6/5, 4/5) -- (0,0);  

\draw[ dotted] (6/5, 4/5) -- (2,0);

\draw[fill=blue] (1,0) circle (1pt); 

 \draw[fill=red] (1/2,1) circle (1 pt);  

\draw[fill=red](1.75,.75) circle (1pt);  

  \draw[fill=red] (1.75,-.75) circle (1pt);  

\draw[blue, thick] (1,0) -- (1.75,.75);
\draw[blue, thick](1,0) -- (1.75, -.75);

\draw[blue, thick] (1/2,1) -- (3/4,1/2);
\draw[blue, dotted] (3/4,1/2) -- (1,0);

\node at (1/4, 1.3){$\tau_1$};
\node at (2.4,1.3){$\tau_2$};
\node at (2.4, -1.3){$\tau_3$};
\node at (1.7, -.15){$\sigma$};
\end{tikzpicture}\end{center}
\vskip-.5cm\centerline{Dual of $\sigma$ in blue}

A {\em sheaf} $B$ on the dual cell decomposition assigns a vector space $B_{\xi}$ to each dual cell
$\xi = D(\sigma)$ and a homomorphism $B_{\xi} \to B_{\eta}$ whenever $\xi < \eta$ in a way that is compatible 
with compositions.  So a sheaf on the dual cell decomposition is the same thing as a cosheaf on the
cell decomposition. 

\subsection{Perversities}  \index{perversity}
\quash{
In \cite{IH1} a perversity is a function $p(c)$ defined on strata,
such that $p(c) \le p(c+1) \le p(c) +1$ and
where $c$ denotes the codimension of the stratum.  In \cite{BBDG} a perversity is a function defined on strata but
generally indexed by the dimension $d$ of that stratum so that $p(d) \ge p(d+1) \ge p(d) -1$.  (At various points in
\cite{BBDG} arbitrary integer valued functions $p$ are allowed, but the above restrictions are required for a
theory that is invariant under refinement of the stratification.)  We shall use the convention of \cite{BBDG} but often
write $p(S)$ rather than $p(d)$ when $S$ is a stratum of dimension $d$.   Since we will be considering the
cells of a simplicial complex to be strata we will also need to include the possibility (as in \cite{BBDG})
that $p(S) \ne 0$ when $S$ is a codimension one cell.  Finally, there is a question of normalization.  In \cite{IH1}
it is assumed that $p(S) = 0$ when the codimension of $S$ is $0$.  In \cite{BBDG} the authors use the
implicit normalization that $p(S) = 0$ when $\dim(S) = 0$.   Since the perversity corresponds to a cutoff at
a certain degree of cohomology, changing the normalization also involves changing the cohomological degrees.
In categories of complexes of sheaves, this is effected by a shift, e.g. $A\b[-1]$.
}
We follow the notation of  \cite{BBDG} by modfying the definition of a perversity $p$ to allow for codimension one strata
and considering $p$ to be a function of the {\em dimension}  of the stratum. 
  This results in a shift of cohomological degree that depends on $p$. 
 We refer to the notation of \S \ref{subsec-numbering}.

Although it is possible to develop the theory of simplicial perverse sheaves for perversities that vary with the simplex, in
these notes we assume for simplicity that the perversity $p(k)\le 0$ is a function only of the  dimension $k$ of the 
simplex\footnote{ Suppose  $K$ is a  triangulation of a stratified space $X$ and $A\b \in D^b(X)$ is constructible 
  with respect to the stratification $\mathcal S$.  Then it is also constructible with respect to the triangulation $K$
  and (by \cite{IH2} \S 4.1, or \cite{BBDG} \S 2.1.14) it is $p$-perverse with respect to the stratification $S$ 
  if and only if it is $p$-perverse with respect to the triangulation.}
 with $p(0) = 0$ and $p(k-1) \ge p(k) \ge p(k-1)$ for all $k$.

Let $C$ be a $k$-dimensional cell of $K$, with closure $\sigma$.
 Following \cite{Polishchuk} we say that $C$ (or its closure $\sigma$ or its 
dimension $k$) is of type  ! if $p(k) = p(k-1)$ and is of type * if $p(k) = p(k-1)-1$.  The number $k=0$ is considered
to be both types. 
  The simple objects in the
category of perverse sheaves are the complexes $\uul{IC}{}^{\bullet}_{\sigma}$. From \S \ref{subsec-singlecell},
and observed  in \cite{Polishchuk}, for any simplex $\sigma$:
\begin{tcolorbox}[colback=cyan!30!white]\index{IC (intersection complex)}\index{intersection complex}
The intersection complex extending the constant sheaf $\uul{\QQ}{}_{\sigma^o}$ on the interior $\sigma^o$ is: 
\begin{equation}\label{eqn-ICtau}
\uul{IC}{}_{\sigma} = \begin{cases}
\uul{\QQ}{}_{\sigma*}[-p(\sigma)] &\text{if }\ \sigma {\text  { is type *} } \\
\uul{\QQ}{}_{\sigma!}[-p(\sigma)] &\text{if }\ \sigma {\text { is type !} }  \end{cases}\end{equation}
\end{tcolorbox}

\subsection{Perverse dimension}  \index{perverse!dimension}  \label{subsec-perverse-dimension}

MacPherson's insight (\cite{MP1, MP2}) was to define the {\em perverse dimension}\index{perverse!dimension}
\index{dimension, perverse} of a $d$-simplex $\sigma$ (or of $C=\sigma^o$)  to be
\[ \delta(\sigma) = \begin{cases}-p(d) \ge 0&\text{if $\sigma$ has type * } \\
-p(d)-d\le 0 &\text{if $\sigma$ has type !}
\end{cases}\]
Therefore $|\delta(\sigma)|$ is the number of  integers $k$ with $1 \le k \le d$ of the same type as $d$.
Then $\delta:\ZZ_{\ge 0} \to \ZZ$ is one to one, and for each $d$
the image $\delta([0,d])$ is a subinterval of integers in $[-d,d]$  containing $0$.

\quash{
Using \S \ref{subsec-singlecell} we find: 
\begin{tcolorbox}[colback=cyan!30!white]
The intersection complex extending the constant sheaf on the interior $\sigma^o$ is:
\[ \uul{IC}{}_{\sigma} = \begin{cases}
\uul{\QQ}{}_{\sigma*}[-p(\sigma)] & \text{ if  $\sigma$ has type }\ *\\
\uul{\QQ}{}_{\sigma!}[-p(\sigma)] & \text{ if $\sigma$  has type }\  !
\end{cases}\]
\end{tcolorbox}
}


\begin{tikzpicture}[
dot/.style = {circle, fill, minimum size=#1,
              inner sep=0pt, outer sep=0pt},
dot/.default = 6pt  
                    ] 
\foreach \i in {7,...,10}
  {\draw [dotted](0,\i) -- (6,\i);}; 
\foreach \j in {0,...,6}
  {\draw[line width = 0.01mm, color = blue, dotted] (\j,10) -- (\j,7);};

\node at (0,10) [dot, color=red, radius=.1pt]{};
\node at (1,10)[dot, color=red, radius = .1pt]{};
\node at (2,10)[dot,color=red,radius=.1pt]{};
\node at(3,9)[dot,color=red,radius=.1pt]{};
\node at(4,8)[dot,color=red,radius=.1pt]{};
\node at (5,7)[dot,color=red,radius=.1pt]{};
\node at (6,7)[dot,color=red,radius=.1pt]{};
\draw[color=red, line width = .4mm](0,10) -- (1,10) -- (2,10) -- (3,9) -- (4,8) -- (5,7) -- (6,7);

\foreach \i in {0,...,6}
 {\node at (\i,6.5){\i};  };

\foreach \i in {-3, ..., 0}
  {\node at (-.5,10 + \i) {\i}; };

\def\ls{9} \def\bot{10}  

\draw[fill=gray!30!white,  line width = 0mm] (\ls,\bot) -- (\ls +3, \bot +3) -- (\ls + 6, \bot+3) --
(\ls + 6, \bot -3) -- (\ls+3, \bot -3) -- (\ls,\bot);

\draw [dotted, blue](\ls,\bot) -- (\ls+3,\bot+3); \draw[dotted, blue] (\ls,\bot) -- (\ls+3,\bot-3); 
\quash{
\foreach \i in {0,...,3}
  {\draw[color=gray, dotted] (\ls +\i, \bot + \i) -- (\ls +\i, \bot - \i); 
    \draw[color=gray, dotted] (\ls +\i, \bot + \i) -- (\ls + 6, \bot + \i);
    \draw[color=gray, dotted] (\ls + \i, \bot - \i) -- (\ls + 6, \bot -\i);
  };
\foreach \i in {4, ..., 6}
   \draw[color=gray, dotted] (\ls +\i, \bot + 3) -- (\ls + \i, \bot -3);  
}

\foreach \i in {-3, ..., 3}
  { \draw[dotted] (\ls, \bot + \i) -- (\ls + 6, \bot + \i);
  \node at (\ls-0.5, \bot + \i){\i};  };

  \foreach \i in  {0, ..., 6}
    {\draw[dotted] (\ls + \i, \bot + 3) -- (\ls  + \i, \bot -3);};

   \foreach \i in {0, ..., 6}
    {\node at (\ls + \i, \bot -3.5){\i};};

\draw[color=red, line width = .4mm] (\ls,\bot) -- (\ls+2, \bot-2) --(\ls +3,\bot+1) -- (\ls +4, \bot + 2) -- (\ls+5, \bot+3) -- (\ls+ 6, \bot -3);
\node at (\ls,\bot)[dot,color=red,radius=.1pt]{};
\node at (\ls +1,\bot -1)[dot,color=red,radius=.1pt]{};
\node at (\ls+2, \bot-2)[dot,color=red,radius=.1pt]{};
\node at (\ls+3, \bot+1)[dot,color=red,radius=.1pt]{};
\node at (\ls+4, \bot +2)[dot,color=red,radius=.1pt]{};
\node at(\ls+5,\bot +3)[dot,color=red,radius=.1pt]{};
\node at (\ls+6, \bot-3)[dot,color=red,radius=.1pt]{};

\node at (3,6){graph of $p(d)$};
\node at (\ls +3, \bot -4){graph of $\delta(d)$};

\end{tikzpicture}

\subsection{Partial order}\label{subsec-relations}
There is a unique partial order on the set of simplices of $K$ that is generated
by the following {\em elementary relation}: 
\index{elementary relation}\index{relation, elementary}
 $\sigma \succeq \tau$ if $\sigma \leftrightarrow \tau$ and $\delta(\sigma) = \delta(\tau)+1$,
in other words, two simplices $\sigma, \tau$ satisfy $\sigma \succeq \tau$ if there is a sequence of elementary relations
$\sigma = \sigma_0 \succeq \sigma_1 \cdots \succeq \sigma_r = \tau$.  It is easy to check the following.

\begin{tcolorbox}[colback=cyan!30!white]
Suppose $\sigma \succeq \tau$.  If both $\sigma, \tau$ are type * then $\tau \subset \sigma$.  If both are type
! then $\sigma \subset\tau$.  If they are opposite types then $\sigma$ is type * (so $\delta(\sigma)>0$)
and $\tau$ has type ! (so $\delta(\tau)<0$) and $\sigma \cap \tau \ne \phi$.
\end{tcolorbox}
Hence an arrow $\sigma \succeq \tau$ can extend no farther than adjacent simplices.  This may happen if
$\delta(\sigma)$ and $\delta(\tau)$ have opposite signs, for example, $\dim(\sigma) = 1$,
$\dim(\tau) = 2$ and $\delta(0,1,2) = 0, 1, -1$ respectively.
\begin{center}
\begin{tikzpicture}[
dot/.style = {circle, fill, minimum size=#1,
              inner sep=0pt, outer sep=0pt},
dot/.default = 6pt  
                    ] 
\def\ls{0} \def\bot{0}

\foreach \i in {-1, ..., 1}
  { \draw[dotted] (\ls, \bot + \i) -- (\ls + 2, \bot + \i);
  \node at (\ls-0.5, \bot + \i){\i};  };

  \foreach \i in  {0, ..., 2}
    {\draw[dotted] (\ls + \i, \bot + 1) -- (\ls  + \i, \bot -1);};

   \foreach \i in {0, ..., 2}
    {\node at (\ls + \i, \bot -1.5){\i};};

\node at (\ls, \bot)[dot, color=red, radius = .1pt]{};
\node at (\ls +1, \bot +1)[dot, color = red, radius = .1pt]{};
\node at (\ls +2, \bot -1)[dot, color=red, radius = .1pt]{};

\draw[color=red, line width = .4mm](\ls,\bot) -- (\ls +1, \bot + 1) -- (\ls + 2, \bot -1);

\node at (\ls +1, \bot -2){graph of $\delta(d)$};

\def\ls{5}

\draw[line width = .4mm] (\ls,\bot-1) -- (\ls +4,\bot-1);
\node at (\ls,\bot-1)[dot, radius = .1pt]{};
\node at (\ls+2, \bot-1)[dot]{};
\node at (\ls+4, \bot-1)[dot, radius = .01pt]{}; 
\node at (\ls +3,\bot+1)[dot, radius = .1pt]{};
\shadedraw[line width = .4mm] (\ls +2,\bot-1) -- (\ls +3, \bot +1) -- (\ls + 4, \bot-1);

\node at (\ls + 1, \bot -1.3) {$\sigma$};
\node at (\ls +3, \bot-.5){$\tau$};

\node at (\ls +2, \bot-2){$\sigma \succeq \tau$};
\end{tikzpicture}
\end{center}


\begin{tcolorbox}[colback=cyan!30!white]
Using \S \ref{subsec-singlecell} we find:
if $\sigma, \tau$ are simplices then 
\begin{equation}\label{eqn-ext}
\Ext^i(\uul{IC}{}^{\bullet}_{\sigma}, \uul{IC}{}^{\bullet}_{\tau}) =\begin{cases}
\QQ & \text{ if } \sigma \succeq\tau \text{ and}\ i = \delta(\sigma) - \delta(\tau) \\
0 &\text{ otherwise}. \end{cases}\end{equation}
\end{tcolorbox}

\subsection{Perverse cells}\index{perverse!cell}\index{cell, perverse}
For each simplex $\sigma$ the corresponding {\em perverse cell} ${}^{\delta}\sigma$ 
is the following union of simplices in the barycentric subdivision $K'$:
\[ {}^{\delta}\sigma = \bigcup \large\left\{ \langle \hat{\tau}_0 \hat{\tau}_1 \cdots \hat{\tau}_r\rangle |\ 
\forall  j,\  \delta(\tau_j) \le \delta(\sigma) \text{ and }\exists i \text{ with}\ \tau_i = \sigma
  \large\right\} \]
In this equation $\tau_0<\tau_1<\cdots<\tau_r$ form a partial flag of simplices hence $\delta(\tau_j)$ takes distinct values.
The ``boundary''  $\partial ({}^{\delta}\sigma)$ (also called the perverse link of $\sigma$)\index{perverse!link}
\index{link!perverse} is the union of those simplices $\langle \hat\tau_0 \hat \tau_1 \cdots \hat\tau_r
\rangle$ in ${ {}^{\delta}\sigma}$ such that
$\delta(\tau_j) < \delta (\sigma)$ for all $j$ ($0 \le j \le r$).  (That is, remove the vertex $\hat\sigma$ from
each simplex.) If $\sigma$ is type ! then ${}^{\delta}\sigma \subset D(\sigma)$.  If $\sigma$
is type * then $\sigma \subset \pd\sigma$ and $\partial\sigma \subset \partial \pd\sigma$.
The {\em interior} is $\pd\sigma^o = \pd\sigma - \partial\pd\sigma$.  If   $C = \tau^o$ is the
interior of $\tau$ we sometimes denote $\pd\sigma^o$ by $\pd C$.
The key technical point in the theory is the following, which implies the functor 
$H^r_{\pd C}$ from the category of perverse sheaves to vector spaces is exact.

\begin{tcolorbox}[colback=cyan!30!white]
\begin{prop}\label{prop-perversecell}  
 \cite{MP1,MP2}\cite{Polishchuk}  Let $A$ be a perverse sheaf on the simplicial complex $K$ and let $\sigma$
be a simplex of $K$ with $C = \pd\sigma^o$. The cohomology with support
 in ${}^{\delta}C$ vanishes:
\begin{equation}\label{eqn-perversecell}
 H^r_{ {}^{\delta}C}(K,A) = 0
\ \text{ unless }\ r = -\delta(\sigma)
 \end{equation}
\end{prop}
\end{tcolorbox}

\paragraph{\bf Proof.}\label{subsec-proof-of-vanishing}
It suffices to prove  equation (\ref{eqn-perversecell}) when $A = \uul{IC}{}_{\tau}$ since these are the simple
objects in the Noetherian category of perverse sheaves.   If $\sigma \leftrightarrow\tau$  
and $i: \pd\sigma^o\cap\tau \to K$, we will show that 
\begin{equation}\label{eqn-cellular-IH}
H_{\{\pd\sigma^o\}}^r(\uul{IC}{}_{\tau}) =H^r_{\{ \pd \sigma^o \cap \tau  \}}(\uul{IC}{}_{\tau}) =
H^r( i^!\uul{IC}{}_{\tau}) =0\ \text{ unless  }\ \sigma = \tau \text{ and } r = -\delta(\sigma)\end{equation}
in which case,  $H^{-\delta(\sigma)}_{\pd\sigma^o}(\uul{IC}{}_{\sigma})=\QQ.$

  If $\sigma = \tau$ is type * then $\pd\sigma\cap\tau = \sigma$ and $\uul{IC}{}_{\sigma}
= \uul{\QQ}{}_{\sigma}[\delta(\sigma)]$ so
 $H^*_{ \{\sigma\} }(\uul{IC}{}_{\sigma}) = \QQ[-p(\sigma)]=\QQ[\delta(\sigma)].$
 
If $\sigma = \tau$ is type ! then $\pd\sigma \cap \tau = \hat\sigma$ and by \S \ref{subsec-singlecell},
 \[ H^r_{ \{\hat\sigma\} }(\uul{IC}{}_{\sigma}) = H^r_{ \{\hat\sigma\} }(\QQ_{\sigma!}[-p(\sigma)]) =
 \begin{cases} \QQ & \text{ if } r = \dim(\sigma)+p(\sigma) \\
 0 & \text{ otherwise} \end{cases}\]

Next, suppose $\sigma \ne \tau$ and $\tau$ is type *.
 If $\sigma < \tau$ then $\delta(\sigma) < \delta(\tau)$ so $\pd \sigma^o \subset
 \pd\sigma \subset \partial\tau$ since $\hat\tau$.  Let
 $j:\tau^o \to K$ denote the inlusion.  Then
 $ H^*(i^!Rj_*\uul{\QQ}{}_{\tau^o})$ is dual to
 $ H^*(i^*Rj_!\uul{\QQ}{}_{\tau^o}) = 0.$  If $\sigma > \tau$ then $\pd\sigma^o \cap \tau = \phi$ so 
 $H^*(i^!\uul{IC}{}_{\tau}) = 0$.
 

 Now suppose that $\tau$ is type !. Then $H^*(i^!\uul{IC}{}_{\tau}) = H^*(i^!\uul{\QQ}{}_{\tau!})$ is dual to
 \[ H^{*}_c(i^*\uul{\QQ}{}_{\tau*})=H^*_c(\pd\sigma^o \cap \tau) = H^{*}(\pd\sigma \cap \tau, \partial \pd\sigma \cap \tau).\]
 If $\sigma$ is type *, or if $\sigma$ is type ! and $\sigma < \tau$ then $\delta(\sigma)  > \delta (\tau)$ so every
 maximal simplex in  $\pd \sigma \cap \tau$ and every maximal simplex in $\partial \pd\sigma \cap \tau$ 
 contains the vertex $\hat\tau$ hence both are contractible and the relative cohomology is trivial.
In the remaining case, $\tau$ is type !, $\sigma$ is type ! and $\sigma > \tau$ so $\pd\sigma\cap\tau = \phi$.  \qed

\subsection{Perverse skeleton}  \index{perverse!skeleton}\index{skeleton, perverse}
MacPherson defines the $r$-th perverse skeleton:
\[ \pd K_r = \bigcup \left\{ \pd\sigma |\ \delta(\sigma) \le r \right\} .\]
It is a closed subcomplex\footnote{after appropriate shifts and translation of perversity conventions it coincides with
the ``basic set'' $Q_r^{\bar p}$ of \cite{IH1}.} 
of the barycentric subdivision of $K$ and it is the span 
of the barycenters $\hat\sigma$ of simplices with $\delta(\sigma) \le r$.  Its interior is the disjoint union,
 \[\pd K_r^o = \pd K_r-\pd K_{r-1} = \coprod_{\delta(\sigma)=r} \pd \sigma^o.\] 
Let $U_r = K-\pd K_{r-1}$.  This gives a filtration of $K$ by open sets,
\begin{equation}\label{eqn-filtration} \cdots \subset U_{r+1} \subset U_r \subset U_{r-1} \subset \cdots\end{equation}
and $\pd K_r^o = U_r-U_{r+1}$ is closed in $U_r$ with inclusions
\[ \begin{CD}
\pd K_r^o @>>{i_r}> U_r @<<{j_r}< U_{r+1}
\end{CD}\]
The long exact sequence on cohomology (\S \ref{prop-pair-sequence}) for $A|U_r$ together with 
Proposition \ref{prop-perversecell} says that the relative cohomology lives in a single degree:
\begin{tcolorbox}[colback=cyan!30!white]
For any perverse sheaf $A$ we have:
\begin{equation}\label{eqn-Asigma} H^t(U_r, U_{r+1};A) = H^t(i_r^!(A|U_r))=
 \begin{cases}
\underset{\delta(\sigma) = r}{\bigoplus}H^r_{\pd \sigma^o}(K;A) &\ \text{ if } t = -r \\ 0 &\ \text{ otherwise}.  \end{cases}\end{equation}
\end{tcolorbox}
For each simplex $\sigma$ let $A_{\sigma} = H^{-r}_{\pd \sigma^o}(K;A)$ where $r = \delta(\sigma)$. 
From the exact sequences for the triple $(U_{r-1}, U_r, U_{r+1})$
we obtain a chain complex (meaning $d\circ d = 0$) whose cohomology is $H^*(K;A)$:
\begin{equation}\label{eqn-key}\begin{diagram}[size=2em]
 & \rTo^{d} &  H^{-r-1}(U_{r+1}, U_{r+2};A) & \rTo^{d} &H^{-r}(U_r, U_{r+1};A) & \rTo^{d} &  H^{-r+1}(U_{r-1}, U_{r};A) &\rTo^{d}  \\
   &                 &         ||                           &       &      ||    &   &  || &\\
   &                 &  \bigoplus_{\delta(\sigma) = r+1}A_{\sigma} &\rTo^{d}& \bigoplus_{\delta(\sigma)=r}A_{\sigma} 
   & \rTo^{d}& \bigoplus_{\delta(\sigma)= r-1}A_{\sigma} &
\end{diagram}\end{equation}
 Writing $d= \oplus_{\sigma,\tau}s_{\sigma\tau}$  leads to the following ``combinatorial'' definition.
\begin{tcolorbox}[colback=yellow!30!white]
\begin{defn}\index{perverse!sheaf!cellular}\index{cellular perverse sheaf}
A cellular perverse sheaf $S$ on the simplicial complex $K$ is a rule that assigns to each simplex $\sigma$
a $\QQ$-vector space $S_{\sigma}$ and ``attaching'' \index{attaching homomorphism}
homomorphism $s_{\sigma\tau}:S_{\sigma} \to S_{\tau}$
whenever $\sigma \leftrightarrow \tau$ and $\delta(\sigma) = \delta(\tau) +1$ such that 
$d\circ d = 0$ in the resulting sequence
\begin{equation}\label{eqn-dd0}\begin{diagram}[size=2em]
&\rTo^d& \bigoplus_{\delta(\sigma) = r+1}S_{\sigma} &\rTo^{d}& \bigoplus_{\delta(\sigma)=r}S_{\sigma} 
   & \rTo^{d}& \bigoplus_{\delta(\sigma)= r-1}S_{\sigma} &\rTo^d
\end{diagram}\end{equation}\end{defn}
\end{tcolorbox}
The cohomology of the cellular perverse sheaf $A$ is the cohomology of this sequence.
The chain complex condition is equivalent to the statement that whenever $\delta(\sigma) = r+1$ and $\delta(\tau)=r-1$ then 
the sum over $\theta$ vanishes:
\begin{equation}\label{eqn-relations}
 \underset{\begin{smallmatrix}{\delta(\theta)=r}\\{\sigma \leftrightarrow \theta \leftrightarrow\tau}\end{smallmatrix}}
{\sum}   s_{\theta \tau} \circ s_{\sigma\theta} = 0\end{equation}
A morphism $A \to B$ of cellular perverse sheaves is a homomorphism $A_{\sigma} \to B_{\sigma}$ for every simplex
$\sigma$, that commutes with the boundary maps $s_{\sigma \tau}$.  The category of cellular perverse sheaves
is denoted $\mathcal P^{\Delta} (K)$.  By (\ref{eqn-cellular-IH}) the complex $\uul{IC}{}_{\sigma}$
corresponds to the single non-zero assignment $S_{\sigma} = \QQ$.

\begin{tcolorbox}[colback=cyan!30!white]
\begin{prop} \cite{MP1, Polishchuk}\label{prop-equivalence}
The cohomology functor $T=\bigoplus_{\sigma} H^{-\delta(\sigma)}_{\pd\sigma^o}$ is exact and it
defines an equivalence of categories 
\[\begin{CD}\mathcal P(K) @>{T}>> \mathcal P^{\Delta}(K) \end{CD}\]
between the category of perverse sheaves that are constructible with respect to the triangulation
 and the category of cellular perverse sheaves (with the same perversity). 
\end{prop}\end{tcolorbox}
\paragraph{\bf Proof:}
The proof in \cite{Polishchuk} uses Koszul duality (see \S \ref{subsec-quadratic-dual}).  A direct proof is tedious but
the main point is that by (\ref{eqn-Asigma}) the spectral sequence for the filtration (\ref{eqn-filtration}) 
collapses after the $E_1$ page with
$E_1^{t,r} = H^t(U_r, U_{r+1};A)$ leaving only the diagonal terms, $E_2^{-r,r}$ on the $E_2$ page.  These
 are precisely the cohomology  groups of the  sequence (\ref{eqn-key}) (in degree $-r$).  Therefore
the inclusion $U_r \to K$ induces an isomorphism
\[ H^{-r}(\oplus_{\sigma}S_{\sigma}^{\bullet}) \cong H^{-r}(K;A)\]
 of the cohomology in degree $-r$ of the above complex (\ref{eqn-dd0})
with the cohomology of the perverse sheaf $A$. 
 The same remark applies to subcomplexes of $K$.  However a perverse
sheaf on $K$ is determined by the cohomology of every subcomplex of $K$.  Similarly, a homomorphism
between perverse sheaves $A \to B$ is a global section of the sheaf $\uul{\RHom}{}^{\bullet}(A,B)$
and is determined by its restriction to every subcomplex. \qed
\medskip

\paragraph{\bf Remark.}  Rather than the ``chain complex'' relations (\ref{eqn-dd0}) one might consider objects $T$
which assign to each simplex $\sigma$ a vector space $T_{\sigma}$ and  attaching
homomorphisms $T_{\sigma} \to T_{\tau}$ whenever $\sigma \leftrightarrow \tau$ with $\delta(\sigma)
= \delta(\tau) +1$, that are compatible with compositions.  The resulting category turns out to be
the {\em Koszul dual} category, see \S \ref{prop-Amodel}.

\subsection{Verdier duality} It is easy to check that\index{Verdier duality}
 Verdier duality in the category of sheaves on $K$ constructible with respect to the triangulation exchanges
the category of perverse sheaves with perversity $p$ and associated function $\delta$ with the category
of perverse sheaves with perversity $q = t-p$ and associated function $-\delta$ (where $t(d) = -d$ is the
``top'' perversity). It reverses the arrows of the quiver $Q$ associated to $\delta$ (\S \ref{subsec-quiver}).

\subsection{Remarks}\label{subsec-remarks-on-cellular}
\paragraph{\bf 1}
If $p= 0$ is the zero perversity then every cell is type ! and $\pd\sigma =D(\sigma)$ is the dual of the simplex
$\sigma$.  The perverse link $\partial \pd\sigma$ is the ``usual'' link of $\sigma$.  A perverse sheaf $A$ is a simplicial
sheaf in the usual sense, except that signs must be added (which is always possible in this case)
to the  mappings $s_{\sigma\tau}$ so as to
obtain $s_{\theta\tau}s_{\sigma\theta} = s_{\sigma\tau}$.  The complex (\ref{eqn-dd0}) computes the cellular
cohomology of $A$. cf. \S \ref{subsec-simplicialsheaf}.  

If $p$ is the top perversity then $\pd\sigma = \sigma$
and a perverse sheaf $A$ is a simplicial cosheaf (\S \ref{subsec-simplicial-cosheaf})\index{cosheaf!simplicial}
\index{simplicial!cosheaf} in the usual sense, up to the same issue involving signs.  The complex (\ref{eqn-dd0})
computes its homology.

\medskip\noindent
\paragraph{\bf 2}
One might attempt to mimic Proposition \ref{prop-equivalence} using the filtration by closed subsets,
\[ \cdots \subset K_{r-1} \subset K_r \subset K_{r+1} \subset \cdots\]
 in order to obtain a chain complex with terms
\[ H^t(K_r, K_{r-1};A) \cong \bigoplus_{\delta(\sigma) = r} H^t(\pd\sigma, \partial\pd\sigma, A).\]
Unfortunately the relative cohomology $H^*(\pd \sigma, \partial \pd\sigma;A)$ may live in several degrees as the
perverse sheaf $A$ varies.  Let $K = \tau$
be the 3-simplex, let $p(3) = -1$ and $p(d) = 0$ otherwise, so that $\delta(3) = 1$ and $\delta(d) = -d$ for
$d = 0, 1, 2$.  Let $\sigma$ be a vertex of $\tau$.  Then 
$H^0(\pd\sigma, \partial \pd\sigma; \uul{IC}{}_{\sigma}) = \QQ$ but
$H^1(\pd\sigma,\partial\pd\sigma; \uul{IC}{}_{\tau}) = \QQ$.

\medskip\noindent
\paragraph{\bf 3}
In the example of \S \ref{subsec-relations} the sum (\ref{eqn-dd0}) contains a single term, so that $s_{\sigma\tau}=0$.
  In general, if $\sigma, \tau$ are simplices of $K$ such that
   $\delta(\sigma)$ and $\delta(\tau)$ have opposite signs, and $\sigma \succeq \tau$ but
   $\sigma \not\leftrightarrow \tau$ then any composition
\[ s(\sigma_r,\tau)\circ s(\sigma_{r-1},\sigma_r) \circ \cdots \circ s(\sigma_1,\sigma_2)\circ s(\sigma,\sigma_1)...=0\]
of attaching homomorphisms from $A_{\sigma}$ to $A_{\tau}$
will vanish because such a composition eventually passes through a triple $\theta_1 \succeq \theta_0 \succeq\theta_{-1}$
with $\delta= +1,0,-1$ respectively and $\theta_1 \not\leftrightarrow \theta_{-1}$.
Then $\theta_0$ is a vertex and one of the other two is a 1-simplex which meets the remaining simplex
in a single point, $\{\theta_0\}$. The chain complex
requirement for this triple has only one term: 
$s(\theta_0,\theta_{-1})s(\theta_1,\theta_0) = 0$.

\medskip\noindent
\paragraph{\bf 4} In the more general situation when $K$ is a regular cell complex the preceding argument fails.
In this case the definition of a cellular perverse sheaf requires the additional {\em locality axiom} (\cite{MP1}):
If $\sigma \succeq \tau$ and $\sigma \not\leftrightarrow\tau$ then any composition of morphisms from
$A_{\sigma}$ to $A_{\tau}$ vanishes.

\quash{
\subsection{Proof of Proposition \ref{prop-perversecell}}\label{subsec-proof-of-vanishing}
It suffices to prove  equation (\ref{eqn-perversecell}) when $A = \uul{IC}{}_{\tau}$ since these are the simple
objects in the Noetherian category of perverse sheaves.  We may assume that $\sigma\leftrightarrow\tau$.
Let $i: \pd\sigma^o\cap\tau \to K$  Following \cite{Polishchuk} we show that 
\[H_{\{\pd\sigma^o\}}^*(\uul{IC}{}_{\tau}) =H^*_{\{ \pd \sigma^o \cap \tau  \}}(\uul{IC}{}_{\tau}) =
H^*( i^!\uul{IC}{}_{\tau}) =0\ \text{ unless  }\ \sigma = \tau\ \text{ and }\]
 \[H^r_{\pd\sigma^o}(\uul{IC}{}_{\sigma}) = \begin{cases} \QQ & \text{ if } r = -\delta(\sigma)\\
0 &\text{ otherwise}\end{cases}\] 
  If $\hat\tau$ does not
occur in ${}^{\delta}\sigma$ (for example, if $\sigma < \tau$) then ${}^{\delta}\sigma \cap \tau \subset \partial \tau$.  
If $j:\tau^o \to \tau$ is
the inclusion then $H^*(i^!Rj_*(\uul{\QQ{}_{\tau^o}}))$ is (Verdier) dual to $H^*(i^*Rj_!(\uul{\QQ}{}_{\tau^o})) = 0$
which verifies the first equation in (\ref{eqn-perversecell}).  
}
\quash{

\subsection{Partial order}
  If $\sigma, \tau$ are simplices of $K$ write $\sigma \leftrightarrow \tau$ in
case {\em either} $\sigma \le \tau$ or $\tau \le \sigma$.  \index{perverse!partial order}\index{partial order, perverse}
There is a unique partial order on the set of simplices of $K$ that is generated
by the following {\em elementary relation}:  $\sigma \succeq \tau$ if $\sigma \leftrightarrow \tau$ and $\delta(\sigma) = \delta(\tau)+1$,
in other words, two simplices $\sigma, \tau$ satisfy $\sigma \succeq \tau$ if there is a sequence of elementary relations
$\sigma = \sigma_0 \succeq \sigma_1 \cdots \succeq \sigma_r = \tau$.  It is easy to check the following.

\begin{tcolorbox}[colback=cyan!30!white]
Suppose $\sigma \succeq \tau$.  If both $\sigma, \tau$ are type * then $\tau \subset \sigma$.  If both are type
! then $\sigma \subset\tau$.  If they are opposite types then $\sigma$ is type * (so $\delta(\sigma)>0$)
and $\tau$ has type ! (so $\delta(\tau)<0$) and $\sigma \cap \tau \ne \phi$.
\end{tcolorbox}
Hence an arrow $\sigma \succeq \tau$ can extend no farther than adjacent simplices.  The simple objects in the
category of perverse sheaves are the complexes $\uul{IC}{}^{\bullet}_{\sigma}$.
Using \S \ref{subsec-singlecell} we find:
\begin{tcolorbox}[colback=cyan!30!white]
If $\sigma, \tau$ are simplices then $\RHom(\uul{IC}{}^{\bullet}_{\sigma}, \uul{IC}{}^{\bullet}_{\tau}) = 0$
unless $\sigma \succeq\tau$.
\end{tcolorbox}

}

\begin{tcolorbox}[colback=cyan!30! white]
\begin{prop} \cite{Polishchuk}\label{prop-Polish-equivalence}
The natural functor $\mathcal P(K) \to D^b_c(K)$ induces an equivalence of derived categories $D^b(\mathcal P(K))
\overset{R}{\longrightarrow} D^b_c(K)$. 
\end{prop}
\end{tcolorbox}
The proof in \cite{Polishchuk} uses Koszul duality (cf.~\S \ref{subsec-quadratic-dual}, \S\ref{subsec-Koszulderived})
 together with the results of \cite{BBDG}  (p. 84-85).   
 
 Every complex $A\b$ in $D^b_c(K)$ is quasi-isomorphic to a complex of perverse sheaves.  For, if $\sigma$ is
 a simplex of $K$ of type * then the elementary injective sheaf $\uul{\QQ}{}_{\sigma*}$ is isomorphic to
 the perverse sheaf $\uul{IC}{}_{\sigma}$.  If it has type ! then there is a triangle
 \begin{diagram}[size=2em]
 \uul{IC}{}_{\sigma}= \uul{\QQ}{}_{\sigma!} & & \rTo & & \uul{\QQ}{}_{\sigma*}  \\
  & \luTo_{d} && \ldTo &  \\
  && B\b &&
  \end{diagram}
    The complex $B\b$ is supported on the boundary, $\partial \sigma$, so by induction it is quasi-isomorphic to a complex
    $\to \cdots \to C^{r-1} \to C^r$ of perverse sheaves.  Then $\cdots C^{r-1} \to C^r \overset{d}{\to} \uul{\QQ}{}_{\sigma!}$ is
    a complex of perverse sheaves quasi-isomorphic to the elementary injective sheaf $\uul{\QQ}{}_{\sigma*}$ so the
    functor $R$ above is essentially surjective.  The content of \cite{Polishchuk} is showing
 that $R$ induces an isomorphism on morphisms,
 that is, $\Ext^i_{\mathcal P(K)}(A,B) \cong \Ext^i(A,B)$.  (The Yoneda $\Ext$ does not require enough injectives
in the category $\mathcal P(K)$ for its definition while the $\Ext^i(A,B)$ is defined using injective resolutions in the 
full derived $D^b_c(K)$.)

\newcommand{\CPS}{\mathcal C}
\newcommand{\PA}{{\rm F}}
\section{The path algebra\label{sec-Koszul}}\index{path algebra}
\subsection{}  Aside from standard definitions the material in this section is due to Vybornov 
\cite{Vybornov, Vybornovthesis2, Vybornov2, Vybornov3} and Polishchuk \cite{Polishchuk}.
Throughout this section we retain the notation of \S \ref{sec-cellular-perverse-sheaves}:  $K$ is a
finite simplicial complex, $p$ is a perversity with associated function $\delta$, 
For simplices $\sigma,\tau$ in $K$ we write $\sigma \leftrightarrow \tau$ if either $\sigma < \tau$ or $\tau < \sigma$.
We denote by $\CPS$ the category of cellular perverse sheave of $\QQ$-vector spaces on $K$.

\subsection{}\label{subsec-quiver}
Associated to $K, \delta$ there is a quiver $Q$, or directed graph without oriented cycles, 
whose vertices are the simplices of $K$ and whose arrows correspond
to elementary relations $\sigma \succ \tau$ with $\delta(\sigma) = \delta(\tau)+1$.   There is a
graded {\em path algebra} $\PA=\bigoplus_{j\ge 0}\PA_j$, which is generated by elements of degree $\le 1$.
 The vector space $\PA_0$ has a basis element
$[\sigma]$ for each simplex $\sigma$ in $K$ (the ``trivial paths'), and $\PA_1$ has a (canonical) 
basis consisting of the elementary relations, that is, 
$[s_{\sigma\tau}]$ where $\sigma \leftrightarrow \tau$ and $\delta(\sigma) = \delta(\tau)+1$, subject to the following
relations \begin{enumerate}
\item $[\sigma]^2 = [\sigma]$ 
\item  $[\sigma].[\tau] = 0$ unless $\tau = \sigma$.
\item $ [\tau].[s_{\sigma\tau}] = [s_{\sigma\tau}].[\sigma] = [s_{\sigma\tau}] $
\end{enumerate}
These relations already imply that $[s_{\sigma\tau}]. [s_{\alpha\beta}] = 0$ unless $\sigma = \beta$.

If $\sigma \succeq \tau$ and $\delta(\sigma) = \delta(\tau)+r $ (with $r \ge 1$) then there is a sequence 
of elementary relations  $\sigma =\sigma_0\succeq \sigma_1\succeq \cdots \succeq \sigma_r = \tau$ 
connecting them.  It follows that the vector space $\PA_r$ has a (canonical) basis whose elements are 
 paths of length $r$ with decreasing $\delta$.  Products in $\PA$ corresponds to concatenation of paths. Consequently
 a representation of the quiver $Q$ (that is, a vector space $A_{\sigma}$ for each simplex and a homomorphism
 $A_{\sigma} \to A_{\tau}$ for each elementary relation $\sigma \succ\tau$) is the same thing as a module over
 the path algebra $\PA$.  Given a module $M$ the associated vector spaces are given by $M_{\sigma} = [\sigma].M$.

 Not every representation of $Q$ corresponds to a perverse sheaf:  the relations (\ref{eqn-relations}) are required to hold.
From these relations we immediately conclude
\begin{tcolorbox}[colback=cyan!30!white]
\begin{prop}\label{prop-representation}
The category $\mathcal P(K)$ of perverse sheaves, constructible with respect to the triangulation of $K$ is naturally equivalent to
the category $Mod-B$ of modules over the ring $B=\PA/J$ where $J\subset \PA$ is the two-sided homogeneous ideal generated by
the vector subspace $E \subset \PA_2$ spanned by the elements  
$u_{\sigma\tau}$ whenever  $\delta(\sigma) = r+1$
and $\delta(\tau) = r-1$ defined by
\begin{equation}\label{eqn-relations1}
u_{\sigma\tau}= \underset{\begin{smallmatrix}{\delta(\theta)=r}\\{\sigma \leftrightarrow \theta \leftrightarrow\tau}\end{smallmatrix}}
{\sum}   [s_{\theta \tau}] . [s_{\sigma\theta}]\end{equation}
 
\end{prop}\end{tcolorbox}

\subsection{ The $\Ext$ algebra} \index{Ext!algebra}
 The ring $B$ is a quadratic algebra.  The Koszul dual of such an algebra is
sometimes defined to be the $\Ext$ algebra of the simple $B$-modules.
Following \cite{Polishchuk} we consider the graded algebra
\[A= \bigoplus_{\sigma,\tau} \Ext^*(\uul{IC}{}_{\sigma}, \uul{IC}{}_{\tau}).\]
(As in \S \ref{subsec-singlecell},  $\Ext^i$ is computed in the full derived category $D^b_c(K)$.)
By (\ref{eqn-ext}) we know that $\Ext^j(\uul{IC}{}_{\sigma}, \uul{IC}_{\tau}) = 0$ unless $\sigma \succeq \tau$ and
$\delta(\sigma) = \delta(\tau)+j$, in which case this group is canonically isomorphic to $\QQ$ and it has a
canonical generator $t_{\sigma\tau} \in  \Ext^{\delta(\sigma)-\delta(\tau)}(\uul{IC}{}_{\sigma},
\uul{IC}{}_{\tau})$ that corresponds to the identity mapping $\QQ \to \QQ$ between the stalk cohomology of
$\uul{IC}{}_{\sigma}$ and $\uul{IC}{}_{\tau}$ at points in the interior.

Similarly by considering the composition of two morphisms it is easy to see that   
the product $t_{\alpha\beta}t_{\sigma\tau}$ is zero unless $\alpha = \tau$ in which case
it is $t_{\sigma\beta}$  (\cite{Polishchuk}).  Equivalently,
  if $\delta(\sigma) = r+1$,
$\delta(\tau) = r-1$ and if $\sigma \leftrightarrow \theta \leftrightarrow \tau$ and
$\sigma \leftrightarrow \theta' \leftrightarrow \tau$ with $\delta(\theta) = \delta(\theta') = r$ then
$ t_{\theta\tau}.t_{\sigma\theta} = t_{\theta'\tau}.t_{\sigma\theta'}$.  (Compare \S \ref{subsec-remarks-on-cellular} (3).)
 We may conclude:

\begin{tcolorbox}[colback=cyan!30!white]
\begin{prop}\label{prop-Amodel}
A module $T$ over the Ext algebra $A$ is a vector space 
$T_{\sigma}$ attached to each simplex $\sigma$, and  ``attaching'' morphisms $t_{\sigma\tau}:
T_{\sigma} \to T_{\tau}$ whenever $\sigma \succeq \tau$  which are compatible with composition. 
Equivalently, the algebra $A$ is isomorphic to the quotient $\PA/I$ where $I$ is the two-sided homogeneous ideal
generated by the vector space $D \subset \PA_2$ spanned by the following elements:
\begin{equation}\label{eqn-trelations}
 [s_{\theta\tau}].[s_{\sigma\theta}] - [s_{\theta'\tau}].[s_{\sigma\theta'}]\end{equation}
for every $\sigma, \tau$ with $\delta(\sigma) =\delta(\tau)+2$ and 
every $\theta, \theta'$ with $\delta(\theta) = \delta(\theta') = \delta(\tau)+1$ and $\sigma \leftrightarrow \theta \leftrightarrow \tau$
and $\sigma \leftrightarrow \theta' \leftrightarrow \tau$. 

 \qed
\end{prop}\end{tcolorbox}

\subsection{} \label{subsec-quadratic-dual}
  The canonical basis $\{[s_{\sigma\tau}]\}$ of $\PA_1$ determines an inner product
$\langle [s_{\sigma\tau}], [s_{\alpha\beta}]\rangle = \delta_{\sigma\alpha}\delta_{\tau\beta}$ on $\PA_1$
and hence also on each $\PA_r$.
Any vector subsapce $X \subset \PA_2$ generates a homogeneous ideal $(X) \subset \PA$ with quotient
algebra $Y = \PA/(X)$. The {\em quadratic dual} \index{quadratic dual}\index{dual!quadratic}
algebra is $Y^! = \PA/V$ where $V = (X^{\perp})$ is the
homogeneous ideal generated by the complementary subspace $X^{\perp} \subset \PA_2$.

\begin{tcolorbox}[colback=cyan!30!white]
\begin{prop} \cite{Vybornovthesis} The ring $B$ and the $Ext$ algebra $A$ are quadratic duals. 
\end{prop}\end{tcolorbox}

\paragraph{\bf Proof.}  In fact the subspaces $D, E\subset \PA_2$ above are orthogonal complements.   
The vector space $F_2$ in the path algebra is the orthogonal direct sum of the subspaces $V_{\sigma\tau}$ 
spanned by paths of length two, that is,
\begin{equation}\label{eqn-sigmatau}
 V_{\sigma\tau} =\text{span}  \left\{ [s_{x\tau}].[s_{\sigma x}] |\ \sigma \succeq x \succeq \tau, \delta(\sigma) = \delta(x)+1 =
\delta(\tau)+2 \right\}.\end{equation}
We claim that within each $V_{\sigma\tau}$ the corresponding subspaces $D_{\sigma\tau},
E_{\sigma\tau}\subset\PA_2$ are 
orthogonal complements. Fix $\sigma,\tau$ as in (\ref{eqn-sigmatau}) and  first suppose they are both of type *.  
This implies that $\sigma \leftrightarrow \tau$ and so, given $\sigma, \tau$ there are two or more 
choices for $x$.  The resulting vectors  (\ref{eqn-trelations}) span the subspace $D_{\sigma\tau}$ that is orthogonal to the 
vector $u_{\sigma\tau}$ of (\ref{eqn-relations1}) which spans $E_{\sigma\tau}$.  The
same holds if $\sigma,\tau$ are both of type !.  If $\sigma,\tau$ are different types then $\delta(\sigma) = 1$,
$\delta(\tau) = -1$, $\delta(x) = 0$ and either $\sigma$ or $\tau$ is one-dimensional, cf.~the example in
\S \ref{subsec-relations}.   Then $\{x\} = \sigma \cap \tau$ is the unique point in the
intersection so the sum (\ref{eqn-relations1})  has a  single term: it is the monomial $[s_{x\tau}].[s_{\sigma x}]$.
Hence $V_{\sigma\tau} = E_{\sigma\tau}$ is one dimensional, while $D_{\sigma\tau} = 0$.  \qed.

\subsection{} \label{subsec-Koszulderived} \index{injective!object}\index{elementary injective}
 Injective objects in the category of modules over the $\Ext$ algebra $A$ are sums of {\em elementary
injectives} $I_{\sigma}$ where
\[ I_{\sigma}(\tau) = \begin{cases} \QQ & \text{ if } \tau \succeq \sigma\\
0 & \text{ otherwise}\end{cases} \]
with identity attaching morphisms.  This category has enough injectives and it is not difficult to see that the
derived category $D^b(Mod-A)$ is naturally equivalent to the full derived category $D^b_c(K)$.  The general theory
of Koszul duality \cite{BeilinsonGS} guarantees an equivalence of categories $D^b(Mod-A) \sim D^b(Mod-B)$ with
the category of modules over the Koszul dual algebra.  This is the main idea in the proof of
Proposition \ref{prop-Polish-equivalence}, that is,
\[ D^b_c(K) \simeq D^b(Mod-A) \simeq D^b(Mod-B) \simeq D^b(\mathcal P(K)) \simeq D^b(\mathcal P^{\Delta}(K)).\]

\part{Appendices}
\begin{appendix}
\section{Transversality}\label{appendix-transversality}\index{transversality}
Let $A,B\subset M$ be a smooth submanifolds of a smooth manifold $M$.  Let $\{ v_1, \cdots, v_r\}$
be vector fields on $M$ which span the tangent space at each point\footnote{In fact we only need to assume
that they span the normal space to $B$ at each point in $B$.}  Let $V$ be the $r$-dimensional vector space of
formal linear combinations of these vector fields.  Let $\Phi_v$ denote the time $=1$ flow of the vector field $v$.

\begin{thm}  The set of elements $v \in V$ such that $\Phi_v$ fails to take $A$ transversally to $B$ has measure zero.
  \end{thm}
  
  \begin{proof}
  The method of proof (using Sard's theorem, of course) is due originally to Marston Morse.  Let
  $\Phi:V \times M \to M$ be the time $=1$ flow mapping.  The assumptions on the vector fields imply that
  $\Phi$ takes $V \times A$ transversally to $B$.  (In fact for each $a\in A$ the partial map $\Phi_a:V \to M$ is a
  submersion.)  So $\Phi^{-1}(B) $ is a submanifold of $V \times M$ and it is transverse
to $V \times A$.  Let $\pi:\Phi^{-1}(B)\cap (V \times A) \to V$ denote the
  projection.  We claim the following are equivalent: \begin{enumerate}
  \item $v \in V$ is a regular value of $\pi$
  \item $\{v\} \times A$ is transverse to $\Phi^{-1}(B)$ in $V \times M$
  \item $\Phi_v$ takes $A$ transversally to $B$.
  \end{enumerate}
Equivalence of (2) and (3) is a short calculation.
  Equivalence of (1) and (2) is a dimension count.  
First note that $(V\times A)$ is transverse to $\phi^{-1}(B)$ in $V \times M$,
by the assumption on the vector fields.
Fix $(w,a) \in (V\times A) \cap \Phi^{-1}(B)$.  Let $\dot{A} = T_aA$.  Let $\dot{M} = T_aM$,
let $\dot{V} = T_wV$ and let $\dot{B} = T_{w,a}\Phi^{-1}(B)$.  
 So we have a diagram of vector spaces
\[ \begin{diagram}[size=2em]
 & & \dot{V} \times \dot{M} && \\
& & \uInto & &&\\
&& (\dot{V} \times \dot{A}) & \cap \dot{B} &
\end{diagram}\]

Let $r$ be the rank of
the linear map $\pi:\dot{B} \cap(\dot{V} \times \dot{A}) \to \dot{V}$; its kernel is 
$\dot{B} \cap (\{0\} \times A$.  Therefore
\[ \dim(\dot{B} \cap (\dot{V} \times \dot{A})) = r + \dim(\dot{B} \cap (\{0\} \times \dot{A})\]
or
\[ b+(v+a)-m = r + (b+a - \dim(\dot{B}+\{0\}\times \dot{A}\]
where $b = \dim(\dot{B})$, $a = \dim(\dot{A})$, $v = \dim(\dot{V})$.  Therefore
\[ v-r = m - \dim(\dot{B} + \{0\} \times \dot{A}).\]
Therefore the map $\pi$ is surjective ($r = v$) if and only if $\dot{B} + \{0\}\times \dot{A}$ spans
$\dot M$.

  Now, by Sard's theorem the set of critical values $v \in V$ has measure zero.  But this is exactly the
  set of elements such that $\Phi_v$ fails to take $A$ transversally to $B$.

  \end{proof}
\subsection{}  Exactly the same argument shows, for example, that two submanifolds of Euclidean space
can be made transverse by an arbitrarily small {\em translation}.  Similarly, two stratified subvarieties of
projective space can be made transverse by an arbitrarily small projective transformation. 
The same method also applies to
transversality of maps:  if $F:M' \to M$ is a smooth mapping, if
  $A\subset M'$ and if $B \subset M$ then by composing $F$ with the time $=1$ flow of an arbitrarily
  small smooth vector field we can guarantee that the resulting map $\tilde{F}:M' \to M$ takes
  $A$ transversally to $B$.

\section{Proof of Theorem \ref{thm-equivalence}}
The following lemma  provides lifts of morphisms in the derived category, see \cite{IH2}
\begin{tcolorbox}[colback=cyan!30!white]
\begin{lem}\label{lem-lifting}
Let $A\b, B\b$ be a objects in the derived category.  Suppose $\uul{H}^r(A\b) = 0$ for all
$r >p$ and suppose that $\uul{H}^r(B\b) = 0 $ for all $r<p$.  Then the natural map
\[ \hHom_{D^b_c(X)}(A\b, B\b) \to \hHom_{Sh(X)}(\uul{H}^p(A\b), \uul{H}^p(B\b))\]
is an isomorphism.
\end{lem}\end{tcolorbox}
\begin{proof}  When we wrote IH II, Verdier (who was one of the referees) showed us how to replace our
4 page proof with the following simple proof.  Up to quasi-isomorphism it is possible to replace the complexes
$A\b$, $B\b$ with complexes
\[
\begin{CD}
 \cdots @>>> A^{p-1} @>{d_A}>> A^p @>>> 0 @>>>0 @>>> \cdots \\
 \cdots @>>> 0    @>>>  I^p @>>{d_B}>  I^{p+1}@>>> I^{p+2} @>>> \cdots
 \end{CD}\]
where $I^r$ are injective.  This means that a morphism in the derived category is represented by an honest
morphism between these complexes, that is, a mapping
\[\phi: \uul{H}^p(A\b) = \coker(d_A) \to \ker(d_b) = \uul{H}^p(B\b).\qedhere\]
\end{proof}

\subsection{Proof of Theorem \ref{thm-equivalence}}\label{ProofOfEquivalence}
We have a Whitney stratification of $W$ and inclusions
\[ \begin{CD}
U_2 @>>{j_2}> U_3 @>>{j_3}>  \cdots @>>{j_{n-1}}> U_n @>>{j_n}> U_{n+1}=W \end{CD}\]
Let us suppose that $A\b$ is constructible with respect to this stratification and that it satisfies the
support (but not necessarily the co-support) conditions, that is
\[ H^r(A\b)_x = 0 \text{ for } r \ge p(c)+1\]
whenever $x \in X^{n- c}$ lies in a stratum of codimension $c$.  Fix $k \ge 2$ and consider the situation
 \[ \begin{CD}
U_k @>>{j_k}> U_{k+1} @<<{i_k}< X^{n-k} \end{CD}\]
where $X^{n-k}$ is the union of the codimension $k$ strata.  Let $A^{\bullet}_k = A\b|U_k$.
Let $\bar q$ be the complementary perversity, $q(c) = c-2-p(c)$.
The following proposition says that the vanishing of the stalk cohomology with
compact supports $H^r(i_x^!A\b)$ is equivalent to the condition that the attaching map is an isomorphism:
\begin{tcolorbox}[colback=cyan!30!white]
\begin{prop} \label{prop-attachingequivalent} The following statements are equivalent.
\begin{enumerate}
\item $A^{\bullet}_{k+1} \cong \tau_{\le p(k)}Rj_{k*}A^{\bullet}_k$
\item $\uul{H}{}^r(A^{\bullet}_{k+1})_x \to \uul{H}{}^r(Rj_{k*}A^{\bullet}_k)_x$ is an isomorphism for all $x \in X^{n-k}$
\item $H^r(i_k^!A^{\bullet}_{k+1}) =0$ for all $r \le p(k)+1$
\item $H^r(i_x^!A^{\bullet}_{k+1}) = 0$ for all $r < n-q(k)$ for all $x \in X^{n-k}$
\end{enumerate}
\end{prop}\end{tcolorbox}
Since the sheaf $\uul{IC}{}^{\bullet}_{\bar p}[-n]$ satisfies these conditions, this proposition (and  induction)
gives another proof of the isomorphism $\uul{IC}{}^{\bullet}_{\bar p}[-n] \cong  \uul{P}{}^{\bullet}_{\bar p}$ 
of Theorem \ref{thm-ICbytruncation}.
\begin{proof}
Items (1) and (2) are equivalent because there is a canonical morphism
\[ A^{\bullet}_{k+1} \to Rj_{k*}j_k^*A\b = Rj_{k*}A_k\]
truncation $\tau_{\le p(k)}$ leaves an isomorphism in degrees $\le p(k)$.
Items (3) and (4) are equivalent because $i_x:\{x\} \to X^{n-k}$ is the inclusion into a manifold so $i_x^! = i_k^![n-k]$,
and because $r <p(k)+2 + (n-k) = n -(k-2-p(k)) = n-q(k)$.
Items (2) and (3) are equivalent because there is a distinguished triangle,
\begin{diagram}[size=2em]
Ri_{k*}i_k^!({A}{}^{\bullet}_{k+1}) & & \rTo &&{A}{}^{\bullet}_{k+1}\\
& \luTo && \ldTo_{\alpha} &\\
&&Rj_{k*}j_k^*({A}{}^{\bullet}_{k+1})&&
\end{diagram}
and therefore an exact sequence on stalk cohomology as follows:
\newcommand{\Ako}{A^{\bullet}_{k+1}} \newcommand{\Ak}{A^{\bullet}_k}
\[  \begin{diagram}[size=2em]
\uul{H}^{p+2}(i_k^!\Ako)_x & \rTo & \highlight{\uul{H}^{p+2}(\Ako)_x} &\rTo& \uul{H}^{p+2}(Rj_{k*}\Ak)_x & & \\
&\luTo(4,2)&\\
\uul{H}^{p+1}(i_k^!\Ako)_x & \rTo & \highlight{\uul{H}^{p+1}(\Ako)_x} &\rTo& \uul{H}^{p+1}(Rj_{k*}\Ak)_x & & \\
&\luTo(4,2)&\\
 \uul{H}^{p}(i_k^!\Ako)_x & \rTo & \uul{H}^{p}(\Ako)_x &
 \rTo_{\highlight[green]{\alpha}}& \uul{H}^{p}(Rj_{k*}\Ak)_x & & \\
&\luTo(4,2)&\\
 \uul{H}^{p-1}(i_k^!\Ako)_x & \rTo & \uul{H}^{p-1}(\Ako)_x &
 \rTo_{\highlight[green]{\alpha}}& \uul{H}^{p-1}(Rj_{k*}\Ak)_x & & 
\end{diagram}\]
Now use the fact that the yellow highlighted terms are zero and the green highlighted morphisms are isomorphisms
to conclude the proof of the Proposition.
\end{proof}

\subsection{Continuation of the proof of Theorem \ref{thm-equivalence}}
 Now let us show that if $\mathcal E_1, \mathcal E_2$ are local systems on $U_2$ and if
$A\b = \uul{P}{}^{\bullet}_{\bar p}(\mathcal E_1)$ and if $B\b = \uul{P}{}^{\bullet}_{\bar p}(\mathcal E_2)$ then
we have an isomorphism
\[ \hHom_{Sh}(\mathcal E_1,\mathcal E_2) \cong \hHom_{D^b_c(X)}(A\b, B\b).\]
As before, let $A^{\bullet}_{k+1} = A\b|U_{k+1} = \tau_{\le p(k)}Rj_{k*}A^{\bullet}_k$.
Assume by induction that we have established an isomorphism
\[ \hHom(\mathcal L_1, \mathcal L_2) \cong \hHom_{D^b_c(U_k)}(A^{\bullet}_k, B^{\bullet}_k).
\]

Using the above triangle for $B\b$ we get an exact triangle of $\uul{\RHom}$ sheaves,
\[
\begin{diagram}[size=2em]
  \uul{\RHom}^{\bullet}(A^{\bullet}_{k+1} , Ri_{k*}i_k^!B^{\bullet}_{k+1}) && \rTo &&\uul{\RHom}^{\bullet}
  (A^{\bullet}_{k+1}, B^{\bullet}_{k+1}) \\
  & \luTo_{[1]} && \ldTo_{\alpha} && \\
&&  \uul{\RHom}^{\bullet}(A^{\bullet}_{k+1}, Rj_{k*}B^{\bullet}_k) &&
\end{diagram}\]
By Lemma \ref{lem-lifting} and the support conditions, we see that $\alpha$ is an isomorphism in degree zero,
\[ \hHom_{D^b_cX}(A^{\bullet}_{k+1}, B^{\bullet}_{k+1}) = H^0(U_{k+1}; \uul{\RHom}^{\bullet}(A^{\bullet}_{k+1},
B^{\bullet}_{k+1}) \cong H^0(U_{k+1}; \uul{\RHom}^{\bullet}(A^{\bullet}_{k+1}, Rj_{k*}B^{\bullet}_k)).\]
Moreover, 
\[Rj_{k*}\uul{\RHom}^{\bullet}(A^{\bullet}_k, B^{\bullet}_k) \cong
Rj_{k*}\uul{\RHom}^{\bullet}(j_k^*A^{\bullet}_{k+1}, B^{\bullet}_k) \cong
\uul{\RHom}^{\bullet}(A^{\bullet}_{k+1}, Rj_{k*}B^{\bullet}_k)\]
by the standard identites (above), whose cohomology is
\[ H^0(U_{k+1}; Rj_{k*}\uul{\RHom}^{\bullet}(A^{\bullet}_k, B^{\bullet}_k)) \cong
H^0(U_k; \uul{\RHom}^{\bullet}(A^{\bullet}_k, B^{\bullet}_k)) \cong
\hHom_{D^b_c(U_k)}(A^{\bullet}_k, B^{\bullet}_k).\]
So, putting these together we have a canonical isomorphism
\[ \hHom_{D^b_c(U_k)}(A^{\bullet}_k, B^{\bullet}_k) \cong
\hHom_{D^b_c(U_{k+1}}(A^{\bullet}_{k+1}, B^{\bullet}_{k+1})\]
which was canonically isomorphic to $\hHom_{Sh}(\mathcal E_1, \mathcal E_2)$ by induction.
This completes the proof of the theorem, but the main point is that the depth of the argument is the moment in which
Lemma \ref{lem-lifting} was used in order to lift a morphism $A^{\bullet}_{k+1} \to Rj_{k*} B^{\bullet}_k$  to a
morphism $A^{\bullet}_{k+1} \to B^{\bullet}_{k+1}$.  \qed

\section{Complex Morse theory of perverse sheaves}\label{sec-Morseperverse}
\subsection{The complex link}\index{link!complex}\index{complex link}
 Details for this section may be found in \cite{SMT} part II \S 2.
(See also \cite{Schuermann} chapt. 5.)
 Throughout this section $W$ denotes a complex analytically stratified complex analytic variety  of
complex dimenion $n$ contained in some complex manifold $M$, say, of dimension $m$. 
Let $X$ be a stratum of $W$, let $p\in X$
and let $N = \widetilde{N} \cap W$ be a normal slice at $p$, that is,
the intersection with a complex submanifold $\widetilde{N}$ that
is transverse (and complementary) to $X$ with $\widetilde \cap X = \{p\}$.  

Suppose $\xi \in T^*_{X,p}M$ is a nondegenerate covector. Since
\[ T^*_pM = \hHom_{\CC}(T_pM,\CC) \cong \hHom_{\RR}(T_pM,\RR)\]
we may choose local coordinates near $p$  that  identify $\xi$ with the linear projection $\CC^m \to \CC$
with $p = 0 \in \CC^m$ and $\pi(p) = 0$.
The restriction $\pi|N:N \to \CC$ takes the zero dimensional singular point to the origin.  Every other
stratum of $N$ is taken sumersively to $\CC$ (over a neighborhood of $0$), so $N \to \CC$ is a
stratified fiber bundle except at the point $p$.  Let $r(z)$ denote the square of the distance from $p$ in
$\CC^m$.  As in Lemma \ref{lem-epsilon} there is a region $0 < \delta << \epsilon \subset \RR^2$ such that for any
pair $(\delta,\epsilon)$ in this region the following holds:
\begin{enumerate}
\item $\partial B_{\epsilon}(p)$  is transverse to every stratum of $N$ (where $B_{\epsilon}(p)$ denotes the
ball of radius $\epsilon$).
\item  For each stratum $Y$ of $N$ (except $Y = \{p\}$) the restriction $\pi|Y$ has no critical values in
the disk $D_{\delta}(0) \subset \CC$
\item For each stratum $Y$ of $N$ with $\dim_{\CC}(Y) \ge 2$,
 and for each point $z \in Y \cap \partial B_{\epsilon}(p)$ if $|\pi(z)| \le \delta$ then
the complex linear map
\[ (dr(z), d\pi(z)):T_zY \to \CC^2\]
has rank $2$.
\end{enumerate}
Identify $\delta = \delta + 0i \in \CC$ and define the {\em complex link} (\cite{SMT} II \S 2):
\[ \mathcal L = \pi^{-1}(\delta) \cap N \cap B_{\epsilon}(0), \ \partial\mathcal L = \pi^{-1}(\delta)
\cap \partial B_{\epsilon}(0).\]

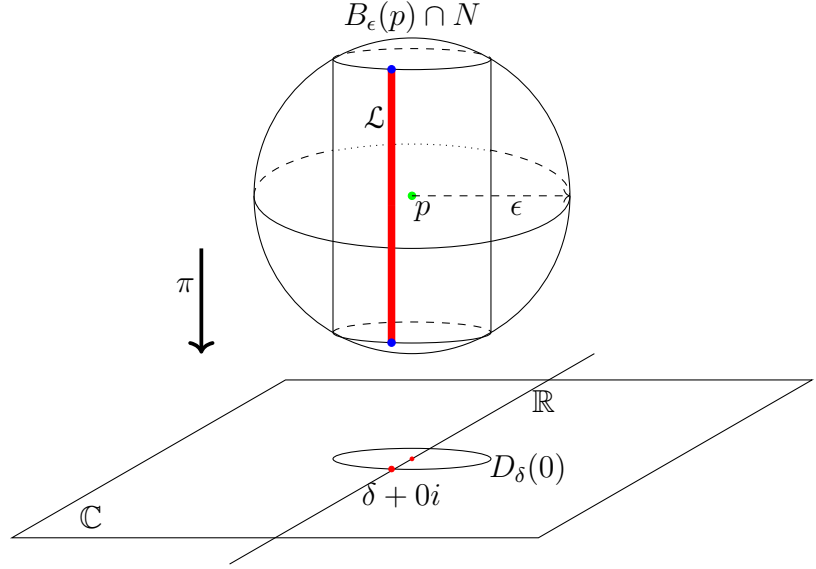
\begin{figure}[H]\centering
\begin{tikzpicture}[scale=.7]
\draw (0,0) circle (3cm);
\draw (-3,0) arc (180:360:3cm and 1cm);  
\begin{scope}
\clip (-3,0) rectangle ++(1.5,1);
\draw[dashed] (-3,0) arc (180:0:3cm and 1cm);  
\end{scope}
\begin{scope}
\clip(3,0) rectangle ++(-1.5,1);
\draw[dashed] (-3,0) arc (180:0:3cm and 1cm);  
\end{scope}
\pgfmathsetmacro{\x}{sqrt(3)};\pgfmathsetmacro{\y}{3*sqrt(3)/2};

\draw[dotted] (1.5,.85) arc (60:120:3 cm and 1 cm);
\draw (-1.5,\y) arc (180:360: 1.5 cm and .2 cm);  
\draw[dashed] (-1.5,\y) arc (180:0: 1.5 cm and .2 cm);  
\draw (-1.5, -\y) arc (180:360: 1.5 cm and .2 cm);  
\draw[dashed] (-1.5,-\y) arc (180:0: 1.5 cm and .2 cm);  
\draw (-1.5,\y) -- (-1.5,-\y); \draw (1.5,\y) -- (1.5,-\y);  
\pgfmathsetmacro{\z}{-5};

\pgfmathsetmacro{\angl}{30};
\draw(-5,\z) -- ++(\angl:3cm) -- ++(10,0) -- ++(\angl:-6cm) -- ++(-10,0) -- cycle;
\draw(0,\z) -- ++(\angl:4cm); \draw(0,\z) -- ++ (\angl:-4cm);  
\draw(0,\z) ellipse (1.5 cm and .2cm);
\draw [red, fill=red] (0,\z) circle (1pt); 

\coordinate (P) at (255:1.5 cm and .2 cm); 
\coordinate (Q) at ($(P) + (0,\z)$);  
\draw [red, fill=red] (Q) circle (1.5pt); 
\draw[line width = 1mm, red] (P) -- ++(0,\y); \draw[line width = 1 mm, red] (P) -- ++(0,-\y);
\draw[blue,fill=blue] ( $(P) + (0,\y)$) circle (2pt); 
\draw[blue, fill=blue] ( $(P) + (0,-\y)$) circle (2pt);
\draw [green, fill=green] (0,0) circle (2pt);
\draw [dashed,->] (0,0) -- ++(0:3cm); 

\node at ({-\x/2.5},1.5) () {$\mathcal L$};
\node at (-6.1,{\z-1.1}) () {$\mathbb C$};
\node at (2.5,{\z+1.1}) () {$\mathbb R$};
\node at (-.2, {\z -.7}) () {$\delta + 0i$};
\node at (-7:2cm) () {$\epsilon$};
\node at (0,{\y + .8}) () {$B_{\epsilon}(p)\cap N$};
\node at (2.2, {\z-.2} ) () {$D_{\delta}(0)$};
\node at (.2,-.3) () {$p$};

\draw [line width = .5mm, ->] (-4,-1) -- (-4,-3);
\node at (-4.3,-1.7) () {$\pi$};
\end{tikzpicture}
\caption{{Normal slice with complex link, $\epsilon$-ball and $\delta$-disk}}
\label{fig-complexlink}
\end{figure}

It is a single fiber of the (stratified) fiber bundle over the circle $S^1=\partial D_{\delta}\subset \CC$,
\[ \mathcal E =\pi^{-1}(\delta e^{i \theta}) \cap N \cap B_{\epsilon}(0),\
 \partial\mathcal E = \pi^{-1}(\delta e^{i \theta}) \cap N \cap \partial B_{\epsilon}(0)\]
\[(\mathcal E, \partial \mathcal E) \to 
S^1=\partial D_{\delta} = \{ \delta e^{i\theta}|\ 0 \le \theta \le 2\pi\} \] 
and the boundary $\partial \mathcal E$ is a trivial bundle over $S^1 = \partial D_{\delta}$.
See Figure \ref{fig-complexlink}.  The stratified homeomorphism type of the spaces $\mathcal L, \partial\mathcal L$
is independent of the choice of $\epsilon,\delta$.
\begin{thm}\label{theorem-link-structure} \cite{SMT}
The bundle $\mathcal E$ is stratified-homeomorphic to the mapping cylinder of a
(stratified) monodromy  homeomorphism
\[ \mu: (\mathcal L, \partial \mathcal L) \to (\mathcal L, \partial \mathcal L)\]
that is the identity on $\partial \mathcal L$ and  is well defined up to stratum preserving isotopy.  
The link $L_X(p)$   is homeomorphic to the ``cylinder with caps'', 
\begin{equation}
\label{equation-link} L_X(p) \cong \mathcal E \cup_{\partial \mathcal E} \left(\partial \mathcal L \times D_{\delta}\right).\end{equation}
\end{thm}

\subsection{The braid diagram}\label{subsec-braid}\index{braid diagram}
All these facts fit together into a single commutative diagram that explains their cohomological interrelationships.
Throughout this section we fix a constructible complex of sheaves $A\b \in D^b_c(W)$ of $K$-vector spaces, where
$K$ is a (coefficient) field.  Fix a point $x \in W$ in some stratum
$X$ with normal slice $N$.  Fix a nondegenerate covector $\xi \in T^*_{X,x}M$ with resulting  
complex link $\mathcal L$.  As in Theorem \ref{theorem-MorseSheaf}
 we may identify $\xi = d\phi$ with the
differential of a locally defined (near $x$) smooth function $\phi:M \to \RR$
such that $\phi(x)=0$ and $d\phi|T_xX = 0$.  Let 
\[ H^r(\xi,A\b) = H^r(N_{\le \delta}, N_{<0}; A\b) = H^r(R\Gamma_ZA\b)\]
denote the Morse group for $\xi$ in degree $r$ for $\delta >0$ sufficiently small.  
Let $L = L_X(x)$ denote the link of the stratum $X$ at $x$.

The homeomorphisms
described in the Theorem \ref{theorem-link-structure} are stratum preserving so they induce isomorphisms
on cohomology with coefficients in $A\b$ and they allow us to interpret these cohomology groups,
\begin{align*}
H^r(N-x; A\b)  &\cong  H^r(L_X(x);A\b) & H^r(N,N-x;A\b) &\cong  H^r(i_x^!A\b)\\
H^r(N;A\b) &\cong H^r(i_x^*A\b) &H^r(N,N_{<0};A\b) &\cong H^r(\xi;A\b)  \\
H^r(N_{<0};A\b) &\cong H^r(\mathcal L; A\b)&
H^{r+1}(N-x, N_{<0};A\b)&\cong H^{r}(\mathcal L,\partial\mathcal L; A\b)
\end{align*}
By (\ref{equation-link})  the ``variation" map $I-\mu: H^r(\mathcal L;A\b) \to 
H^r(\mathcal L, \partial \mathcal L; A\b)$ may be
identified with the connecting homomorphism in
the long exact sequence for the pair $(N-x, N_{<0})$,  cf. \cite{GM2}.
\quash{
\begin{diagram}[size=2em]
&\rTo & H^r(N_{<0};A\b) &\rTo^{\delta} & H^{r+1}(N-x, N_{<0};A\b) & \rTo\end{diagram}


The relation between the Morse group, stalk cohomology and stalk cohomology with
compact supports are best expressed by the three exact sequences that correspond to the triple
of spaces}
As in \cite{SMT} p. 215 the three  long exact sequences for the triple of spaces
\[ N_{<0} \subset N-\{x\} \subset N\]
may be assembled into a braid diagram with 
exact sinusoidal rows.  
(cf. \cite{Schuermann} \S 6.1 where the same sequences are considered separately):

\newcommand{\sst}{}
\newcommand{\AZ}{\sst H^{r}(i_x^*A\b)}
\newcommand{\BZ}{\sst H^r(\mathcal L, A\b)}
\newcommand{\CZ}{\sst H^{r}(\mathcal L,\partial\mathcal L;A\b)}
\newcommand{\DZ}{\sst H^{r+2}(i_x^!A\b)}
\newcommand{\XZ}{\sst H^r(L;A\b)}
\newcommand{\YZ}{\sst H^{r+1}(\xi; A\b)}
\newcommand{\VZ}{\sst H^{r+1}(L;A\b)}
\newcommand{\WZ}{\sst H^{r+1}(\xi, A\b)}
\newcommand{\aZ}{\sst H^{r-1}(\mathcal L,\partial\mathcal L;A\b)}
\newcommand{\bZ}{\sst H^{r+1}(i_x^!A\b)}
\newcommand{\cZ}{\sst H^{r+1}(i_x^*A\b)}
\newcommand{\dZ}{\sst H^{r+1}(\mathcal L;A\b)}

\begin{figure}[h!]
\begin{tikzcd}[column sep=-4ex]\centering
\AZ \arrow[bend left]{rr}   \arrow[dr] & &   \BZ \arrow[bend left,"var"]{rr} \arrow[dr] & & 
\CZ \arrow[bend left]{rr} \arrow[dr]&& \DZ \\
& \XZ \arrow[ur]\arrow[dr] && \YZ \arrow[ur]\arrow[dr]&& \VZ \arrow[ur]\arrow[dr] \\
\aZ \arrow[bend right]{rr}   \arrow[ur] & &   \bZ \arrow[bend right]{rr} \arrow[ur] & & 
\cZ \arrow[bend right]{rr} \arrow[ur]&& \dZ 
\end{tikzcd}
\caption{Braid diagram}\end{figure}
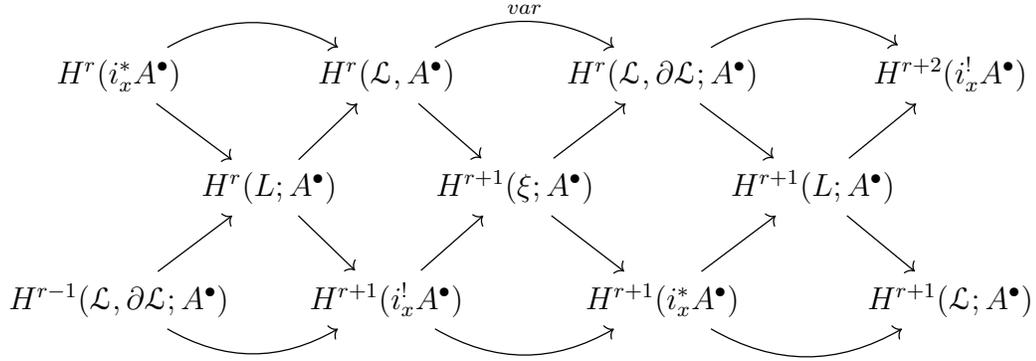

\subsection{Levi form} \label{subsec-Levi}\index{Levi form}
 Theorem \ref{thm-SMT} reduces the Morse theory of $W$ to that of $N$, a variety with
a zero dimensional stratum.  Theorem \ref{theorem-link-structure}  reduces the analysis of the singularities of 
$N$ to an analysis of the complex link $\mathcal L$.  This in turn can be analyzed with Morse theory.
For generic choice of $\delta$ (in the definition of $\mathcal L$)
the $(\text{distance})^2$ from the singular point $p$ is a Morse function on $\mathcal L$ with a 
constant maximum value on $\partial L$. The  Levi form $\partial^2f(x)/(\partial z_i \partial \bar z_j)$ of
$f(z) = \sum_i z_i \bar z_i$ is the identity $I_{m\times m}$ and the  Levi form of $-f$ is $-I_{m \times m}$.
 Let $S$ be a stratum of $\mathcal L$, let $y \in S$ be a
critical point of $f|S$.
The Levi form of $f|S$ equals the restriction of the Levi form of $f$ to $T_yS$. From this, and a
few facts from linear algebra (\cite{SMT} \S 4.A) it follows that for each (nondegenerate) critical point
$y \in S$,  {\em the Morse index of $f|S$ is $\le \dim_{\CC}S$} and
{\em the Morse index of $-f|S$ is $\ge \dim_{\CC}S$}.

\subsection{Vanishing conditions}\label{subsec-vanishing}
The braid diagram, together with induction and the estimates in \S \ref{subsec-Levi},
 may be used to prove 
many vanishing theorems and Lefschetz-type theorems in sheaf cohomology, see \cite{Grothendieck},
\cite{GM2}, \cite{SMT}, \cite{Schuermann},
\cite{Kashiwara}, \cite{Hamm}, \cite{Hamm0}.  
Using the convex function distance${}^2$ from the point $\{x\}$ below, the braid diagram, and induction, one finds 
the following two results (cf.~\cite{Artin}) which may be proven together
 (since the statement for one becomes the inductive step for the other):
\begin{thm}\label{prop-Hestimate}
Suppose $A\b\in D^b_c(W)$ is a complex of sheaves of $K$-vector spaces on
a complex analytic set $W$ such that for each stratum $X$ and for
each point $x \in X$, with $i_x:\{x\} \to W$, the stalk cohomology vanishes: 
\begin{equation}\label{eqn-Pplus}H^r(i_x^*A\b) = 0 \text{ whenever }\
r > {\rm codim}_{\CC}(X).\end{equation}  
Then $H^r(\mathcal L; A\b) = 0$ for all $r > \ell = \dim_{\CC}(\mathcal L)$.
If the Verdier dual $\DD_W(A\b)$ satisfies (\ref{eqn-Pplus}), or equivalently, 
(if $W$ has pure (complex) dimension $n$ and)
\begin{equation}\label{eqn-Pminus}
H^r(i_x^!A\b) = 0  \text{ whenever }\ r < n+\dim_{\CC}(X),\end{equation}  
then $H^r(\mathcal L, \partial \mathcal L; A\b) = 0$ for all $r < \ell$.
\end{thm}

\begin{thm}\label{thm-Stein}
Let $W$ be a Stein space or an affine complex algebraic variety of dimension $n$ and let
$A\b\in D^b_c(W)$ be  a complex of sheaves on $W$ that satisfies (\ref{eqn-Pplus}).  Then $H^r(W;A\b) = 0$ for
all $r > n$.  Let $W$ be a projective variety and let $H$ be a hyperplane that is transverse
to each stratum of some Whitney stratification of $W$.  Let $A\b$ be a complex of sheaves on $W$ that
satisfies (\ref{eqn-Pminus}).  Then $H^r(W, W\cap H;A\b) = 0$ for all $r< n$. 
\end{thm}
By examining the braid diagram again it is easy to see that the statements at the end of Theorems
\ref{prop-Hestimate} and \ref{thm-Stein} provide the following Morse-theoretic characterization:  
a sheaf is perverse if and only if the Morse groups lie in a single degree.  To be explicit,

\begin{tcolorbox}[colback=cyan!30!white]
\begin{cor}
The complex $A\b$ is  in ${}^pD_c(W)^{\le 0}$ (resp. ${}^pD_c(W)^{\ge 0}$) if and only if  for every stratum 
$X$ of $W$ and every point $x \in X$ and every
nondegenerate covector $\xi \in T_{X,x}^*M$ the Morse group $M^j(\xi;A\b)$ vanishes for all
$j > \cod_{\CC}(X)$ (resp. for all $j < \cod_{\CC}(X)$). 
\end{cor}\end{tcolorbox}
  
\end{appendix}

\printindex

\end{document}